\documentclass[11pt]{article}

\usepackage[numbers]{natbib}
\usepackage[T1]{fontenc}
\usepackage{lmodern}
\usepackage{microtype}
\usepackage{amsmath}
\usepackage{mathtools}
\usepackage{amssymb}
\usepackage{amsthm}
\usepackage{xcolor}

\usepackage{array}
\usepackage{booktabs}
\usepackage{pgfplots}
\usetikzlibrary{arrows.meta,calc}
\definecolor{overviewblue}{RGB}{31,77,116}
\definecolor{overviewteal}{RGB}{28,107,97}
\pgfplotsset{compat=1.18}
\usepackage{tabularx}
\usepackage{hyperref}
\usepackage{algorithm}
\makeatletter
\providecommand*{\toclevel@algorithm}{0}
\makeatother
\usepackage{algpseudocode}
\usepackage{enumitem}
\usepackage[margin=1in]{geometry}

\hypersetup{colorlinks=true,citecolor=blue!45!black,linkcolor=blue!45!black, urlcolor=blue!45!black,hypertexnames=false}
\setlist[itemize]{topsep=4pt,itemsep=3pt,parsep=0pt}
\numberwithin{equation}{section}
\allowdisplaybreaks

\newcommand{\paperTitle}{A Matrix Lee--Sidford Barrier and Faster Semidefinite Programming}
\newcommand{\paperAuthor}{Zhao Song\thanks{\texttt{magic.linuxkde@gmail.com}. Independent Researcher.} \and Ruizhe Zhang\thanks{\texttt{rzzhang@purdue.edu}. Purdue University.}}

\DeclareMathOperator{\tr}{tr}
\newcommand{\E}{\mathbb E}
\newcommand{\R}{\mathbb R}
\newcommand{\Tmat}{{\cal T}_{\mathrm{mat}}}
\renewcommand{\d}{\mathrm{d}}

\theoremstyle{plain}
\newtheorem{theorem}{Theorem}[section]
\newtheorem{lemma}[theorem]{Lemma}
\newtheorem{fact}[theorem]{Fact}
\newtheorem{definition}[theorem]{Definition}

\newtheorem{corollary}[theorem]{Corollary}
\theoremstyle{definition}

\newtheorem{problem}[theorem]{Problem}
\newtheorem{remark}[theorem]{Remark}

\begin{document}

\date{October 31, 2026}
\title{\paperTitle}
\author{\paperAuthor}
\maketitle

\begin{abstract}
We construct a matrix analog of the Lee--Sidford barrier for semidefinite programming (SDP) with $m$ equality constraints and $n$-by-$n$ matrices. The barrier has self-concordance parameter $\widetilde{O}(m)$, so short-step interior-point methods require only $\widetilde O(\sqrt m)$ iterations. We also show that our matrix Lee--Sidford barrier can be implemented efficiently by developing a randomized algorithm that solves SDPs with $m\leq n$ to accuracy $\epsilon$ in $(n^{\omega+o(1)}+n^2m^{3/2+o(1)})\operatorname{polylog}(1/\epsilon)$ time, where $\omega\leq 9/4$ is the matrix multiplication exponent, achieving the fastest known running time in this regime.
\end{abstract}

{\hypersetup{linkcolor=black}
\tableofcontents
}
\clearpage
\section{Introduction}
\label{sec:introduction}

Semidefinite programming (SDP) optimizes a linear function of a positive semidefinite matrix subject to linear constraints, and it is one of the most powerful algorithmic tools in theoretical computer science. The relaxation of Goemans and Williamson~\cite{gw95} gives the best-known approximation ratio for Max-Cut, and SDP relaxations also yield approximation algorithms for graph coloring~\cite{kms98} and sparsest cut~\cite{arv09}. The Lov\'asz theta function~\cite{l79} is the value of an SDP, leading to polynomial-time algorithms for the maximum stable set and the chromatic number of perfect graphs~\cite{gls88}. Hierarchies of SDP relaxations, such as the sum-of-squares hierarchy~\cite{par03,lass01}, turn polynomial optimization into a sequence of SDPs with increasingly strong certificates, and they underlie many recent algorithms for average-case problems~\cite{bs14}.

For a broad class of problems, SDP is in fact optimal. Under the Unique Games Conjecture of Khot~\cite{k02}, the Goemans--Williamson ratio is optimal for Max-Cut~\cite{kkmo07,moo10}, and Raghavendra~\cite{r08} showed that a basic SDP relaxation achieves the best polynomial-time approximation ratio for every constraint satisfaction problem. With the recent Lean-formalized resolution of the Unique Games Conjecture~\cite{o26c}, these optimality results now hold under $\mathsf P\neq\mathsf{NP}$.

SDP is also a standard modeling tool across many applied fields, including operations research and combinatorial optimization~\cite{wsv00}, control theory~\cite{begfb94}, statistics and machine learning~\cite{lcbgj04,dgjl07,cr09}, quantum information~\cite{npa08,ffs24}, quantum chemistry~\cite{m04}, and theoretical physics~\cite{psv12,sd15,prv19}.

Let $\mathbb S^n$ denote the set of real symmetric $n\times n$ matrices, with trace inner product $\langle U,V\rangle:=\tr[U^\top V]$. Given $A_1,\ldots,A_m,C\in\mathbb S^n$ and $b\in\R^m$, the primal SDP is defined as
\begin{equation}\label{eq:primal_sdp}
 \max_{X\succeq0}\ \langle C,X\rangle
 \quad\text{subject to}\quad
 \langle A_i,X\rangle=b_i\quad\forall i\in[m],
\end{equation}
and its dual is
\begin{equation}\label{eq:intro_dual_sdp}
 \min_{y\in\R^m}\ b^\top y
 \quad\text{subject to}\quad
 S(y):=\sum_{i=1}^m y_iA_i-C\succeq0.
\end{equation}
The matrix $S(y)$ is the dual slack, $m$ is both the number of primal equality constraints and the dimension of the dual variable, and $n$ is the matrix dimension. At feasible points, the duality gap $b^\top y-\langle C,X\rangle=\langle S(y),X\rangle\geq0\geq 0$. If the dual problem is strictly feasible and has a finite optimal value, then the two optimal values agree~\cite[Theorem~3.1]{vb96}.

The SDP~\eqref{eq:primal_sdp}--\eqref{eq:intro_dual_sdp} is highly expressive. When $C,A_1,\ldots,A_m$ are diagonal, it reduces to a linear program (LP) with $m$ dual variables and $n$ inequality constraints. However, this expressiveness comes at a computational price. High-accuracy solvers\footnote{A \emph{high-accuracy} solver is an algorithm whose running time to output an $\epsilon$-accurate solution depends polylogarithmically on $1/\epsilon$.} for dense SDPs rely on two main frameworks: interior-point methods (IPMs) and cutting-plane methods (CPMs). For dense instances, the running time of an IPM-based SDP solver has two components: the number of iterations and the computational cost per iteration. The classical IPM solvers take $\widetilde O(\sqrt n)$ iterations, and each iteration forms and solves an $m\times m$ Newton system~\cite{nn92,nn94,ali95}. Jiang, Kathuria, Lee, Padmanabhan, and Song~\cite{jklps20} maintain approximate Hessians across iterations, thereby reducing the total cost to $$\widetilde O(\sqrt n(mn^2+m^\omega+n^\omega)),$$ where $\omega\leq 2.25$ is the matrix multiplication exponent. In the large-constraint regime $m=\Omega(n^2)$, a robust IPM framework~\cite{hjstrz22} solves SDPs in $m^{\omega+o(1)}$ time. By contrast, a CPM-based SDP solver makes $\widetilde O(m)$ queries to a separation oracle for the dual feasible region and solves an SDP in $$\widetilde O(m(mn^2+n^\omega+m^2))$$ time~\cite{lsw15,jlsw20}, which is advantageous when $m$ is very small.

Notice that when $m\ll n$, existing IPM algorithms for SDP still require $\widetilde O(\sqrt n)$ iterations, whereas for linear programming this bottleneck has already been overcome. An LP with $n$ variables and $O(n)$ constraints can be solved in matrix multiplication time~\cite{cls21,b20,jswz21,song2026solving}, and a tall dense LP with $m$ dual variables and $n\gg m$ inequality constraints can be solved in nearly linear time~\cite{blss20,van2021minimum}. The key ingredient behind the latter results is the Lee--Sidford barrier~\cite{ls14,ls19}, which reduces the number of IPM iterations from $\widetilde O(\sqrt n)$ to $\widetilde O(\sqrt m)$ by lowering the barrier parameter.

Specifically, an IPM for the dual problem~\eqref{eq:intro_dual_sdp} chooses a smooth barrier function $F$ that tends to infinity at the boundary of the feasible region $\Omega:=\{y:S(y)\succ0\}$ and tracks the central path $y(t):=\arg\min_y\{tb^\top y+F(y)\}$ as $t\to\infty$~\cite{k84,r88,nn94}. In the framework of Nesterov and Nemirovski~\cite{nn94}, if $F$ is a \emph{$\nu$-self-concordant barrier}, which means that its second and third derivatives are suitably controlled in the local norm induced by $\nabla^2 F(y)$ (Definition~\ref{def:tech_ov_sc_barrier}), then the central path can be followed to high accuracy in $\widetilde O(\sqrt\nu)$ Newton steps~\cite{nn94,renegar2001mathematical}.

The classical logarithmic barrier for an LP with $m$ variables and $n$ inequality constraints has barrier parameter $\nu=n$. While every convex domain in $\R^m$ admits a universal or entropic barrier with optimal parameter $\nu=O(m)$~\cite{nn94,be19,ly18}, these barriers are defined through high-dimensional integrals that are expensive to evaluate. Lee and Sidford~\cite{ls14,ls19} bridged this gap for LPs using $\ell_p$-Lewis weights~\cite{lewis78,cp15} to construct an efficiently computable barrier with $\nu=\widetilde O(m)$. This enables any LP over $\R^m$ to be solved in $\widetilde O(\sqrt m)$ IPM iterations regardless of $n$. %

For a general SDP~\eqref{eq:intro_dual_sdp}, the analogous logarithmic barrier $-\log\det S(y)$ has barrier parameter $\nu=n$, which underlies the $\widetilde O(\sqrt n)$-iteration IPM solvers of~\cite{jklps20,hjstrz22}. The hybrid volumetric--logarithmic barrier~\cite{nn94,ans00} reduces the barrier parameter to $O(\sqrt{mn})$ when $m<n$, giving $\widetilde O((mn)^{1/4})$ iterations~\cite{ans00,hjstrz22}. An SDP analog of the Lee--Sidford barrier would reduce the iteration count all the way to $\widetilde O(\sqrt m)$. However, since each IPM iteration requires computing (or approximating) the gradient and Hessian of the barrier, reducing the barrier parameter helps only if these derivatives are cheap to evaluate, making the geometric question inseparable from an algorithmic one:
\begin{quote}
\emph{Can the Lee--Sidford barrier be extended to semidefinite programs, with a barrier parameter nearly linear in $m$, and with derivatives that can be evaluated efficiently?}
\end{quote}
Huang, Jiang, Song, Tao, and Zhang~\cite[Section~5]{hjstrz22} raised this question as a possible route to improving the state-of-the-art runtime of SDP solvers. In this paper, we answer it affirmatively.

\subsection{Our results}
\label{sec:intro_results}
This paper establishes two main results. On the geometric side, we construct a matrix Lee--Sidford barrier for SDPs whose self-concordance parameter is nearly linear in the dual dimension $m$ (the number of primal constraints). On the algorithmic side, we develop an IPM-based SDP algorithm that achieves the fastest known running time in the regime $m\leq n$, improving upon all previous high-accuracy algorithms.

\paragraph{Main result I:} Our matrix Lee--Sidford barrier function is defined as the minimum value of a strictly convex optimization problem over $m\times m$ positive-definite matrices, whose assumptions and barrier parameter bound are stated in the following theorem.  In Section~\ref{sec:technical_overview}, we present the geometric intuition behind this construction and sketch the proof of the theorem.

\begin{theorem}[Matrix Lee--Sidford barrier; informal version of Theorem~\ref{thm:ml_matrix_sc}]
\label{thm:geo_informal_barrier}
Assume that $\Omega=\{y:S(y)\succ0\}$ is nonempty, the coefficient matrices $A_1,\ldots,A_m$ are linearly independent, and $\bigcap_i\ker A_i=\{0\}$. The matrix Lee--Sidford barrier is a $\nu$-self-concordant barrier for $\Omega$ with
\begin{align*}
 \nu=O(m\log^2(2mn)).
\end{align*}
\end{theorem}
Theorem~\ref{thm:geo_informal_barrier} implies that, starting from a suitable central point, an IPM algorithm using this barrier can solve SDPs in $\widetilde O(\sqrt m)$ iterations. Note that the last two assumptions for the $A_i$'s hold without loss of generality (Section~\ref{sec:ml_assumptions}). We also remark that when $C,A_1,\ldots,A_m$ are diagonal, this barrier reduces to the Lee--Sidford barrier of the corresponding LP, up to scaling. For LPs,  Theorem~\ref{thm:geo_informal_barrier} thus improves the barrier parameter bound $O(m\log^5n)$ in~\cite[Theorem~3]{ls19} to $O(m\log^2(2mn))$.

\paragraph{Main result II:} Our SDP algorithm implements the matrix Lee--Sidford barrier in IPM iterations via randomized sketching and amortized maintenance. In addition to the assumptions in Theorem~\ref{thm:geo_informal_barrier}, we assume the algorithm is given numbers $B\geq1$ and $\rho>0$ such that there exist a dual optimum $y_*$ and a strictly feasible point $y_{\rm f}$ satisfying
\begin{align}\label{eq:intro_SDP_assumptions}
 \|y_*\|_\infty,\|y_{\rm f}\|_\infty\leq B,
 \qquad S(y_{\rm f})\succeq\rho I.
\end{align}
Bounds of this kind are standard for high-accuracy SDP solvers~\cite{jklps20,hjstrz22}. The point $y_{\rm f}$ is not supplied to the algorithm. Its existence ensures strict dual feasibility, so the primal and dual optimal values agree even when every primal feasible matrix is singular.

\begin{theorem}[Randomized SDP algorithm; informal version of Theorem~\ref{thm:more_complete_sdp}]\label{thm:intro_randomized}
Suppose that $1\leq m\leq n$ and the assumptions of Theorem~\ref{thm:geo_informal_barrier} and~\eqref{eq:intro_SDP_assumptions} hold, and let $\mathrm{OPT}$ be the common optimal value of~\eqref{eq:primal_sdp} and~\eqref{eq:intro_dual_sdp}. For every fixed $\zeta>0$, accuracy $0<\epsilon\leq1$, there is a randomized algorithm that, with high probability, outputs $y\in\R^m$ and $X\succeq0$ satisfying
\begin{equation}\label{eq:intro_output_accuracy}
 S(y)\succ0,\qquad
 0\leq b^\top y-\mathrm{OPT}\leq\epsilon,\qquad
 \sum_{i=1}^m|\langle A_i,X\rangle-b_i|\leq\epsilon,\qquad
 \langle C,X\rangle\geq\mathrm{OPT}-\epsilon
\end{equation}
in time
\begin{align*}
 \widetilde O(
 n^{\omega+\zeta}+n^2m^{3/2+\zeta}
 ).
\end{align*}
\end{theorem}
Here, $\widetilde O$ hides polylogarithmic factors in $n$, $1/\epsilon$, $B$, $1/\rho$, and the input magnitudes, and $\omega$ denotes the fast matrix multiplication exponent (currently, $\omega \leq 9/4$~\cite{o26b,p26,amanov2026multiplying,Shp26}). Section~\ref{sec:runtime_comparison} compares our runtime with that of previous high-accuracy SDP solvers.

In Section~\ref{sec:algorithmic_overview}, we first explain how to implement an IPM with the matrix Lee--Sidford barrier, in which the algorithm must solve two nested optimization problems: the SDP itself and the convex optimization that defines the barrier function at the current iterate. We then introduce the algorithmic techniques that accelerate each iteration to achieve the running time of Theorem~\ref{thm:intro_randomized}.

\subsection{Comparison of running times}\label{sec:runtime_comparison}
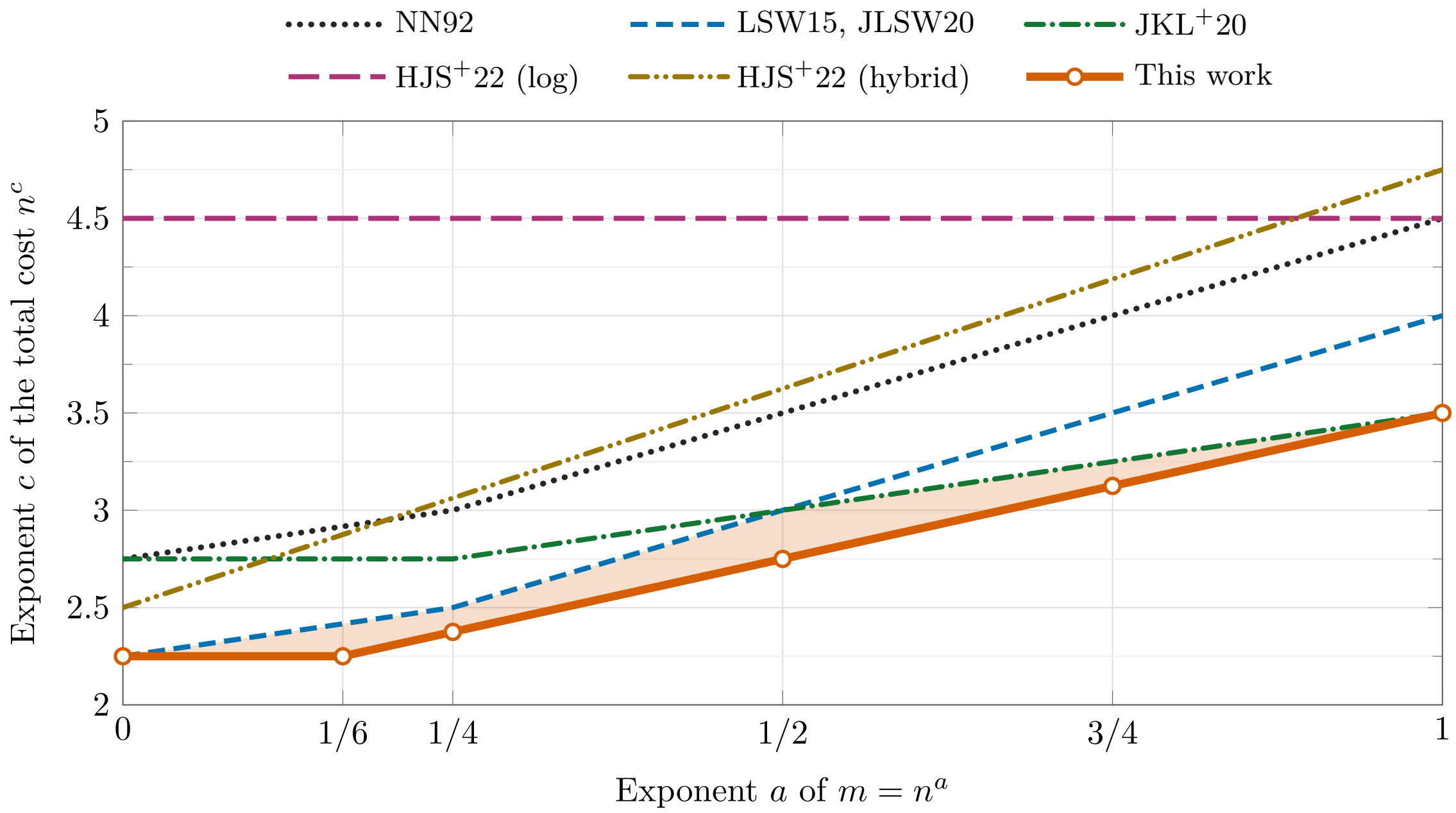
\begin{figure}[!t]
\centering
\definecolor{cmpBlack}{gray}{0.15}
\definecolor{cmpBlue}{HTML}{0072B2}
\definecolor{cmpGreen}{HTML}{117733}
\definecolor{cmpPurple}{HTML}{AA3377}
\definecolor{cmpGold}{HTML}{997700}
\definecolor{cmpOurs}{HTML}{D55E00}
\begin{tikzpicture}
\begin{axis}[
 width=0.92\linewidth,
 height=7.6cm,
 xmin=0, xmax=1,
 ymin=2, ymax=5,
 xtick={0,0.1666667,0.25,0.5,0.75,1},
 xticklabels={$0$,$1/6$,$1/4$,$1/2$,$3/4$,$1$},
 ytick={2,2.5,3,3.5,4,4.5,5},
 minor ytick={2.25,2.75,3.25,3.75,4.25,4.75},
 xlabel={Exponent $a$ of $m=n^a$},
 ylabel={Exponent $c$ of the total cost $n^c$},
 grid=both,
 major grid style={black!12},
 minor grid style={black!6},
 axis line style={black!60},
 tick style={black!60},
 tick label style={font=\small},
 label style={font=\small},
 legend columns=3,
 legend style={font=\footnotesize,at={(0.5,1.03)},anchor=south,
  draw=none,fill=none,inner sep=1pt,row sep=2pt,
  /tikz/every even column/.append style={column sep=12pt}},
 legend cell align=left,
 legend image code/.code={\draw[mark repeat=2,mark phase=2,#1]
  plot coordinates {(0cm,0cm) (0.5cm,0cm) (1cm,0cm)};},
]
\addplot[draw=none,fill=cmpOurs,fill opacity=0.2,forget plot]
 coordinates {(0,2.25) (0.25,2.5) (0.5,3) (1,3.5) (0.1666667,2.25)} -- cycle;
\addplot[cmpBlack,line width=1.7pt,line cap=round,dash pattern=on 0pt off 3.2pt]
 coordinates {(0,2.75) (0.25,3) (1,4.5)};
\addlegendentry{NN92}
\addplot[cmpBlue,line width=1.5pt,dash pattern=on 5pt off 2.5pt]
 coordinates {(0,2.25) (0.25,2.5) (1,4)};
\addlegendentry{LSW15, JLSW20}
\addplot[cmpGreen,line width=1.5pt,line cap=round,
 dash pattern=on 5pt off 2.6pt on 0pt off 2.6pt]
 coordinates {(0,2.75) (0.25,2.75) (1,3.5)};
\addlegendentry{JKL$^{+}$20}
\addplot[cmpPurple,line width=1.5pt,dash pattern=on 9pt off 3pt]
 coordinates {(0,4.5) (1,4.5)};
\addlegendentry{HJS$^{+}$22 (log)}
\addplot[cmpGold,line width=1.5pt,line cap=round,
 dash pattern=on 5pt off 2.6pt on 0pt off 2.6pt on 0pt off 2.6pt]
 coordinates {(0,2.5) (1,4.75)};
\addlegendentry{HJS$^{+}$22 (hybrid)}
\addplot[cmpOurs,line width=2.4pt,mark=*,mark size=2.2pt,
 mark options={solid,fill=white,line width=1pt}]
 coordinates {(0,2.25) (0.1666667,2.25) (0.25,2.375) (0.5,2.75) (0.75,3.125) (1,3.5)};
\addlegendentry{This work}
\end{axis}
\end{tikzpicture}
\caption{The total runtime of SDP solvers as exponents $c$ of $n^c$ for $m=n^a$, using $\omega=9/4$. The shaded region is the improvement over the best previous bound.}
\label{fig:sdp_runtime_exponents}
\end{figure}

We compare the total running time of high-accuracy SDP solvers in the regime $m\leq O(n)$, where the matrix Lee--Sidford barrier (Theorem~\ref{thm:geo_informal_barrier}) takes no more iterations than the logarithmic barrier (up to polylogarithmic factors). Table~\ref{tab:sdp_runtime_history} shows the iteration counts and total running times of previous algorithms and our result, with all bounds simplified under $m\leq n$ and $2<\omega\leq9/4$. Figure~\ref{fig:sdp_runtime_exponents} plots the base-$n$ exponent $c$ of the total running time $\widetilde O(n^c)$ as a function of $a\in[0,1]$ for $m=n^a$.

\begin{table}[!t]
\centering
\begin{tabular}{@{}llcc@{}}
\toprule
Year & Reference & Iterations & Total cost for $1\leq m\leq n$ \\
\midrule
1992 & Nesterov--Nemirovski~\cite{nn92} & $\sqrt n$ & $\sqrt n(m^2n^2+mn^\omega)$ \\
2015 & Lee--Sidford--Wong~\cite{lsw15} (CPM) & $m$ & $m(mn^2+n^\omega)$ \\
2020 & Jiang--Lee--Song--Wong~\cite{jlsw20} (CPM) & $m$ & $m(mn^2+n^\omega)$ \\
2020 & Jiang et al.~\cite{jklps20} & $\sqrt n$ & $\sqrt n(mn^2+n^\omega)$ \\
2022 & Huang et al.~\cite{hjstrz22} (log barrier) & $\sqrt n$ & $n^{9/2}$ \\
2022 & Huang et al.~\cite{hjstrz22} (hybrid barrier) & $(mn)^{1/4}$ & $(mn)^{1/4}m^2n^\omega$ \\
\midrule
This work & Theorem~\ref{thm:more_complete_sdp} & $\sqrt m$ & $n^\omega+n^2m^{3/2}$ \\
\bottomrule
\end{tabular}
\caption{Comparison of high-accuracy SDP solvers in the regime $m\leq n$, where $n$ is the matrix dimension and $m$ is the number of primal equality constraints (and the dual dimension).}
\label{tab:sdp_runtime_history}
\end{table}

As shown in Figure~\ref{fig:sdp_runtime_exponents}, the runtime exponent of our algorithm matches the best previous bounds at $m=\Theta(1)$ and $m=\Theta(n)$, and is strictly smaller for every $m^a$ with $0<a<1$. The largest improvement occurs at $m=n^{1/2}$, where the running time is reduced from $\widetilde O(n^3)$ to $\widetilde O(n^{11/4})$.

\subsection{Technical overview for the matrix Lee--Sidford barrier}
\label{sec:technical_overview}

The matrix Lee--Sidford barrier of Theorem~\ref{thm:geo_informal_barrier} has parameter $\nu=O(m\log^2(2mn))$, whereas the logarithmic barrier $-\log\det S(y)$ has $\nu=n$. This section explains the geometric motivation for this result and sketches its proof. After reviewing self-concordant barriers (Section~\ref{sec:geo_sc}), we recall how Lee and Sidford~\cite{ls19} build a barrier for linear programs from inscribed ellipsoids (Section~\ref{sec:geo_lp}), and we extend their construction to semidefinite constraints through a matrix sum-of-squares certificate (Section~\ref{sec:geo_sdp}). The main technical difficulty lies in the calculus of the barrier (Section~\ref{sec:geo_proof}). The scalar weights of the Lee--Sidford barrier can be differentiated explicitly, whereas the derivative of the matrix weight involves noncommuting Fr\'echet derivatives without a tractable closed form. We avoid this derivative by variational principles along exponential curves of matrices.

\subsubsection{Self-concordant barriers}
\label{sec:geo_sc}

An interior-point method (IPM) minimizes a linear objective $c^\top x$ over a convex feasible region $\mathcal X\subseteq\R^d$. It enforces feasibility by a smooth convex barrier function $F:\mathcal X\to\R$ that tends to $+\infty$ at the boundary $\partial\mathcal X$. When $F$ is strictly convex and the penalized objectives $f_t(x):=F(x)+tc^\top x$ attain their minima, the unique minimizers $x(t)$ form the \emph{central path}. For the matrix Lee--Sidford barrier, Lemma~\ref{lem:ml_central_minimizer} shows that $x(t)$ exists for every $t>0$ under primal and dual strict feasibility. A path-following IPM starts near $x(t_0)$ for some $t_0>0$, increases $t$ geometrically, and keeps its iterates close to $x(t)$. The self-concordance parameter $\nu$ of $F$ governs the number of steps needed to reach an $\epsilon$-accurate solution. %

\begin{definition}[Self-concordant barrier~\cite{nn94}]~\label{def:tech_ov_sc_barrier}
A convex $C^3$ function $F:\mathcal X\to\R$ with positive definite Hessian is a \emph{$\nu$-self-concordant barrier} if, for every $x\in\mathcal X$ and $h\in\R^d$,
\begin{enumerate}[label=(\alph*),leftmargin=*]
\item $F(x_k)\to+\infty$ whenever $x_k\in\mathcal X$ converges to a boundary point of $\mathcal X$.
\item $|\nabla^3F(x)[h,h,h]|\leq2\,(\nabla^2F(x)[h,h])^{3/2}$.
\item $|\nabla F(x)[h]|^2\leq\nu\,\nabla^2F(x)[h,h]$.
\end{enumerate}
\end{definition}

Condition~(a) keeps the central path inside $\mathcal X$. Condition~(b) bounds the change of the Hessian in the \emph{local norm} $\|h\|_x:=(h^\top\nabla^2F(x)h)^{1/2}$. Hence Newton's method for $f_t$ converges quadratically from every \emph{near-central} point, that is, from every point at which the Newton step is small in $\|\cdot\|_x$.

Condition~(c) bounds the gradient in the dual local norm $\|v\|_{x,*}:=(v^\top\nabla^2F(x)^{-1}v)^{1/2}$, which limits how far $x(t)$ moves in $\|\cdot\|_x$ when $t$ increases. Each outer iteration can thus multiply $t$ by $1+\Theta(\nu^{-1/2})$ and restore near-centrality in $O(1)$ Newton steps. Since $x(t)$ has suboptimality gap at most $\nu/t$, a near-central start at $t=t_0$ yields an $\epsilon$-accurate solution in $O(\sqrt\nu\log(2+\nu/(t_0\epsilon)))$ Newton steps.

The parameter $\nu$ also has a geometric characterization that guides our construction. The unit balls of the local norm, called \emph{Dikin ellipsoids}, capture the local geometry seen by Newton's method. If the barrier attains its minimum at $x_F=\arg\min_x F(x)$, the Dikin ellipsoid $E_{x_F}:=\{h:\|h\|_{x_F}\leq1\}$ rounds the feasible region ${\cal X}$ up to a factor controlled by $\nu$~\cite{nn94}, in the sense that
\begin{align}\label{eq:tech_ov_barrier_rounding}
 x_F+E_{x_F}\subseteq\overline{\mathcal X}\subseteq x_F+(\nu+2\sqrt\nu)\,E_{x_F}.
\end{align}
In particular, a barrier with $\nu=\widetilde O(m)$ must induce a local metric whose unit ellipsoid approximates the feasible region within a factor that depends mainly on the dimension $m$ of the domain.

\subsubsection{The Lee--Sidford barrier through inscribed ellipsoids}
\label{sec:geo_lp}

Consider a linear program (LP) in dual form with the $n$ constraints
\begin{align*}
    a_k^\top y-b_k> 0,\quad a_k\in \R^m,~ b_k\in \R\quad \forall k\in [n],\quad y\in\R^m.
\end{align*}
This LP is an SDP whose constraint matrices are diagonal. The logarithmic barrier $-\sum_{k=1}^n\log(a_k^\top y-b_k)$ has one term for each inequality, and its parameter is $\nu=n$. Duplicating a constraint increases $\nu$, although the feasible polytope stays the same. To make the parameter depend on $m$ instead of $n$, Lee and Sidford~\cite{ls19} build a barrier from the maximum-volume inscribed ellipsoid. This ellipsoid depends only on the feasible region, and by John's theorem it approximates every symmetric convex body in $\R^m$ within a factor of $\sqrt m$, which is consistent with the rounding requirement~\eqref{eq:tech_ov_barrier_rounding} for $\nu=\widetilde O(m)$.

More specifically, we normalize each constraint by its slack at $y$ and define the symmetrized local feasible region
\begin{align*}
 K_y:=\left\{h\in\R^m:|r_k^\top h|\leq1~\forall k\in [n]\right\},\qquad \text{where}\quad  r_k=r_k(y):=\frac{a_k}{a_k^\top y-b_k}\in \R^m.
\end{align*}
A step $h$ lies in $K_y$ if and only if both $y+h$ and $y-h$ satisfy all nonstrict inequalities. For a positive definite matrix $G\in \R^{m\times m}$, the centered ellipsoid $\mathcal E_G:=\{G^{1/2}z:\|z\|_2\leq1\}$ has volume proportional to $\sqrt{\det G}$, and its squared width along $r_k$ is $\max_{h\in\mathcal E_G}(r_k^\top h)^2=r_k^\top Gr_k$. The maximum-volume centered ellipsoid contained in $K_y$ is therefore the solution of
\begin{align*}
 \max_{G\succ0}\ \log\det G\quad\text{subject to}\quad r_k^\top Gr_k\leq1\quad\forall k\in [n].
\end{align*}
Because this optimization problem depends only on $K_y$, redundant or duplicated constraints do not affect it.

The maximum-volume ellipsoid need not vary smoothly with $y$, because the set of active constraints can change abruptly. For instance, on the interval $(-1,1)\subseteq\R$, the largest inscribed interval centered at $y$ has radius $1-|y|$, and its negative log-volume has a kink at $y=0$. Lee and Sidford~\cite{ls19} therefore replace the hard width constraints $r_k^\top Gr_k\leq1$ by an $\ell_s$-penalty and define the potential
\begin{equation}\label{eq:geo_lp_potential}
 \Phi^{\rm LP}(y):=\frac12\min_{G\succ0}\left\{\frac1s\sum_{k\in [n]}(r_k^\top Gr_k)^s-\log\det G-m\right\}.
\end{equation}
Here $-\log\det G$ rewards a large ellipsoid, while the power sum penalizes violations of the width constraints $r_k^\top Gr_k\leq1$. At $s=1$, the potential~\eqref{eq:geo_lp_potential} recovers Vaidya's volumetric barrier~\cite{vai96}, and as $s\to\infty$, its minimizer converges to the John ellipsoid of $K_y$. Since $r_k^\top G r_k\geq 0$, the power sum satisfies
\begin{align*}
    \max_{k\in [n]}r_k^\top G r_k\leq \Bigl(\sum_{k\in [n]}(r_k^\top G r_k)^s\Bigr)^{1/s}\leq n^{1/s}\cdot \max_{k\in [n]}r_k^\top G r_k.
\end{align*}
The choice $s=\Theta(\log(2n))$ thus approximates the maximum width within a constant factor and keeps the higher derivatives bounded.

The optimality condition of~\eqref{eq:geo_lp_potential} assigns a weight to each constraint. Let $G_*$ be the minimizer in~\eqref{eq:geo_lp_potential}, and let $w_k:=(r_k^\top G_*r_k)^s$. Setting the derivative with respect to $G$ to zero, and then taking the trace inner product with $G_*$, gives
\begin{align*}
 G_*^{-1}=\sum_{k\in [n]}w_k^{1-1/s}r_kr_k^\top,\qquad\sum_{k\in [n]}w_k=\tr\left[G_*\sum_{k\in [n]}w_k^{1-1/s}r_kr_k^\top\right]=m.
\end{align*}
The scalars $w_k$ are the $\ell_{2s}$ Lewis weights of the matrix with rows $r_k^\top$~\cite{lewis78,cp15}. Indeed, by the first identity, $w_k$ is the leverage score of row $k$ after every row $r_j^\top$ is rescaled by $w_j^{(s-1)/(2s)}$. Geometrically, $w_k$ measures how strongly the $k$th constraint restricts the ellipsoid. Identical constraints receive equal weights, and the total weight equals the dimension $m$. A constraint with $r_k^\top G_*r_k\leq1/2$ has weight $w_k\leq2^{-s}$, so loose constraints receive almost no weight. Lee and Sidford originally formulated their barrier as a maximization over these weights,
\begin{align*}
 \Phi^{\rm LP}(y)=\frac12\max_{w\geq0}\Bigg\{\log\det\left(\sum_kw_k^{1-1/s}r_kr_k^\top\right)-\left(1-\frac1s\right)\sum_kw_k\Bigg\}.
\end{align*}
To prove that the self-concordance parameter is $\widetilde O(m)$, we bound the derivatives of $\Phi^{\rm LP}$. By the envelope theorem at the stationary point $G_*$, and since $\nabla r_k[h]=-(r_k^\top h)r_k$, the gradient of $\Phi^{\rm LP}$ equals that of the weighted logarithmic barrier $L_w(y):=-\sum_k w_k\log(a_k^\top y-b_k)$ with the weights $w$ held fixed,
\begin{align*}
 \nabla\Phi^{\rm LP}(y)[h]=-\sum_{k\in [n]}w_kr_k^\top h=\nabla L_w(y)[h].
\end{align*}
The Hessians differ by a term that involves the motion of the weights,
\begin{align*}
    \nabla^2\Phi^{\rm LP}(y)[h,h]=Q_y^{\rm LP}[h,h]-\sum_{k}\nabla w_k(y)[h](r_k^\top h),
\end{align*}
where $Q_y^{\rm LP}:=\nabla^2L_w(y)=\sum_k w_k r_kr_k^\top$ serves as a proxy Hessian for $\nabla^2\Phi^{\rm LP}$. The Cauchy--Schwarz inequality gives the gradient bound in the $Q_y^{\rm LP}$-norm,
\begin{align*}
 \left|\nabla\Phi^{\rm LP}(y)[h]\right|^2=\Bigg(\sum_{k\in [n]}w_k r_k^\top h\Bigg)^2\leq\Bigg(\sum_{k\in [n]}w_k\Bigg)\Bigg(\sum_{k\in [n]}w_k(r_k^\top h)^2\Bigg)=m\,Q_y^{\rm LP}[h,h].
\end{align*}
The identity $\sum_kw_k=m$ is the source of the factor $m$ in $\nu$. It remains to prove the spectral comparison $Q_y^{\rm LP}\preceq\nabla^2\Phi^{\rm LP}$ and to bound $\nabla^3\Phi^{\rm LP}$. Lee and Sidford~\cite{ls19} prove both by differentiating the weight map $y\mapsto w(y)$ explicitly. Since the constraints commute, the Jacobian of $w(y)$ reduces to entrywise (Hadamard) products of projection matrices, and direct matrix algebra then gives $Q_y^{\rm LP}\preceq\nabla^2\Phi^{\rm LP}\preceq(2s+1)Q_y^{\rm LP}$ together with the third-derivative bound.

\begin{figure}[!t]
\centering
\begin{tikzpicture}[x=1.36cm,y=1.36cm,font=\scriptsize]
 \fill[black!7] (1.0000,0.3333) -- (0.5333,0.8000) -- (-1.0000,0.8000) -- (-1.0000,-0.3333) -- (-0.5333,-0.8000) -- (1.0000,-0.8000) -- cycle;
\draw[black!75,line width=.45pt] (0.965,-1.300) -- (0.965,1.300);
\draw[black!75,line width=.45pt] (1.000,-1.300) -- (1.000,1.300);
\draw[black!75,line width=.45pt] (1.035,-1.300) -- (1.035,1.300);
\draw[black!75,line width=.6pt] (-1.750,0.800) -- (1.750,0.800);
\draw[black!75,line width=.6pt] (0.033,1.300) -- (1.750,-0.417);
\draw[black!75,line width=.6pt,densely dotted] (0.778,-1.300) -- (1.750,-0.425);
\draw[black!75,line width=.45pt] (-1.035,-1.300) -- (-1.035,1.300);
\draw[black!75,line width=.45pt] (-1.000,-1.300) -- (-1.000,1.300);
\draw[black!75,line width=.45pt] (-0.965,-1.300) -- (-0.965,1.300);
\draw[black!75,line width=.6pt] (-1.750,-0.800) -- (1.750,-0.800);
\draw[black!75,line width=.6pt] (-1.750,0.417) -- (-0.033,-1.300);
\draw[black!75,line width=.6pt,densely dotted] (-1.750,0.425) -- (-0.778,1.300);
 \fill (0,0) circle (0.035);
 \node[anchor=west] at (0,0) {$y$};
 \node[anchor=south] at (1.0,1.3) {$1+1+1$};
 \node[anchor=south] at (-0.2,0.8) {$1$};
 \node[anchor=west] at (1.42,0.06) {$1$};
 \node[anchor=north west] at (1.25,-0.98) {$1$};
 \node[anchor=north,font=\small] at (0,-1.38) {(a) logarithmic barrier: $\nu=6$};
\end{tikzpicture}\hfill
\begin{tikzpicture}[x=1.36cm,y=1.36cm,font=\scriptsize]
 \fill[black!7] (1.0000,0.3333) -- (0.5333,0.8000) -- (-1.0000,0.8000) -- (-1.0000,-0.3333) -- (-0.5333,-0.8000) -- (1.0000,-0.8000) -- cycle;
\draw[black!35,line width=.45pt] (0.965,-1.300) -- (0.965,1.300);
\draw[black!35,line width=.45pt] (1.000,-1.300) -- (1.000,1.300);
\draw[black!35,line width=.45pt] (1.035,-1.300) -- (1.035,1.300);
\draw[black!35,line width=.6pt] (-1.750,0.800) -- (1.750,0.800);
\draw[black!35,line width=.6pt] (0.033,1.300) -- (1.750,-0.417);
\draw[black!35,line width=.6pt,densely dotted] (0.778,-1.300) -- (1.750,-0.425);
\draw[black!35,line width=.45pt] (-1.035,-1.300) -- (-1.035,1.300);
\draw[black!35,line width=.45pt] (-1.000,-1.300) -- (-1.000,1.300);
\draw[black!35,line width=.45pt] (-0.965,-1.300) -- (-0.965,1.300);
\draw[black!35,line width=.6pt] (-1.750,-0.800) -- (1.750,-0.800);
\draw[black!35,line width=.6pt] (-1.750,0.417) -- (-0.033,-1.300);
\draw[black!35,line width=.6pt,densely dotted] (-1.750,0.425) -- (-0.778,1.300);
 \filldraw[fill=overviewteal!22,draw=overviewteal!90!black,line width=.8pt] plot coordinates {(-0.1673,-0.7370) (-0.2127,-0.7256) (-0.2576,-0.7123) (-0.3018,-0.6969) (-0.3452,-0.6797) (-0.3876,-0.6606) (-0.4290,-0.6397) (-0.4692,-0.6170) (-0.5081,-0.5927) (-0.5456,-0.5667) (-0.5816,-0.5392) (-0.6160,-0.5102) (-0.6487,-0.4798) (-0.6797,-0.4480) (-0.7088,-0.4151) (-0.7359,-0.3810) (-0.7610,-0.3459) (-0.7841,-0.3098) (-0.8050,-0.2729) (-0.8237,-0.2352) (-0.8401,-0.1968) (-0.8542,-0.1580) (-0.8660,-0.1187) (-0.8754,-0.0791) (-0.8824,-0.0392) (-0.8870,0.0007) (-0.8892,0.0407) (-0.8889,0.0805) (-0.8862,0.1201) (-0.8810,0.1594) (-0.8735,0.1982) (-0.8635,0.2365) (-0.8512,0.2742) (-0.8366,0.3111) (-0.8196,0.3471) (-0.8004,0.3822) (-0.7790,0.4163) (-0.7555,0.4492) (-0.7299,0.4809) (-0.7023,0.5112) (-0.6728,0.5402) (-0.6415,0.5677) (-0.6084,0.5936) (-0.5736,0.6179) (-0.5372,0.6405) (-0.4994,0.6613) (-0.4602,0.6804) (-0.4197,0.6975) (-0.3781,0.7128) (-0.3355,0.7261) (-0.2919,0.7374) (-0.2475,0.7467) (-0.2025,0.7539) (-0.1569,0.7591) (-0.1109,0.7622) (-0.0645,0.7632) (-0.0180,0.7621) (0.0286,0.7590) (0.0750,0.7537) (0.1213,0.7464) (0.1673,0.7370) (0.2127,0.7256) (0.2576,0.7123) (0.3018,0.6969) (0.3452,0.6797) (0.3876,0.6606) (0.4290,0.6397) (0.4692,0.6170) (0.5081,0.5927) (0.5456,0.5667) (0.5816,0.5392) (0.6160,0.5102) (0.6487,0.4798) (0.6797,0.4480) (0.7088,0.4151) (0.7359,0.3810) (0.7610,0.3459) (0.7841,0.3098) (0.8050,0.2729) (0.8237,0.2352) (0.8401,0.1968) (0.8542,0.1580) (0.8660,0.1187) (0.8754,0.0791) (0.8824,0.0392) (0.8870,-0.0007) (0.8892,-0.0407) (0.8889,-0.0805) (0.8862,-0.1201) (0.8810,-0.1594) (0.8735,-0.1982) (0.8635,-0.2365) (0.8512,-0.2742) (0.8366,-0.3111) (0.8196,-0.3471) (0.8004,-0.3822) (0.7790,-0.4163) (0.7555,-0.4492) (0.7299,-0.4809) (0.7023,-0.5112) (0.6728,-0.5402) (0.6415,-0.5677) (0.6084,-0.5936) (0.5736,-0.6179) (0.5372,-0.6405) (0.4994,-0.6613) (0.4602,-0.6804) (0.4197,-0.6975) (0.3781,-0.7128) (0.3355,-0.7261) (0.2919,-0.7374) (0.2475,-0.7467) (0.2025,-0.7539) (0.1569,-0.7591) (0.1109,-0.7622) (0.0645,-0.7632) (0.0180,-0.7621) (-0.0286,-0.7590) (-0.0750,-0.7537) (-0.1213,-0.7464) (-0.1673,-0.7370)} -- cycle;
 \fill (0,0) circle (0.035);
 \node[anchor=west] at (0,0) {$y$};
 \node[anchor=south,text=overviewteal!80!black] at (1.0,1.3) {$3\times 0.31$};
 \node[anchor=south,text=overviewteal!80!black] at (-0.2,0.8) {$0.94$};
 \node[anchor=west,text=overviewteal!80!black] at (1.42,0.06) {$0.14$};
 \node[anchor=north west,text=overviewteal!80!black] at (1.25,-0.98) {$0.00$};
 \node[anchor=north,font=\small] at (0,-1.38) {(b) Lee--Sidford: $\sum_kw_k=m=2$};
\end{tikzpicture}\hfill
\begin{tikzpicture}[x=1.36cm,y=1.36cm,font=\scriptsize]
\draw[black!22,line width=.35pt] (0.706,1.550) -- (1.296,-1.550);
\draw[black!22,line width=.35pt] (-1.296,1.550) -- (-0.706,-1.550);
\draw[black!22,line width=.35pt] (0.527,1.550) -- (1.523,-1.550);
\draw[black!22,line width=.35pt] (-1.523,1.550) -- (-0.527,-1.550);
\draw[black!22,line width=.35pt] (0.346,1.550) -- (1.750,-1.394);
\draw[black!22,line width=.35pt] (-1.750,1.394) -- (-0.346,-1.550);
\draw[black!22,line width=.35pt] (0.132,1.550) -- (1.750,-0.758);
\draw[black!22,line width=.35pt] (-1.750,0.758) -- (-0.132,-1.550);
\draw[black!22,line width=.35pt] (-0.158,1.550) -- (1.750,-0.222);
\draw[black!22,line width=.35pt] (-1.750,0.222) -- (0.158,-1.550);
\draw[black!22,line width=.35pt] (-0.566,1.550) -- (1.750,0.220);
\draw[black!22,line width=.35pt] (-1.750,-0.220) -- (0.566,-1.550);
\draw[black!22,line width=.35pt] (-1.154,1.550) -- (1.750,0.559);
\draw[black!22,line width=.35pt] (-1.750,-0.559) -- (1.154,-1.550);
\draw[black!22,line width=.35pt] (-1.750,1.472) -- (1.750,0.824);
\draw[black!22,line width=.35pt] (-1.750,-0.824) -- (1.750,-1.472);
\draw[black!22,line width=.35pt] (-1.750,1.289) -- (1.750,1.058);
\draw[black!22,line width=.35pt] (-1.750,-1.058) -- (1.750,-1.289);
\draw[black!22,line width=.35pt] (0.046,-1.550) -- (1.750,-0.273);
\draw[black!22,line width=.35pt] (-1.750,0.273) -- (-0.046,1.550);
\draw[black!22,line width=.35pt] (0.380,-1.550) -- (1.750,0.405);
\draw[black!22,line width=.35pt] (-1.750,-0.405) -- (-0.380,1.550);
\draw[black!22,line width=.35pt] (0.654,-1.550) -- (1.702,1.550);
\draw[black!22,line width=.35pt] (-1.702,-1.550) -- (-0.654,1.550);
\draw[black!22,line width=.35pt] (0.887,-1.550) -- (1.221,1.550);
\draw[black!22,line width=.35pt] (-1.221,-1.550) -- (-0.887,1.550);
\draw[black!22,line width=.35pt] (0.920,1.550) -- (1.095,-1.550);
\draw[black!22,line width=.35pt] (-1.095,1.550) -- (-0.920,-1.550);
 \filldraw[fill=overviewblue!9,draw=overviewblue,line width=.8pt] plot coordinates {(1.0000,0.0000) (0.9959,0.0261) (0.9914,0.0520) (0.9866,0.0776) (0.9814,0.1032) (0.9759,0.1285) (0.9701,0.1536) (0.9639,0.1786) (0.9573,0.2035) (0.9505,0.2282) (0.9432,0.2527) (0.9357,0.2772) (0.9278,0.3015) (0.9195,0.3256) (0.9109,0.3497) (0.9019,0.3736) (0.8925,0.3974) (0.8828,0.4211) (0.8727,0.4446) (0.8621,0.4681) (0.8511,0.4914) (0.8398,0.5146) (0.8279,0.5377) (0.8156,0.5606) (0.8029,0.5833) (0.7896,0.6059) (0.7759,0.6283) (0.7616,0.6505) (0.7468,0.6725) (0.7315,0.6942) (0.7157,0.7157) (0.6992,0.7368) (0.6822,0.7577) (0.6647,0.7782) (0.6465,0.7984) (0.6278,0.8182) (0.6085,0.8375) (0.5886,0.8565) (0.5682,0.8749) (0.5472,0.8929) (0.5256,0.9103) (0.5035,0.9273) (0.4808,0.9437) (0.4577,0.9595) (0.4340,0.9748) (0.4099,0.9895) (0.3853,1.0037) (0.3602,1.0173) (0.3348,1.0303) (0.3089,1.0428) (0.2826,1.0547) (0.2559,1.0661) (0.2289,1.0769) (0.2015,1.0873) (0.1738,1.0971) (0.1457,1.1065) (0.1172,1.1153) (0.0884,1.1237) (0.0593,1.1317) (0.0298,1.1392) (0.0000,1.1462) (-0.0302,1.1529) (-0.0607,1.1591) (-0.0917,1.1649) (-0.1230,1.1703) (-0.1547,1.1753) (-0.1869,1.1799) (-0.2195,1.1841) (-0.2525,1.1879) (-0.2860,1.1913) (-0.3200,1.1943) (-0.3546,1.1970) (-0.3896,1.1991) (-0.4253,1.2009) (-0.4615,1.2023) (-0.4955,1.1961) (-0.5254,1.1800) (-0.5547,1.1630) (-0.5834,1.1449) (-0.6113,1.1259) (-0.6386,1.1060) (-0.6650,1.0852) (-0.6906,1.0634) (-0.7153,1.0408) (-0.7391,1.0173) (-0.7620,0.9930) (-0.7838,0.9679) (-0.8047,0.9422) (-0.8245,0.9157) (-0.8433,0.8887) (-0.8611,0.8611) (-0.8777,0.8329) (-0.8934,0.8044) (-0.9079,0.7755) (-0.9215,0.7462) (-0.9340,0.7167) (-0.9455,0.6870) (-0.9561,0.6571) (-0.9657,0.6271) (-0.9744,0.5971) (-0.9822,0.5671) (-0.9891,0.5370) (-0.9952,0.5071) (-1.0006,0.4772) (-1.0051,0.4475) (-1.0090,0.4179) (-1.0122,0.3885) (-1.0147,0.3593) (-1.0166,0.3303) (-1.0179,0.3015) (-1.0187,0.2730) (-1.0189,0.2446) (-1.0186,0.2165) (-1.0178,0.1886) (-1.0165,0.1610) (-1.0148,0.1336) (-1.0126,0.1064) (-1.0101,0.0795) (-1.0071,0.0528) (-1.0037,0.0263) (-1.0000,0.0000) (-0.9959,-0.0261) (-0.9914,-0.0520) (-0.9866,-0.0776) (-0.9814,-0.1032) (-0.9759,-0.1285) (-0.9701,-0.1536) (-0.9639,-0.1786) (-0.9573,-0.2035) (-0.9505,-0.2282) (-0.9432,-0.2527) (-0.9357,-0.2772) (-0.9278,-0.3015) (-0.9195,-0.3256) (-0.9109,-0.3497) (-0.9019,-0.3736) (-0.8925,-0.3974) (-0.8828,-0.4211) (-0.8727,-0.4446) (-0.8621,-0.4681) (-0.8511,-0.4914) (-0.8398,-0.5146) (-0.8279,-0.5377) (-0.8156,-0.5606) (-0.8029,-0.5833) (-0.7896,-0.6059) (-0.7759,-0.6283) (-0.7616,-0.6505) (-0.7468,-0.6725) (-0.7315,-0.6942) (-0.7157,-0.7157) (-0.6992,-0.7368) (-0.6822,-0.7577) (-0.6647,-0.7782) (-0.6465,-0.7984) (-0.6278,-0.8182) (-0.6085,-0.8375) (-0.5886,-0.8565) (-0.5682,-0.8749) (-0.5472,-0.8929) (-0.5256,-0.9103) (-0.5035,-0.9273) (-0.4808,-0.9437) (-0.4577,-0.9595) (-0.4340,-0.9748) (-0.4099,-0.9895) (-0.3853,-1.0037) (-0.3602,-1.0173) (-0.3348,-1.0303) (-0.3089,-1.0428) (-0.2826,-1.0547) (-0.2559,-1.0661) (-0.2289,-1.0769) (-0.2015,-1.0873) (-0.1738,-1.0971) (-0.1457,-1.1065) (-0.1172,-1.1153) (-0.0884,-1.1237) (-0.0593,-1.1317) (-0.0298,-1.1392) (-0.0000,-1.1462) (0.0302,-1.1529) (0.0607,-1.1591) (0.0917,-1.1649) (0.1230,-1.1703) (0.1547,-1.1753) (0.1869,-1.1799) (0.2195,-1.1841) (0.2525,-1.1879) (0.2860,-1.1913) (0.3200,-1.1943) (0.3546,-1.1970) (0.3896,-1.1991) (0.4253,-1.2009) (0.4615,-1.2023) (0.4955,-1.1961) (0.5254,-1.1800) (0.5547,-1.1630) (0.5834,-1.1449) (0.6113,-1.1259) (0.6386,-1.1060) (0.6650,-1.0852) (0.6906,-1.0634) (0.7153,-1.0408) (0.7391,-1.0173) (0.7620,-0.9930) (0.7838,-0.9679) (0.8047,-0.9422) (0.8245,-0.9157) (0.8433,-0.8887) (0.8611,-0.8611) (0.8777,-0.8329) (0.8934,-0.8044) (0.9079,-0.7755) (0.9215,-0.7462) (0.9340,-0.7167) (0.9455,-0.6870) (0.9561,-0.6571) (0.9657,-0.6271) (0.9744,-0.5971) (0.9822,-0.5671) (0.9891,-0.5370) (0.9952,-0.5071) (1.0006,-0.4772) (1.0051,-0.4475) (1.0090,-0.4179) (1.0122,-0.3885) (1.0147,-0.3593) (1.0166,-0.3303) (1.0179,-0.3015) (1.0187,-0.2730) (1.0189,-0.2446) (1.0186,-0.2165) (1.0178,-0.1886) (1.0165,-0.1610) (1.0148,-0.1336) (1.0126,-0.1064) (1.0101,-0.0795) (1.0071,-0.0528) (1.0037,-0.0263)} -- cycle;
 \filldraw[fill=overviewteal!22,draw=overviewteal!90!black,line width=.8pt] plot coordinates {(-0.6585,-0.3031) (-0.6783,-0.2579) (-0.6961,-0.2120) (-0.7121,-0.1654) (-0.7261,-0.1185) (-0.7381,-0.0712) (-0.7481,-0.0237) (-0.7560,0.0238) (-0.7619,0.0713) (-0.7656,0.1186) (-0.7673,0.1656) (-0.7669,0.2121) (-0.7644,0.2580) (-0.7598,0.3032) (-0.7530,0.3476) (-0.7443,0.3910) (-0.7335,0.4334) (-0.7207,0.4746) (-0.7059,0.5144) (-0.6891,0.5529) (-0.6705,0.5899) (-0.6501,0.6252) (-0.6278,0.6588) (-0.6039,0.6906) (-0.5782,0.7205) (-0.5511,0.7485) (-0.5223,0.7744) (-0.4922,0.7982) (-0.4607,0.8197) (-0.4280,0.8391) (-0.3940,0.8561) (-0.3590,0.8708) (-0.3230,0.8831) (-0.2862,0.8930) (-0.2485,0.9004) (-0.2102,0.9054) (-0.1713,0.9079) (-0.1319,0.9079) (-0.0921,0.9054) (-0.0521,0.9004) (-0.0120,0.8930) (0.0282,0.8831) (0.0683,0.8708) (0.1082,0.8561) (0.1478,0.8390) (0.1870,0.8197) (0.2257,0.7981) (0.2638,0.7743) (0.3012,0.7484) (0.3377,0.7205) (0.3733,0.6905) (0.4079,0.6587) (0.4413,0.6251) (0.4736,0.5898) (0.5045,0.5528) (0.5341,0.5143) (0.5622,0.4745) (0.5888,0.4333) (0.6137,0.3909) (0.6370,0.3475) (0.6585,0.3031) (0.6783,0.2579) (0.6961,0.2120) (0.7121,0.1654) (0.7261,0.1185) (0.7381,0.0712) (0.7481,0.0237) (0.7560,-0.0238) (0.7619,-0.0713) (0.7656,-0.1186) (0.7673,-0.1656) (0.7669,-0.2121) (0.7644,-0.2580) (0.7598,-0.3032) (0.7530,-0.3476) (0.7443,-0.3910) (0.7335,-0.4334) (0.7207,-0.4746) (0.7059,-0.5144) (0.6891,-0.5529) (0.6705,-0.5899) (0.6501,-0.6252) (0.6278,-0.6588) (0.6039,-0.6906) (0.5782,-0.7205) (0.5511,-0.7485) (0.5223,-0.7744) (0.4922,-0.7982) (0.4607,-0.8197) (0.4280,-0.8391) (0.3940,-0.8561) (0.3590,-0.8708) (0.3230,-0.8831) (0.2862,-0.8930) (0.2485,-0.9004) (0.2102,-0.9054) (0.1713,-0.9079) (0.1319,-0.9079) (0.0921,-0.9054) (0.0521,-0.9004) (0.0120,-0.8930) (-0.0282,-0.8831) (-0.0683,-0.8708) (-0.1082,-0.8561) (-0.1478,-0.8390) (-0.1870,-0.8197) (-0.2257,-0.7981) (-0.2638,-0.7743) (-0.3012,-0.7484) (-0.3377,-0.7205) (-0.3733,-0.6905) (-0.4079,-0.6587) (-0.4413,-0.6251) (-0.4736,-0.5898) (-0.5045,-0.5528) (-0.5341,-0.5143) (-0.5622,-0.4745) (-0.5888,-0.4333) (-0.6137,-0.3909) (-0.6370,-0.3475) (-0.6585,-0.3031)} -- cycle;
 \fill (0,0) circle (0.035);
 \node[anchor=west] at (0,0) {$y$};
 \node[text=overviewblue] at (-0.32,1.06) {$K_y$};
 \node[anchor=north,font=\small] at (0,-1.6) {(c) matrix Lee--Sidford: $\|B_h\|_{\rm op}\le1$};
\end{tikzpicture}
\caption{Barrier weights on a local region. (a)~The symmetrized local region $K_y$ (gray) of an LP with $n=6$ inequalities and $m=2$ variables. The six inequalities come from four distinct constraints, one of which appears three times (triple line) and one of which is redundant (dotted). The logarithmic barrier gives every inequality weight $1$. (b)~The ellipsoid of the optimizer of~\eqref{eq:geo_lp_potential} with $s=8$, scaled by $m^{-1/(2s)}$ so that it lies in $K_y$. The numbers are the Lewis weights rounded to two decimal places. The redundant constraint has a positive weight that rounds to $0.00$, and the unrounded weights sum to $m=2$. (c)~For an SDP with $m=2$ and $n=3$, the region $K_y$ (blue) is the intersection of the slabs $|r(u)^\top h|\leq1$ over all unit vectors $u$. The certified ellipse $m^{-1/(2s)}\mathcal E_{G_*}$ lies inside $K_y$.}
\label{fig:geo_ellipsoids}
\end{figure}
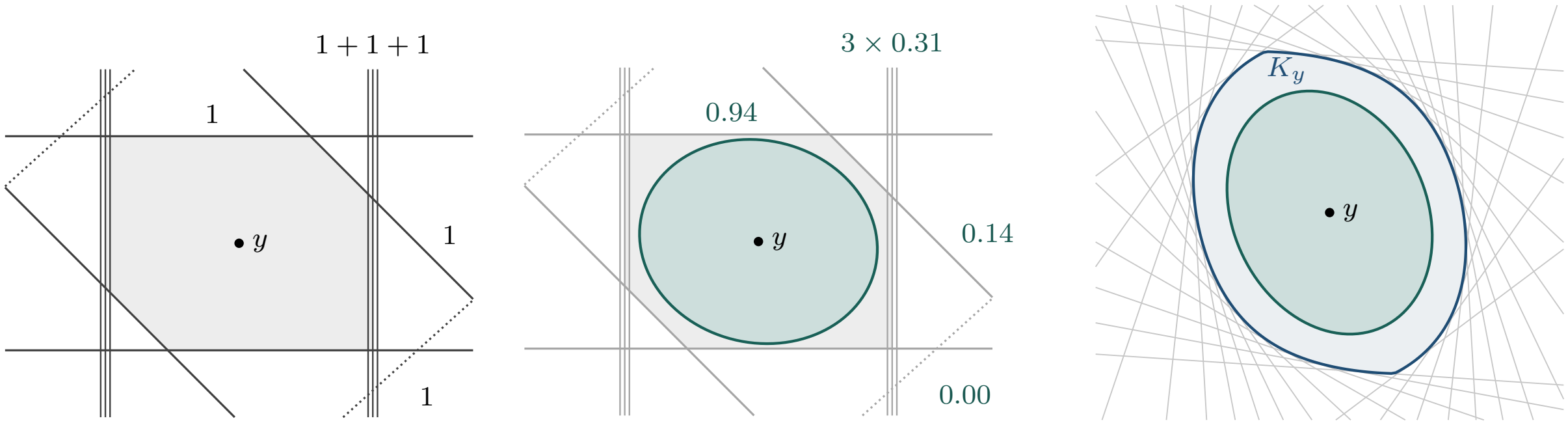

\subsubsection{From linear to semidefinite constraints}
\label{sec:geo_sdp}
An SDP can be viewed as an LP with infinitely many constraints, since $S(y)\succ0$ holds if and only if $u^\top S(y)u>0$ for every unit vector $u\in\R^n$. We normalize by the current slack $S(y)$ and define the normalized constraint matrices
\begin{align*}
 B_i:=S(y)^{-1/2}A_iS(y)^{-1/2}\in \R^{n\times n}~~~\forall i\in [m],\qquad B_h:=\sum_{i\in [m]}h_iB_i.
\end{align*}
The symmetrized local region around $y$ is then
\begin{align*}
 K_y:=\Bigl\{h\in\R^m:\|B_h\|_{\rm op}\leq1\Bigr\}=\Bigl\{h\in\R^m:|r(u)^\top h|\leq1\text{ for every unit vector }u\Bigr\},
\end{align*}
where $r(u)\in \R^m$ has the entries $r(u)_i:=u^\top B_iu$ for $i\in [m]$.

Two obstacles prevent a direct extension of the Lee--Sidford construction. First, the constraints form an infinite family. When the $B_i$ do not commute, the perturbations $B_h$ have different eigenbases, and no fixed list of $n$ scalar inequalities describes all of them. Second, the exact containment $\mathcal E_G\subseteq K_y$ requires $\sup_{\|u\|_2=1}r(u)^\top Gr(u)\leq1$, which is an intractable quartic maximization over the unit sphere.

A sum-of-squares relaxation resolves both obstacles by bounding this quartic by a quadratic form. To test whether the ellipsoid $\mathcal E_G=\{G^{1/2}z:\|z\|_2\leq1\}$ lies in $K_y$, we write the slack perturbation for $h=G^{1/2}z$ as
\begin{align*}
    B_h=\sum_{i\in [m]} (G^{1/2}z)_iB_i=\sum_{a\in [m]}z_a \underbrace{\Bigg(\sum_{i\in [m]}(G^{1/2})_{ia}B_i\Bigg)}_{:=V_a}=\sum_{a\in [m]} z_a V_a.
\end{align*}
The squared width in direction $u$ is then
\begin{align*}
    r(u)^\top Gr(u)=\sum_{a\in [m]} (u^\top V_au)^2=\sum_{a\in [m]}u^\top V_a (uu^\top)V_a u.
\end{align*}
For every unit vector $u$, we have $uu^\top\preceq I$. Replacing the inner $uu^\top$ by $I$, which amounts to the Cauchy--Schwarz bound $(u^\top V_au)^2\leq\|V_au\|_2^2$, relaxes the quartic width to the quadratic form of a single matrix sum of squares,
\begin{equation}\label{eq:geo_sos}
 r(u)^\top Gr(u)\leq u^\top \mathcal L_y(G)u,\qquad \mathcal L_y(G):=\sum_{a\in [m]} V_a^2=\sum_{i,j\in [m]}G_{ij}B_iB_j\in \R^{n\times n}.
\end{equation}
Consequently, the single operator inequality $\mathcal L_y(G)\preceq I$ certifies $\sup_{\|u\|_2=1}r(u)^\top Gr(u)\leq1$, and hence $\mathcal E_G\subseteq K_y$ (Lemma~\ref{lem:ml_sos}). The eigenvalues of $\mathcal L_y(G)$ thus play the role of the squared widths of an LP. When $C,A_1,\ldots,A_m$ are diagonal in a common basis, $\mathcal L_y(G)=\operatorname{Diag}(r_k^\top Gr_k)$ with $(r_k)_i:=(B_i)_{kk}$. For a general pencil, $\mathcal L_y(G)$ is defined independently of any basis, and it is linear in $G$. %

We define our potential by replacing the power sum of squared widths in~\eqref{eq:geo_lp_potential} with the trace power $\tr[\mathcal L_y(G)^s]$. The index $p:=2s$ matches the order of the $\ell_{2s}$ Lewis weights in Section~\ref{sec:geo_lp}.
\begin{equation}\label{eq:geo_potential}
 \Phi_p(y):=\frac12\min_{G\succ0}\left\{\frac1s\tr\left[\mathcal L_y(G)^s\right]-\log\det G-m\right\},\qquad p:=2s.
\end{equation}
Since $\|\mathcal L_y(G)\|_{\rm op}\leq\tr[\mathcal L_y(G)^s]^{1/s}\leq n^{1/s}\|\mathcal L_y(G)\|_{\rm op}$, the trace power approximates $\|\mathcal L_y(G)\|_{\rm op}$, and as for LPs, we take the even integer $s:=2\lceil\log(2mn)\rceil$. For a diagonal pencil, the potential~\eqref{eq:geo_potential} reduces to the Lee--Sidford potential~\eqref{eq:geo_lp_potential}. %

The Lewis-weight structure carries over to the matrix setting. Let $G_*=G_*(y)$ be the unique minimizer in~\eqref{eq:geo_potential}, let $W:=\mathcal L_y(G_*)^{s-1}$, and define the map $\mathcal L_y^*:\mathbb S^n\to\mathbb S^m$ by $(\mathcal L_y^*(Z))_{ij}:=\frac12\tr[Z(B_iB_j+B_jB_i)]$ for $i,j\in [m]$. Setting the derivative with respect to $G$ to zero gives the optimality condition and the mass identity (Lemma~\ref{lem:ml_optimizer}),
\begin{equation}\label{eq:geo_equilibrium}
 G_*^{-1}=\mathcal L_y^*(W),\qquad\tr[\mathcal L_y(G_*)^s]=m,
\end{equation}
which generalize the LP identities $G_*^{-1}=\sum_kw_k^{1-1/s}r_kr_k^\top$ and $\sum_kw_k=m$. Convex duality gives a formulation of~\eqref{eq:geo_potential} over slack weights (Lemma~\ref{lem:ml_optimizer}),
\begin{equation}\label{eq:geo_weight_form}
 \Phi_p(y)=\frac12\max_{W\succeq0}\Bigl\{\log\det \mathcal L_y^*(W)-\frac{s-1}s\tr\bigl[W^{s/(s-1)}\bigr]\Bigr\},
\end{equation}
and the maximum is attained at $W$. Here $W$ acts as a matrix weight on the $n$-dimensional slack space, and $G$ acts as a weight on the $m$ coefficient matrices. By~\eqref{eq:geo_equilibrium}, $G_*$ and $W$ determine each other, and our algorithm maintains a matrix in $\mathbb S^m_{++}$ directly. We therefore call every $G\succ0$ a \emph{weight} and $G_*$ the \emph{optimal weight}. The logarithmic barrier $-\log\det S(y)$ has Hessian $\mathcal L_y^*(I)$, so it weights all $n$ slack directions equally with total mass $\tr[I]=n$. Our potential instead chooses $W$ to maximize $\log\det \mathcal L_y^*(W)$ under a Schatten-norm penalty of order $s/(s-1)$, which reduces the total mass to $\tr[W^{s/(s-1)}]=\tr[\mathcal L_y(G_*)^s]=m$.

The sum-of-squares condition $\mathcal L_y(G)\preceq I$ can be stricter than the exact containment $\mathcal E_G\subseteq K_y$. For example, for
\begin{align*}
 B_1=\begin{pmatrix}1&0\\*0&-1\end{pmatrix},\qquad B_2=\begin{pmatrix}0&1\\*1&0\end{pmatrix},
\end{align*}
the region $K_y={\cal E}_I$ is the unit disk in $\R^2$, whereas $\mathcal L_y(I)=2I\not\preceq I$. A certified ellipsoid suffices for the barrier, and a rescaling provides one. At the optimum $G_*$, the mass identity in~\eqref{eq:geo_equilibrium} gives $\|\mathcal L_y(G_*)\|_{\rm op}\leq(\tr[\mathcal L_y(G_*)^s])^{1/s}=m^{1/s}$. Hence the rescaled ellipsoid $m^{-1/(2s)}\mathcal E_{G_*}$ is certified to lie in $K_y$, and
\begin{align*}
 \Phi_p(y)=\underbrace{-\frac12\log\det G_*+\frac m{2s}}_{-\log \operatorname{vol}(m^{-1/(2s)}{\cal E}_{G_*})-c_{m,s}}-\frac m2,
\end{align*}
where $c_{m,s}$ depends only on $m$ and $s$. Thus, up to an additive constant, $\Phi_p(y)$ is the negative log-volume of a certified inscribed ellipsoid (Figure~\ref{fig:geo_ellipsoids}(c)). The next subsection proves the barrier properties of $\Phi_p$ directly.

\subsubsection{Barrier properties of our matrix potential}
\label{sec:geo_proof}

We now sketch why $\Phi_p$, after scaling, satisfies conditions~(a)--(c) of Definition~\ref{def:tech_ov_sc_barrier}. We write the potential~\eqref{eq:geo_potential} as
\begin{align*}
    \Phi_p(y)=\frac12\min_{G\succ0}f_y(G),\qquad f_y(G):=\frac1s\tr[\mathcal L_y(G)^s]-\log\det G-m,
\end{align*}
and $G_*(y)$ is the minimizer of $f_y$. As for $\Phi^{\rm LP}$ in Section~\ref{sec:geo_lp}, the envelope theorem gives the gradient
\begin{align}\label{eq:tech_ov_slack_gradient}
\nabla\Phi_p(y)[h]=-\tr[M_yB_h],\qquad M_y:=\frac12\Bigl(\mathcal L_y(G_*)^s+\sum_{i,j\in[m]}(G_*)_{ij}B_iWB_j\Bigr)\succ0.
\end{align}
For a diagonal pencil, $M_y=\operatorname{Diag}(w_k)$ with the scalar Lewis weights $w_k$. In general, $\tr[M_y]=m$ by the mass identity in~\eqref{eq:geo_equilibrium}, and the Cauchy--Schwarz inequality, applied as in the LP analysis, gives
\begin{equation}\label{eq:geo_budget}
 |\nabla\Phi_p(y)[h]|^2\leq\tr[M_y]\,\tr\left[M_yB_h^2\right]=m\,Q_y[h,h],
\end{equation}
where the \emph{proxy Hessian} $Q_y[h,h]:=\tr\left[M_yB_h^2\right]$ plays the role of $Q_y^{\rm LP}$. It remains to prove the spectral comparison $Q_y\preceq\nabla^2\Phi_p(y)$ and to bound $\nabla^3\Phi_p$.

Unlike the gradient, the second and third derivatives of $\Phi_p$ depend on how $G_*(y)$ moves with $y$. For LPs, the constraints commute, and Lee and Sidford~\cite{ls19} differentiate the weights $w(y)$ explicitly through entrywise products of projection matrices. For SDPs, the matrices $B_i$ need not commute, and differentiating the matrix power $\mathcal L_y(G_*)^{s-1}$ in~\eqref{eq:geo_equilibrium} gives no tractable formula for the derivative of $G_*(y)$. We avoid this derivative and use only the stationarity of $G_*(y)$, through the following standard fact.

\begin{fact}[Derivatives of an optimal value; informal version of Lemma~\ref{lem:ml_contact}]
\label{fact:geo_value}
Let $g_*(t):=\min_zg(t,z)$ for a smooth function $g$ whose minimizer $z_*(t)$ is nondegenerate, and estimate $g_*(t)$ by $g_z(t):=g(t,z(t))$ along a trial curve $z(t)$ with $z(0)=z_*(0)$. Then the following statements hold.
\begin{enumerate}[label=(\roman*)]
\item $g_z''(0)\geq g''_*(0)$, with equality attained if and only if $z'(0)=z_*'(0)$.
\item If $z'(0)=z_*'(0)$, then $g_z'''(0)= g'''_*(0)$.
\end{enumerate}
\end{fact}
For $\Phi_p$, this fact means that the Hessian and the third derivative in a direction $h$ require only the velocity $\frac{\d}{\d t}G_*(y+th)|_{t=0}$ of the optimal weight. By part~(i), this velocity minimizes the second derivative over the velocities of trial curves, which avoids differentiating~\eqref{eq:geo_equilibrium}. The fact requires curves along which the objective can be differentiated explicitly. Step~1 below replaces straight lines by exponential curves, which have this property, and Steps~2 and~3 apply parts~(i) and~(ii) to the Hessian and to the third derivative.

\paragraph{Step 1: exponential curves make the objective explicit.} Fix a base point $y$ and a direction $h$. We view the potential as a function $\varphi(S)$ of a positive definite slack $S$, with $A_1,\ldots,A_m$ fixed and with the optimal weight $G_*(S)$, so that $\Phi_p(y)=\varphi(S(y))$. A congruence of $S$ and a change of basis of the $A_i$ normalize the base point to $S(y)=I$ and $G_*(S(y))=I$, and they change $\varphi$ only by an additive constant (Section~\ref{sec:ml_normalization}). In these coordinates, the optimality condition~\eqref{eq:geo_equilibrium} reads $\mathcal L_y^*(W)=I$. Along the line $y+th$, the slack is $S(y+th)=I+tT$ with $T:=B_h$, and hence $\Phi_p(y+th)=\varphi(I+tT)$.

Differentiating $\varphi(I+tT)$ along the straight line $I+tT$, however, requires the derivatives of the matrix square root $(I+tT)^{-1/2}$ that normalizes the constraints. We instead replace the line $I+tT$ by the exponential curve $S(t)=e^{tT}$, and we move the trial weight along $G(t)=e^{tV}$ for a symmetric matrix $V\in\mathbb S^m$. Along these curves, $-\log\det G(t)=-t\tr[V]$ is linear in $t$, and the matrix of~\eqref{eq:geo_sos} for the weight $G(t)$ at the slack $S(t)$ factorizes as $\mathsf B(t)\mathsf B(t)^\top$, where
\begin{align*}
 \mathsf B(t)=&e^{t\mathsf L}\mathsf Be^{t\mathsf R}\in \R^{n\times mn},\\
 {\sf B}:=\begin{bmatrix}
     B_1 & \cdots & B_m
 \end{bmatrix}\in \R^{n\times mn},\qquad \mathsf L:=&-\frac12T\in \mathbb{S}^n,\qquad \mathsf R:=\frac12(V\otimes I_n-I_m\otimes T)\in \mathbb{S}^{mn}.
\end{align*}
The objective in~\eqref{eq:geo_potential} at the slack $S(t)$ and the weight $G(t)$ is then
\begin{equation}\label{eq:geo_explicit}
 u_{T,V}(t)-t\tr[V]-m,\qquad u_{T,V}(t):=\frac1s\tr\left[(\mathsf B(t)\mathsf B(t)^\top)^s\right].
\end{equation}
Its second and third derivatives along $(e^{tT},e^{tV})$ therefore come from the single trace power $u_{T,V}$, which involves only the fixed matrix $\mathsf B$ and the constant generators $\mathsf L$ and $\mathsf R$.

The exponential curve $S(t)=e^{tT}$ agrees with the line $I+tT$ only to first order, so we relate the derivatives of $\psi(t):=\varphi(e^{tT})$ to those of $\Phi_p(y+th)=\varphi(I+tT)$. Since $S'(0)=T$, $S''(0)=T^2$ and $S'''(0)=T^3$, the chain rule and the identity $\nabla\varphi(I)[Z]=-\tr[M_yZ]$ from~\eqref{eq:tech_ov_slack_gradient} give
\begin{equation}\label{eq:geo_line_curve}
\begin{aligned}
 \nabla^2\Phi_p(y)[h,h]=&~ \nabla^2\varphi(I)[T,T]=\psi''(0)+Q_y[h,h],\\
 \nabla^3\Phi_p(y)[h,h,h]=&~ \nabla^3\varphi(I)[T,T,T]=\psi'''(0)-3\nabla^2\varphi(I)[T,T^2]+\tr[M_yT^3].
\end{aligned}
\end{equation}
\eqref{eq:geo_line_curve} is the only place where we pass between the line and the exponential curve, and its correction terms involve only the first two derivatives of $\varphi$. In particular, the spectral comparison $Q_y\preceq\nabla^2\Phi_p(y)\preceq(2s+1)Q_y$ is equivalent to
\begin{align}\label{eq:tech_ov_step1_goal}
    0\leq\psi''(0)\leq2s\,Q_y[h,h]\qquad \forall h\in \R^{m},
\end{align}
and the lower bound $\psi''(0)\geq0$ states that $\varphi$ is convex along exponential curves. %

\paragraph{Step 2: variational Hessian comparison.} Our goal is~\eqref{eq:tech_ov_step1_goal}. The trial curves $e^{tV}$ have initial velocity $\frac{\d}{\d t}e^{tV}|_{t=0}=V$, which ranges over all $V\in\mathbb S^m$, and the linear term $-t\tr[V]$ in~\eqref{eq:geo_explicit} vanishes after two derivatives. Hence Fact~\ref{fact:geo_value}~(i) gives the variational formula
\begin{align*}
 \psi''(0)=\frac12\min_{V\in\mathbb S^m}u_{T,V}''(0),
\end{align*}
and the minimum is attained at the velocity $V_*:=\frac{\d}{\d t}G_*(e^{tT})\big|_{t=0}$ of the optimal weight. Since $e^{tT}$ and $I+tT$ have the same initial velocity $T$, the velocity $V_*$ also equals $\frac{\d}{\d t}G_*(S(y+th))\big|_{t=0}$. We call $u_{T,V}''(0)$ the \emph{energy} of $V$, and it is a sum of squares (Lemma~\ref{lem:ml_exponential_curvature}),
\begin{align}
 u_{T,V}''(0)=&~ 4\sum_{a=0}^{s-1}\bigl\|({\sf B}{\sf B}^\top)^{(s-1-a)/2}\mathsf X(\mathsf B^\top\mathsf B)^{a/2}\bigr\|_F^2\geq4\sum_{i=1}^m\tr[WX_i^2],\label{eq:tech_ov_energy_T_V}\\
 \mathsf X:=&~ \mathsf B'(0)=[X_1\ \cdots\ X_m],\qquad X_i=\frac12\sum_{j=1}^mV_{ji}B_j-\frac12(TB_i+B_iT).\label{eq:tech_ov_X_Xi}
\end{align}
Combining this expression with~\eqref{eq:geo_line_curve} gives the following lemma.

\begin{lemma}[Variational Hessian; informal version of Lemmas~\ref{lem:ml_second_order} and~\ref{lem:ml_directional_response}]
\label{lem:geo_hessian}
Let $h$ be a direction, let $T=B_h$ in the normalized coordinates, and let $v:=\nabla^2\Phi_p(y)[h,h]$. Then
\begin{equation}\label{eq:geo_identity}
 v=Q_y[h,h]+\frac12\min_{V\in\mathbb S^m}u_{T,V}''(0),
\end{equation}
and the minimum is attained at the velocity $V_*$ of the optimal weight. Consequently,
\begin{align*}
 Q_y\preceq\nabla^2\Phi_p(y)\preceq(2s+1)\,Q_y,\qquad|\nabla\Phi_p(y)[h]|^2\leq mv,\qquad\|V_*\|_F\leq\sqrt{6v}.
\end{align*}
\end{lemma}
\begin{proof}[Proof sketch]
We read the two spectral bounds off the energy, and we bound $V_*$ through the orthonormality of the $B_i$ under $W$. The energy is nonnegative by~\eqref{eq:tech_ov_energy_T_V}, which gives the lower bound $Q_y\preceq\nabla^2\Phi_p(y)$. Together with~\eqref{eq:geo_budget}, this lower bound implies $|\nabla\Phi_p(y)[h]|^2\leq mv$. For the upper bound, the trial velocity $V=0$ has energy $u_{T,0}''(0)\leq4s\,Q_y[h,h]$ by Lemma~\ref{lem:ml_exponential_curvature}, and hence $\nabla^2\Phi_p(y)\preceq(2s+1)Q_y$.

For the bound on $\|V_*\|_F$, consider the norm $\|(Z_i)\|_W^2:=\sum_{i=1}^m\tr[WZ_i^2]$ on $m$-tuples of symmetric matrices. By~\eqref{eq:tech_ov_energy_T_V}, $\|(X_i)\|_W\leq \frac{1}{2}(u_{T,V}''(0))^{1/2}$, and the triangle inequality applied to~\eqref{eq:tech_ov_X_Xi} gives
\begin{align*}
    \frac{1}{2}\Bigl\|\Bigl(\sum_jV_{ji}B_j\Bigr)_i\Bigr\|_W\leq \|(X_i)\|_W + \frac{1}{2}\|(TB_i+B_iT)_i\|_W.
\end{align*}
The optimality condition $\mathcal L_y^*(W)=I$ states that $B_1,\ldots,B_m$ are orthonormal for the inner product $\langle Z,Z'\rangle_W:=\frac12\tr[W(ZZ'+Z'Z)]$. It follows that the tuple $(\sum_jV_{ji}B_j)_i$ has $W$-norm $\|V\|_F$, and the tuple $(TB_i+B_iT)_i$ has $W$-norm at most $2\sqrt{Q_y[h,h]}$. At $V=V_*$, these bounds give $\frac12\|V_*\|_F\leq\frac12(u_{T,V_*}''(0))^{1/2}+\sqrt{Q_y[h,h]}$, and the identity $u_{T,V_*}''(0)=2(v-Q_y[h,h])$ from~\eqref{eq:geo_identity} yields $\|V_*\|_F\leq\sqrt{6v}$.
\end{proof}

\paragraph{Step 3: third derivatives from the first-order velocity.} It remains to prove
\begin{align*}
    |\nabla^3\Phi_p(y)[h,h,h]|\leq150s\,v^{3/2}.
\end{align*}
By~\eqref{eq:geo_line_curve}, the main term is $\psi'''(0)$. Differentiating~\eqref{eq:geo_identity} would involve the second derivative of $G_*$. Instead, the trial curve $e^{tV_*}$ has velocity $V_*$ at $t=0$, so Fact~\ref{fact:geo_value}~(ii) gives $\psi'''(0)=\frac12u_{T,V_*}'''(0)$. It therefore suffices to bound the third derivative of the explicit trace power $u_{T,V_*}$ by its second derivative $u_{T,V_*}''(0)=2(v-Q_y[h,h])\leq2v$.

When $T$, $V$ and all $B_i$ are diagonal, $u_{T,V}$ is an exponential sum $u(t)=\sum_jc_je^{\lambda_jt}$ with $c_j\geq0$ and $|\lambda_j|\leq4s\max\{\|\mathsf L\|_{\rm op},\|\mathsf R\|_{\rm op}\}$, and a term-by-term comparison gives $|u'''|\leq(\max_j|\lambda_j|)\,u''$. The following lemma extends this inequality to noncommuting matrices.

\begin{lemma}[Exponential trace inequality; informal version of Lemma~\ref{lem:ml_exponential_curvature}]
\label{lem:geo_trace}
Let $\mathsf B(t)=e^{t\mathsf L}\mathsf Be^{t\mathsf R}$ for a matrix $\mathsf B$ and symmetric matrices $\mathsf L,\mathsf R$, and let $u(t):=\frac1s\tr[(\mathsf B(t)\mathsf B(t)^\top)^s]$. Then $u''(0)\geq0$ and
\begin{equation}\label{eq:geo_trace_inequality}
 |u'''(0)|\leq40s\max\{\|\mathsf L\|_{\rm op},\|\mathsf R\|_{\rm op}\}\,u''(0).
\end{equation}
\end{lemma}

It remains to bound $\|\mathsf L\|_{\rm op}=\frac12\|T\|_{\rm op}$ and $\|\mathsf R\|_{\rm op}\leq\frac12(\|V_*\|_{\rm op}+\|T\|_{\rm op})$ by $O(\sqrt v)$. Lemma~\ref{lem:geo_hessian} gives $\|V_*\|_{\rm op}\leq\|V_*\|_F\leq\sqrt{6v}$. For $T=B_h$, we must bound $\|B_h\|_{\rm op}$ by a multiple of $\sqrt v=(\nabla^2\Phi_p(y)[h,h])^{1/2}$. Since $K_y$ is the unit ball of the norm $h\mapsto\|B_h\|_{\rm op}$ and the Dikin ellipsoid $\{h:\nabla^2\Phi_p(y)[h,h]\leq1\}$ is the unit ball of $h\mapsto(\nabla^2\Phi_p(y)[h,h])^{1/2}$, a bound $\|B_h\|_{\rm op}\leq C\sqrt v$ for every $h$ is equivalent to the containment of the Dikin ellipsoid in $C\,K_y$. The following lemma proves this rounding property.

\begin{lemma}[Rounding; informal version of Lemma~\ref{lem:ml_rounding}]
\label{lem:geo_rounding}
For every direction $h$, $\|B_h\|_{\rm op}^2\leq2(mn)^{1/(s-1)}Q_y[h,h]$. In particular, $\|B_h\|_{\rm op}^2\leq9\,\nabla^2\Phi_p(y)[h,h]$ when $s\geq1+\log(mn)$.
\end{lemma}

\begin{proof}[Proof sketch]
We compare $B_h$ with the certified ellipsoid and then remove the matrix $W$ by H\"older's inequality. In the normalized coordinates of Step~1, $G_*=I$, so the certified ellipsoid $m^{-1/(2s)}\mathcal E_{G_*}$ of Section~\ref{sec:geo_sdp} is the Euclidean ball of radius $m^{-1/(2s)}$. Its containment in $K_y$ and the optimality condition $\mathcal L_y^*(W)=I$ give
\begin{align*}
    \|B_h\|_{\rm op}^2\leq m^{1/s}\|h\|_2^2= m^{1/s}\tr[WB_h^2]=m^{1/s}\tr[\mathcal L_y(G_*)^{s-1}B_h^2].
\end{align*}
Dropping the second term of $M_y$ in~\eqref{eq:tech_ov_slack_gradient} gives $Q_y[h,h]\geq \frac{1}{2}\tr[\mathcal L_y(G_*)^sB_h^2]$, and H\"older's inequality then yields
\begin{align*}
    \tr[\mathcal L_y(G_*)^{s-1}B_h^2] \leq (\tr[\mathcal L_y(G_*)^sB_h^2])^{1-1/s}(\tr[B_h^2])^{1/s}\leq (2Q_y[h,h])^{1-1/s} (n\|B_h\|_{\rm op}^2)^{1/s}.
\end{align*}
Combining the two displays, we obtain
\begin{align*}
 \|B_h\|_{\rm op}^2\leq m^{1/s}\,(2Q_y[h,h])^{1-1/s}\,(n\|B_h\|_{\rm op}^2)^{1/s},
\end{align*}
and solving for $\|B_h\|_{\rm op}^2$ proves the first bound $\|B_h\|_{\rm op}^2\leq2(mn)^{1/(s-1)}Q_y[h,h]$. The second bound follows from $Q_y\preceq\nabla^2\Phi_p(y)$ and $2(mn)^{1/(s-1)}\leq2e<9$.
\end{proof}

Combining the bounds on $\|V_*\|_{\rm op}$ and $\|T\|_{\rm op}$ gives $\|\mathsf L\|_{\rm op},\|\mathsf R\|_{\rm op}\leq3\sqrt v$, and Lemma~\ref{lem:geo_trace} then yields $|\psi'''(0)|\leq120s\,v^{3/2}$. The two correction terms in~\eqref{eq:geo_line_curve} are $O(\sqrt s)\,v^{3/2}$ by the spectral comparison and Lemma~\ref{lem:geo_rounding}, which completes the bound on the third derivative (Theorem~\ref{thm:ml_matrix_sc}).

\paragraph{Completing the proof.} The lower bound $\nabla^2\Phi_p(y)\succeq Q_y\succ0$ gives strict convexity. For condition~(a) of Definition~\ref{def:tech_ov_sc_barrier}, the formulation~\eqref{eq:geo_weight_form} at $W=I$ gives $2\Phi_p(y)\geq\log\det \mathcal L_y^*(I)-n$, which tends to $+\infty$ as $y\to\partial\Omega$ (Theorem~\ref{thm:ml_matrix_sc}). Finally, scaling by $\alpha=10^6s^2$ reduces the constant in the third-derivative bound to $150s/\sqrt\alpha<2$ and multiplies the gradient constant $m$ by $\alpha$. Hence the matrix Lee--Sidford barrier $F:=\alpha\Phi_p$ is $\nu$-self-concordant with
\begin{align*}
 \nu=\alpha m=O(m\log^2(2mn)).
\end{align*}

\subsection{Technical overview for the efficient SDP solver}
\label{sec:algorithmic_overview}
The matrix Lee--Sidford barrier $F=\alpha\Phi_p$ of Section~\ref{sec:technical_overview} has parameter $\nu=\widetilde O(m)$, so a short-step interior-point method needs only $\widetilde O(\sqrt m)$ iterations. This section explains how to implement each iteration efficiently. The target running time of Theorem~\ref{thm:intro_randomized} decomposes as
\begin{equation}\label{eq:overview_anatomy}
 \underbrace{\widetilde O(\sqrt m)}_{\text{queries}}
 \cdot\underbrace{\widetilde O(mn^2)}_{\text{amortized work per query}}
 +\underbrace{\widetilde O(n^\omega)}_{\text{slack factorizations}}
 =\widetilde O(n^2m^{3/2}+n^\omega),
\end{equation}
where $mn^2$ is the size of the input $A_1,\ldots,A_m$ and $n^\omega$ is the cost of factoring one $n\times n$ slack matrix. We first explain why a direct implementation misses this target by a factor $n^{\omega-2}$ (Section~\ref{sec:overview_direct}). We then describe the three data structures that remove its bottlenecks (Sections~\ref{sec:overview_sketch}--\ref{sec:overview_slack}), and Section~\ref{sec:overview_interfaces} adds preprocessing and postprocessing and combines all components into the complete algorithm. Throughout this overview, we assume $m\leq n$ and use the multiplication bounds of Fact~\ref{fact:rectangular_thresholds}. The notation $\widetilde O$ also hides polynomial factors in $s$ and $\alpha$, logarithmic factors in the accuracy and the input magnitudes, and the arbitrarily small exponent losses of Theorem~\ref{thm:intro_randomized}.

\subsubsection{A direct implementation and its bottlenecks}
\label{sec:overview_direct}
A short-step method tracks the central path of $F(y)+tb^\top y$. Each of its iterations, which we call \emph{phases}, multiplies $t$ by $1+\Theta(\nu^{-1/2})$ and recenters the point by a few Newton-type steps. The matrix Lee--Sidford barrier is defined by the inner minimization $\Phi_p(y)=\frac12\min_{G\succ0}f_y(G)$ of Section~\ref{sec:geo_proof}, so its derivatives at $y$ depend on the optimal weight $G_*(y)$. An algorithm can compute $G_*(y)$ only approximately, and it evaluates the formulas of Section~\ref{sec:geo_proof} at a nearby weight $G\succ0$. With $\mathcal L_y$ as in~\eqref{eq:geo_sos}, let
\begin{gather*}
 M_G:=\frac12\Bigl(\mathcal L_y(G)^s+\sum_{i,j\in[m]}G_{ij}B_i\mathcal L_y(G)^{s-1}B_j\Bigr),\\
 (g_G)_i:=-\alpha\tr[M_GB_i],\qquad
 Q_G[h,h]:=\tr[M_GB_h^2],
\end{gather*}
where $g_G$ is the \emph{fixed-weight gradient} since it is the gradient of $y\mapsto\frac\alpha2f_y(G)$ with $G$ held fixed, and $Q_G$ is the \emph{proxy Hessian}. At the optimal weight $G_*(y)$,~\eqref{eq:tech_ov_slack_gradient} and Lemma~\ref{lem:geo_hessian} give
\begin{equation}\label{eq:overview_derivatives}
 g_{G_*}=\nabla F(y),\qquad
 \alpha Q_{G_*}\preceq\nabla^2F(y)\preceq(2s+1)\,\alpha Q_{G_*}.
\end{equation}
Because differentiating $G_*(y)$ is expensive, the algorithm preconditions each step with the scaled proxy Hessian $\alpha Q_G$ instead of $\nabla^2F(y)$, taking damped steps of the form $y\gets y-\gamma_y(\alpha Q_G)^{-1}(g_G+tb)$, which we call \emph{$y$-steps}. Since $\alpha Q_{G_*}$ approximates $\nabla^2F(y)$ within a factor $O(s)$ by~\eqref{eq:overview_derivatives}, and $\alpha Q_G$ stays within a constant factor of $\alpha Q_{G_*}$ when $G$ is close to $G_*(y)$ (see~\eqref{eq:overview_bias} below), preconditioned gradient descent with step size $\gamma_y=\Theta(1/s)$ contracts the gap between $F(y)+tb^\top y$ and its minimum by a factor $1-\Theta(1/s)$ per step, up to an additive term from the gradient error (Lemma~\ref{lem:random_direct_centering}). Thus, $O(s)$ such steps suffice to replace the single exact Newton step of a textbook short-step IPM.

We track the quality of a weight $G=LL^\top$ through two quantities. The first is the \emph{weight gap} $\Delta_y(G):=f_y(G)-f_y(G_*(y))\geq0$, which controls how closely $(g_G,Q_G)$ approximate the true barrier derivatives in~\eqref{eq:overview_derivatives}. Whenever $s\Delta_y(G)$ is at most a small constant, Theorem~\ref{thm:weight_interface} gives
\begin{equation}\label{eq:overview_bias}
\begin{gathered}
 \|g_G-\nabla F(y)\|_{\nabla^2F(y)^{-1}}\leq3\sqrt{\alpha s\Delta_y(G)},\\
 \Bigl(1-3\sqrt{s\Delta_y(G)}\Bigr)^2Q_{G_*}\preceq Q_G
 \preceq\Bigl(1+3\sqrt{s\Delta_y(G)}\Bigr)^2Q_{G_*}.
\end{gathered}
\end{equation}
Because each $y$-step requires gradient error $O(s^{-1/2})$ and $\alpha=\Theta(s^2)$, we need $\Delta_y(G)\leq\eta=\Theta(s^{-4})$, and we denote such $G$ an \emph{accurate} weight. The second quantity is the \emph{stationarity matrix} $R:=L^\top\mathcal L_y^*(\mathcal L_y(G)^{s-1})L$, which measures the multiplicative violation of the optimality condition~\eqref{eq:geo_equilibrium} ($R=I$ if and only if $G=G_*(y)$). For a larger gap $\Delta_y(G)\leq\eta_{\rm w}=\Theta(1/s)$, it still holds that $\frac12I\preceq R\preceq2I$, and we denote such $G$ a \emph{warm} weight.%

Algorithm~\ref{alg:overview_direct} shows one phase of a direct implementation. Each $y$-step forms the $m$ normalized constraint matrices $B_i$ explicitly and recomputes an accurate weight, which costs $\widetilde O(mn^\omega)$. Over the $\widetilde O(\sqrt m)$ phases, the total is $\widetilde O(m^{3/2}n^\omega)$, a factor $n^{\omega-2}$ above the target~\eqref{eq:overview_anatomy}.

\begin{algorithm}[!ht]
\caption{One phase of a direct implementation with the matrix Lee--Sidford barrier}
\label{alg:overview_direct}
\begin{algorithmic}[1]
\State $t\gets(1+\gamma_t)\,t$, where $\gamma_t=\Theta(\nu^{-1/2})$
\For{each of the $N_y=O(s)$ $y$-steps} \Comment{one query each}
  \State form $S(y)=\sum_{i\in[m]}y_iA_i-C$ \Comment{$O(mn^2)$}
  \State factor $S(y)$ and form $B_i=S(y)^{-1/2}A_iS(y)^{-1/2}$ for all $i\in[m]$ \Comment{$O(mn^\omega)$}
  \State compute an accurate weight $G$ by $\widetilde O(1)$ dense iterations \Comment{$\widetilde O(mn^\omega)$}
  \State form $M_G$, $g_G$ and $Q_G$ \Comment{$O(mn^\omega)$}
  \State $y\gets y-\gamma_y(\alpha Q_G)^{-1}(g_G+tb)$ \Comment{$O(m^\omega)$}
\EndFor
\end{algorithmic}
\end{algorithm}

The excess cost has three sources, and our algorithm removes each of them with one data structure.
\begin{itemize}
\item Forming the $m$ normalized constraint matrices costs $O(mn^\omega)$. The gradient/Hessian estimation data structure $\mathcal D_{\rm est}$ never forms them. It applies them to thin blocks of vectors and estimates $g_G$, $R$ and $Q_G$ from random sketches (Section~\ref{sec:overview_sketch}).
\item The weight is recomputed at every step. The weight maintenance data structure $\mathcal D_{\rm weight}$ maintains the weight $\mathcal D_{\rm weight}.G$ between steps and corrects it by a damped update driven by an estimate of $R$ (Section~\ref{sec:overview_client}).
\item The slack is factored at every step. When $m\leq n^{1/4}$, this cost can dominate, and the slack maintenance data structure $\mathcal D_{\rm slack}$ reuses one factorization while $t$ changes by a factor of two (Section~\ref{sec:overview_slack}).
\end{itemize}
When $m>\sqrt{n}$, even sketches of width $m$ are too expensive. We therefore implement $\mathcal D_{\rm est}$ in two ways. For $m\leq \sqrt{n}$, the data structure $\mathcal D^{\rm f}_{\rm est}$ draws fresh estimates at every query. For $m>\sqrt{n}$, the data structure $\mathcal D^{\rm m}_{\rm est}$ maintains its estimates between queries (Section~\ref{sec:overview_gradient}). Part~\ref{part:algorithm} states the guarantee of each data structure in Section~\ref{sec:data_structures}. Table~\ref{tab:overview_ladder} previews how the cost of the main loop decreases as these data structures are added.

\begin{table}[!ht]
\centering
\resizebox{\textwidth}{!}{%
\small
\begin{tabular}{@{}llllccc@{}}
\toprule
 & Implementations & Per query & Total & $m=n^{1/6}$ & $m=n^{3/8}$ & $m=n^{3/4}$\\
\midrule
(A) & forms $B_i$, $G$, $M_G$, $Q_G$ densely & $mn^\omega$ & $m^{3/2}n^\omega$ & $n^{5/2}$ & $n^{45/16}$ & $n^{27/8}$\\
(B) & maintains $G$, sketches at width $m$ & $n^\omega+\mathcal C(m)$ & $\sqrt m\,(n^\omega+\mathcal C(m))$ & $n^{7/3}$ & $\boldsymbol{n^{41/16}}$ & $n^{103/32}$\\
(C) & maintains estimates ($m>\sqrt n$) & $n^\omega+mn^2$ & $\sqrt m\,n^\omega+n^2m^{3/2}$ & -- & -- & $\boldsymbol{n^{25/8}}$\\
(D) & reuses factorizations ($m\leq n^{1/4}$) & $mn^2$ & $n^\omega+n^2m^{3/2}$ & $\boldsymbol{n^{9/4}}$ & -- & --\\
\bottomrule
\end{tabular}
}
\caption{The cost of different implementations of the main loop over $\widetilde O(\sqrt m)$ queries as the data structures are added. Row~(A) is Algorithm~\ref{alg:overview_direct}. Row~(B) uses block operations and fresh sketches of width $m$ (Section~\ref{sec:overview_sketch}) and maintains the weight between queries (Section~\ref{sec:overview_client}). Here, $\mathcal C(m)=\widetilde O(mn^2)$ if $m\leq \sqrt{n}$, and $\mathcal C(m)=\widetilde O(n^\omega+n^2m^{9/8})$ in general. This implementation already achieves the target runtime of Theorem~\ref{thm:intro_randomized} for $n^{1/4}\leq m\leq \sqrt{n}$. Row~(C) applies when $m>\sqrt{n}$ (Section~\ref{sec:overview_gradient}), and row~(D) applies when $m\leq n^{1/4}$ (Section~\ref{sec:overview_slack}). The last three columns evaluate the total time with $\omega=9/4$ in three regimes of $m$, with the target runtime of Theorem~\ref{thm:intro_randomized} in bold.}
\label{tab:overview_ladder}
\end{table}

\subsubsection{Gradient/Hessian estimation}
\label{sec:overview_sketch}
To avoid the $O(mn^\omega)$ cost of forming the normalized constraint matrices $B_i$ explicitly, the estimation data structure $\mathcal D_{\rm est}$ accesses them only implicitly through matrix-block products. Specifically, let $J\in \R^{n\times n}$ be a \emph{whitener} of the slack $S(y)$ such that $JS(y)J^\top=I$, for example the inverse of a Cholesky factor of $S(y)$. Every such $J$ equals $OS(y)^{-1/2}$ for an orthogonal matrix $O$. The substitution $B_i\mapsto JA_iJ^\top=OB_iO^\top$ conjugates the slack-space matrices $\mathcal L_y(G)$ and $M_G$ by $O$, and it leaves $f_y$ and the coefficient-space quantities $g_G\in\R^m$ and $Q_G,R\in\mathbb S^m$ unchanged, since these depend on the $B_i$ only through traces such as $\tr[(OM_GO^\top)(OB_iO^\top)]=\tr[M_GB_i]$. We therefore use $B_i=JA_iJ^\top$ from now on, and the slack-space matrices below are formed with these $B_i$. Let $G=LL^\top$ be the Cholesky factorization of the current weight. We define the weighted combinations $V_a:=\sum_{i=1}^m L_{ia}B_i$ as in Section~\ref{sec:geo_sdp} (with $L$ in place of $G^{1/2}$), so that $\mathcal L_y(G)=\sum_{a=1}^m V_a^2$. Instead of forming $V_1,\ldots,V_m$ explicitly, $\mathcal D_{\rm est}$ applies them only through \emph{block operations} of \emph{width} $q$, which map a matrix $X\in\R^{n\times q}$ to $[V_1X\ \cdots\ V_mX]$ or apply the adjoint map $(X_1,\ldots,X_m)\mapsto\sum_{a=1}^m V_aX_a$. Since $V_aX=\sum_{i=1}^m L_{ia}JA_iJ^\top X$, the forward operation multiplies $X$ by $J^\top$, multiplies the result by the stacked matrix $A_{\rm stack}:=[A_1\ \cdots\ A_m]^\top\in\R^{mn\times n}$, multiplies the $m$ resulting blocks by $J$, and mixes them by $L$. Its cost is
\begin{equation}\label{eq:overview_block_cost}
 \mathcal C(q):=\underbrace{\Tmat(mn,n,q)}_{\text{multiply by }A_{\rm stack}}
 +\underbrace{\Tmat(n,n,mq)}_{\text{whiten $m$ blocks}}
 +\underbrace{\Tmat(m,m,nq)}_{\text{mix by }L}
 +\underbrace{mnq}_{\text{read and write}},
\end{equation}
and the adjoint map has the same cost by transposition. Composing the two maps applies $\mathcal L_y(G)$ to an $n\times q$ block. Fast rectangular multiplication gives
\begin{equation}\label{eq:overview_width}
\begin{aligned}
 \mathcal C(m)=&~\widetilde O(n^\omega+n^2m^{9/8}),\quad
 \mathcal C(m^{2/3})=\widetilde O(n^\omega+mn^2)\quad\forall m\leq n,\\
 \mathcal C(m)=&~\widetilde O(mn^2)\quad\forall m\leq \sqrt{n}.
\end{aligned}
\end{equation}
However, exact evaluation of $M_G$ and $R$ needs width $n$ and costs $\mathcal C(n)=O(mn^\omega)$, which is no cheaper than the direct implementation. The data structure $\mathcal D_{\rm est}$ reduces the width by sketching the following factors.

Let $\mathcal V:=[V_1\ \cdots\ V_m]\in\R^{n\times mn}$, so that $\mathcal V\mathcal V^\top=\mathcal L_y(G)$. The integer-power factors
\begin{equation}\label{eq:overview_factors}
 \mathcal K:=\mathcal L_y(G)^{s/2-1}\mathcal V,\qquad F_0:=2^{-1/2}\mathcal L_y(G)^{s/2},\qquad
 F_a:=2^{-1/2}V_a\mathcal K\quad\forall a\in[m]
\end{equation}
satisfy $\mathcal K\mathcal K^\top=\mathcal L_y(G)^{s-1}$, $M_G=\sum_{a=0}^mF_aF_a^\top$ and $R_{ab}=2\langle F_a,F_b\rangle$. After padding by zero columns, the factors are the block rows of an even power of a symmetric dilation,
\begin{equation}\label{eq:overview_dilation}
 \mathcal F_G:=\begin{pmatrix}F_0&0\\0&\mathsf F\end{pmatrix}
 =\frac1{\sqrt2}\,\mathcal Z_G^s,\qquad
 \mathcal Z_G:=\begin{pmatrix}0&\mathcal V\\\mathcal V^\top&0\end{pmatrix},
\end{equation}
where $\mathsf F=[F_1^\top\ \cdots\ F_m^\top]^\top\in \R^{mn\times mn}$. A trace-power inequality gives $\|\mathcal F_G-\mathcal F_{G_*}\|_F^2\leq s\Delta_y(G)$ (Lemma~\ref{lem:weight_factor_distance}). This single bound converts the weight gap into the accuracy of every derived quantity, and it implies the error bounds for the gradient and the proxy Hessian in~\eqref{eq:overview_bias}.

The data structure $\mathcal D_{\rm est}$ sketches $F_0$ and $F_1,\ldots,F_m$ in their respective input spaces, with independent sketching matrices $\Xi_0,\Xi_1$ that satisfy $\E[\Xi_b\Xi_b^\top]=I$. The sketched factors $Y_0=F_0\Xi_0$ and $Y_a=F_a\Xi_1$ give
\begin{equation}\label{eq:overview_sampler}
 \widehat g_i=-\alpha\tr\Bigl[B_i\sum_{a=0}^mY_aY_a^\top\Bigr],\qquad
 \widehat R_{ab}=2\langle Y_a,Y_b\rangle,
\end{equation}
whose expectations are $g_G$ and $R$. Evaluating the sketched factors uses only integer powers of $\mathcal L_y(G)$ and block operations, at cost $\widetilde O(n^\omega+\mathcal C(q))$ for $q$ columns.

A plain Gaussian sketch of width $q$ has mean squared error $\widetilde O(m^2/q)$ in the norm $\|v\|_{G/\alpha}:=(v^\top Gv/\alpha)^{1/2}$. Indeed, $\widehat g$ is a vector of Gaussian quadratic forms whose mean squared error in this norm scales with $\|\mathcal K\|_{\rm op}^2\|\mathcal K\|_F^2/q$. We have $\|\mathcal K\|_F^2=\tr[\mathcal L_y(G)^{s-1}]=O(m)$, but the high power $s-1$ can concentrate most of it onto the top eigenvector of $\mathcal L_y(G)$, so $\|\mathcal K\|_{\rm op}^2$ can also be of order $m$. To reduce the sketching error (or equivalently, the required width $q$), we follow Hutch++~\cite{mmmw21} by computing the contribution of an approximate top-$\widetilde O(q)$ singular subspace of $\mathcal K$ exactly and sketching only the residual. Because at most $q$ squared singular values of $\mathcal K$ can exceed $\|\mathcal K\|_F^2/q$, the squared operator norm of the residual is at most $O(m/q)$, and the mean squared errors of $\widehat g$ and $\widehat R$ drop to $\widetilde O(m^2/q^2)$. Moreover,  $\mathcal L_y(G)\preceq3I$ at warm weights, so $\mathsf F^\top\mathsf F=\frac12\mathcal K^\top\mathcal L_y(G)\mathcal K \preceq\frac32\mathcal K^\top\mathcal K$ and $F_0F_0^\top=\frac12\mathcal L_y(G)^s\preceq\frac32\mathcal K\mathcal K^\top$. Thus, a single subspace approximation on $\mathcal K$ works for both $F_1,\ldots,F_m$ and $F_0$ (Lemma~\ref{lem:common_gram_deflation}). Therefore, width $\widetilde O(m)$ sketching matrices suffice for a full estimate, which are sampled by $\mathcal D^{\rm f}_{\rm est}$ at every query when $m\leq\sqrt n$.

The same sketched factors $Y_0,\ldots,Y_m$ also yield an estimate of the proxy Hessian $Q_G$. Because $Y:=[Y_0\ \cdots\ Y_m]$ has $O(mq)$ columns in total, we first compress it from the right using an independent Gaussian sketching matrix $\Omega$ of width $\widetilde O(\sqrt m)$, and then compute $\widehat Q_{ij}:=\langle B_iY\Omega,B_jY\Omega\rangle$ using only one block operation of width $\widetilde O(\sqrt m)$. For every width $q$ used by the algorithm, $\widehat Q$ approximates $Q_G$ within a constant factor (Lemma~\ref{lem_more_metric}).

\subsubsection{The path-following method and weight maintenance}
\label{sec:overview_client}
\label{sec:overview_weight}

Our algorithm has the same phases as the direct method (Algorithm~\ref{alg:overview_direct}), but avoids recomputing the preconditioner and optimal weight from scratch at every step. At the first query of each phase, we obtain $\widehat Q$ from $\mathcal D_{\rm est}$ and factor the \emph{phase preconditioner} $P:=\alpha\widehat Q$. Since $\widehat Q$ approximates $Q_G$ within a constant factor,~\eqref{eq:overview_derivatives} and~\eqref{eq:overview_bias} show that $P$ approximates $\nabla^2F(y)$ within a factor $O(s)$ at the first point of the phase. The phase stays within a small local ball, and self-concordance of $F$ keeps this factor $O(s)$ throughout the phase. The iterate then advances via $N_y=O(s)$ damped $y$-steps
\begin{equation}\label{eq:overview_point_step}
 y^+=y-\gamma_yP^{-1}(\widehat g+tb),\qquad\gamma_y=\Theta(1/s),
\end{equation}
each of which moves $y$ by a small constant distance in the local norm.

\begin{figure}[!ht]
\centering
\begin{tikzpicture}[x=0.82cm,y=0.9cm,font=\footnotesize,
  rowlab/.style={anchor=east,font=\small},
  note/.style={anchor=west,align=left,text width=4.6cm},
  seg/.style={rounded corners=1pt,line width=.5pt}]
 \def\len{8}
 \node[rowlab] at (-0.15,3) {queries};
 \node[rowlab] at (-0.15,2) {phases};
 \node[rowlab] at (-0.15,1) {epochs ($m>\sqrt{n}$)};
 \node[rowlab] at (-0.15,0) {anchors ($m\leq n^{1/4}$)};
 \draw[black!25] (0,3) -- (\len,3);
 \foreach \i in {0,...,31} \fill[black!40] ({0.125+0.25*\i},3) circle (0.035);
 \foreach \i in {0,4,...,28} \fill[black!85] ({0.125+0.25*\i},3) circle (0.075);
 \foreach \i in {0,4}
   \draw[seg,draw=overviewblue!70,fill=overviewblue!14] (\i,1.8) rectangle ++(2,0.4);
 \foreach \i in {2,6}
   \draw[seg,draw=overviewblue!70,fill=overviewblue!4] (\i,1.8) rectangle ++(2,0.4);
 \draw[seg,draw=overviewteal!80,fill=overviewteal!14] (0,0.8) rectangle (4,1.2);
 \draw[seg,draw=overviewteal!80,fill=overviewteal!4] (4,0.8) rectangle (\len,1.2);
 \foreach \x in {0,4}
   \node[anchor=west,inner sep=1.5pt,text=overviewteal!60!black] at (\x,1) {full estimate};
 \draw[seg,draw=black!55,fill=black!5] (0,-0.2) rectangle (6,0.2);
 \draw[seg,draw=black!55,fill=black!13] (6,-0.2) rectangle (\len,0.2);
 \foreach \x in {0,6}
   \node[anchor=west,inner sep=1.5pt] at (\x,0) {factor $S$};
 \node[note] at (\len+0.3,3) {one $y$-query, then $N_G$ $G$-queries at the same point};
 \node[note] at (\len+0.3,2) {raise $t$, compress $\widehat Q$ and factor $P$, then $N_y$ $y$-steps};
 \node[note] at (\len+0.3,1) {$H$ phases, which start with a full estimate};
 \node[note] at (\len+0.3,0) {factor the slack when $t$ has doubled};
\end{tikzpicture}
\caption{The schedule of the main loop on four nested time scales (not to scale). In the first row, each large black dot represents a $y$-query, followed by $N_G=O(s\log s)$ $G$-queries at the same point $y$ (the small grey dots). The second row shows the phases, each consisting of $N_y=O(s)$ $y$-steps. The third and fourth rows correspond to the large-$m$ and small-$m$ regimes: when $m>\sqrt n$, phases are grouped into epochs of length $H\approx m^{1/6}$ that start with a full estimate; when $m\leq n^{1/4}$, an \emph{anchor} slack $S(y_{\rm anc})$ is Cholesky factorized only at the first phase after $t$ doubles (every $\widetilde O(\sqrt m)$ phases), reusing it at all intermediate queries, whereas for $m>n^{1/4}$, $S(y)$ is factored at every query.}
\label{fig:overview_clocks}
\end{figure}
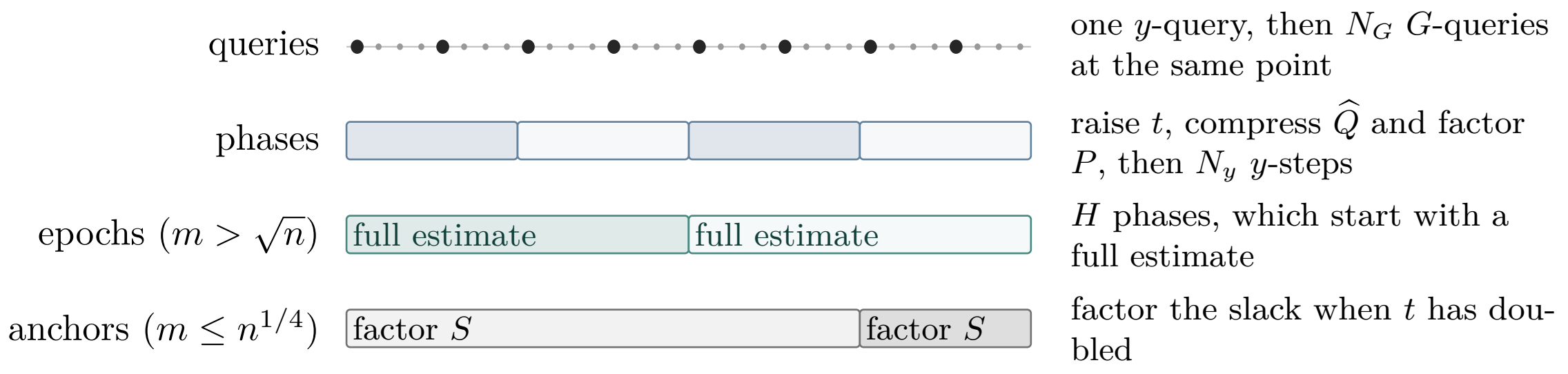

\paragraph{Maintaining the weight via $\mathcal D_{\rm weight}$.} To keep the weight $G$ accurate before every $y$-step without solving the inner minimization from scratch, the data structure $\mathcal D_{\rm weight}$ stores the current weight $\mathcal D_{\rm weight}.G=LL^\top$, and the algorithm follows each $y$-step~\eqref{eq:overview_point_step} by $N_G=O(s\log s)$ damped multiplicative updates of this weight at the new point, which we call \emph{$G$-steps}:
\begin{equation}\label{eq:overview_weight_step}
 G^+=L\bigl((1-\gamma_G)I+2\gamma_G(I+\widehat R)^{-1}\bigr)L^\top,\qquad\gamma_G=\Theta(1/s).
\end{equation}
Each $G$-step~\eqref{eq:overview_weight_step} uses an estimate $\widehat R$ of the stationarity matrix $R$ from $\mathcal D_{\rm est}$ and takes $O(m^\omega)$ time to invert $I+\widehat R$ and update the Cholesky factor $L$. Theorem~\ref{thm:weight_interface} establishes two properties that make this alternation work:
\begin{itemize}%
\item \emph{Bounded drift during a $y$-step:} When $y$ moves and $G$ is held fixed, $(\Delta_y(G))^{1/2}$ changes by at most $O(\sqrt{s/\alpha})=O(s^{-1/2})$ per unit of local distance. Consequently, a $y$-step~\eqref{eq:overview_point_step} of small constant local length starting from an accurate weight ($\Delta_y(G)\leq\eta=\Theta(s^{-4})$) increases the gap to at most $\eta_{\rm w}/4$, leaving $G$ warm ($\Delta_{y^+}(G)\leq\eta_{\rm w}=\Theta(1/s)$) at the new point $y^+$.
\item \emph{Linear contraction during $G$-steps:} At a fixed point, whenever $G$ is warm ($\Delta_y(G)\leq\eta_{\rm w}$) and $\|\widehat R-R\|_F\leq\tau:=\sqrt\eta/32$, the $G$-step~\eqref{eq:overview_weight_step} contracts the gap according to
\begin{equation}\label{eq:overview_contraction}
 \Delta_y(G^+)\leq(1-\gamma_G/16)\,\Delta_y(G)+6\gamma_G\tau^2.
\end{equation}
Because the additive error terms in~\eqref{eq:overview_contraction} accumulate to at most $96\tau^2<\eta/10$, applying $N_G:=\lceil(16/\gamma_G)\log(4\eta_{\rm w}/\eta)\rceil=O(s\log s)$ $G$-steps at the fixed point $y^+$ restores an accurate weight ($\Delta_{y^+}(G)\leq\eta$) before the next $y$-step (Theorem~\ref{thm:single_loop_centering}).
\end{itemize}
Thus, each phase performs $N_y=O(s)$ $y$-steps and $(N_G+1)N_y=O(s^2\log s)$ queries in total. We call the query that precedes a $y$-step a \emph{$y$-query} and the $N_G$ queries that drive the subsequent $G$-steps \emph{$G$-queries}. Figure~\ref{fig:overview_clocks} illustrates how queries and phases nest within the epochs and slack anchors of Sections~\ref{sec:overview_gradient} and~\ref{sec:overview_slack}.

In Part~\ref{part:algorithm}, we formalize the guarantees of the three data structures as modular interface theorems:
\begin{itemize}%
\item $\mathcal D_{\rm slack}$ returns an exact or nearly exact whitener of $S(y)$ (Theorem~\ref{thm:slack_interface});
\item $\mathcal D_{\rm weight}$ computes an accurate initial weight and contracts the weight gap via~\eqref{eq:overview_weight_step} (Theorem~\ref{thm:weight_interface}); and
\item $\mathcal D_{\rm est}$ returns $\widehat g$ within the gradient error tolerance of a $y$-step, $\widehat R$ within the Frobenius tolerance $\tau$ of a $G$-step, and $\widehat Q$ within a constant spectral factor of $Q_G$ (Theorem~\ref{thm:estimator_interface}).
\end{itemize}
An induction over the queries then shows that every query remains within the valid operating regime of all three data structures throughout the run (Theorem~\ref{thm:centering_correctness}).

\subsubsection{Gradient/Hessian estimate maintenance when \texorpdfstring{$m>\sqrt{n}$}{m>sqrt n}}
\label{sec:overview_gradient}
When $m>\sqrt n$, computing a fresh sketch of width $\widetilde O(m)$ at every query costs $\mathcal C(m)=\widetilde O(n^\omega+n^2m^{9/8})$, which exceeds our target per-query budget of $\widetilde O(n^\omega+mn^2)$. To avoid paying for a width-$m$ sketch at every step, the data structure $\mathcal D^{\rm m}_{\rm est}$ computes a full estimate only periodically and updates it across the queries in between by sketching the increments at smaller widths. We first formulate this task as a standalone problem, then specialize it to the queries in the path-following loop.

\begin{problem}[Online Gram estimation]\label{prob:gram_estimate}
A sequence of matrices $\mathcal F_1,\mathcal F_2,\ldots$ arrives one at a time, with $\|\mathcal F_j\|_F^2=O(m)$ and $\|\mathcal F_j-\mathcal F_{j-1}\|_F^2=\widetilde O(1)$. At every step $j$, we need to estimate $g_j\in \R^m$, a linear function of the Gram matrix $\mathcal F_j\mathcal F_j^\top$, within distance $1$ in a given norm $\|x\|_{Z_j}:=(x^\top Z_j x)^{1/2}$. To do so, we may invoke a width-$q$ sketching primitive at cost $\widetilde O(n^\omega+\mathcal C(q))$ in two ways:
\begin{itemize}%
\item \emph{Full estimate:} sketch $\mathcal F_j$ alone and obtain an unbiased estimate of $g_j$ with mean squared error $\widetilde O(m^2/q^2)$, so error $O(1)$ needs width $q=\widetilde O(m)$.
\item \emph{Coupled pair estimate:} sketch $\mathcal F_{j-1}$ and $\mathcal F_j$ with the same sketching matrix and obtain a coupled pair $(\widehat g_{j-1}^{\rm pair},\widehat g_j^{\rm pair})$ whose difference is an unbiased estimate of $g_j-g_{j-1}$ with mean squared error $\widetilde O(m/q^2)$, since the shared sketch makes the two sketching errors nearly cancel.
\end{itemize}
How can we combine full estimates and coupled pair estimates to minimize the amortized cost per step?
\end{problem}

When the norm is fixed ($Z_j\equiv Z$), our solution follows the variance-reduction framework of SARAH~\cite{nglst17} and SPIDER~\cite{fllz18}, originally developed for minimizing a finite sum $f=\sum_i f_i$ by maintaining a running gradient estimate $v^+:=v+\nabla f_i(x^+)-\nabla f_i(x)$ with a shared random index $i$ and periodically resetting $v$ to a full gradient. In Problem~\ref{prob:gram_estimate}, the shared sketching matrix serves as the random index $i$. We partition the steps into \emph{epochs} of $N$ consecutive steps, where the first step computes a full estimate of width $\widetilde O(m)$, and every subsequent step $j$ in the epoch draws a coupled pair estimate of width $q$ and updates
\begin{equation}\label{eq:overview_plain_recursion}
 \widehat g_j:=\widehat g_{j-1}+\widehat g_j^{\rm pair}-\widehat g_{j-1}^{\rm pair}.
\end{equation}
Since $\widehat g_j^{\rm pair}-\widehat g_{j-1}^{\rm pair}$ is an unbiased estimator of $g_j-g_{j-1}$ conditioned on the past (even if $\mathcal F_j$ depends on earlier estimates), and
\begin{align*}
   \widehat g_j-g_j = (\widehat g_{j-1}-g_{j-1})+\left((\widehat g_j^{\rm pair}-\widehat g_{j-1}^{\rm pair}) - (g_j-g_{j-1})\right),
\end{align*}
the error $\widehat g_j-g_j$ is a martingale within an epoch, and the squared norms of its increments add in expectation. After $k\leq N$ steps in an epoch, the mean squared error is $\widetilde O(1+km/q^2)$, and choosing $q=\widetilde O(\sqrt{Nm})$ keeps it $O(1)$ throughout the epoch. Moreover, by a Hilbert-space martingale concentration inequality~\cite{pinelis94,pinelis1999correction}, this accuracy can be achieved for all steps simultaneously with high probability. The amortized cost of a step is $\widetilde O(n^\omega+\mathcal C(\sqrt{Nm})+\mathcal C(m)/N)$.

When the norm $Z_j$ changes from step to step, an error $\widehat g_{j-1}-g_{j-1}$ that was small in $\|\cdot\|_{Z_{j-1}}$ can become large in $\|\cdot\|_{Z_j}$. To prevent past errors from expanding, the update transports the previous estimate through a linear map $\Gamma_j$ chosen before a new pair $(\widehat g_{j-1}^{\rm pair},\widehat g_j^{\rm pair})$ is drawn:
\begin{equation}\label{eq:overview_recursion}
 \widehat g_j:=\Gamma_j\widehat g_{j-1}
       +\widehat g_j^{\rm pair}-\Gamma_j\widehat g_{j-1}^{\rm pair},
 \qquad \Gamma_j^\top Z_j\Gamma_j\preceq Z_{j-1}.
\end{equation}
The spectral condition for $\Gamma_j$ ensures $\|\Gamma_j x\|_{Z_j}\leq\|x\|_{Z_{j-1}}$ for all $x\in\R^m$, so the transported error $\Gamma_j(\widehat g_{j-1} - g_{j-1})$ never grows. Meanwhile, we can rewrite the correction term in~\eqref{eq:overview_recursion} as
\begin{align*}
    \widehat g_j^{\rm pair}-\Gamma_j\widehat g_{j-1}^{\rm pair} = \Gamma_j(\widehat g_j^{\rm pair}-\widehat g_{j-1}^{\rm pair})+(I-\Gamma_j)\widehat g_j^{\rm pair},
\end{align*}
which decomposes it into a contracted coupled increment with mean squared error $\widetilde O(m/q^2)$ and a fresh full-estimate component along $I-\Gamma_j$. We use the \emph{refresh charge} $r_j\geq0$ as an amortization parameter that measures the size of this second component (with $r_j=0$ when $\Gamma_j=I$), giving a per-step mean squared error of $\widetilde O\bigl(m(1+r_j)/q^2\bigr)$ and requiring width $q=\widetilde O((m(N+\sum_j r_j))^{1/2})$ over an $N$-step epoch. We choose $\Gamma_j$ to act as $I$ except on subspaces where $Z_j$ grows relative to $Z_{j-1}$, which keeps the refresh charges small. It therefore remains to instantiate $\mathcal F_j$, $Z_j$, and $\Gamma_j$ for the path-following queries and to bound the total refresh charge $\sum_jr_j$ of an epoch.

\paragraph{The data structure $\mathcal D^{\rm m}_{\rm est}$.} We now instantiate the online Gram estimation framework where each step $j$ corresponds to a query $(y_j,G_j)$ of the path-following loop, $\mathcal F_j=\mathcal F_{G_j}$ is the dilated factor from~\eqref{eq:overview_dilation}, and each epoch consists of $H$ consecutive phases. Lemmas~\ref{lem_more_projected_gradient} and~\ref{lem_more_gradient_increments} establish the full-estimate and coupled-increment error bounds required by the framework for our two estimation targets: the stationarity matrix $R_{G_j}(y_j)$ and the fixed-weight gradient $g_{G_j}(y_j)$. Because $R_{G_j}(y_j)$ is measured in the fixed Frobenius norm, $\mathcal D^{\rm m}_{\rm est}$ maintains $\widehat R$ directly via the fixed-norm update~\eqref{eq:overview_plain_recursion}. On the other hand, the gradient estimate $\widehat g_j$ must be accurate in a local norm that tracks the evolving weight and Hessian, so $\mathcal D^{\rm m}_{\rm est}$ maintains a reference metric $\mathcal D^{\rm m}_{\rm est}.Z$ and updates $\widehat g_j$ via~\eqref{eq:overview_recursion}. Within each phase, $\mathcal D^{\rm m}_{\rm est}.Z$ is held fixed, so $\Gamma_j=I$ and no refresh charge is incurred ($r_j=0$). At the start of a new phase, $\mathcal D^{\rm m}_{\rm est}$ updates the previous metric $Z_{j-1}$ in two spectral clipping steps, which give the new metric an upper and a lower bound.
\begin{enumerate}
    \item Clipping $Z_{j-1}$ from above relative to the current weight $G_j$ gives $Z_j^-\preceq Z_{j-1}$ with $Z_j^-\preceq2G_j$. The sketching errors are controlled in the norm $\|\cdot\|_{G_j/\alpha}$ (Section~\ref{sec:overview_sketch}), and this upper bound keeps them controlled in the new metric when the weight shrinks in some directions.
    \item Clipping $Z_j^-$ from below relative to $c_Z\widehat Q_j^{-1}$, for a small constant $c_Z$, gives the new metric $Z_j\succeq Z_j^-$ with $c_Z\widehat Q_j^{-1}\preceq Z_j\preceq3G_j$. Since $\nabla^2F(y)\succeq\Omega(\alpha\widehat Q_j)$ throughout the phase, this lower bound turns accuracy in $\|\cdot\|_{Z_j/\alpha}$ into the gradient accuracy required by the $y$-step.
\end{enumerate}
The first step only shrinks the metric, so errors that were small in the old metric stay small at no cost. The second step can enlarge the metric, and the map $\Gamma_j:=Z_j^{-1}Z_j^-$ satisfies the spectral condition in~\eqref{eq:overview_recursion}, because $Z_j^-\preceq Z_j$ gives $\Gamma_j^\top Z_j\Gamma_j=Z_j^-Z_j^{-1}Z_j^-\preceq Z_j^-\preceq Z_{j-1}$. Each spectral clipping uses a constant-accuracy diagonalization of an $m\times m$ matrix, which takes $\widetilde O(m^\omega)\leq\widetilde O(mn^2)$ time since $m\leq n$ (Lemma~\ref{lem:approximate_diagonalization}).

\begin{figure}[!ht]
\centering
\begin{tikzpicture}
\begin{axis}[
 width=0.72\linewidth, height=5.6cm,
 xmin=0, xmax=0.5, ymin=0.8, ymax=1.19,
 xlabel={$H$},
 ylabel={cost per query},
 xtick={0,0.125,0.1666667,0.3333333,0.5},
 xticklabels={$1$,$m^{1/8}$,$m^{1/6}$,$m^{1/3}$,$m^{1/2}$},
 ytick={0.875,1,1.125},
 yticklabels={$m^{7/8}n^2$,$mn^2$,$m^{9/8}n^2$},
 tick align=outside,
 axis line style={black!60},
 clip=false]
 \fill[black!10] (axis cs:0.125,0.8) rectangle
                        (axis cs:0.1666667,1.19);
 \addplot[black!45,dashed,line width=0.8pt] coordinates {(0,1) (0.5,1)};
 \addplot[overviewblue,line width=1.4pt]
   coordinates {(0,1) (0.1666667,1) (0.5,1.125)};
 \addplot[overviewteal!80,line width=1.4pt]
   coordinates {(0,1.125) (0.325,0.8)};
 \node[font=\small,anchor=south east,text=overviewblue]
   at (axis cs:0.49,1.127) {coupled pairs $\mathcal C(H\sqrt m)$};
 \node[font=\small,anchor=west,text=overviewteal!85!black]
   at (axis cs:0.272,0.885) {full estimates $\mathcal C(m)/H$};
 \node[font=\small,anchor=north east,text=black!60]
   at (axis cs:0.5,0.99) {target cost};
 \node[font=\small,anchor=west]
   at (axis cs:0.012,1.152) {$H=1$: full estimates dominate};
\end{axis}
\end{tikzpicture}
\caption{Choosing the epoch length $H$ when $m > \sqrt{n}$. The curves show the worst-case upper bounds on the amortized cost per query. The coupled pair estimate bound (rising blue curve) stays within the target $\widetilde O(mn^2)$ budget for $H \leq m^{1/6}$, while the amortized full-estimate bound (falling green curve) fits once $H \geq m^{1/8}$. Any choice of $H$ in the shaded window $[m^{1/8}, m^{1/6}]$ balances these two bounds to achieve the target per-query cost of $\widetilde O(n^\omega + mn^2)$. }
\label{fig:overview_window}
\end{figure}
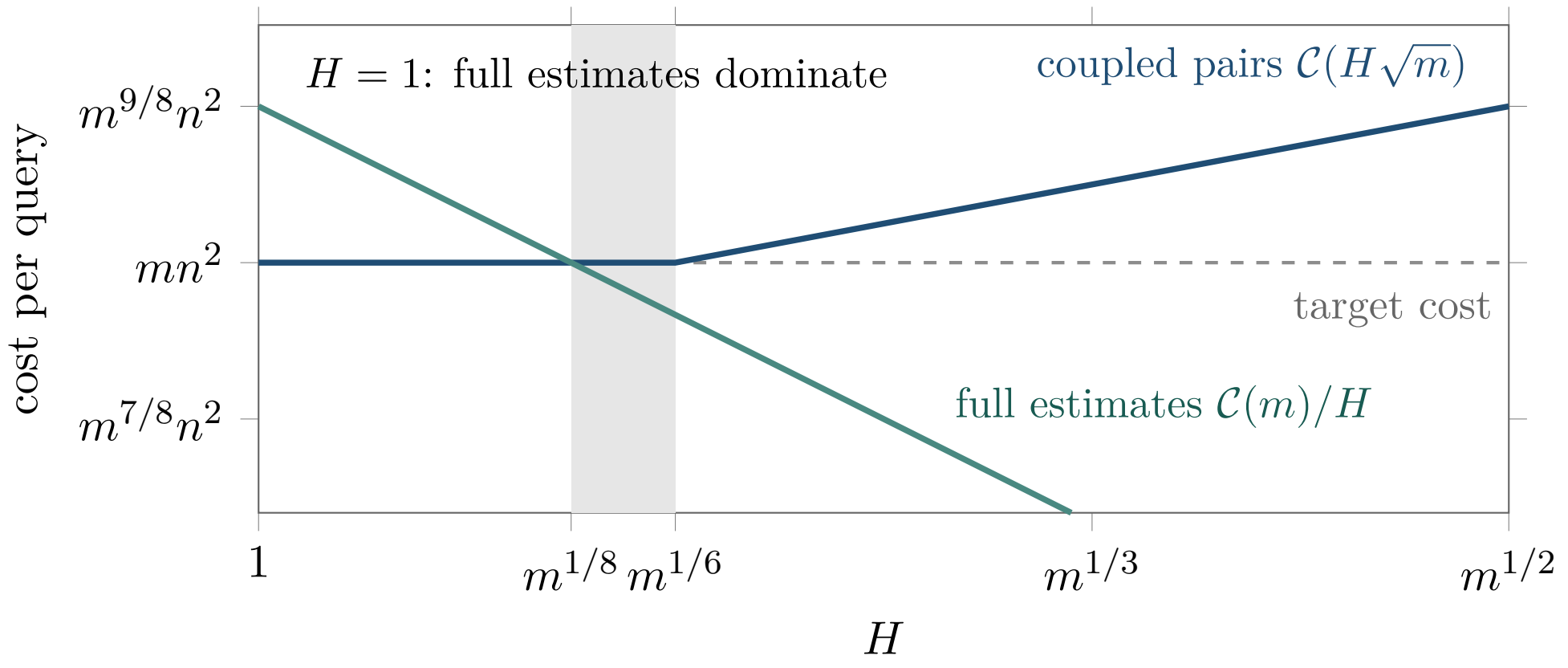

To bound the cumulative refresh charge over an $H$-phase epoch, we track the log-determinant divergence $\tr[G_j^{-1}Z_j]-\log\det(G_j^{-1}Z_j)-m$ between the reference metric and the weight: the slow drift of $G$ across $H$ phases increases this potential by at most $\widetilde O(H^2)$, whereas the two spectral clippings together decrease it by at least $r_j/2$, implying $\sum_j r_j=\widetilde O(H^2)$ (Corollary~\ref{cor:phase_reference_metric}). Since each phase makes $O(s^2\log s)=\widetilde O(1)$ queries, an epoch makes $N=\widetilde O(H)$ queries, and its accumulated mean squared error is
\begin{align*}
    \widetilde{O}\Big(m\Big(N+\sum_jr_j\Big)/q^2\Big)=\widetilde{O}\left(mH^2/q^2\right).
\end{align*}
Hence, the coupled pair estimates only need width $q=\widetilde O(H\sqrt m)$. Combined with the width-$\widetilde O(m)$ full estimate at the start of each epoch, this choice guarantees accuracy across all queries simultaneously, even under adaptive query sequences (Section~\ref{sec:oracle_maintenance}). The amortized per-query cost is
\begin{equation}\label{eq:overview_query_cost}
 \widetilde O\left(n^\omega+\mathcal C(H\sqrt m)
                                  +\mathcal C(m)/H\right).
\end{equation}
Since $\sqrt{n}<m\leq n$, by~\eqref{eq:overview_width}, we have $\mathcal C(H\sqrt m)=\widetilde O(n^\omega+mn^2)$ if $H\sqrt m\leq m^{2/3}$, and $\mathcal C(m)/H=\widetilde O(n^\omega/H+n^2m^{9/8}/H)=\widetilde O(mn^2)$ if $H\geq m^{1/8}$. Hence, any epoch length $H\in[m^{1/8},m^{1/6}]$ achieves $\widetilde O(n^\omega+mn^2)$ per query (Figure~\ref{fig:overview_window}), and we take $H:=\lfloor m^{1/6}\rfloor$.

\subsubsection{Slack maintenance when \texorpdfstring{$m\leq \sqrt{n}$}{m<=sqrt n}}
\label{sec:overview_slack}

When $m\leq\sqrt n$, a block operation of width $\widetilde O(m)$ costs $\mathcal C(m)=\widetilde O(mn^2)$ by~\eqref{eq:overview_width}, which fits the per-query budget of~\eqref{eq:overview_anatomy}, so the algorithm can use the fresh estimator $\mathcal D^{\rm f}_{\rm est}$ at every query instead of maintaining estimates across queries. The only remaining bottleneck is applying a whitener of the slack $S(y)\in \R^{n\times n}$ during each block operation (Section~\ref{sec:overview_sketch}). Recall that a whitener $J$ of $S(y)$ satisfies $JS(y)J^\top=I$. In the intermediate regime $n^{1/4}<m\leq\sqrt n$, we have $n^\omega\leq n^{9/4}<mn^2$, so $\mathcal D_{\rm slack}$ simply computes an exact Cholesky factorization $S(y)=L_SL_S^\top$ at every query in $O(n^\omega)$ time and applies the exact whitener $J=L_S^{-1}$. When $m\leq n^{1/4}$, a single Cholesky factorization costs $O(n^\omega)$, which can exceed the per-query budget $\widetilde O(mn^2)$. To overcome this barrier, we use polynomial approximation and amortization techniques discussed below.

\paragraph{Approximate whitening via polynomial approximation.} The data structure $\mathcal D_{\rm slack}$ applies an approximate whitener $\widehat J=\mathsf K J\in \R^{n\times n}$ for the block operations, where $J$ is an exact whitener for $S(y)$ and $\mathsf K$ is a symmetric matrix satisfying $\|\mathsf K-I\|_{\rm op}\leq\varepsilon_{\rm w}$ for an error tolerance $\varepsilon_{\rm w}=1/\operatorname{poly}(mn)$ fixed in~\eqref{eq:bias_budget}. To construct $\widehat J$ cheaply across many queries, $\mathcal D_{\rm slack}$ saves an \emph{anchor slack} $S_{\rm anc}:=S(y_{\rm anc})$ together with its Cholesky factor $L_{\rm anc}$, so that $J_{\rm anc}:=L_{\rm anc}^{-1}$ satisfies $J_{\rm anc}S_{\rm anc}J_{\rm anc}^\top=I$. This anchor slack will be reused until the path parameter $t$ doubles. Throughout this doubling interval, the current slack $S(y)$ remains spectrally close to $S_{\rm anc}$ (Lemma~\ref{thm:slack_window}):
\begin{align*}
 \Upsilon^{-1}S_{\rm anc}\preceq S(y)\preceq\Upsilon S_{\rm anc}
 \qquad\text{for}\qquad \Upsilon=O(sm)=\widetilde O(m).
\end{align*}
Thus, the relative slack $\Sigma:=J_{\rm anc}S(y)J_{\rm anc}^\top$ has all its eigenvalues in $[\Upsilon^{-1},\Upsilon]$. To compute a whitener for $\Sigma$, we apply a matrix polynomial that approximates $\Sigma^{-1/2}$. More specifically, it is a standard result in approximation theory that there exists an explicit polynomial $\chi_\Upsilon:\R\rightarrow \R$ of degree $d_{\rm w}=O\bigl(\Upsilon\log(\Upsilon/\varepsilon_{\rm w})\bigr)=\widetilde O(m)$ that approximates $x\mapsto x^{-1/2}$ in the sense that $|\sqrt x\,\chi_\Upsilon(x)-1|\leq\varepsilon_{\rm w}$ on $[\Upsilon^{-1},\Upsilon]$ (Lemma~\ref{lem:inverse_sqrt_polynomial}). Then, we can define an approximate whitener for $S(y)$:
\begin{align*}
    \widehat J:=\chi_\Upsilon(\Sigma)J_{\rm anc},\qquad \mathsf K:=\chi_\Upsilon(\Sigma)\Sigma^{1/2}.
\end{align*}
We can verify that $\widehat J = {\sf K} J$, where $J:=\Sigma^{-1/2}J_{\rm anc}$ is an exact whitener for $S(y)$. Indeed, $JS(y)J^\top=\Sigma^{-1/2}(J_{\rm anc}S(y)J_{\rm anc}^\top)\Sigma^{-1/2}=\Sigma^{-1/2}\Sigma \Sigma^{-1/2}=I$. We can also verify that $\|\mathsf K-I\|_{\rm op}\leq\varepsilon_{\rm w}$, which follows from the polynomial approximation guarantee of $\chi_{\Upsilon}$ and the spectrum of $\Sigma$.

\paragraph{Amortized slack whitening cost.} Recall from Section~\ref{sec:overview_sketch} that each block operation in~\eqref{eq:overview_block_cost} applies the whitener twice: first multiplying $X\in\R^{n\times q}$ by $\widehat J^\top$, and then multiplying the stacked matrix $\bigl[A_1\widehat J^\top X\ \cdots\ A_m\widehat J^\top X\bigr]\in\R^{n\times mq}$ by $\widehat J$, where $q=\widetilde{O}(m)$ in this regime. To apply $\widehat J=\chi_\Upsilon(\Sigma)J_{\rm anc}$ (or $\widehat J^\top$) to an $n\times mq$ block without forming the $n\times n$ matrices $\Sigma$ or $\widehat J$ explicitly, we evaluate the degree-$d_{\rm w}$ matrix polynomial $\chi_\Upsilon(\Sigma)$ using Horner's method, which reduces to $d_{\rm w}=\widetilde O(m)$ multiplications of $\Sigma=J_{\rm anc}S(y)J_{\rm anc}^\top$ with $n\times mq$ matrices. Each multiplication by $\Sigma$ consists of a triangular solve with $L_{\rm anc}^\top$, a product with $S(y)$, and a triangular solve with $L_{\rm anc}$, which takes $O\bigl(\Tmat(n,n,mq)\bigr)$ time (Lemma~\ref{lem:fast_factorizations}). Hence, the total cost for the approximate whitener $\widehat{J}$ is
\begin{align*}
    O(d_{\rm w}\Tmat(n,n,m q))=\widetilde{O}(m\Tmat(n,n,m^2))=\widetilde{O}(mn^2),
\end{align*}
where the last step uses $m\leq n^{1/4}$ and the fast rectangular matrix multiplication in Fact~\ref{fact:rectangular_thresholds}. Meanwhile, because the path parameter $t$ doubles only $O(\log(2+\nu/(t_0\epsilon)))=\widetilde O(1)$ times over all $\widetilde O(\sqrt m)$ phases of the algorithm, $\mathcal D_{\rm slack}$ computes an exact anchor whitener $J_{\rm anc}$ in $O(n^\omega)$ time only $\widetilde O(1)$ times in total, contributing $\widetilde O(n^\omega/\sqrt m)$ amortized time per query. Therefore, the total runtime for $\mathcal D_{\rm slack}$ is $\widetilde{O}(n^\omega + m^{3/2}n^2)$.

Finally, the whitening error can bias each estimate. Lemma~\ref{lem:whitening_consequences} bounds the bias of $\widehat g$ and $\widehat R$ by $(2mn)^{O(1)}\varepsilon_{\rm w}$, and the choice of $\varepsilon_{\rm w}$ in~\eqref{eq:bias_budget} keeps it below a quarter of the accuracy required by Theorem~\ref{thm:estimator_interface}. Moreover, since $m\leq n^{1/4}$, the approximate whitener is only used by $\mathcal D^{\rm f}_{\rm est}$, which recomputes its estimates from fresh sketching matrices at every query, so the bias of an estimate depends only on the current query and will not accumulate.

\subsubsection{Preprocessing and postprocessing}
\label{sec:overview_interfaces}
Algorithm~\ref{alg:overview_main_loop} combines the data structures of Sections~\ref{sec:overview_direct}--\ref{sec:overview_slack} into the main path-following loop, which takes a near-central starting point and tracks the central path to a near-optimal dual solution. To obtain a complete primal-dual SDP solver from the bounds $B$ and $\rho$, we formulate an auxiliary SDP with a known interior point and run three stages on it: (i) \emph{preprocessing}, which uses reverse path following to turn the known interior point into a near-central starting point for the target objective; (ii) the \emph{main loop}, which follows the central path forward to a near-optimal dual solution $y$; and (iii) \emph{postprocessing}, which converts the final dual iterate into an explicit primal matrix $X\succeq0$. Because all three stages operate on the same auxiliary SDP and reuse the components of the main loop, the entire solver has the same runtime bound as the main loop.

\begin{algorithm}[!htbp]
\caption{A simplified main loop of our SDP solver. The full versions are Algorithms~\ref{alg:more_sdp}, \ref{alg:maintained_center}, \ref{alg:fresh_query}, \ref{alg:weight_inverse_step}, \ref{alg:more_gradient} and~\ref{alg:implicit_slack}.}
\label{alg:overview_main_loop}
\begin{algorithmic}[1]
\Procedure{AmortizedSDP}{$y,t$} \Comment{$y$ is near-central for $tb$}
 \State Initialize $\mathcal D_{\rm weight}.G$ to be an accurate weight at $y$
 \For{each of the $\widetilde O(\sqrt m)$ phases} \Comment{Section~\ref{sec:overview_direct}}
  \State $t\gets(1+\gamma_t)t$
  \For{$k=1,\ldots,N_y$} \Comment{$N_y=\widetilde{O}(1)$, Section~\ref{sec:overview_client}}\label{ln:tech_ov_phase_start}
   \State $(\widehat g,\widehat R,\widehat Q)\gets\Call{BarrierOracle}{y,t,\mathbf 1_{\{k=1\}}}$ \Comment{$y$-query}
   \State if $k=1$, then $P\gets\alpha\widehat Q$, factorize $P$ \Comment{phase preconditioner}
   \State $y\gets y-\gamma_yP^{-1}(\widehat g+tb)$ \Comment{$y$-step}
   \For{$\ell=1,\ldots,N_G$}\Comment{$N_G=\widetilde{O}(1)$}
    \State $(\widehat g,\widehat R,\widehat Q)\gets\Call{BarrierOracle}{y,t,0}$ \Comment{$G$-query}
    \State update $\mathcal D_{\rm weight}.G$ using the $G$-step~\eqref{eq:overview_weight_step} with $\widehat R$
   \EndFor
  \EndFor\label{ln:tech_ov_phase_end}
 \EndFor
 \State \Return $(y,t)$
\EndProcedure
\Procedure{BarrierOracle}{$y,t,a$} \Comment{$a=1$ at the first query of a phase}
 \If{$m\leq n^{1/4}$} \Comment{$\mathcal D_{\rm slack}$, Section~\ref{sec:overview_slack}}
    \If{this is the first query or $t$ has doubled since the last anchor}
        \State factorize a new anchor $S(y)=L_{\rm anc}L_{\rm anc}^\top$, and $J_{\rm anc}\gets L_{\rm anc}^{-1}$
    \EndIf
  \State $\widehat J\gets\chi_\Upsilon(\Sigma)J_{\rm anc}$ \Comment{$\Sigma=J_{\rm anc}S(y)J_{\rm anc}^\top$}
 \Else
  \State factorize $S(y)=L_SL_S^\top$, and $\widehat J\gets L_S^{-1}$
 \EndIf
 \If{$m\leq\sqrt n$} \Comment{$\mathcal D^{\rm f}_{\rm est}$, Section~\ref{sec:overview_sketch}}
  \State sketch ${\cal F}$ of width $\widetilde O(m)$ and obtain sketched factors $Y$\Comment{full estimate}\label{ln:tech_ov_Y}
  \State compute $(\widehat g,\widehat R)$ by~\eqref{eq:overview_sampler}\label{ln:tech_ov_g_R}
  \State if $a=1$, then sketch $Y$ on the right and compute $\widehat Q$\label{ln:tech_ov_Q}
 \Else\Comment{$\mathcal D^{\rm m}_{\rm est}$, Section~\ref{sec:overview_gradient}}
    \If{the query starts an epoch of $H$ phases}
        \State compute $Y$, $(\widehat g,\widehat R)$, and $\widehat{Q}$ as Lines~\ref{ln:tech_ov_Y}--\ref{ln:tech_ov_Q}\Comment{full estimate}
        \State $Z\gets\mathcal D_{\rm weight}.G$
    \Else\Comment{coupled pair estimate}
        \State $\Gamma\gets I$
        \State if $a=1$, then compute $\widehat{Q}$ using $\mathcal D^{\rm m}_{\rm est}.Y$ (for the previous ${\cal F}$), clip $Z$, $\Gamma\gets Z^{-1}Z^-$
        \State replace $\widehat J$ by a whitener of $S(y)$ aligned with $\mathcal D^{\rm m}_{\rm est}.J$\Comment{Lemma~\ref{lem_more_local_geometry}}
        \State sketch the previous and the current ${\cal F}$ with the same sketching matrices of width $\widetilde O(H\sqrt m)$, obtain $(\widehat g^{\rm pair}_{\rm prev},\widehat g^{\rm pair})$ and $(\widehat R^{\rm pair}_{\rm prev},\widehat R^{\rm pair})$, and store the current sketched factors in $Y$
        \State update $\widehat g$ by~\eqref{eq:overview_recursion} and $\widehat R$ by~\eqref{eq:overview_plain_recursion}
    \EndIf
    \State $\mathcal D^{\rm m}_{\rm est}.y\gets y$, $\mathcal D^{\rm m}_{\rm est}.G\gets\mathcal D_{\rm weight}.G$, $\mathcal D^{\rm m}_{\rm est}.J\gets\widehat J$, $\mathcal D^{\rm m}_{\rm est}.Y\gets Y$
 \EndIf
 \State \Return $(\widehat g,\widehat R,\widehat Q)$
\EndProcedure
\end{algorithmic}
\end{algorithm}

\paragraph{The auxiliary SDP.} Recall from~\eqref{eq:intro_SDP_assumptions} that the input SDP instance is guaranteed to have some dual optimum $y_*$ and some strictly feasible point $y_{\rm f}$ that satisfy $\|y_*\|_\infty,\|y_{\rm f}\|_\infty\leq B$ and $S(y_{\rm f})\succeq\rho I$, and we are only given the bounds $B\geq1$ and $\rho>0$, but not $y_{\rm f}$ itself. Moreover, under these assumptions, the original dual domain $\Omega$ may be unbounded. We therefore augment $y\in\R^m$ with a scalar shift $z\in\R$ and solve an \emph{auxiliary SDP} that minimizes $c^\top u$ with $c:=(b,K)$ over $u=(y,z)\in\R^{m+1}$ subject to
\begin{equation}\label{eq:overview_auxiliary_slack}
 \overline S(u):=\operatorname{diag}\bigl(S(y)+zI,\ 2B\mathbf 1+y,\
 2B\mathbf 1-y,\ z_{\max}-z,\ \delta+z\bigr)\succ0,
\end{equation}
where $z_{\max}:=2(\|C\|_F+\rho+1)$, $K:=1+(4B\|b\|_1+1)/\rho$, and $\delta:=\epsilon/(16K)$. This formulation has two key properties (Lemma~\ref{lem:penalty_pencil}):
\begin{itemize}%
\item \emph{Bounded domain and explicit interior point:} The scalar diagonal blocks restrict $u$ to the box $\|y\|_\infty<2B$ and $-\delta<z<z_{\max}$. Furthermore, at an explicit point $u_0:=(0,z_{\max}/2)$, we have $\overline S(u_0)\succeq I$, so that the auxiliary domain contains a Euclidean ball around $u_0$ of inverse-polynomial radius.
\item \emph{Exact recovery of a strictly feasible dual solution:} For any auxiliary feasible point $u=(y,z)$ with objective value within $\epsilon/64$ of the auxiliary optimum, it must satisfy $-\delta<z\leq-\delta/4$, which implies that $y$ must be strictly feasible for the original dual SDP with $b^\top y\leq {\rm OPT}+3\epsilon/64$.
\end{itemize}
The auxiliary SDP satisfies the assumptions of Theorem~\ref{thm:geo_informal_barrier}, and its primal and dual problems are strictly feasible (Lemma~\ref{lem:penalty_pencil}), even when every feasible matrix of the original primal problem is singular. Since it has $m+1$ variables and matrix dimension $n+2m+2$, Theorem~\ref{thm:geo_informal_barrier} and all data structures from Sections~\ref{sec:overview_sketch}--\ref{sec:overview_slack} apply directly with $(m,n)$ replaced by $(O(m),O(n))$.

\paragraph{Preprocessing: reverse path following from $u_0$.} Although $u_0$ is an explicit interior point of the auxiliary domain, it may lie far from the central path of the target objective $c$. To initialize the main loop, we first construct an artificial objective $c_0$ for which $u_0$ is already near-central. Querying $\mathcal D_{\rm est}$ once at $u_0$ returns a gradient estimate $\widehat g_0$ satisfying $\|\widehat g_0-\nabla F(u_0)\|_{\nabla^2F(u_0)^{-1}}\leq1/100$, which implies that $u_0$ is near-central for $c_0:=-\widehat g_0$. Starting from $(u_0,c_0)$, preprocessing follows this auxiliary central path \emph{backward} toward the barrier center (the minimizer of $F$) by multiplying the auxiliary objective by $1-\gamma_t$ in each phase and recentering via the same procedure as the main loop (Lines~\ref{ln:tech_ov_phase_start}--\ref{ln:tech_ov_phase_end} in Algorithm~\ref{alg:overview_main_loop}). Because the auxiliary domain is bounded and contains an explicit ball around $u_0$, the initial objective $c_0$ is polynomially bounded (Lemma~\ref{lem:cost_bounded_metric}), so $N_{\rm aux}=\widetilde O(\sqrt m)$ phases suffice to shrink the objective until it is negligible in the dual local norm. Two final phases then transition from this tiny objective to the zero objective and from zero to $t_0 c$ for a polynomially small $t_0>0$, yielding a near-central starting point for the main loop. In total, preprocessing makes $\widetilde O(\sqrt m)$ queries at the same per-query cost as the main loop (Theorem~\ref{thm:preprocessing}).

\paragraph{Postprocessing: primal recovery by a modified Newton iteration.} From the starting point produced by preprocessing, the main loop (Algorithm~\ref{alg:overview_main_loop}) runs $\widetilde O(\sqrt m)$ phases on the auxiliary SDP to reach a near-central point $u=(y,z)$ within $\epsilon/64$ of the auxiliary optimum, giving the dual output $y$. To recover an explicit primal solution $X\succeq0$, let $\overline S(u)=\sum_{i=1}^{m+1}u_i\overline A_i-\overline C$, $J$ be an exact whitener of $\overline S(u)$, and $B_i:=J\overline A_iJ^\top$. For any matrix $M\succeq0$, the matrix $\overline X:=\frac{\alpha}{t}J^\top MJ\succeq 0$ satisfies $\langle\overline S(u),\overline X\rangle=(\alpha/t)\tr[M]$, and it is auxiliary primal feasible if and only if $\tr[B_iM]=tc_i/\alpha$ for all $i\in[m+1]$. Because $\nabla F(u)_i=-\alpha\tr[M_uB_i]\approx-tc_i$ by near-centrality and $\tr[M_u]=m+1$, the gradient matrix $M_u$ of~\eqref{eq:tech_ov_slack_gradient} nearly satisfies these conditions, and the corresponding $\overline{X}$ has complementarity $\nu/t=O(\epsilon)$.

Since $M_u$ involves the optimal weight, postprocessing approximates it by sketched factors. It computes an accurate weight at $u$, draws a fresh full estimate with new sketching matrices, and keeps its sketched factors $Y=[Y_0\ \cdots\ Y_{m+1}]$ (Section~\ref{sec:overview_sketch}). The matrix $M_{\rm init}:=YY^\top$ is positive semidefinite by construction, $-\alpha\tr[M_{\rm init}B_i]$ form a gradient estimate within $1/100$ of $\nabla F(u)$ in the dual local norm, and $\tr[M_{\rm init}]=O(m)$ (Lemma~\ref{lem:more_recovery_density}). Its residual is therefore small in a local norm, whereas the primal accuracy in~\eqref{eq:intro_output_accuracy} requires a residual that is even smaller by a factor polynomial in the input magnitudes. Since subtracting a correction from $M_{\rm init}$ could destroy positive semidefiniteness, we instead refine $M_{\rm init}$ using a congruence parametrization and apply a modified Newton iteration, starting from $w=0$:
\begin{equation}\label{eq:overview_recovery}
 M(w):=(I+B_w)M_{\rm init}(I+B_w),\qquad
 f(w)_i:=\tr[B_iM(w)]-\tfrac t\alpha c_i,\qquad
 w\gets w-\tfrac12\widehat Q^{-1}f(w),
\end{equation}
where $B_w:=\sum_{i=1}^{m+1}w_iB_i$. This parametrization ensures that $M(w)\succeq0$ throughout the iteration. The residual map $f$ is quadratic, and its Jacobian at the origin satisfies $0.9\,\partial f(0)\preceq2\widehat Q\preceq1.1\,\partial f(0)$. Lemma~\ref{lem:more_congruence_correction} shows that the iterates remain in a small neighborhood of the origin, where the update is a contraction and $\|f(w)\|_{\widehat Q^{-1}}$ decays geometrically. Thus $N_{\rm rec}=\widetilde O(1)$ iterations suffice to reduce the feasibility residual to the required inverse-polynomial accuracy.

Finally, let $X$ be the leading $n\times n$ block of $\overline X:=\frac{\alpha}{t}J^\top M(w)J$. Then $X\succeq0$, and by bounding the contributions of the auxiliary scalar blocks, Theorem~\ref{thm:postprocessing} shows that $X$ satisfies the primal accuracy guarantees of~\eqref{eq:intro_output_accuracy}. This entire postprocessing requires one width-$\widetilde O(m)$ sketch of cost $\widetilde O(n^\omega+\mathcal C(m))=\widetilde O(n^\omega+n^2m^{\omega-1})$ and $N_{\rm rec}=\widetilde O(1)$ updates of cost $O(n^\omega+mn^2)$ each. Hence, postprocessing takes $\widetilde O(n^\omega+n^2m^{\omega-1})\leq\widetilde O(n^\omega+n^2m^{3/2})$ time, since $\omega-1\leq5/4<3/2$.

\subsection{Organization}
\label{sec:organization}
Part~\ref{part:barrier} develops the matrix Lee--Sidford barrier, and Part~\ref{part:algorithm} gives the randomized algorithm. Each part opens with its main result and has its own preliminaries, and the shared notation appears once in Section~\ref{sec:barrier_preliminaries}.

In Part~\ref{part:barrier}, Section~\ref{sec:barrier_overview} states our first main result (Theorem~\ref{thm:ml_matrix_sc}). Section~\ref{sec:barrier_preliminaries} also reviews the self-concordance barrier theory and matrix calculus. Section~\ref{sec:matrix_lewis_barrier} constructs the barrier from inscribed ellipsoids. Section~\ref{sec:barrier_geometry} establishes the barrier properties and completes the proof of Theorem~\ref{thm:ml_matrix_sc}.

In Part~\ref{part:algorithm}, Section~\ref{sec:algorithm_part_overview} states our central algorithmic result of solving SDP with a near-central starting point. Section~\ref{sec:algorithm_preliminaries} gives the preliminaries for the efficient algorithms. Section~\ref{sec:algorithm_description} states the algorithm and the guarantees of its three data structures.  Section~\ref{sec:path_following} proves the correctness of path following from these guarantees. Sections~\ref{sec:randomized_oracles}--\ref{sec:reused_factorization} implement the data structures $\mathcal D^{\rm f}_{\rm est}$, $\mathcal D_{\rm weight}$, $\mathcal D^{\rm m}_{\rm est}$, and $\mathcal D_{\rm slack}$. Section~\ref{sec:randomized_solver} combines these costs and completes the proof of the main algorithmic result. Finally, Sections~\ref{sec:preprocessing} and~\ref{sec:postprocessing} develop the preprocessing and postprocessing algorithms, and Section~\ref{sec:solver_wrapup} combines them together and proves our second main result (Theorem~\ref{thm:more_complete_sdp}).

\clearpage
\part{The matrix Lee--Sidford barrier}
\label{part:barrier}
\section{The main geometric result}
\label{sec:barrier_overview}
We construct the matrix Lee--Sidford barrier, a self-concordant barrier for the dual feasible region $\Omega:=\{y\in\R^m:S(y)\succ0\}$ with parameter $O(m\log^2(2mn))$. This parameter is nearly linear in the number $m$ of dual variables, while the logarithmic barrier $-\log\det S(y)$ has parameter $n$. The construction extends the Lee--Sidford barrier for linear programs~\cite{ls19}. The value of the barrier at a point is determined by an optimal weight, the minimizer of a strictly convex problem over $m\times m$ positive definite matrices. Theorem~\ref{thm:ml_matrix_sc} also describes the geometry of the optimal weight, expresses the first two derivatives of the barrier through it, and bounds how fast it moves. Part~\ref{part:algorithm} evaluates the barrier through these formulas.

With $S(y)=\sum_iy_iA_i-C$ as in~\eqref{eq:intro_dual_sdp}, define at a point $y\in\Omega$
\begin{align*}
 B_i:=S(y)^{-1/2}A_iS(y)^{-1/2},\qquad B_h:=\sum_ih_iB_i,\qquad
 \mathcal L_y(G):=\sum_{i,j}G_{ij}B_iB_j
\end{align*}
for $h\in\R^m$ and $G\in\mathbb S^m$. Since $S(y\pm h)=S(y)^{1/2}(I\pm B_h)S(y)^{1/2}$, both $y+h$ and $y-h$ lie in $\Omega$ if and only if $\|B_h\|_{\rm op}<1$. We call $K_y:=\{h\in\R^m:\|B_h\|_{\rm op}\leq1\}$ the symmetrized local region. For $G\succ0$, the ellipsoid $\mathcal E_G:=\{G^{1/2}z:\|z\|_2\leq1\}$ lies in $K_y$ whenever $\mathcal L_y(G)\preceq I$ (Lemma~\ref{lem:ml_sos}). The matrix Lewis-weight potential is defined by maximizing the volume of $\mathcal E_G$, with a trace-power penalty in place of this condition. For an integer $s\geq2$, it is given by
\begin{equation}\label{eq:main_potential}
 \Phi_p(y):=\frac12\min_{G\succ0}\left\{\frac1s\tr[\mathcal L_y(G)^s]
 -\log\det G-m\right\},\qquad p:=2s.
\end{equation}
When $C,A_1,\ldots,A_m$ are diagonal, $\Phi_p$ is the Lee--Sidford potential of the corresponding linear program.

\begin{theorem}[Matrix Lee--Sidford barrier]
\label{thm:ml_matrix_sc}
Let $A_1,\ldots,A_m,C\in\mathbb S^n$, and assume that $\Omega=\{y\in\R^m:S(y)\succ0\}$ is nonempty, the $A_i$ are linearly independent, and $\bigcap_i\ker A_i=\{0\}$. Let
\begin{align*}
 s:=2\lceil\log(2mn)\rceil,\qquad \alpha:=10^6s^2,\qquad
 \nu:=\alpha m=O(m\log^2(2mn)),\qquad F:=\alpha\Phi_p.
\end{align*}
Then the following hold for every $y\in\Omega$ and $h\in\R^m$, with $B_i$, $B_h$, $\mathcal L_y$ and $K_y$ defined at $y$.
\begin{itemize}
\item[(a)] The minimum in~\eqref{eq:main_potential} is attained at a unique matrix $G_*=G_*(y)\succ0$, called the optimal weight, which depends smoothly on $y$. It satisfies
\begin{align*}
 \tr[\mathcal L_y(G_*)^s]=m,\qquad
 \Phi_p(y)=-\frac12\log\det G_*+\frac m{2s}-\frac m2.
\end{align*}
\item[(b)] The function $F$ is a $\nu$-self-concordant barrier for $\Omega$. Its Hessian is positive definite, $F$ diverges at every boundary point of $\Omega$, and
\begin{align*}
 |\nabla^3F(y)[h,h,h]|\leq2\bigl(\nabla^2F(y)[h,h]\bigr)^{3/2},\qquad
 |\nabla F(y)[h]|^2\leq\nu\,\nabla^2F(y)[h,h].
\end{align*}
\item[(c)] Let $W_y:=\mathcal L_y(G_*)^{s-1}$. The gradient matrix $M_y$ and the proxy Hessian $Q_y$, defined by
\begin{align*}
 M_y:=\frac12\Bigl(\mathcal L_y(G_*)^s
 +\sum_{i,j}(G_*)_{ij}B_iW_yB_j\Bigr),\qquad
 Q_y[h,h]:=\tr[M_yB_h^2],
\end{align*}
satisfy $M_y\succ0$, $\tr[M_y]=m$, and
\begin{align*}
 (\nabla F(y))_i=-\alpha\tr[M_yB_i]\quad(i\in[m]),\qquad
 \alpha Q_y\preceq\nabla^2F(y)\preceq\alpha(2s+1)Q_y.
\end{align*}
\item[(d)] The symmetrized local region satisfies
\begin{align*}
 m^{-1/(2s)}\mathcal E_{G_*}\subseteq K_y
 \subseteq\bigl(m^{1-1/s}n^{1/s}\bigr)^{1/2}\mathcal E_{G_*},
\end{align*}
and $\|B_h\|_{\rm op}^2\leq9\,Q_y[h,h]$.
\item[(e)] The derivative of the optimal weight satisfies
\begin{align*}
 \bigl\|G_*^{-1/2}\,\nabla G_*(y)[h]\,G_*^{-1/2}\bigr\|_F
 \leq\sqrt{6\,\nabla^2F(y)[h,h]/\alpha}.
\end{align*}
\end{itemize}
\end{theorem}

Part~(a) shows that the barrier is well defined. Up to an additive constant that depends only on $m$ and $s$, the potential $\Phi_p(y)$ is the negative log-volume of the ellipsoid $m^{-1/(2s)}\mathcal E_{G_*}$, which lies in $K_y$. Part~(b) is the main conclusion. The factor $m$ in $\nu$ comes from the trace identity of Part~(a) through $\tr[M_y]=m$. The scaling $\alpha=\Theta(s^2)$ reduces the constant in the third-derivative bound from order $s$ for $\Phi_p$ to at most $2$ for $F$. Part~(c) gives the gradient and a spectral approximation of the Hessian in terms of the value of $G_*$, whereas the Hessian itself involves the derivative of $G_*$. By Part~(d), the ellipsoid of the optimal weight approximates $K_y$ within the factor $(mn^{1/s})^{1/2}=O(\sqrt m)$, and the ellipsoid $\{h:Q_y[h,h]\leq1/9\}$ also lies in $K_y$. Finally, Part~(e) shows that the optimal weight moves slowly. Relative to $G_*$, its derivative in the direction $h$ is at most $\sqrt{6/\alpha}=O(1/s)$ times the local norm $(\nabla^2F(y)[h,h])^{1/2}$.

Section~\ref{sec:technical_overview} sketches the proof of Theorem~\ref{thm:ml_matrix_sc}, and Sections~\ref{sec:barrier_preliminaries}--\ref{sec:barrier_geometry} give a self-contained proof. Section~\ref{sec:barrier_preliminaries} collects notation, self-concordance, and matrix calculus. Section~\ref{sec:matrix_lewis_barrier} derives the potential from inscribed ellipsoids, explains why the last two assumptions of the theorem hold without loss of generality, and proves Part~(a). Section~\ref{sec:barrier_geometry} differentiates the potential along exponential curves and proves Parts~(c)--(e). It completes the proof with Part~(b) in Section~\ref{sec:ml_third_derivative}.
\section{Preliminaries for the barrier}
\label{sec:barrier_preliminaries}
Both parts use the notation and local-norm conventions below. This section also recalls the standard self-concordance estimates used in Part~\ref{part:algorithm} and the matrix calculus used in the barrier construction. %

\subsection{Notation}
For every positive integer $m$, define $[m]:=\{1,\ldots,m\}$. For a vector $x\in\R^d$, define its Euclidean norm by $\|x\|_2:=(\sum_{i=1}^d x_i^2)^{1/2}$. For a real matrix $A\in\R^{a\times b}$, define its operator norm and Frobenius norm by
\begin{align*}
 \|A\|_{\mathrm{op}}:=\sup_{\|x\|_2=1}\|Ax\|_2,
 \qquad
 \|A\|_F:=\sqrt{\tr[A^\top A]}
             =(\sum_{i=1}^a\sum_{j=1}^b A_{ij}^2)^{1/2}.
\end{align*}
For matrices of the same dimensions, define the trace inner product by $\langle A,B\rangle:=\tr[A^\top B]$, so that $\|A\|_F^2=\langle A,A\rangle$. Let $\mathbb S^n$ denote the real symmetric $n\times n$ matrices. For $A,B\in\mathbb S^n$, write $A\succeq B$ when $A-B$ is positive semidefinite and $A\succ B$ when $A-B$ is positive definite. Write $I$ for the identity matrix of the dimension determined by context. The symbol $\otimes$ denotes the Kronecker (tensor) product. For $A\in\R^{a\times b}$ and $B\in\R^{c\times d}$, the matrix $A\otimes B\in\R^{ac\times bd}$ has $(i,j)$ block $A_{ij}B$.

\subsection{Local norms and self-concordance}
\label{sec:sc_preliminaries}
For a positive definite matrix $H$, write $\|h\|_H:=(h^\top Hh)^{1/2}$ and $\|g\|_{H^{-1}}:=(g^\top H^{-1}g)^{1/2}$ for its primal and dual norms. For a smooth map $f$, $\nabla^kf(x)[h_1,\ldots,h_k]$ denotes its $k$th derivative at $x$ applied to the directions $h_1,\ldots,h_k$. For a scalar function, $\nabla f(x)$ and $\nabla^2f(x)$ are the gradient and the Hessian, so $\nabla f(x)[h]=\nabla f(x)^\top h$ and $\nabla^2f(x)[h,h]=h^\top\nabla^2f(x)h$. For a vector-valued map, $\nabla f(x)$ is the Jacobian matrix. When the function has several arguments, a subscript names the variable of differentiation, as in $\nabla_yf_y(G)[h]$. For a quantity $X$ that depends on a scalar parameter $t$, such as a matrix along a curve, $\partial_tX$ denotes its derivative with respect to $t$. For a function of one variable written with its argument, we use primes, as in $y'(t)$.

We use the standard definition of a self-concordant barrier~\cite{nn94}.

\begin{definition}[Self-concordant barrier]
\label{def:sc_barrier}
A convex $C^3$ function $F$ on an open convex set $\Omega\subseteq\R^m$ is a $\nu$-self-concordant barrier if its Hessian is positive definite, $F$ diverges at every boundary point of $\Omega$, and
\begin{align*}
 |\nabla^3F(y)[h,h,h]|\leq2\bigl(\nabla^2F(y)[h,h]\bigr)^{3/2},\qquad
 |\nabla F(y)[h]|^2\leq\nu \nabla^2F(y)[h,h]
\end{align*}
for every $y\in\Omega$ and $h\in\R^m$.
\end{definition}
The definition implies the following standard estimates, whose proofs can be found, for example, in~\cite[Chapter~5]{nesterov18}.

\begin{lemma}[Local self-concordance estimates]
\label{lem:sc_local_estimates}
Let $F$ be a self-concordant barrier on $\Omega$, and let $y\in\Omega$. If $d:=\|z-y\|_{\nabla^2F(y)}<1$, then $z\in\Omega$ and
\begin{equation}\label{eq:sc_hessian_comparison}
 (1-d)^2\nabla^2F(y)\preceq\nabla^2F(z)\preceq(1-d)^{-2}\nabla^2F(y).
\end{equation}
For $c\in\R^m$, let $f(x):=F(x)+c^\top x$ and $\lambda:=\|\nabla f(y)\|_{\nabla^2F(y)^{-1}}$. If $\lambda<1$, then $f$ has a unique minimizer $y_c$, and
\begin{equation}\label{eq:sc_center_estimates}
 \|y-y_c\|_{\nabla^2F(y)}\leq\frac\lambda{1-\lambda},\qquad
 \lambda-\log(1+\lambda)
 \leq f(y)-f(y_c)\leq-\lambda-\log(1-\lambda).
\end{equation}
If $f$ has a minimizer $y_c$, then the lower bound in~\eqref{eq:sc_center_estimates} holds for every $\lambda\geq0$.
\end{lemma}

\subsection{Matrix powers and trace calculus}
\label{sec:matrix_calculus_preliminaries}
For a real symmetric matrix $D=U\operatorname{diag}(\lambda_i)U^\top$ and a scalar function $f$ defined on its spectrum, write $f(D):=U\operatorname{diag}(f(\lambda_i))U^\top$. In particular, positive definite matrices have well-defined real powers, logarithms, and positive square roots. For symmetric $V$, $e^V\succ0$ and $\log\det(e^V)=\tr[V]$.

We use the cyclicity of the trace and the Cauchy--Schwarz inequality for the trace inner product,
\begin{align*}
 \tr[AB]&=\tr[BA],\qquad |\tr[X^\top Y]|\leq\|X\|_F\|Y\|_F.
\end{align*}
The first identity holds whenever both products are square. Since congruence preserves the semidefinite order, $A\preceq B$ and $C\succeq0$ imply $\tr[CA]\leq\tr[CB]$.

For a differentiable positive definite curve $D(t)$ and a positive integer $s$, differentiation gives
\begin{align*}
 \partial_t\log\det D&=\tr[D^{-1}\partial_tD],\qquad \partial_t(D^{-1})=-D^{-1}(\partial_tD)D^{-1},\\
 \partial_t(D^s)&=\sum_{a=0}^{s-1}D^a(\partial_tD)D^{s-1-a},\qquad \partial_t\tr[D^s]=s\tr[D^{s-1}\partial_tD].
\end{align*}
The $s$ terms in the derivative of $D^s$ are different matrices in general, and they have equal traces by cyclicity.
\section{Construction of the matrix Lee--Sidford barrier}
\label{sec:matrix_lewis_barrier}

We construct the matrix Lee--Sidford barrier by fitting an ellipsoid to the local feasible region. For a linear program, an ellipsoid lies in this region if and only if its widths along finitely many constraint normals are at most one. The Lee--Sidford potential penalizes these widths by a power sum. For a semidefinite program, there is one such constraint for each unit vector. A matrix sum of squares bounds all of them at once, and a trace power of this matrix gives the penalty.

\subsection{Local ellipsoids for linear programs}
\label{sec:ml_lp_ellipsoids}
Consider a linear program with $n$ inequality constraints on $y\in\R^m$. At a strictly feasible $y$, let $r_k\in\R^m$ be the $k$th constraint normal divided by its slack. The symmetrized local region is $K_y:=\{h:|r_k^\top h|\leq1\text{ for all }k\}$. Both $y+h$ and $y-h$ are strictly feasible if and only if $h$ lies in the interior of $K_y$.

For $G\succeq0$, the ellipsoid $\mathcal E_G:=\{G^{1/2}z:\|z\|_2\leq1\}$ has volume proportional to $\sqrt{\det G}$. Its squared width along a vector $r\in\R^m$ is
\begin{equation}\label{eq:ml_width}
 \max_{h\in\mathcal E_G}(r^\top h)^2=\max_{\|z\|_2\leq1}(z^\top G^{1/2}r)^2=r^\top Gr.
\end{equation}
Hence $\mathcal E_G\subseteq K_y$ if and only if $r_k^\top Gr_k\leq1$ for every $k$. For an integer $s\geq2$, replacing these hard constraints by a power-sum penalty gives the Lee--Sidford potential~\cite{ls19} in the form
\begin{equation}\label{eq:scalar_shape_potential}
 \Phi^{\rm LP}(y):=\frac12\min_{G\succ0}
 \left\{\frac1s\sum_k(r_k^\top G r_k)^s-\log\det G-m\right\}.
\end{equation}
When the $r_k$ span $\R^m$, the minimizer $G_*$ exists, and the numbers $(r_k^\top G_*r_k)^s$ are the $\ell_{2s}$ Lewis weights of the matrix with rows $r_k^\top$~\cite{lewis78,cp15}. The matrix Lewis-weight potential defined below replaces the power sum by a trace power, and the two potentials coincide when the SDP data $C,A_1,\ldots,A_m$ are diagonal.

\subsection{A sum-of-squares constraint for semidefinite programs}
\label{sec:ml_assumptions}
We now turn to the SDP~\eqref{eq:primal_sdp}--\eqref{eq:intro_dual_sdp} and its affine pencil $y\mapsto S(y)=\mathcal A^*(y)-C$, where $\mathcal A^*(h):=\sum_ih_iA_i$. Throughout the barrier analysis, we assume
\begin{equation}\label{eq:ml_assumptions}
 \Omega:=\{y:S(y)\succ0\}\ne\varnothing,\qquad
 A_1,\ldots,A_m\text{ are linearly independent},\qquad
 \bigcap_i\ker A_i=\{0\}.
\end{equation}
The last two conditions hold without loss of generality. If the $A_i$ are linearly dependent, then $\mathcal A^*$ has a nonzero kernel. The slack $S(y)$ is constant along this kernel, and so is $b^\top y$ if the primal problem is feasible. Thus we may restrict $y$ to a complementary subspace. If the $A_i$ have a common nonzero kernel, then $S(y)v=-Cv$ for every $v$ in this kernel, so the block of $S(y)$ on this kernel is constant. It is also positive definite because $\Omega\neq\varnothing$. Replacing $S(y)$ by the Schur complement of this block gives a pencil of smaller order with the same domain $\Omega$.

Fix $y\in\Omega$. The normalized constraint matrices are
\begin{equation}\label{eq:ml_whitened_coefficients}
 B_i:=S(y)^{-1/2}A_iS(y)^{-1/2},\qquad B_h:=\sum_i h_iB_i.
\end{equation}
For every $h\in\R^m$,
\begin{align*}
 S(y\pm h)=S(y)^{1/2}(I\pm B_h)S(y)^{1/2},
\end{align*}
so both $y+h$ and $y-h$ lie in $\Omega$ if and only if $\|B_h\|_{\rm op}<1$. We therefore define the symmetrized local region $K_y:=\{h:\|B_h\|_{\rm op}\leq1\}$. We call the pencil diagonal if $C$ and all $A_i$ are diagonal. The SDP is then a linear program, $(B_i)_{kk}=(r_k)_i$, and $K_y$ is the region of Section~\ref{sec:ml_lp_ellipsoids}.

As in the linear case, $K_y$ is an intersection of slabs. For a unit vector $u\in\R^n$, let $r(u)\in\R^m$ have entries $r(u)_i:=u^\top B_iu$. Since $B_h$ is symmetric,
\begin{align*}
 \|B_h\|_{\rm op}=\max_{\|u\|_2=1}|u^\top B_hu|=\max_{\|u\|_2=1}|r(u)^\top h|.
\end{align*}
Consequently, $K_y$ is the intersection of the slabs $\{h:|r(u)^\top h|\leq1\}$ over all unit vectors $u$. By~\eqref{eq:ml_width}, the ellipsoid $\mathcal E_G$ lies in $K_y$ if and only if
\begin{align*}
 \sup_{\|u\|_2=1}r(u)^\top Gr(u)\leq1.
\end{align*}
The left side is the maximum of a quartic polynomial over the unit sphere, and the constraints form an infinite family. Instead, we use a single semidefinite inequality that implies this condition. Fix $G\succeq0$. Every $h\in\mathcal E_G$ has the form $h=G^{1/2}z$ with $\|z\|_2\leq1$, and then
\begin{align*}
 B_h=\sum_az_aV_a,\qquad V_a:=\sum_i(G^{1/2})_{ia}B_i\quad(a\in[m]).
\end{align*}
The next lemma bounds the square of this matrix by the linear map
\begin{equation}\label{eq:ml_map}
 \mathcal L_y(G):=\sum_{i,j}G_{ij}B_iB_j,\qquad G\in\mathbb S^m.
\end{equation}

\begin{lemma}[Sum-of-squares bound]
\label{lem:ml_sos}
Let $y\in\Omega$ and $G\succeq0$, and let $V_a$ be as above. Then $\mathcal L_y(G)=\sum_aV_a^2\succeq0$, and every $z\in\R^m$ satisfies
\begin{align*}
 \left(\sum_a z_aV_a\right)^2\preceq\|z\|_2^2\mathcal L_y(G).
\end{align*}
In particular, $\mathcal L_y(G)\preceq I$ implies $\mathcal E_G\subseteq K_y$. If $G\succ0$, then $\mathcal L_y(G)\succ0$.
\end{lemma}
\begin{proof}
We first verify the sum-of-squares identity and then derive the matrix inequality from the Cauchy--Schwarz inequality. The last two claims follow from these two facts.

Since $G^{1/2}$ is symmetric and $G^{1/2}G^{1/2}=G$,
\begin{align*}
 \sum_aV_a^2=\sum_{i,j}\Bigl(\sum_a(G^{1/2})_{ia}(G^{1/2})_{aj}\Bigr)B_iB_j=\sum_{i,j}G_{ij}B_iB_j=\mathcal L_y(G).
\end{align*}
For every vector $u$, the triangle and Cauchy--Schwarz inequalities give
\begin{align*}
 \left\|\sum_a z_aV_au\right\|_2^2
 \leq\left(\sum_a z_a^2\right)\sum_a\|V_au\|_2^2
 =\|z\|_2^2u^\top\mathcal L_y(G)u,
\end{align*}
which is the stated matrix inequality.

If $\mathcal L_y(G)\preceq I$ and $h=G^{1/2}z$ with $\|z\|_2\leq1$, then $B_h=\sum_az_aV_a$, and the matrix inequality gives $B_h^2\preceq\mathcal L_y(G)\preceq I$, so $h\in K_y$.

Now let $G\succ0$ and $\mathcal L_y(G)u=0$. Then
\begin{align*}
 \sum_a\|V_au\|_2^2=u^\top\mathcal L_y(G)u=0,
\end{align*}
so $V_au=0$ for every $a$. Since $G^{1/2}$ is invertible, each $B_i$ is a linear combination of the $V_a$. Therefore $B_iu=0$, that is, $A_iS(y)^{-1/2}u=0$, for every $i$. The assumption $\bigcap_i\ker A_i=\{0\}$ gives $S(y)^{-1/2}u=0$, and so $u=0$.
\end{proof}

\begin{remark}[Tightness of the sum-of-squares condition]
\label{rem:ml_sos_tightness}
For diagonal pencils, the condition $\mathcal L_y(G)\preceq I$ is equivalent to $\mathcal E_G\subseteq K_y$. In this case $(B_i)_{kk}=(r_k)_i$ gives
\begin{align*}
 \mathcal L_y(G)=\operatorname{diag}(r_k^\top G r_k),
\end{align*}
so $\mathcal L_y(G)\preceq I$ if and only if $r_k^\top Gr_k\leq1$ for all $k$. For general coefficients, the condition can be stricter. For example, let $m=n=2$ and consider the normalized constraint matrices
\begin{align*}
 B_1:=\begin{pmatrix}1&0\\0&-1\end{pmatrix},\qquad
 B_2:=\begin{pmatrix}0&1\\1&0\end{pmatrix}.
\end{align*}
Since $B_1^2=B_2^2=I$ and $B_1B_2=-B_2B_1$,
\begin{align*}
 (h_1B_1+h_2B_2)^2=\|h\|_2^2I,\qquad \mathcal L_y(I)=B_1^2+B_2^2=2I.
\end{align*}
It follows that $K_y=\mathcal E_I$ is the unit disk, while $\mathcal L_y(cI)=2cI\preceq I$ holds if and only if $c\leq1/2$. We nevertheless use $\mathcal L_y(G)\preceq I$ in place of exact containment, because it is a single semidefinite inequality and it is exact for diagonal pencils.
\end{remark}

\subsection{The potential and the optimal weight}
\label{sec:ml_variational_construction}
As for linear programs, we replace the hard constraint $\mathcal L_y(G)\preceq I$ by a penalty. For $G\succeq0$, the eigenvalues of $\mathcal L_y(G)$ give
\begin{align*}
 \|\mathcal L_y(G)\|_{\rm op}\leq\tr[\mathcal L_y(G)^s]^{1/s}\leq n^{1/s}\|\mathcal L_y(G)\|_{\rm op}.
\end{align*}
Thus, for $s\geq\log n$, the trace power $\tr[\mathcal L_y(G)^s]^{1/s}$ approximates $\|\mathcal L_y(G)\|_{\rm op}$ within a factor $e$.

\begin{definition}[Matrix Lewis-weight potential]
\label{def:ml_potential}
Let $s\geq2$ be an integer, and let $p:=2s$. For $y\in\Omega$, with $B_i$ and $\mathcal L_y$ as in~\eqref{eq:ml_whitened_coefficients} and~\eqref{eq:ml_map}, define
\begin{equation}\label{eq:ml_dual_weight}
 f_y(G):=\frac1s\tr[\mathcal L_y(G)^s]-\log\det G-m,
 \qquad
 \Phi_p(y):=\frac12\min_{G\succ0} f_y(G).
\end{equation}
By Lemma~\ref{lem:ml_optimizer} below, $f_y$ has a unique minimizer $G_*(y)$, which we call the \emph{optimal weight}. The adjoint of $\mathcal L_y$ with respect to the trace inner product is the map $\mathcal L_y^*:\mathbb S^n\to\mathbb S^m$ with
\begin{align*}
 (\mathcal L_y^*(W))_{ij}:=\tfrac12\tr[W(B_iB_j+B_jB_i)],
 \qquad h^\top \mathcal L_y^*(W)h=\tr[WB_h^2].
\end{align*}
\end{definition}

Besides existence and uniqueness, the next lemma shows that $G_*(y)$ depends smoothly on $y$, which is needed to differentiate $\Phi_p$. It also gives the stationarity conditions and a dual formulation of $\Phi_p$. We call the trace identity $\tr[\mathcal L_y(G_*)^s]=m$ in~\eqref{eq:ml_mass} the mass identity. It is the source of the factor $m$ in the barrier parameter.

\begin{lemma}[Existence and stationarity]
\label{lem:ml_optimizer}\label{lem:ml_dual_weight}
Let $y\in\Omega$. The function $f_y$ has a unique minimizer $G_*=G_*(y)$ over $G\succ0$, and $G_*(y)$ depends smoothly on $y$. Let
\begin{equation}\label{eq:ml_stationarity}
 W_y:=\mathcal L_y(G_*)^{s-1}.
\end{equation}
Then
\begin{equation}\label{eq:ml_mass}
 \mathcal L_y^*(W_y)=G_*^{-1},\qquad \tr[\mathcal L_y(G_*)^s]=\tr\bigl[W_y^{s/(s-1)}\bigr]=m.
\end{equation}
Moreover,
\begin{equation}\label{eq:ml_weight_maximum}
 \Phi_p(y)=\frac12\max_{W\succeq0}
 \Bigl\{\log\det \mathcal L_y^*(W)-\frac{s-1}s\tr\bigl[W^{s/(s-1)}\bigr]\Bigr\},
\end{equation}
where the objective is $-\infty$ when $\mathcal L_y^*(W)$ is singular, and $W_y$ is the unique maximizer.
\end{lemma}
\begin{proof}
We first show that $f_y$ is strictly convex and that $f_y(G)\to\infty$ when $G$ grows without bound or approaches a singular matrix. These two facts give a unique minimizer, and the implicit function theorem then gives its smooth dependence on $y$. For the dual formulation, we bound the objective in~\eqref{eq:ml_weight_maximum} by $f_y$ using the inequality $\log x\leq x-1$ and Young's inequality, and we check that both are tight at $(G_*,W_y)$.

\emph{Convexity.} For $P\succ0$ and symmetric $Z$, differentiation gives
\begin{equation}
 \nabla_P^2\left(\frac1s\tr[P^s]\right)[Z,Z]
 =\sum_{a=0}^{s-2}\tr[P^aZP^{s-2-a}Z]
 =\sum_{a=0}^{s-2}\|P^{a/2}ZP^{(s-2-a)/2}\|_F^2.\label{eq:ml_power_energy}
\end{equation}
For $G\succ0$, the matrix $\mathcal L_y(G)$ is positive definite by Lemma~\ref{lem:ml_sos}. Applying~\eqref{eq:ml_power_energy} with $P=\mathcal L_y(G)$ and $Z=\mathcal L_y(V)$, and adding the Hessian of $-\log\det G$, gives for every $V\in\mathbb S^m$
\begin{equation}\label{eq:ml_weight_hessian_energy}
 \nabla^2f_y(G)[V,V]
 =\|G^{-1/2}VG^{-1/2}\|_F^2
  +\sum_{a=0}^{s-2}
    \|\mathcal L_y(G)^{a/2}\mathcal L_y(V)\mathcal L_y(G)^{(s-2-a)/2}\|_F^2.
\end{equation}
The first term is positive for $V\neq0$, so $f_y$ is strictly convex.

\emph{Existence and uniqueness.} Since the $A_i$ are linearly independent, so are the $B_i$, and their Gram matrix $\mathcal L_y^*(I)=(\tr[B_iB_j])_{ij}$ is positive definite. Let $c_y:=\lambda_{\min}(\mathcal L_y^*(I))>0$. Since
\begin{align*}
 \tr[\mathcal L_y(G)]=\tr[G\mathcal L_y^*(I)]\geq c_y\tr[G],
\end{align*}
H\"older's inequality gives
\begin{align*}
 \tr[\mathcal L_y(G)^s]
 \geq n^{1-s}\bigl(\tr[\mathcal L_y(G)]\bigr)^s
 \geq n^{1-s}(c_y\tr[G])^s.
\end{align*}
By the AM--GM inequality,
\begin{align*}
 \log\det G\leq m\log(\tr[G]/m).
\end{align*}
Together, these bounds give $f_y(G)\to\infty$ as $\tr[G]\to\infty$. If $\tr[G]$ stays bounded and $G$ approaches a singular matrix, then $f_y(G)\to\infty$ as well, because the trace power is nonnegative and $-\log\det G\to\infty$. Therefore every sublevel set of $f_y$ is a compact subset of the positive definite cone, and $f_y$ attains its minimum. By strict convexity, the minimizer is unique.

\emph{Stationarity and smoothness.} The gradient of $f_y$ is
\begin{align*}
 \nabla f_y(G)=\mathcal L_y^*(\mathcal L_y(G)^{s-1})-G^{-1},
\end{align*}
so the condition $\nabla f_y(G_*)=0$ is $\mathcal L_y^*(W_y)=G_*^{-1}$. Taking the trace inner product with $G_*$ and using $W_y=\mathcal L_y(G_*)^{s-1}$, we obtain
\begin{align*}
 m=\tr[G_*\mathcal L_y^*(W_y)]=\tr[W_y\mathcal L_y(G_*)]=\tr[\mathcal L_y(G_*)^s]=\tr\bigl[W_y^{s/(s-1)}\bigr].
\end{align*}
This proves~\eqref{eq:ml_mass}. The $B_i$ depend smoothly on $y\in\Omega$, and so does $\nabla f_y(G)$. Moreover, the Hessian~\eqref{eq:ml_weight_hessian_energy} at $G_*$ is invertible. The implicit function theorem applied to $\nabla f_y(G)=0$ shows that $G_*(y)$ is smooth.

\emph{Dual formulation.} Let $G\succ0$ and $W\succeq0$. If $\mathcal L_y^*(W)$ is singular, the objective in~\eqref{eq:ml_weight_maximum} is $-\infty$. Otherwise, we use the two inequalities
\begin{align*}
 \log\det \mathcal L_y^*(W)&\leq\tr[G\mathcal L_y^*(W)]-\log\det G-m,\\
 \tr[W\mathcal L_y(G)]&\leq\frac1s\tr[\mathcal L_y(G)^s]+\frac{s-1}s\tr\bigl[W^{s/(s-1)}\bigr].
\end{align*}
For the first, let $\lambda_1,\ldots,\lambda_m>0$ be the eigenvalues of $G^{1/2}\mathcal L_y^*(W)G^{1/2}$. Since $\log x\leq x-1$,
\begin{align*}
 \log\det\mathcal L_y^*(W)+\log\det G
 =\sum_{k=1}^m\log\lambda_k
 \leq\sum_{k=1}^m(\lambda_k-1)
 =\tr[G\mathcal L_y^*(W)]-m,
\end{align*}
with equality if $\mathcal L_y^*(W)=G^{-1}$. For the second, write
\begin{align*}
 \mathcal L_y(G)=\sum_id_iu_iu_i^\top,\qquad W=\sum_jw_jv_jv_j^\top
\end{align*}
in orthonormal eigenbases, and let $\pi_{ij}:=(u_i^\top v_j)^2$. The matrix $(\pi_{ij})$ is doubly stochastic, so Young's inequality $d_iw_j\leq d_i^s/s+(s-1)w_j^{s/(s-1)}/s$ gives
\begin{align*}
 \tr[W\mathcal L_y(G)]=\sum_{i,j}\pi_{ij}d_iw_j
 \leq\sum_{i,j}\pi_{ij}\left(\frac{d_i^s}{s}+\frac{s-1}s w_j^{s/(s-1)}\right)
 =\frac1s\tr[\mathcal L_y(G)^s]+\frac{s-1}s\tr\bigl[W^{s/(s-1)}\bigr].
\end{align*}
Since $\tr[G\mathcal L_y^*(W)]=\tr[W\mathcal L_y(G)]$, adding the two inequalities gives
\begin{align*}
 \log\det\mathcal L_y^*(W)-\frac{s-1}s\tr\bigl[W^{s/(s-1)}\bigr]
 \leq\frac1s\tr[\mathcal L_y(G)^s]-\log\det G-m=f_y(G).
\end{align*}
Both inequalities are equalities at $(G,W)=(G_*,W_y)$. The first is an equality because $\mathcal L_y^*(W_y)=G_*^{-1}$. The second is an equality because $W_y\mathcal L_y(G_*)=\mathcal L_y(G_*)^s=W_y^{s/(s-1)}$, so both of its sides equal $\tr[\mathcal L_y(G_*)^s]$. Consequently, the maximum in~\eqref{eq:ml_weight_maximum} equals $f_y(G_*)=2\Phi_p(y)$ and is attained at $W_y$.

\emph{Uniqueness of the maximizer.} Let $W$ be any maximizer. Its objective value equals $f_y(G_*)$, so both inequalities hold with equality for $G=G_*$. In particular, $d_iw_j=d_i^s/s+(s-1)w_j^{s/(s-1)}/s$ whenever $\pi_{ij}>0$, and equality in Young's inequality forces $w_j=d_i^{s-1}$ for each such pair. Hence each eigenvector $v_j$ of $W$ is a linear combination of eigenvectors $u_i$ of $\mathcal L_y(G_*)$ with $d_i^{s-1}=w_j$. Thus $\mathcal L_y(G_*)^{s-1}v_j=w_jv_j$ for every $j$, and $W=\mathcal L_y(G_*)^{s-1}=W_y$.
\end{proof}

The next corollary returns to the geometry of Section~\ref{sec:ml_lp_ellipsoids}. Up to an additive constant, $\Phi_p(y)$ is the negative log-volume of an ellipsoid that lies in $K_y$, as for linear programs.

\begin{corollary}[Volume interpretation]
\label{cor:ml_volume}
Let $y\in\Omega$. Then
\begin{align*}
 \Phi_p(y)=-\tfrac12\log\det G_*(y)+\tfrac{m}{2s}-\tfrac m2,
 \qquad m^{-1/(2s)}\mathcal E_{G_*(y)}\subseteq K_y.
\end{align*}
\end{corollary}
\begin{proof}
Both claims follow from the mass identity. Since $\Phi_p(y)=\frac12f_y(G_*)$ and $\frac1s\tr[\mathcal L_y(G_*)^s]=\frac ms$, we obtain the formula for $\Phi_p(y)$. For the containment, the mass identity gives
\begin{align*}
 \|\mathcal L_y(G_*)\|_{\rm op}\leq\tr[\mathcal L_y(G_*)^s]^{1/s}=m^{1/s},\qquad\text{so}\qquad
 \mathcal L_y(m^{-1/s}G_*)=m^{-1/s}\mathcal L_y(G_*)\preceq I.
\end{align*}
Since $\mathcal E_{m^{-1/s}G_*}=m^{-1/(2s)}\mathcal E_{G_*}$, Lemma~\ref{lem:ml_sos} applied to $m^{-1/s}G_*$ gives $m^{-1/(2s)}\mathcal E_{G_*}\subseteq K_y$.
\end{proof}
\section{Self-concordance of the matrix Lee--Sidford barrier}
\label{sec:barrier_geometry}

This section completes the proof of Theorem~\ref{thm:ml_matrix_sc}. Since $\Phi_p$ is defined by a minimization, its derivatives depend on how the optimal weight $G_*$ moves with the slack. We compute them along exponential curves in the slack and in the weight. Along such curves, the normalized constraint matrices are multiplied by matrix exponentials and $\log\det G$ is linear in time, so the objective becomes an explicit trace power plus a linear term. A single curvature bound for such trace powers then controls the Hessian of the potential, the derivative of the optimal weight, and the third derivative required for self-concordance.

We fix the trace-power exponent and the scaling of Theorem~\ref{thm:ml_matrix_sc} for the rest of the paper. These choices make the third-derivative bound hold while keeping the barrier parameter nearly linear in $m$.

\begin{definition}[Barrier parameters]\label{def:ml_s_choice}
For the remainder of the paper, let
\begin{align*}
 s:=2\lceil\log(2mn)\rceil,\qquad \alpha:=10^6s^2,\qquad F:=\alpha\Phi_p,\qquad \nu:=\alpha m,
\end{align*}
with $p=2s$ as in Definition~\ref{def:ml_potential}. When the barrier is applied to another pencil, $m$ and $n$ denote the number of variables and the matrix order of that pencil.
\end{definition}

\subsection{Normalized coordinates and exponential curves}
\label{sec:ml_normalization}
Throughout this section, the coefficient matrices $A_1,\ldots,A_m$ are fixed, and we regard the potential as a function of the slack. For $S\succ0$, let $\mathcal L_S$ and $f_S$ be the map~\eqref{eq:ml_map} and the objective~\eqref{eq:ml_dual_weight} with the normalized constraint matrices $S^{-1/2}A_iS^{-1/2}$ in place of $B_i$, and let
\begin{align*}
 \varphi(S):=\frac12\min_{G\succ0}f_S(G).
\end{align*}
Then $\Phi_p(y)=\varphi(S(y))$. The proof of Lemma~\ref{lem:ml_optimizer} applies to every $S\succ0$ under~\eqref{eq:ml_assumptions}, so $f_S$ has a unique minimizer $G_*(S)$, which depends smoothly on $S$. As $S(y)$ is affine in $y$, the derivatives of $\Phi_p$ are derivatives of $\varphi$ along straight lines. The extension to all $S\succ0$ allows us to differentiate $\varphi$ along other curves as well.

Two changes of coordinates simplify the calculation at a base point $S$. First, let $P$ be an invertible matrix. The congruence $S\mapsto PSP^\top$, $A_i\mapsto PA_iP^\top$ replaces each normalized constraint matrix $S^{-1/2}A_iS^{-1/2}$ by $U(S^{-1/2}A_iS^{-1/2})U^\top$, where
\begin{align*}
 U:=(PSP^\top)^{-1/2}PS^{1/2},\qquad UU^\top=I.
\end{align*}
This replaces $\mathcal L_S(G)$ by $U\mathcal L_S(G)U^\top$, so $\varphi_{\rm new}(PSP^\top)=\varphi(S)$ and $G_{*,\rm new}(PSP^\top)=G_*(S)$. Second, let $\Gamma\in\R^{m\times m}$ be invertible. The coefficient change $A_i\mapsto\sum_j\Gamma_{ji}A_j$ replaces $\mathcal L_S(G)$ by $\mathcal L_S(\Gamma G\Gamma^\top)$. Substituting $\Gamma G\Gamma^\top$ for $G$ in the minimum gives
\begin{align*}
 \varphi_{\rm new}(S)=\varphi(S)+\log|\det\Gamma|,\qquad G_{*,\rm new}(S)=\Gamma^{-1}G_*(S)\Gamma^{-\top}.
\end{align*}
Both changes are independent of the point of differentiation, so derivatives of $\varphi$ can be computed in the new coordinates. Taking $P:=S^{-1/2}$ and then $\Gamma:=G_*(I)^{1/2}$, we may assume that the base point satisfies
\begin{equation}\label{eq:ml_fixed_normalization}
\begin{aligned}
 &S=I,\qquad G_*=\mathcal L_I^*(W)=I,\\
 &\mathsf B:=[B_1\ \cdots\ B_m],\qquad
 W:=(\mathsf B\mathsf B^\top)^{s-1},
\end{aligned}
\end{equation}
where $B_i=A_i$ are the coefficients in these coordinates, so that $\mathsf B\mathsf B^\top=\sum_iB_i^2=\mathcal L_I(I)$. The identity $\mathcal L_I^*(W)=I$ is the stationarity condition of Lemma~\ref{lem:ml_optimizer}.

For $T\in\mathbb S^n$ and $V\in\mathbb S^m$, consider the curve $e^{tT}$ in the slack and the curve $e^{tV}$ in the weight. Both pass through the base point at $t=0$, and
\begin{align*}
 \mathcal L_{e^{tT}}(e^{tV})=\mathsf B_{T,V}(t)\mathsf B_{T,V}(t)^\top,
\end{align*}
where
\begin{equation}\label{eq:ml_exponential_factor}
 \mathsf B_{T,V}(t)
 :=e^{-tT/2}\mathsf B
   \exp\left(\frac t2(V\otimes I_n-I_m\otimes T)\right).
\end{equation}
The two terms in the exponent commute, so the exponential equals $e^{tV/2}\otimes e^{-tT/2}$. The $a$th block of $\mathsf B_{T,V}(t)$ is
\begin{align*}
 \sum_i(e^{tV/2})_{ia}\,e^{-tT/2}B_ie^{-tT/2},
\end{align*}
and the factorization is the sum-of-squares identity of Lemma~\ref{lem:ml_sos} for the weight $e^{tV}$ and the normalized constraint matrices $e^{-tT/2}B_ie^{-tT/2}$. Since $\log\det e^{tV}=t\tr[V]$, the objective along these curves is
\begin{equation}\label{eq:ml_orbit_objective}
 f_{e^{tT}}(e^{tV})=u_{T,V}(t)-t\tr[V]-m,\qquad
 u_{T,V}(t):=\frac1s\tr\bigl[(\mathsf B_{T,V}(t)\mathsf B_{T,V}(t)^\top)^s\bigr].
\end{equation}

\subsection{Trace powers along exponential curves}
\label{sec:ml_exponential_curvature}
By~\eqref{eq:ml_exponential_factor}, $u_{T,V}$ is a trace power of a matrix of the form $e^{t\mathsf L}\mathsf Be^{t\mathsf R}$ with symmetric $\mathsf L$ and $\mathsf R$. This subsection bounds the second and third derivatives of such trace powers. We first treat a symmetric matrix that moves along $e^{t\mathsf J/2}Ze^{t\mathsf J/2}$, and then deduce the rectangular case from a symmetric dilation.

\begin{lemma}[Trace powers along a symmetric exponential curve]
\label{lem:ml_symmetric_curvature}
Let $Z$ and $\mathsf J$ be real symmetric matrices of the same order, let $s\geq1$ be an integer, and define
\begin{align*}
 Z(t):=e^{t\mathsf J/2}Ze^{t\mathsf J/2},\qquad Y:=\tfrac12(\mathsf JZ+Z\mathsf J),\qquad
 g(t):=\frac1{2s}\tr[Z(t)^{2s}].
\end{align*}
Then
\begin{equation}\label{eq:orbit_dilation_second}
 g''(0)=2\sum_{a=0}^{s-1}\|Z^aYZ^{s-1-a}\|_F^2\leq2s\tr[Z^{2s}\mathsf J^2]
\end{equation}
and
\begin{equation}\label{eq:orbit_dilation_third}
 |g'''(0)|\leq20s\|\mathsf J\|_{\rm op}\,g''(0).
\end{equation}
\end{lemma}
\begin{proof}
We first expand $g''(0)$ in an eigenbasis of $Z$. A polynomial factorization turns it into the sum of squares in~\eqref{eq:orbit_dilation_second}, and the weighted AM--GM inequality gives the upper bound. For the third derivative, we differentiate this sum of squares and bound every resulting term by the squares themselves, using the cyclicity of the trace and the Cauchy--Schwarz inequality.

Let
\begin{equation}\label{eq:orbit_square_energies}
 \sigma_a:=\|Z^aYZ^{s-1-a}\|_F^2=\tr[Z^{2a}YZ^{2s-2-2a}Y]\quad(0\leq a\leq s-1),
 \qquad \mathcal S:=\sum_{a=0}^{s-1}\sigma_a.
\end{equation}

\emph{The second derivative.} Since $Z'(t)=\frac12(\mathsf JZ(t)+Z(t)\mathsf J)$, the cyclicity of the trace gives $g'(t)=\tr[\mathsf JZ(t)^{2s}]$ for every $t$. Differentiating once more,
\begin{align*}
 g''(0)=\sum_{k=0}^{2s-1}\tr[\mathsf JZ^kYZ^{2s-1-k}].
\end{align*}
We work in an orthonormal eigenbasis of $Z$, with eigenvalues $\lambda_1,\lambda_2,\dots$. Then $Y_{ij}=\frac12(\lambda_i+\lambda_j)\mathsf J_{ij}$. Together with the factorization
\begin{align*}
 \sum_{k=0}^{2s-1}\xi^k\eta^{2s-1-k}
 =(\xi+\eta)\sum_{a=0}^{s-1}\xi^{2a}\eta^{2s-2-2a},
\end{align*}
this gives
\begin{align*}
 g''(0)
 =\frac12\sum_{i,j}\mathsf J_{ij}^2(\lambda_i+\lambda_j)^2
       \sum_{a=0}^{s-1}\lambda_i^{2a}\lambda_j^{2s-2-2a}.
\end{align*}
In the same basis,
\begin{align*}
 \sigma_a=\frac14\sum_{i,j}\mathsf J_{ij}^2(\lambda_i+\lambda_j)^2\lambda_i^{2a}\lambda_j^{2s-2-2a},
\end{align*}
and so $g''(0)=2\mathcal S$. This proves the identity in~\eqref{eq:orbit_dilation_second}.

For the upper bound, let $\xi,\eta\geq0$. Expanding $(\xi+\eta)^2$ and applying the weighted AM--GM inequality to each of the three resulting monomials gives
\begin{align*}
 \tfrac14\xi^{2a}\eta^{2s-2-2a}(\xi+\eta)^2
 \leq\frac{2a+1}{2s}\xi^{2s}+\frac{2s-1-2a}{2s}\eta^{2s}.
\end{align*}
Since $(\lambda_i+\lambda_j)^2\leq(|\lambda_i|+|\lambda_j|)^2$, we may apply this with $\xi=|\lambda_i|$ and $\eta=|\lambda_j|$ and sum over $i,j$. Since $\mathsf J_{ij}^2=\mathsf J_{ji}^2$, the two resulting sums coincide, and their coefficients add up to $\bigl((2a+1)+(2s-1-2a)\bigr)/(2s)=1$. This gives
\begin{align*}
 \sigma_a\leq\sum_{i,j}\mathsf J_{ij}^2\lambda_i^{2s}=\tr[Z^{2s}\mathsf J^2].
\end{align*}
Summing over the $s$ values of $a$ proves the upper bound in~\eqref{eq:orbit_dilation_second}.

\emph{The third derivative.} The identity $g''=2\mathcal S$ holds at every time, because $Z(t+\tau)=e^{\tau \mathsf J/2}Z(t)e^{\tau \mathsf J/2}$. Thus $g'''(0)=2\mathcal S'(0)$, and~\eqref{eq:orbit_dilation_third} follows from
\begin{equation}\label{eq:orbit_energy_derivative}
 |\mathcal S'(0)|\leq20s\rho\,\mathcal S,\qquad \rho:=\|\mathsf J\|_{\rm op}.
\end{equation}
If $s=1$, then $\mathcal S=\tr[Y^2]$ and $\mathcal S'=2\tr[\mathsf JY^2]$, so $|\mathcal S'|\leq2\rho\mathcal S$. Assume now that $s\geq2$. The basic estimate is that, whenever $a+b+c=2s-2$ and $0\leq a,c\leq s-1$,
\begin{equation}\label{eq:orbit_two_factor_bound}
 |\tr[\mathsf JZ^aYZ^bYZ^c]|\leq\rho\sqrt{\sigma_a\sigma_c}.
\end{equation}
Here the product of the matrices following $\mathsf J$ equals $(Z^aYZ^{s-1-a})(Z^{s-1-c}YZ^c)$. By the Cauchy--Schwarz inequality for the Frobenius inner product, the trace is at most $\rho$ times the product of the Frobenius norms of the two factors, and these norms are $\sqrt{\sigma_a}$ and $\sqrt{\sigma_c}$.

Since $Z'=Y$ and $Y'=\frac12(\mathsf JY+Y\mathsf J)$, differentiating~\eqref{eq:orbit_square_energies} produces two kinds of terms. Define
\begin{align*}
 \mathcal T:=\sum_{a+b+c=2s-3}|\tr[Z^aYZ^bYZ^cY]|,
 \qquad
 \mathcal U:=\sum_{a=0}^{s-1}|\tr[\mathsf JZ^{2a}YZ^{2s-2-2a}Y]|,
\end{align*}
where the first sum runs over all ordered triples of nonnegative integers. Differentiating the factor $Z^{2a}$ of $\sigma_a$ gives traces with exponent triples $(j,2a-1-j,2s-2-2a)$, and the third and first entries determine $a$ and $j$. Differentiating the factor $Z^{2s-2-2a}$ gives the triples $(2a,j,2s-3-2a-j)$, which are again distinct. Each of these two families contributes at most $\mathcal T$. Differentiating the two factors $Y$ contributes at most $2\mathcal U$, after cyclic rotation and transposition. Altogether,
\begin{equation}\label{eq:orbit_differentiated_sum}
 |\mathcal S'|\leq2\mathcal T+2\mathcal U.
\end{equation}

We first bound $\mathcal T$. A cyclic rotation of a trace in $\mathcal T$ permutes its exponent triple cyclically. Every term of $\mathcal T$ equals, up to rotation, a term whose largest exponent is in the middle, and each such term arises from at most three terms of $\mathcal T$. For such a term, the first and last exponents satisfy $a,c\leq s-2$, because $a+b+c=2s-3$ and $b\geq\max\{a,c\}$. Substituting $Y=\frac12(\mathsf JZ+Z\mathsf J)$ for the last factor and applying~\eqref{eq:orbit_two_factor_bound} to the two resulting traces gives
\begin{align*}
 |\tr[Z^aYZ^bYZ^cY]|
 \leq\frac\rho2\bigl(\sqrt{\sigma_{a+1}\sigma_c}+\sqrt{\sigma_a\sigma_{c+1}}\bigr).
\end{align*}
For each pair $(a,c)$ there is at most one $b$. Extending the sum to all $0\leq a,c\leq s-2$ and applying the Cauchy--Schwarz inequality, we obtain
\begin{equation}\label{eq:orbit_triple_bound}
 \mathcal T
 \leq3\rho\Bigl(\sum_{a=1}^{s-1}\sqrt{\sigma_a}\Bigr)\Bigl(\sum_{c=0}^{s-2}\sqrt{\sigma_c}\Bigr)
 \leq3(s-1)\rho\,\mathcal S.
\end{equation}

Next we bound $\mathcal U$. Consider its summand $\tr[\mathsf JZ^{2a}YZ^{2s-2-2a}Y]$. If $2a\leq s-1$, then~\eqref{eq:orbit_two_factor_bound} with $c=0$ bounds it by $\rho\sqrt{\sigma_{2a}\sigma_0}\leq\rho\mathcal S$. If $2a>s-1$, let $k:=2a-s+1$. Iterating the identity $\mathsf JZ=2Y-Z\mathsf J$ gives
\begin{align*}
 \mathsf JZ^{2a}=(-1)^kZ^k\mathsf JZ^{s-1}+2\sum_{r=0}^{k-1}(-1)^rZ^rYZ^{2a-1-r}.
\end{align*}
We multiply on the right by $YZ^{2s-2-2a}Y$ and take traces. After cyclic rotation, the first term becomes $\pm\tr[\mathsf JZ^{s-1}YZ^{2s-2-2a}YZ^k]$, whose first and last exponents are $s-1$ and $k\leq s-1$. By~\eqref{eq:orbit_two_factor_bound}, its absolute value is at most $\rho\sqrt{\sigma_{s-1}\sigma_k}\leq\rho\mathcal S$. Each remaining term is twice a term of $\mathcal T$ with exponent triple $(r,2a-1-r,2s-2-2a)$, and distinct pairs $(a,r)$ give distinct triples. Since $\mathcal U$ has $s$ summands, it follows that
\begin{align*}
 \mathcal U\leq s\rho\mathcal S+2\mathcal T.
\end{align*}
Combining this with~\eqref{eq:orbit_differentiated_sum} and~\eqref{eq:orbit_triple_bound} gives
\begin{align*}
 |\mathcal S'|\leq6\mathcal T+2s\rho\mathcal S
 \leq(20s-18)\rho\mathcal S\leq20s\rho\mathcal S,
\end{align*}
which proves~\eqref{eq:orbit_energy_derivative}.
\end{proof}

We now transfer the bounds of Lemma~\ref{lem:ml_symmetric_curvature} to rectangular matrices by a symmetric dilation. This yields the derivative estimates for the trace power $u_{T,V}$ in~\eqref{eq:ml_orbit_objective}.

\begin{lemma}[Curvature along a matrix exponential curve]
\label{lem:ml_exponential_curvature}
Let $B$ be a rectangular matrix, let $\mathsf L$ and $\mathsf R$ be real symmetric matrices of compatible orders, and let $s\geq1$ be an integer. Define
\begin{align*}
 B(t)&:=e^{t\mathsf L}Be^{t\mathsf R},\qquad u(t):=\frac1s\tr[(B(t)B(t)^\top)^s],\\
 D&:=BB^\top,\qquad E:=B^\top B,\qquad X:=\mathsf LB+B\mathsf R.
\end{align*}
We call $\mathsf L$ and $\mathsf R$ the left and right generators. Then
\begin{align}
 u''(0)&=4\sum_{a=0}^{s-1}\tr[D^{s-1-a}XE^aX^\top]
                \geq4\tr[D^{s-1}XX^\top],
 \label{eq:ml_orbit_curvature}\\
 u''(0)&\leq8s\bigl(\tr[D^s\mathsf L^2]+\tr[E^s\mathsf R^2]\bigr),
 \label{eq:ml_orbit_upper}\\
 |u'''(0)|&\leq40s\max\{\|\mathsf L\|_{\rm op},\|\mathsf R\|_{\rm op}\}\,u''(0).
 \label{eq:ml_orbit_cubic}
\end{align}
The same bounds hold at every time $t$, with $B$ replaced by $B(t)$.
\end{lemma}
\begin{proof}
We apply Lemma~\ref{lem:ml_symmetric_curvature} to the symmetric dilation of $B$. Let
\begin{align*}
 Z:=\begin{pmatrix}0&B\\B^\top&0\end{pmatrix},\qquad
 \mathsf J:=2\begin{pmatrix}\mathsf L&0\\0&\mathsf R\end{pmatrix}.
\end{align*}
Then $e^{t\mathsf J/2}Ze^{t\mathsf J/2}$ is the dilation of $B(t)$, and its square is $\operatorname{diag}(B(t)B(t)^\top,B(t)^\top B(t))$. Hence the function $g$ of Lemma~\ref{lem:ml_symmetric_curvature} equals $u$, and $\|\mathsf J\|_{\rm op}=2\max\{\|\mathsf L\|_{\rm op},\|\mathsf R\|_{\rm op}\}$. The matrix $Y=\frac12(\mathsf JZ+Z\mathsf J)$ is the dilation of $X$, and $Z^{2a}=\operatorname{diag}(D^a,E^a)$. It follows that
\begin{align*}
 \|Z^aYZ^{s-1-a}\|_F^2=\tr[Z^{2a}YZ^{2s-2-2a}Y]
 =\tr[D^aXE^{s-1-a}X^\top]+\tr[E^aX^\top D^{s-1-a}X].
\end{align*}
We sum over $a$, rotate the second trace cyclically, and replace $a$ by $s-1-a$ in the first sum. The identity in~\eqref{eq:orbit_dilation_second} then gives the identity in~\eqref{eq:ml_orbit_curvature}. Every term of that sum is a squared Frobenius norm, which gives the lower bound. Since $\tr[Z^{2s}\mathsf J^2]=4(\tr[D^s\mathsf L^2]+\tr[E^s\mathsf R^2])$, the upper bound in~\eqref{eq:orbit_dilation_second} gives~\eqref{eq:ml_orbit_upper}, and~\eqref{eq:orbit_dilation_third} gives~\eqref{eq:ml_orbit_cubic}. Finally, $B(t+\tau)=e^{\tau \mathsf L}B(t)e^{\tau \mathsf R}$, so the bounds also hold at time $t$.
\end{proof}

\subsection{Derivatives of the optimal value}
\label{sec:ml_contact}
To differentiate $\varphi$ along $e^{tT}$ up to third order, it suffices to know the optimal weight to first order. Since the weight is stationary, an error of order $t^2$ in the weight changes the objective by $O(t^4)$. The next lemma makes this precise for a general parametric minimization.

\begin{lemma}[First-order approximation of a minimizer]
\label{lem:ml_contact}
Let $f(t,z)$ be a smooth function of $t\in\R$ and $z$ in an open subset of a Euclidean space. Suppose that, for every $t$ near $0$, the point $z_*(t)$ minimizes $f(t,\cdot)$, that $z_*$ is smooth, and that $H:=\nabla_z^2f(0,z_*(0))$ is positive definite. Let $z(t)$ be a smooth curve with $z(0)=z_*(0)$.
\begin{itemize}
\item[(a)] If $z'(0)=z_*'(0)$, then $f(t,z(t))-f(t,z_*(t))=O(t^4)$.
\item[(b)] The second derivatives at $t=0$ satisfy
\begin{align*}
 \frac{\d^2}{\d t^2}f(t,z(t))\Big|_{t=0}\geq\frac{\d^2}{\d t^2}f(t,z_*(t))\Big|_{t=0},
\end{align*}
with equality if and only if $z'(0)=z_*'(0)$.
\end{itemize}
\end{lemma}
\begin{proof}
Both claims follow from a Taylor expansion of $f(t,\cdot)$ around its minimizer $z_*(t)$. Let $\Delta(t):=z(t)-z_*(t)$. Since $\nabla_zf(t,z_*(t))=0$, Taylor's formula with integral remainder gives
\begin{align*}
 f(t,z(t))-f(t,z_*(t))
 =\int_0^1(1-a)\,
   \langle\Delta(t),\nabla_z^2f(t,z_*(t)+a\Delta(t))\Delta(t)\rangle\,\d a.
\end{align*}
If $z'(0)=z_*'(0)$, then $\Delta(t)=O(t^2)$, and the right side is $O(t^4)$. This proves part~(a). In general, $\Delta(t)=tw+O(t^2)$ with $w:=z'(0)-z_*'(0)$, so the right side equals $\frac12t^2\langle w,Hw\rangle+O(t^3)$. Comparing the second derivatives at $t=0$ shows that the left side of part~(b) exceeds the right side by $\langle w,Hw\rangle$. This quantity is nonnegative and vanishes if and only if $w=0$.
\end{proof}

We apply Lemma~\ref{lem:ml_contact} with $z=G\in\mathbb S^m$, $f(t,G):=f_{e^{tT}}(G)$, and $z_*(t):=G_*(e^{tT})$. Lemma~\ref{lem:ml_optimizer} shows that $z_*$ is smooth, and~\eqref{eq:ml_weight_hessian_energy} shows that the Hessian $H$ is positive definite. Let
\begin{align*}
 \psi(t):=\varphi(e^{tT}),\qquad V_*:=\frac{\d}{\d t}G_*(e^{tT})\Big|_{t=0}.
\end{align*}
For every $V\in\mathbb S^m$, the curve $e^{tV}$ starts at $G_*=I$ with velocity $V$, and $e^{tV_*}$ has the same velocity as $G_*(e^{tT})$. By~\eqref{eq:ml_orbit_objective} and Lemma~\ref{lem:ml_contact},
\begin{equation}\label{eq:ml_optimized_orbit}
 2\psi''(0)=\min_{V\in\mathbb S^m}u_{T,V}''(0),
 \qquad 2\psi'''(0)=u_{T,V_*}'''(0),
\end{equation}
and the minimum is attained at $V=V_*$. Thus the derivatives of the optimal value reduce to derivatives of the explicit trace powers $u_{T,V}$, and the motion of the optimal weight enters only through the choice $V=V_*$.

\subsection{Gradient and Hessian}
\label{sec:ml_rounding_gradient}\label{sec:ml_convexity_gradient}
We now compute the gradient of $\Phi_p$ and compare its Hessian with a simpler quadratic form. Both are expressed through the following matrix.

\begin{definition}[Gradient matrix and proxy Hessian]
\label{def:ml_certificate_matrices}
Let $y\in\Omega$, and let $G_*$ and $W_y$ be as in Lemma~\ref{lem:ml_optimizer}. Define
\begin{equation}\label{eq:ml_gradient_matrix}
 M_y:=\frac12\Bigl(\mathcal L_y(G_*)^s+\sum_{i,j}(G_*)_{ij}B_iW_yB_j\Bigr),\qquad
 Q_y[h,h]:=\tr[M_yB_h^2].
\end{equation}
We call $Q_y$ the proxy Hessian.
\end{definition}

The matrix $M_y$ is positive definite and has trace $m$. Both properties are unchanged by the coordinate changes of Section~\ref{sec:ml_normalization}, so we may check them in the normalization~\eqref{eq:ml_fixed_normalization}, where
\begin{align*}
 M_y=\frac12\Bigl((\mathsf B\mathsf B^\top)^s+\sum_iB_iWB_i\Bigr).
\end{align*}
Each $B_iWB_i$ is positive semidefinite, and $(\mathsf B\mathsf B^\top)^s\succ0$ by Lemma~\ref{lem:ml_sos}. By the mass identity,
\begin{align*}
 \tr\Bigl[\sum_iB_iWB_i\Bigr]=\tr[W\mathsf B\mathsf B^\top]=\tr[(\mathsf B\mathsf B^\top)^s]=m,
\end{align*}
and so $\tr[M_y]=m$.

\begin{lemma}[Gradient and proxy Hessian]
\label{lem:ml_gradient_certificates}\label{lem:ml_second_order}
For every $y\in\Omega$ and $h\in\R^m$,
\begin{equation}\label{eq:ml_gradient}
 \nabla\Phi_p(y)[h]=-\tr[M_yB_h].
\end{equation}
Moreover,
\begin{equation}\label{eq:ml_metric_comparison}
 Q_y\preceq\nabla^2\Phi_p(y)\preceq(2s+1)Q_y,
 \qquad |\nabla\Phi_p(y)[h]|^2\leq m\,\nabla^2\Phi_p(y)[h,h].
\end{equation}
In the normalization~\eqref{eq:ml_fixed_normalization}, every direction $T\in\mathbb S^n$ of the slack satisfies
\begin{align*}
 \nabla\varphi(I)[T]=-\tr[M_yT],\qquad Q(T)\leq\nabla^2\varphi(I)[T,T]\leq(2s+1)Q(T),\qquad Q(T):=\tr[M_yT^2].
\end{align*}
\end{lemma}
\begin{proof}
We obtain the gradient from the first derivative of $u_{T,V}$ and the stationarity of $G_*$. For the Hessian, the chain rule along $e^{tT}$ and~\eqref{eq:ml_optimized_orbit} express $\nabla^2\varphi(I)[T,T]$ as $Q(T)$ plus half the minimum of $u_{T,V}''(0)$ over $V$, and Lemma~\ref{lem:ml_exponential_curvature} bounds this minimum from both sides. The two coordinate changes that lead to~\eqref{eq:ml_fixed_normalization} leave $M_y$ unchanged and map the slack direction $\mathcal A^*(h)$ to $T=B_h$, so we may work in this normalization.

\emph{The gradient.} Fix $T\in\mathbb S^n$ and $V\in\mathbb S^m$. By~\eqref{eq:ml_exponential_factor},
\begin{equation}\label{eq:ml_block_tangent}
 \mathsf B_{T,V}'(0)=\tfrac12\mathsf B(V\otimes I_n)-\tfrac12\bigl(T\mathsf B+\mathsf B(I_m\otimes T)\bigr).
\end{equation}
Since $(\mathsf B\mathsf B^\top)^{s-1}=W$ and $\mathsf B(I_m\otimes T)\mathsf B^\top=\sum_iB_iTB_i$,
\begin{align*}
 u_{T,V}'(0)=2\tr[W\mathsf B_{T,V}'(0)\mathsf B^\top]
 =\tr[W\mathsf B(V\otimes I_n)\mathsf B^\top]-\tr[WT\mathsf B\mathsf B^\top]-\tr\Bigl[W\sum_iB_iTB_i\Bigr].
\end{align*}
By stationarity, the first term equals
\begin{align*}
 \sum_{i,j}V_{ij}\tr[WB_iB_j]=\tr[V\mathcal L_I^*(W)]=\tr[V].
\end{align*}
Since $W\mathsf B\mathsf B^\top=(\mathsf B\mathsf B^\top)^s$, the other two terms add up to
\begin{align*}
 -\tr\Bigl[\Bigl((\mathsf B\mathsf B^\top)^s+\sum_iB_iWB_i\Bigr)T\Bigr]=-2\tr[M_yT].
\end{align*}
We conclude that
\begin{align*}
 u_{T,V}'(0)=\tr[V]-2\tr[M_yT].
\end{align*}
The term $\tr[V]$ cancels the derivative of $-\log\det e^{tV}$ in~\eqref{eq:ml_orbit_objective}. Since $G_*=I$ is stationary, the first derivative of $f_{e^{tT}}(e^{tV})$ at $t=0$ equals $2\psi'(0)$ for every $V$. Hence
\begin{align*}
 \nabla\varphi(I)[T]=\psi'(0)=-\tr[M_yT].
\end{align*}

\emph{The Hessian.} Let $v:=\nabla^2\varphi(I)[T,T]$. The curve $e^{tT}$ has velocity $T$ and acceleration $T^2$ at $t=0$, and the chain rule gives
\begin{align*}
 \psi''(0)=v+\nabla\varphi(I)[T^2]=v-Q(T).
\end{align*}
Combined with~\eqref{eq:ml_optimized_orbit}, this gives
\begin{equation}\label{eq:ml_curvature_minimum}
 v=Q(T)+\frac12\min_{V\in\mathbb S^m}u_{T,V}''(0).
\end{equation}
Every term in~\eqref{eq:ml_orbit_curvature} is nonnegative, so $v\geq Q(T)$. For the upper bound, take $V=0$. The generators of $\mathsf B_{T,0}$ in~\eqref{eq:ml_exponential_factor} are $\mathsf L=-\frac12T$ and $\mathsf R=-\frac12I_m\otimes T$. Since $(\mathsf B^\top\mathsf B)^s=\mathsf B^\top W\mathsf B$,
\begin{align*}
 \tr[(\mathsf B^\top\mathsf B)^s(I_m\otimes T)^2]
 =\tr[W\mathsf B(I_m\otimes T^2)\mathsf B^\top]
 =\tr\Bigl[\Bigl(\sum_iB_iWB_i\Bigr)T^2\Bigr].
\end{align*}
Consequently,~\eqref{eq:ml_orbit_upper} gives
\begin{align*}
 u_{T,0}''(0)\leq2s\tr\Bigl[\Bigl((\mathsf B\mathsf B^\top)^s+\sum_iB_iWB_i\Bigr)T^2\Bigr]=4s\,Q(T),
\end{align*}
and~\eqref{eq:ml_curvature_minimum} gives $v\leq(2s+1)Q(T)$.

\emph{The gradient bound.} By the Cauchy--Schwarz inequality,
\begin{align*}
 |\tr[M_yT]|^2=|\tr[M_y^{1/2}\cdot M_y^{1/2}T]|^2\leq\tr[M_y]\tr[M_yT^2]=m\,Q(T)\leq m\,v.
\end{align*}
For $T=B_h$, these statements are~\eqref{eq:ml_gradient} and~\eqref{eq:ml_metric_comparison}.
\end{proof}

For $V\in\mathbb S^m$, we call $u_{T,V}''(0)$ the energy of $V$. By~\eqref{eq:ml_optimized_orbit}, the minimum in~\eqref{eq:ml_curvature_minimum} is attained at $V_*$, and so
\begin{align*}
 u_{T,V_*}''(0)=2(v-Q(T)).
\end{align*}
The next lemma uses this energy to bound the derivative of the optimal weight. Lemma~\ref{lem:ml_rounding} then uses $Q(T)$ to bound the operator norm of $T$.

\subsection{Variation of the optimal weight}
\label{sec:ml_directional_control}
We now bound the derivative of the optimal weight by the curvature $v$. This bound is the main input for the third derivative of the potential.

\begin{lemma}[First derivatives of the optimal weight]
\label{lem:ml_directional_response}
In the normalization~\eqref{eq:ml_fixed_normalization}, let $T\in\mathbb S^n$ and $v:=\nabla^2\varphi(I)[T,T]$, and let $V_*:=(\d/\d t)G_*(e^{tT})|_{t=0}$. Then
\begin{equation}\label{eq:ml_gprime_control}
 \|V_*\|_F\leq\sqrt{6v}.
\end{equation}
In arbitrary coordinates, this bound reads
\begin{align*}
 \|G_*^{-1/2}\nabla G_*(y)[h]\,G_*^{-1/2}\|_F\leq\sqrt{6\nabla^2\Phi_p(y)[h,h]}.
\end{align*}
\end{lemma}
\begin{proof}
The bound~\eqref{eq:ml_gprime_control} comes from the energy $u_{T,V_*}''(0)=2(v-Q(T))$ and the formula~\eqref{eq:ml_orbit_curvature} for it. For this bound, we split the tangent $X:=\mathsf B_{T,V_*}'(0)$ into a weight part and a slack part, and we control the weight part by the term $a=0$ of~\eqref{eq:ml_orbit_curvature}. %

On $n\times mn$ matrices, let $\|\cdot\|_W$ be the norm $X\mapsto\tr[WXX^\top]^{1/2}$. The term $a=0$ of~\eqref{eq:ml_orbit_curvature} equals $4\|X\|_W^2$, which gives
\begin{align*}
 4\|X\|_W^2\leq u_{T,V_*}''(0)=2(v-Q(T)).
\end{align*}
By~\eqref{eq:ml_block_tangent}, $X=X_V-X_T$, where
\begin{align*}
 X_V:=\tfrac12\mathsf B(V_*\otimes I_n),\qquad
 X_T:=\tfrac12\bigl(T\mathsf B+\mathsf B(I_m\otimes T)\bigr).
\end{align*}
The stationarity condition $\mathcal L_I^*(W)=I$ gives
\begin{align*}
 \|X_V\|_W^2=\tfrac14\tr[W\mathsf B(V_*^2\otimes I_n)\mathsf B^\top]
 =\tfrac14\tr[V_*^2\mathcal L_I^*(W)]=\tfrac14\|V_*\|_F^2.
\end{align*}
For the slack part, the triangle inequality and $(\xi+\eta)^2\leq2\xi^2+2\eta^2$ give
\begin{align*}
 \|X_T\|_W^2\leq\frac12\Bigl(\tr[WT\mathsf B\mathsf B^\top T]+\tr\Bigl[\Bigl(\sum_iB_iWB_i\Bigr)T^2\Bigr]\Bigr)\leq Q(T).
\end{align*}
The last step uses $\tr[WT\mathsf B\mathsf B^\top T]\leq\tr[(\mathsf B\mathsf B^\top)^sT^2]$. To prove it, let $d_k$ be the eigenvalues of $\mathsf B\mathsf B^\top$ and $T_{kl}$ the entries of $T$ in a corresponding orthonormal eigenbasis. Then
\begin{align*}
 \tr[(\mathsf B\mathsf B^\top)^sT^2]-\tr[WT\mathsf B\mathsf B^\top T]=\sum_{k,l}d_k^{s-1}(d_k-d_l)T_{kl}^2
 =\frac12\sum_{k,l}(d_k^{s-1}-d_l^{s-1})(d_k-d_l)T_{kl}^2\geq0.
\end{align*}
Since $\|\cdot\|_W$ satisfies the triangle inequality,
\begin{align*}
 \tfrac12\|V_*\|_F=\|X_V\|_W\leq\|X\|_W+\|X_T\|_W
 \leq\sqrt{\tfrac12(v-Q(T))}+\sqrt{Q(T)}.
\end{align*}
Multiplying by $2$ and applying the Cauchy--Schwarz inequality to the two terms gives
\begin{align*}
 \|V_*\|_F\leq\sqrt{2(v-Q(T))}+2\sqrt{Q(T)}\leq\sqrt{6v}.
\end{align*}

Finally, since the curves $e^{tT}$ and $I+tT$ have the same velocity at $t=0$, the matrix $V_*$ is the derivative $\nabla G_*(y)[h]$ written in the normalized coordinates, where $T=B_h$. The quantity $\|G_*^{-1/2}(\partial_tG_*)G_*^{-1/2}\|_F$ is invariant under the two coordinate changes of Section~\ref{sec:ml_normalization}, which proves the last assertion.
\end{proof}
\subsection{Ellipsoid rounding}
\label{sec:ml_local_norm}
We now compare the operator norm of $B_h$ with the local norms of $h$. This comparison uses the logarithmic choice of $s$.

\begin{lemma}[Ellipsoid rounding and the local norm]
\label{lem:ml_rounding}
Let $y\in\Omega$ and $h\in\R^m$. Then
\begin{align*}
 \|B_h\|_{\rm op}^2
 \leq m^{1/s}h^\top G_*^{-1}h
 \leq mn^{1/s}\|B_h\|_{\rm op}^2.
\end{align*}
Moreover,
\begin{equation}\label{eq:ml_operator_hessian}
 \|B_h\|_{\rm op}^2
 \leq2(mn)^{1/(s-1)}Q_y[h,h]
 \leq2(mn)^{1/(s-1)}\nabla^2\Phi_p(y)[h,h].
\end{equation}
For $s$ as in Definition~\ref{def:ml_s_choice}, the factor $2(mn)^{1/(s-1)}$ is less than $9$.
\end{lemma}
\begin{proof}
The first inequality restates Corollary~\ref{cor:ml_volume}. The second follows from H\"older's inequality and the mass identity, and~\eqref{eq:ml_operator_hessian} follows from H\"older's inequality and Lemma~\ref{lem:ml_gradient_certificates}. The quantities $\|B_h\|_{\rm op}$, $h^\top G_*^{-1}h$ and $Q_y[h,h]$ are invariant under the coordinate changes of Section~\ref{sec:ml_normalization}, because the coefficient change maps $h$ to $\Gamma^{-1}h$ and $G_*$ to $\Gamma^{-1}G_*\Gamma^{-\top}$. We may therefore work in the normalization~\eqref{eq:ml_fixed_normalization} and write $T:=B_h$. Then $G_*=I$, and
\begin{align*}
 \ell:=h^\top G_*^{-1}h=h^\top\mathcal L_I^*(W)h=\tr[WT^2].
\end{align*}

The containment $m^{-1/(2s)}\mathcal E_{G_*}\subseteq K_y$ of Corollary~\ref{cor:ml_volume} says that $\|B_h\|_{\rm op}\leq1$ whenever $h^\top G_*^{-1}h\leq m^{-1/s}$. By homogeneity, this is the first inequality. For the second inequality, $T^2\preceq\|T\|_{\rm op}^2I$ and H\"older's inequality give
\begin{align*}
 \ell=\tr[WT^2]\leq\tr[W]\,\|T\|_{\rm op}^2
 \leq n^{1/s}\tr[(\mathsf B\mathsf B^\top)^s]^{(s-1)/s}\|T\|_{\rm op}^2
 =n^{1/s}m^{(s-1)/s}\|T\|_{\rm op}^2,
\end{align*}
and we multiply by $m^{1/s}$.

For~\eqref{eq:ml_operator_hessian}, let $d_k$ be the eigenvalues of $\mathsf B\mathsf B^\top$, and let $t_k\geq0$ be the diagonal entries of $T^2$ in a corresponding orthonormal eigenbasis. H\"older's inequality with exponents $s/(s-1)$ and $s$ gives
\begin{align*}
 \ell=\sum_kd_k^{s-1}t_k
 \leq\Bigl(\sum_kd_k^st_k\Bigr)^{(s-1)/s}\Bigl(\sum_kt_k\Bigr)^{1/s}
 =\tr[(\mathsf B\mathsf B^\top)^sT^2]^{(s-1)/s}\tr[T^2]^{1/s}.
\end{align*}
We combine this with $\|T\|_{\rm op}^2\leq m^{1/s}\ell$ and $\tr[T^2]\leq n\|T\|_{\rm op}^2$. If $T\neq0$, dividing by $\|T\|_{\rm op}^{2/s}$ and raising both sides to the power $s/(s-1)$ gives
\begin{align*}
 \|T\|_{\rm op}^2\leq(mn)^{1/(s-1)}\tr[(\mathsf B\mathsf B^\top)^sT^2]\leq2(mn)^{1/(s-1)}Q(T),
\end{align*}
where the last step uses $M_y\succeq\frac12(\mathsf B\mathsf B^\top)^s$. The case $T=0$ is trivial. The comparison $Q_y\preceq\nabla^2\Phi_p(y)$ of Lemma~\ref{lem:ml_gradient_certificates} completes the proof of~\eqref{eq:ml_operator_hessian}. Finally, for $s$ as in Definition~\ref{def:ml_s_choice},
\begin{align*}
 s-1\geq2\log(2mn)-1\geq\log(mn),
\end{align*}
and therefore $2(mn)^{1/(s-1)}\leq2\exp(1)<9$.
\end{proof}

In terms of ellipsoids, the first two inequalities of Lemma~\ref{lem:ml_rounding} state that
\begin{align*}
 m^{-1/(2s)}\mathcal E_{G_*}\subseteq K_y\subseteq\bigl(m^{1-1/s}n^{1/s}\bigr)^{1/2}\mathcal E_{G_*}.
\end{align*}

\subsection{Proof of the main geometric result}
\label{sec:ml_third_derivative}
We now combine Section~\ref{sec:matrix_lewis_barrier} with the bounds of this section to prove Theorem~\ref{thm:ml_matrix_sc}. Its assumptions are~\eqref{eq:ml_assumptions}, its parameters are those of Definition~\ref{def:ml_s_choice}, and~\eqref{eq:main_potential} is the potential of Definition~\ref{def:ml_potential}.

\begin{proof}[Proof of Theorem~\ref{thm:ml_matrix_sc}]
Parts~(a) and~(c)--(e) collect results proved above. For Part~(b), we verify the conditions of Definition~\ref{def:sc_barrier}. The main step is the third-derivative bound, which we first prove along the exponential curve $e^{tT}$ and then transfer to the straight line $I+tT$.

\emph{Parts~(a) and~(c)--(e).} Part~(a) consists of Lemma~\ref{lem:ml_optimizer} and Corollary~\ref{cor:ml_volume}. Since $F=\alpha\Phi_p$, Part~(c) follows from Lemma~\ref{lem:ml_gradient_certificates} and the properties of $M_y$ shown after Definition~\ref{def:ml_certificate_matrices}. The first two inequalities of Lemma~\ref{lem:ml_rounding} give the two containments of Part~(d), and~\eqref{eq:ml_operator_hessian} gives its last bound. Part~(e) is the last assertion of Lemma~\ref{lem:ml_directional_response}, because $\nabla^2F=\alpha\nabla^2\Phi_p$.

\emph{The third derivative along the exponential curve.} Fix $y\in\Omega$ and a nonzero $h\in\R^m$, and work in the normalization~\eqref{eq:ml_fixed_normalization} at $S(y)$. Then $S(y+th)=I+tT$ with $T:=B_h$, so the derivatives of $\Phi_p$ at $y$ in the direction $h$ are the derivatives of $\varphi$ at $I$ in the direction $T$. Let $v:=\nabla^2\varphi(I)[T,T]$, and let $\psi(t):=\varphi(e^{tT})$. By Lemmas~\ref{lem:ml_rounding} and~\ref{lem:ml_directional_response}, together with~\eqref{eq:ml_curvature_minimum} and~\eqref{eq:ml_optimized_orbit},
\begin{align*}
 \|T\|_{\rm op}\leq3\sqrt v,\qquad
 \|V_*\|_{\rm op}\leq\|V_*\|_F\leq\sqrt{6v},\qquad
 0\leq\psi''(0)=v-Q(T)\leq v.
\end{align*}
The generators of $\mathsf B_{T,V_*}$ in~\eqref{eq:ml_exponential_factor} are $\mathsf L=-\frac12T$ and $\mathsf R=\frac12(V_*\otimes I_n-I_m\otimes T)$, so
\begin{align*}
 2\max\{\|\mathsf L\|_{\rm op},\|\mathsf R\|_{\rm op}\}
 \leq\|T\|_{\rm op}+\|V_*\|_{\rm op}\leq(3+\sqrt6)\sqrt v\leq6\sqrt v.
\end{align*}
By~\eqref{eq:ml_optimized_orbit} and~\eqref{eq:ml_orbit_cubic}, and since $u_{T,V_*}''(0)=2\psi''(0)$,
\begin{align*}
 |\psi'''(0)|=\tfrac12|u_{T,V_*}'''(0)|
 \leq20s(\|T\|_{\rm op}+\|V_*\|_{\rm op})\,\psi''(0)
 \leq120s\,v^{3/2}.
\end{align*}

\emph{From the curve to the line.} The curve $e^{tT}$ has derivatives $T$, $T^2$ and $T^3$ at $t=0$. By the chain rule and the gradient formula of Lemma~\ref{lem:ml_gradient_certificates},
\begin{align*}
 \psi'''(0)=\nabla^3\varphi(I)[T,T,T]+3\nabla^2\varphi(I)[T,T^2]-\tr[M_yT^3].
\end{align*}
We bound the last two terms by multiples of $v^{3/2}$. Since $\nabla^2\varphi(I)$ is positive semidefinite, the Cauchy--Schwarz inequality and the upper bound of Lemma~\ref{lem:ml_gradient_certificates} in the direction $T^2$ give
\begin{align*}
 |\nabla^2\varphi(I)[T,T^2]|\leq\sqrt{v\,(2s+1)\tr[M_yT^4]}
 \leq\sqrt{2s+1}\,\|T\|_{\rm op}v,
\end{align*}
where the last step uses $T^4\preceq\|T\|_{\rm op}^2T^2$ and $Q(T)\leq v$. Similarly, since $-\|T\|_{\rm op}T^2\preceq T^3\preceq\|T\|_{\rm op}T^2$,
\begin{align*}
 |\tr[M_yT^3]|\leq\|T\|_{\rm op}Q(T)\leq\|T\|_{\rm op}v.
\end{align*} Combining these bounds with $\|T\|_{\rm op}\leq3\sqrt v$ and $s\geq2$, we obtain
\begin{equation}\label{eq:ml_potential_cubic}
 |\nabla^3\Phi_p(y)[h,h,h]|
 \leq(120s+9\sqrt{2s+1}+3)v^{3/2}
 \leq150s\,v^{3/2}.
\end{equation}

\emph{Scaling and the remaining conditions.} Since $F=\alpha\Phi_p$ and $\alpha^{1/2}=1000s$, the bound~\eqref{eq:ml_potential_cubic} gives
\begin{align*}
 |\nabla^3F(y)[h,h,h]|
 \leq150s\,\alpha^{-1/2}\bigl(\nabla^2F(y)[h,h]\bigr)^{3/2}
 \leq2\bigl(\nabla^2F(y)[h,h]\bigr)^{3/2}.
\end{align*}
By Lemma~\ref{lem:ml_gradient_certificates},
\begin{align*}
 |\nabla F(y)[h]|^2=\alpha^2|\nabla\Phi_p(y)[h]|^2
 \leq\alpha^2m\,\nabla^2\Phi_p(y)[h,h]=\nu\,\nabla^2F(y)[h,h].
\end{align*}
The Hessian is positive definite, because $\nabla^2\Phi_p(y)[h,h]\geq\tr[M_yB_h^2]>0$ for $h\neq0$. Here $M_y\succ0$, and $B_h\neq0$ because the $A_i$ are linearly independent. In particular, $F$ is convex, and it is smooth on the open convex set $\Omega$ by Lemma~\ref{lem:ml_optimizer}.

\emph{Behavior at the boundary.} Taking the identity matrix as the trial point in the dual formulation~\eqref{eq:ml_weight_maximum} gives
\begin{equation}\label{eq:ml_boundary_lower}
 2\Phi_p(y)\geq\log\det \mathcal L_y^*(I)-n(1-1/s).
\end{equation}
Let $y$ tend to a point $\bar y$ on the boundary of $\Omega$. Then $S(y)$ stays bounded, so there is $c>0$ with $S(y)^{-1}\succeq cI$ near $\bar y$. Since the $A_i$ are linearly independent, there is $c'>0$ with $\|\mathcal A^*(h)\|_F^2\geq c'\|h\|_2^2$. Consequently,
\begin{align*}
 h^\top\mathcal L_y^*(I)h=\|S(y)^{-1/2}\mathcal A^*(h)S(y)^{-1/2}\|_F^2\geq c^2c'\|h\|_2^2.
\end{align*}
Thus the smallest eigenvalue of $\mathcal L_y^*(I)$ stays bounded away from zero. On the other hand, $A_0:=\sum_iA_i^2$ is positive definite because $\bigcap_i\ker A_i=\{0\}$, and $A_iS(y)^{-1}A_i\succeq cA_i^2$. Hence
\begin{align*}
 \tr[\mathcal L_y^*(I)]=\sum_i\tr[S(y)^{-1}A_iS(y)^{-1}A_i]
 \geq c\tr[S(y)^{-1}A_0]
 \geq\frac{c\,\lambda_{\min}(A_0)}{\lambda_{\min}(S(y))}\longrightarrow\infty,
\end{align*}
since $S(y)$ becomes singular. Consequently, the largest eigenvalue of $\mathcal L_y^*(I)$ diverges, and so does its log determinant. By~\eqref{eq:ml_boundary_lower}, $F(y)\to\infty$.
\end{proof}

\clearpage
\part{The randomized algorithm}
\label{part:algorithm}
\section{The main algorithmic result}
\label{sec:algorithm_part_overview}
The core of our SDP algorithm is to solve the dual SDP from a given near-central point. The algorithm is a short-step path-following method for the matrix Lee--Sidford barrier, and it accesses the barrier only through three data structures: slack maintenance, weight maintenance, and gradient/Hessian estimation. Its guarantees are stated in the theorem below. The preprocessing and postprocessing modules of our algorithm are presented in Sections~\ref{sec:preprocessing} and~\ref{sec:postprocessing}, which provide a starting point and convert the final dual solution into a primal solution. Section~\ref{sec:solver_wrapup} composes them together and establishes the main algorithmic result.

As throughout the paper, the dual minimizes $b^\top y$ subject to $S(y)\succ0$, and $\operatorname{OPT}$ is the common optimal value of the primal and dual problems.

\begin{theorem}[Randomized dual optimization from a supplied central point]
\label{thm_more_amortized_sdp}
\label{thm:core_dual_sdp}
Let $1\leq m\leq n$, let $A_1,\ldots,A_m,C$ be real symmetric $n\times n$ matrices, and let $b\in\R^m$. Assume that the $A_i$ are linearly independent, $\bigcap_i\ker A_i=\{0\}$, and both primal and dual SDPs are strictly feasible. Use the matrix Lee--Sidford barrier $F=\alpha\Phi_p$ of Theorem~\ref{thm:ml_matrix_sc}, with
\begin{align*}
 s:=2\lceil\log(2mn)\rceil,\qquad \alpha:=10^6s^2,\qquad
 \nu:=\alpha m,\qquad \delta_{\rm c}:=1/100,\qquad \gamma_t:=\frac{\delta_{\rm c}}{4\sqrt\nu}.
\end{align*}
Suppose that the input contains $y_0$ and $t_0>0$ with
\begin{align*}
 S(y_0)\succ0,\qquad
 \lambda_{t_0b}(y_0)
 :=\|\nabla F(y_0)+t_0b\|_{\nabla^2F(y_0)^{-1}}\leq1/10.
\end{align*}
Let $\epsilon>0$, $0<\beta<1/10$, and choose an integer $1\leq H\leq\sqrt m$. For $m^2>n$, each epoch consists of at most $H$ centering phases and begins with a full pair estimate. For $m^2\leq n$, each query uses fresh estimates and $H$ is unused. Let $L_\epsilon:=\log(2+m/(t_0\epsilon))$, and let $\mathcal B_*(H)$ be the per-query cost in~\eqref{eq:effective_query_cost}. If $m^4\leq n$, the algorithm reuses slack factorizations as in Section~\ref{sec:reused_factorization}. Otherwise it factors the slack at every query. Then Algorithm~\ref{alg:more_sdp} has the following guarantees.
\begin{itemize}
\item[(a)] With probability at least $1-\beta$, it returns a pair $(y,t)$ such that
\begin{gather*}
 S(y)\succ0,\qquad 0\leq b^\top y-\operatorname{OPT}\leq\epsilon,
 \qquad \lambda_{tb}(y)\leq\delta_{\rm c},\\
 \max\{t_0,2\nu/\epsilon\}\leq t
 \leq(1+\gamma_t)\max\{t_0,2\nu/\epsilon\}.
\end{gather*}
The multiplier is $t=t_0(1+\gamma_t)^{N_{\rm path}}$, where
\begin{align*}
 N_{\rm path}:=\left\lceil
 \frac{\max\{0,\log(2\nu/(t_0\epsilon))\}}{\log(1+\gamma_t)}
 \right\rceil.
\end{align*}
\item[(b)] The algorithm uses $O\left(s^2\sqrt m\log(2+\nu/(t_0\epsilon))\right)$ $y$-steps and $O\left(s^3\log s\sqrt m\log(2+\nu/(t_0\epsilon))\right)$ oracle queries, including the $G$-queries. Its total cost is
\begin{equation}\label{eq_more_sdp_total}
 \widetilde O\bigl((n^\omega+\sqrt m\mathcal B_*(H))L_\epsilon\bigr)
\end{equation}
operations in the arithmetic model of Section~\ref{sec:algorithm_preliminaries}. On every run, including a failed one, every prefix of $N\geq1$ queries costs
\begin{equation}\label{eq_more_prefix_cost}
 \widetilde O\bigl(L_\epsilon n^\omega+\mathcal C(m)
                         +N\mathcal B_*(H)\bigr).
\end{equation}
These bounds include all initial oracle computations, $G$-steps, proxy Hessian estimates, whitening operations and full pair estimates.
\item[(c)] Under Fact~\ref{fact:rectangular_thresholds}, let $\omega$ be a square exponent bound with $2<\omega\leq9/4$, and choose $H=\lfloor m^{1/6}\rfloor$. For every fixed $\zeta>0$, the total cost is
\begin{equation}\label{eq:reused_core_explicit}
 \widetilde O_\zeta\left(
      (n^{\omega+\zeta}+n^2m^{3/2+\zeta})L_\epsilon\right).
\end{equation}
\end{itemize}
\end{theorem}

The factors suppressed by $\widetilde O$ are specified in Section~\ref{sec:arithmetic_model}. Section~\ref{sec:algorithm_description} gives the algorithm and states the guarantees of its three data structures. Section~\ref{sec:path_following} derives Part~(a) from these guarantees, Sections~\ref{sec:randomized_oracles}--\ref{sec:reused_factorization} implement the data structures, and Section~\ref{sec:randomized_solver} proves Parts~(b) and~(c).
\section{Preliminaries for the randomized algorithm}
\label{sec:algorithm_preliminaries}
Part~\ref{part:algorithm} uses the notation and self-concordance estimates of Section~\ref{sec:barrier_preliminaries}. This section provides preliminaries of the arithmetic model, bounds on matrix multiplication, fast linear algebra, and Gaussian sketching.

\subsection{Arithmetic model and approximation parameters}
\label{sec:arithmetic_model}
The running times count exact real arithmetic operations, comparisons and positive scalar square roots. Randomized procedures also count the draw of each independent standard Gaussian random variable as one operation. The inverse square roots used for slack access are evaluated by polynomial approximations to a prescribed accuracy, and all other matrix powers are integer powers. These bounds count operations in exact real arithmetic and leave the bit precision of a finite-precision implementation as a separate parameter.

The notation $\widetilde O$ suppresses fixed powers of $\log(2mn)$ and of the logarithms of the stated accuracy, magnitude, radius, strict-feasibility and failure parameters. Each oracle states any additional polynomial dependence on its requested local tolerance. A bound that contains a factor $n^{o(1)}$ is used in the following fixed-slack sense. For every fixed $\zeta>0$, the bound holds with $n^\zeta$ in place of $n^{o(1)}$ and with a constant that may depend on $\zeta$. We write $O_\zeta$ and $\widetilde O_\zeta$ for such bounds. The formal solver statements keep the dependence on the accuracy and the failure probability separate from this convention.

\subsection{Rectangular matrix multiplication}
\label{sec:rectangular_matrix_multiplication}
Define $\Tmat(a,b,c)$ to be the arithmetic time to multiply an $a\times b$ matrix by a $b\times c$ matrix. For nonnegative real numbers $a,b,c$, define $\omega(a,b,c)$ to be the infimum of the exponents $\tau$ for which matrices of dimensions $\lceil N^a\rceil\times\lceil N^b\rceil$ and $\lceil N^b\rceil\times\lceil N^c\rceil$ can be multiplied in $O(N^\tau)$ arithmetic operations as $N\to\infty$. An admissible square exponent $\tau>2$ satisfies $\Tmat(d,d,d)=O(d^\tau)$, and every $\tau>\max\{2,\omega(1,1,1)\}$ is admissible. In unspecialized arithmetic bounds, $2<\omega<3$ denotes any fixed admissible square exponent. In numerical specializations with fixed slack $\zeta>0$, $\omega$ instead denotes a bound satisfying $\omega(1,1,1)\leq\omega$. There we apply the unspecialized bounds with the admissible exponent $\tau=\omega+\zeta'$, choosing $\zeta'>0$ sufficiently small that all exponent losses are absorbed by $\zeta$. This convention applies whether or not the infimum is attained. Currently $\omega \leq 9/4$ \cite{o26b,p26}.

\begin{lemma}[Transposition principle {\cite{hm73}}]
\label{lem:matrix_transposition}
With input and output access included in $\Tmat$, every permutation of the dimensions has the same arithmetic cost up to a universal constant. In particular, $\Tmat(p,q,r)=O(\Tmat(p,r,q))$.
\end{lemma}
\begin{proof}
Append the pairing $\tr[XYZ]$ to a bilinear algorithm for $XY$. Reverse accumulation with respect to $Y$ computes $(ZX)^\top$ with constant overhead. The pairing costs $O(pq+qr+rp)$, already included in $\Tmat$; transposition gives the remaining permutations.
\end{proof}

The numerical running-time specializations use one set of matrix-multiplication bounds over $\R$, taken from \cite[Theorem~1.1]{u26}.

\begin{fact}[Rectangular and square multiplication over $\R$]
\label{fact:rectangular_thresholds}
\label{fact:external_multiplication_inputs}
Over $\R$,
\begin{align*}
 \omega(1,1,1/2)=2,\qquad \omega(1,2,2/3)=3,\qquad \omega(2,1,1)\leq25/8,\qquad \omega(1,1,1)\leq9/4.
\end{align*}
Consequently, for every fixed $\zeta>0$,
\begin{align*}
 \Tmat(n,n,\lceil n^{1/2}\rceil)=O_\zeta(n^{2+\zeta}),\qquad \Tmat(n,n^2,\lceil n^{2/3}\rceil)=O_\zeta(n^{3+\zeta}).
\end{align*}
\end{fact}
\begin{proof}
Substitute the four shapes into \cite[Theorem~1.1]{u26}. The size lower bound $\omega(a,b,c)\geq\max\{a+b,a+c,b+c\}$ gives equality for the first two shapes. The arithmetic bounds follow from the fixed-slack convention in Section~\ref{sec:arithmetic_model}.
\end{proof}

\subsection{Fast linear algebra}
\label{sec:linear_algebra_primitives}
We state several tools from previous work.
\begin{lemma}[Fast factorizations {\cite{ddh07}}]
\label{lem:fast_factorizations}
Fix an admissible square matrix-multiplication exponent $\tau>2$. An invertible $d\times d$ matrix can be inverted, and a positive definite $d\times d$ matrix can be Cholesky factored, in $O(d^\tau)$ real arithmetic operations and positive scalar square roots. Triangular solves and products involving the factors have the corresponding matrix-multiplication costs.
\end{lemma}

The reference-metric updates require spectral clipping. The following approximate diagonalization suffices for these updates.

\begin{lemma}[Constant-accuracy diagonalization {\cite{ddh07,bgks20,sobczyk25}}]
\label{lem:approximate_diagonalization}
Given a symmetric $d\times d$ matrix $A$ and a bound $B\geq\|A\|_F$, one can compute an orthogonal matrix $U$ such that
\[
 U^\top AU=\operatorname{diag}(a_1,\ldots,a_d)+E,\qquad
 \operatorname{diag}(E)=0,\qquad \|E\|_F\leq1/4,
\]
using $\widetilde O(d^\omega)$ arithmetic operations, where $\widetilde O$ also hides factors polylogarithmic in $B$.
\end{lemma}

To reuse a slack factorization, we approximate the inverse square root of the relative slack by a polynomial. The following bound controls both the approximation error and the number of matrix products.

\begin{lemma}[Polynomial approximation for inverse square root~{\cite[(18.12.11) and (18.14.1)]{dlmf}}]
\label{lem:inverse_sqrt_polynomial}
\label{lem:legendre_identities}
For every $\Upsilon\geq2$ and $0<\varepsilon\leq1/2$, there is an explicit polynomial $\chi_\Upsilon$ of degree
\begin{equation}\label{eq:polynomial_parameters}
 d=\left\lceil\frac{\Upsilon+1}2\log\frac\Upsilon\varepsilon\right\rceil
\end{equation}
such that
\begin{equation}\label{eq:relative_polynomial_error}
 |\sqrt x\,\chi_\Upsilon(x)-1|\leq\varepsilon
 \qquad\text{for all }x\in[\Upsilon^{-1},\Upsilon].
\end{equation}
Consequently, for every symmetric matrix $\Sigma$ with spectrum in $[\Upsilon^{-1},\Upsilon]$, the matrix $\mathsf K:=\chi_\Upsilon(\Sigma)\Sigma^{1/2}$ is symmetric with $\|\mathsf K-I\|_{\rm op}\leq\varepsilon$, and Horner's rule evaluates $\chi_\Upsilon(\Sigma)X$ with $d$ products by $\Sigma$.
\end{lemma}

\subsection{Gaussian sketches}
\label{sec:gaussian_sketch_preliminaries}
A shared Gaussian right sketch preserves the Gram matrix of a family of matrices. This is the Gaussian case of the tensorized moment bounds in~\cite[Lemma~14]{akkpvwz20} and \cite[Lemmas~A.1--A.2]{ms22}, and we give the needed estimates directly in this subsection.

\begin{lemma}[Gaussian quadratic forms]
\label{lem:gaussian_quadratic_forms}
Let $A$ be a symmetric $d\times d$ matrix, let $z\sim\mathcal N(0,I_d)$, and let $Z:=z^\top Az-\tr[A]$. Then $\E[Z]=0$, $\E[Z^2]=2\|A\|_F^2$, and $\log\E[e^{tZ}]\leq2t^2\|A\|_F^2$ for $|t|\leq\frac1{4\|A\|_{\mathrm{op}}}$. For independent copies $Z_1,\ldots,Z_k$ and $u>0$,
\begin{align*}
 \Pr\Bigl[\Bigl|\frac1k\sum_{j=1}^kZ_j\Bigr|>u\Bigr]
 \leq2\exp\Bigl(-\frac{k}{8}\min\Bigl\{\frac{u^2}{\|A\|_F^2},
 \frac{u}{\|A\|_{\mathrm{op}}}\Bigr\}\Bigr).
\end{align*}
For $A=0$ the error vanishes identically.
\end{lemma}
\begin{proof}
Orthogonal invariance reduces to a diagonal $A$ with eigenvalues $\lambda_i$. Independence and the scalar Gaussian integral give
\begin{align*}
 \E[e^{tZ}]=\prod_i e^{-t\lambda_i}(1-2t\lambda_i)^{-1/2},\qquad
 \log\E[e^{tZ}]=\frac12\sum_i\bigl(-\log(1-2t\lambda_i)-2t\lambda_i\bigr)
 \leq2t^2\sum_i\lambda_i^2.
\end{align*}
The last step uses $-\log(1-x)-x\leq x^2$ for $|x|\leq1/2$. The mean and the variance follow by differentiation at zero. For the upper tail, apply the exponential Markov inequality with $t:=\min\{\frac{u}{4\|A\|_F^2},\frac1{4\|A\|_{\mathrm{op}}}\}$ to $\sum_jZ_j$. If the second choice is active, then $u\geq\|A\|_F^2/\|A\|_{\mathrm{op}}$, which bounds the quadratic term in the exponent by half its linear term. Applying the same argument to $-Z$ and adding the two tail probabilities proves the lemma.
\end{proof}

We next combine Lemma~\ref{lem:gaussian_quadratic_forms} with a finite net to control the sketching error uniformly over a fixed matrix span. This uniform norm bound gives a simultaneous approximation to the Gram matrix of the family.

\begin{lemma}[Gaussian right sketch]
\label{lem:gaussian_right_embedding}
Let $X_1,\ldots,X_m\in\R^{a\times b}$ be fixed, with span dimension $d\leq m$. Let $0<\theta,\rho<1/2$, and draw $\Omega\in\R^{b\times k}$ with independent $\mathcal N(0,1/k)$ entries, where $k\geq C\theta^{-2}(d+\log(2/\rho))$ for a universal constant $C$. Define $H_{ij}:=\langle X_i,X_j\rangle$ and $\widehat H_{ij}:=\langle X_i\Omega,X_j\Omega\rangle$. With probability at least $1-\rho$, $(1-\theta)H\preceq\widehat H\preceq(1+\theta)H$. Equivalently, $\|X\Omega\|_F^2$ approximates $\|X\|_F^2$ to relative error $\theta$ simultaneously for all $X$ in the span. The guarantee also holds conditionally on random $X_i$ drawn before, and independently of, $\Omega$.
\end{lemma}
\begin{proof}
The case $d=0$ is immediate. For a unit-Frobenius matrix $X$ in the span, $A=X^\top X$ has trace one and both $\|A\|_{\rm op}$ and $\|A\|_F$ at most one. Lemma~\ref{lem:gaussian_quadratic_forms} gives
\begin{align*}
 \Pr[|\|X\Omega\|_F^2-1|>u]\leq2e^{-ku^2/8}\qquad(0<u<1).
\end{align*}
A maximal $1/4$-separated subset $\mathcal M$ of the unit sphere has at most $9^d$ points: disjoint balls of radius $1/8$ about its points lie in a ball of radius $9/8$. It is a $1/4$-net by maximality. The union bound with $u=\theta/2$ controls the sketching error on $\mathcal M$ except with probability $2\cdot9^d e^{-k\theta^2/32}\leq\rho$. For a symmetric form $T$, approximating a maximizing unit vector by one in $\mathcal M$ gives $\|T\|_{\rm op}\leq\max_{v\in\mathcal M}|T[v,v]|+\|T\|_{\rm op}/2$. Apply this to $T[X,X]=\|X\Omega\|_F^2-\|X\|_F^2$ to obtain the uniform bound $\theta$. The conditional claim fixes the family before applying the same argument.
\end{proof}
\subsection{Gaussian range approximation and matrix concentration}
\label{sec:random_matrix_preliminaries}
We state a tool from previous work.
\begin{lemma}[Gaussian power iteration {\cite[Corollary~10.10]{hmt11}}]
\label{lem:gaussian_power_method}
Let $A\in\R^{a\times b}$, let $d:=\min\{a,b\}$, and choose integers $2\leq r<d$, $p\geq2$ and $\ell\geq0$ with $r+p\leq d$. Draw a standard Gaussian $b\times(r+p)$ matrix $\Omega$, and let $U$ have orthonormal columns that span the range of $(AA^\top)^\ell A\Omega$. Writing $\sigma_j:=\sigma_j(A)$, one has
\begin{align*}
 \E\|(I-UU^\top)A\|_{\rm op}
 \leq\left[1+\sqrt{\frac r{p-1}}
       +\frac{e\sqrt{r+p}}p\sqrt{d-r}\right]^{1/(2\ell+1)}\sigma_{r+1}.
\end{align*}
The same bounds hold for $r=1$ by the one-row Gaussian calculation. The bound also holds conditionally on a random $A$ drawn before, and independently of, $\Omega$.
\end{lemma}
The proxy Hessian estimate uses the following consequence of the matrix Gaussian-series and Bernstein inequalities \cite[Theorems~4.1.1 and~6.1.1]{tropp15}.

\begin{lemma}[Gaussian Gram concentration]
\label{lem:gaussian_gram_concentration}
Let $X$ be a centered Gaussian $a\times d$ matrix, possibly with correlated entries, and suppose $\E[X^\top X]\preceq I$, $\E[XX^\top]\preceq vI$, and $v\geq1$. For $0<\epsilon,\rho<1/2$, independent copies $X_1,\ldots,X_k$ satisfy
\begin{align*}
 \Bigl\|\frac1k\sum_{j=1}^kX_j^\top X_j-\E[X^\top X]\Bigr\|_{\rm op}
 \leq\epsilon
\end{align*}
with probability at least $1-\rho$, provided $k\geq Cv\epsilon^{-2}\log^2(C(a+d)v/(\epsilon\rho))$. The statement also holds conditionally on the matrices defining $X$.
\end{lemma}
\begin{proof}
Write $Y:=X^\top X$, so that $\|Y\|_{\rm op}=\|X\|_{\rm op}^2$. Writing $X$ as a Gaussian series, the matrix Gaussian-series inequality gives $\Pr[\|Y\|_{\rm op}>u]\leq(a+d)e^{-u/(2v)}$. For the analysis only, truncate at $R_0:=C_1v\log\frac{C_2(a+d)kv}{\epsilon\rho}$ and let $\overline Y:=Y\mathbf1_{\{\|Y\|_{\rm op}\leq R_0\}}$. Integrating the tail gives
\begin{align*}
 \|\E[Y-\overline Y]\|_{\rm op}\leq(a+d)(R_0+2v)e^{-R_0/(2v)}
 \leq\epsilon/4.
\end{align*}
By positivity, $\E[(\overline Y-\E[\overline Y])^2]\preceq\E[\overline Y^2]\preceq R_0\E[\overline Y]\preceq R_0I$, and $\|\overline Y-\E[\overline Y]\|_{\rm op}\leq2R_0$.

The stated number of copies ensures $k\epsilon^2\geq CR_0\log(4d/\rho)$. Indeed, at the lower threshold of $k$ we have $R_0=O(v\log(C(a+d)v/(\epsilon\rho)))$, and the inequality persists as $k$ increases because $k/\log(Ak)$ is increasing for $A>e$. Matrix Bernstein therefore bounds the centered average of the truncated copies by $\epsilon/2$, except with probability $\rho/2$. The tail bound also gives probability at most $\rho/2$ that any of the $k$ copies is truncated. On the complementary event, adding the expectation error proves the assertion for the actual copies.
\end{proof}
\section{The algorithm and its data structures}
\label{sec:algorithm_description}
The algorithm of Theorem~\ref{thm_more_amortized_sdp} is a short-step path-following method that accesses the barrier only through three data structures. The slack maintenance data structure $\mathcal D_{\rm slack}$ supplies whitening actions for the slack matrix. The weight maintenance data structure $\mathcal D_{\rm weight}$ maintains an approximation to the optimal weight. The gradient/Hessian estimation data structure $\mathcal D_{\rm est}$ returns estimates of the gradient, of the stationarity matrix and of the proxy Hessian. It is implemented by $\mathcal D^{\rm f}_{\rm est}$, with fresh estimates at every query, if $m^2\leq n$, and by $\mathcal D^{\rm m}_{\rm est}$, with estimates maintained across queries, if $m^2>n$. This section defines the estimated quantities, gives the algorithm, and states the guarantee of each data structure as a theorem. Section~\ref{sec:path_following} derives Part~(a) of Theorem~\ref{thm_more_amortized_sdp} from these three guarantees alone, and Sections~\ref{sec:randomized_oracles}--\ref{sec:reused_factorization} implement the data structures.

A member $x$ of a data structure $\mathcal D$ is written $\mathcal D.x$, and a procedure $\textsc{Op}$ of $\mathcal D$ is written $\mathcal D.\textsc{Op}$. Thus every symbol stored by the algorithm names its owner. Inside the section that implements a data structure, we drop the prefix after stating it once.

\subsection{Quantities at a query state}
\label{sec:query_quantities}
Throughout Part~\ref{part:algorithm}, $F=\alpha\Phi_p$ is the matrix Lee--Sidford barrier of Theorem~\ref{thm:ml_matrix_sc}, with the parameters $s,\alpha,\nu$ of Definition~\ref{def:ml_s_choice}, and $\mathcal P=(A_1,\ldots,A_m,C)$ is the pencil. For an objective vector $c$, write
\begin{align*}
 \phi_c(y):=F(y)+c^\top y,\qquad
 \lambda_c(y):=\|\nabla\phi_c(y)\|_{\nabla^2F(y)^{-1}}
\end{align*}
for the penalized objective and its Newton decrement. The algorithm evaluates the barrier at a point $y$ through a weight $G$ that approximates the optimal weight $G_*(y)$ of Definition~\ref{def:ml_potential}. We call the pair $(y,G)$ a \emph{query state}.

\begin{definition}[Quantities at a query state]
\label{def:random_implicit_notation}
Let $y$ be feasible, and use $S,B_i,B_h,\mathcal L_y,\mathcal L_y^*$ and $f_y$ from Definition~\ref{def:ml_potential}. For $G\succ0$ with Cholesky factor $L$, let
\begin{align*}
 V_a:=\sum_iL_{ia}B_i,\qquad \mathcal V:=[V_1\ \cdots\ V_m],\qquad
 R:=L^\top\mathcal L_y^*\bigl(\mathcal L_y(G)^{s-1}\bigr)L,
\end{align*}
and define the \emph{weight gap}
\begin{align*}
 \Delta_y(G):=f_y(G)-f_y(G_*(y)).
\end{align*}
The gradient matrix, the fixed-weight gradient and the proxy Hessian at $(y,G)$ are
\begin{align*}
 M:=\frac12\Bigl(\mathcal L_y(G)^s+\sum_aV_a\mathcal L_y(G)^{s-1}V_a\Bigr),\qquad
 g_i:=-\alpha\tr[MB_i],\qquad
 Q_{ij}:=\tfrac12\tr[M(B_iB_j+B_jB_i)].
\end{align*}
We call $R$ the \emph{stationarity matrix}.
\end{definition}

Since $\sum_aV_aXV_a=\sum_{i,j}G_{ij}B_iXB_j$ for every $X\in\mathbb S^n$, the matrix $M$, the vector $g$ and the matrix $Q$ are functions of $G$. The vector $g$ is the gradient of $y\mapsto\frac\alpha2f_y(G)$ with $G$ held fixed (see the proof of Lemma~\ref{lem:weight_gap_motion}), which explains its name. We write $g_G(y),R_G(y),Q_G(y)$ when the state matters, and use the subscript $*$ for the values at $G_*(y)$. At the optimal weight, $M$, $g$ and $Q$ are the matrix $M_y$, the gradient $\nabla F(y)$ and the proxy Hessian $Q_y$ of Definition~\ref{def:ml_certificate_matrices} and~\eqref{eq:ml_gradient}, and stationarity gives $R_*(y)=I$. When the point is fixed, we drop $y$ and write $\Delta(G)$. A weight is \emph{accurate} at $y$ if $\Delta_y(G)\leq\eta$ and \emph{warm} if $\Delta_y(G)\leq\eta_{\rm w}$, where $\eta\leq\eta_{\rm w}$ are fixed in~\eqref{eq:single_loop_parameters}.

The algorithm forms the normalized constraint matrices $B_i=JA_iJ^\top$ for a matrix $J$ with $JS(y)J^\top=I$. We call such a $J$ an \emph{exact whitener} of $S(y)$, for example the inverse of a Cholesky factor of $S(y)$. Every exact whitener equals $OS(y)^{-1/2}$ for an orthogonal $O$. For a fixed coefficient basis, the substitution $B_i\mapsto OB_iO^\top$ conjugates the matrices of Definition~\ref{def:random_implicit_notation} that act on the slack space,
\begin{align*}
 V_a\mapsto OV_aO^\top,\qquad
 \mathcal V\mapsto O\mathcal V(I_m\otimes O^\top),\qquad
 \mathcal L_y(G)\mapsto O\mathcal L_y(G)O^\top,\qquad
 M\mapsto OMO^\top,
\end{align*}
and it leaves $f_y$, the weight gap and the quantities $R$, $g$ and $Q$ on the coefficient space unchanged. Indeed, the latter depend on the slack-space matrices only through traces such as $\tr[(OMO^\top)(OB_iO^\top)]=\tr[MB_i]$. Hence the targets $g$, $R$ and $Q$ of the estimates are the same for every exact whitener, while the sketched factors of Section~\ref{sec:common_gram} rotate with it, and Lemma~\ref{lem_more_local_geometry} tracks this rotation between consecutive queries. The coefficient basis is fixed after the first query, and all outputs are expressed in the original basis.

Most of the work applies the normalized constraint matrices to thin blocks of vectors. We measure its cost as follows.

\begin{definition}[Block-operation cost]
\label{def_more_block_cost}
For $q\geq1$, rounded upward when necessary, let
\begin{align*}
 \mathcal C(q):=\Tmat(mn,n,q)+\Tmat(n,n,mq)+\Tmat(m,m,nq)+mnq.
\end{align*}
\end{definition}

Lemma~\ref{lem:random_implicit_blocks} shows that $\mathcal C(q)$ bounds one application of all normalized constraint matrices to $q$ columns.

\subsection{Path following and centering}
\label{sec:core_algorithm_description}
\label{sec:centering_description}
Algorithm~\ref{alg:more_sdp} follows the central path at the multipliers $t_j=t_0(1+\gamma_t)^j$ for $0\leq j\leq N_{\rm path}$. At each multiplier, a call of \textsc{Center} restores the bound $\lambda_{tb}(y)\leq\delta_{\rm c}$ on the decrement. The final multiplier satisfies $t_{N_{\rm path}}\geq2\nu/\epsilon$, and at this multiplier centrality certifies the dual objective accuracy.

\begin{algorithm}[!ht]
\caption{Amortized SDP optimization from a central start}
\label{alg:more_sdp}
\begin{algorithmic}[1]
\Procedure{AmortizedSDP}{$\mathcal P,b,y_0,t_0,\epsilon,\beta;H$} \Comment{Theorem~\ref{thm_more_amortized_sdp}}
\State Use the parameters of~\eqref{eq:more_centering_parameters} and~\eqref{eq:single_loop_parameters}.
\State $N_{\rm path}\gets\lceil\max\{0,\log(2\nu/(t_0\epsilon))\}/\log(1+\gamma_t)\rceil$ and $N_{\rm tot}\gets(N_{\rm path}+1)N_{\rm c}$.
\State $\mathsf S\gets\Call{InitializeOracleState}{\mathcal P,H,N_{\rm tot},\beta}$. \Comment{Algorithm~\ref{alg:maintained_center}}
\State $y\gets y_0$ and $t\gets t_0$.
\For{$j=0,\ldots,N_{\rm path}$}
  \If{$j>0$}
    \State $t\gets(1+\gamma_t)t$.
  \EndIf
  \State $y\gets\Call{Center}{\mathcal P,tb,y;\mathsf S}$. \Comment{$N_{\rm c}$ queries, Algorithm~\ref{alg:maintained_center}}
\EndFor
\State \Return $(y,t)$.
\EndProcedure
\end{algorithmic}
\end{algorithm}

A call $\textsc{Center}(\mathcal P,c,y;\mathsf S)$ is a \emph{centering phase}. Its first oracle query returns an estimate $\widehat Q$ of the proxy Hessian with its Cholesky factor, and the call keeps the \emph{phase preconditioner} $P:=\alpha\widehat Q$ fixed until it returns. The call makes $N_y$ damped \emph{$y$-steps}
\begin{equation}\label{eq:point_step}
 y^+=y-\gamma_yP^{-1}(\widehat g+c).
\end{equation}
After each $y$-step, it holds the point fixed and makes $N_G$ \emph{$G$-steps}. Each $G$-step draws a new stationarity estimate $\widehat R$ and applies the damped inverse update
\begin{equation}\label{eq:hub_weight_update}
 G^+=L\bigl((1-\gamma_G)I+2\gamma_G(I+\widehat R)^{-1}\bigr)L^\top,
 \qquad G=LL^\top.
\end{equation}
The query before a $y$-step is a \emph{$y$-query}, and the queries before the $G$-steps are \emph{$G$-queries}. A call therefore makes $N_y$ $y$-queries and $N_{\rm c}=(N_G+1)N_y$ queries in total.

Algorithm~\ref{alg:maintained_center} gives the procedures. \textsc{InitializeOracleState} creates the three data structures and a query counter $N_{\rm q}$, which together form the oracle state $\mathsf S$. \textsc{BarrierOracle} answers one query. It obtains the whitening actions from $\mathcal D_{\rm slack}$, computes the weight at the first query, and passes the query state to $\mathcal D_{\rm est}$. The flag $a$ equals one at the first query of a call, where $\mathcal D_{\rm est}$ also estimates the proxy Hessian. A \emph{standalone query} is a call of \textsc{BarrierOracle} with $a=1$ and the zero objective outside \textsc{Center}. Preprocessing uses one before its first centering phase.

\begin{algorithm}[!ht]
\caption{Oracle state, oracle queries and one centering phase}
\label{alg:maintained_center}
\begin{algorithmic}[1]
\Procedure{InitializeOracleState}{$\mathcal P,H,N_{\max},\gamma$}
 \State Create empty data structures $\mathcal D_{\rm slack}$ and $\mathcal D_{\rm weight}$.
 \State Create $\mathcal D_{\rm est}$ with $H$, $N_{\max}$ and failure budget $2\gamma/3$, as $\mathcal D^{\rm f}_{\rm est}$ if $m^2\leq n$ and as $\mathcal D^{\rm m}_{\rm est}$ otherwise.
 \State $N_{\rm q}\gets0$ and \Return $\mathsf S:=(\mathcal D_{\rm slack},\mathcal D_{\rm weight},\mathcal D_{\rm est},N_{\rm q})$.
\EndProcedure
\Statex
\Procedure{BarrierOracle}{$\mathcal P,c,y,a;\mathsf S$}
 \State $N_{\rm q}\gets N_{\rm q}+1$, and report failure if $N_{\rm q}>N_{\max}$.
 \State $\widehat J\gets\mathcal D_{\rm slack}.\textsc{Whiten}(y,c)$. \Comment{$O(mn^2)$, and $O(n^\omega)$ if it factors, Theorem~\ref{thm:slack_interface}}
 \If{$\mathcal D_{\rm weight}$ is empty}
  \State $\mathcal D_{\rm weight}.\textsc{Init}(y,\widehat J,\gamma/3)$. \Comment{$\widetilde O(n^\omega+\mathcal C(m))$ once, Theorem~\ref{thm:weight_interface}(e)}
 \EndIf
 \State \Return $\mathcal D_{\rm est}.\textsc{Query}(y,\mathcal D_{\rm weight}.G,\widehat J,a)$. \Comment{$\widetilde O(\mathcal B_*(H))$ amortized, Theorem~\ref{thm:estimator_interface}}
\EndProcedure
\Statex
\Procedure{Center}{$\mathcal P,c,y;\mathsf S$}
\For{$k=0,\ldots,N_y-1$}
  \State $(\widehat g,\widehat R,\widehat Q)\gets\Call{BarrierOracle}{\mathcal P,c,y,\mathbf 1_{\{k=0\}};\mathsf S}$. \Comment{$y$-query}
  \If{$k=0$}
    \State $P\gets\alpha\widehat Q$, and keep $P$ and its Cholesky factor for this call. \Comment{$O(m^2)$}
       \label{line:more_phase_metric}
  \EndIf
  \State $y\gets y-\gamma_yP^{-1}(\widehat g+c)$, solving with the Cholesky factor of $P$. \Comment{$O(m^2)$}\label{line:more_center_step}
  \For{$\ell=1,\ldots,N_G$}
    \State $(\widehat g,\widehat R,\widehat Q)\gets \Call{BarrierOracle}{\mathcal P,c,y,0;\mathsf S}$. \Comment{$G$-query}
    \State $\mathcal D_{\rm weight}.\textsc{Update}(\widehat R)$. \Comment{$O(m^\omega)$}
  \EndFor
\EndFor
\State \Return $y$.
\EndProcedure
\end{algorithmic}
\end{algorithm}

\subsection{Guarantees of the data structures}
\label{sec:data_structures}
Algorithm~\ref{alg:data_structures} lists the parameters, members, invariants and procedures of $\mathcal D_{\rm slack}$, $\mathcal D_{\rm weight}$, $\mathcal D^{\rm f}_{\rm est}$ and $\mathcal D^{\rm m}_{\rm est}$, and Table~\ref{tab:procedures} in Section~\ref{sec:parameters} lists every procedure with its cost per call. The guarantees below hold at query states that the algorithm can reach on a successful run. We describe these states first.

\begin{algorithm}[!ht]
\caption{The data structures with their parameters, members, invariants and procedures}
\label{alg:data_structures}
\begin{algorithmic}[1]
\State \textbf{data structure} $\mathcal D_{\rm slack}$ \Comment{Theorem~\ref{thm:slack_interface}, Section~\ref{sec:reused_factorization}}
\State \hspace{1.5em}\textbf{members}, used if $m^4\leq n$
\State \hspace{3em}$c\in\R^m$ \Comment{objective of the anchor}
\State \hspace{3em}$L\in\R^{n\times n}$ \Comment{Cholesky factor of the anchor slack}
\State \hspace{3em}$S\in\mathbb S^n$ \Comment{slack of the current query}
\State \hspace{1.5em}\textbf{invariant} the objective of the latest query continues the anchor $(c,L)$
\State \hspace{1.5em}\textbf{procedure} $\textsc{Whiten}(y,c)$ \Comment{Algorithm~\ref{alg:implicit_slack}}
\Statex
\State \textbf{data structure} $\mathcal D_{\rm weight}$ \Comment{Theorem~\ref{thm:weight_interface}, Section~\ref{sec:cost_per_iter}}
\State \hspace{1.5em}\textbf{member}
\State \hspace{3em}$G\in\mathbb S^m_{++}$ \Comment{weight in the original coefficient basis}
\State \hspace{1.5em}\textbf{invariant} on the event of Theorem~\ref{thm:centering_correctness}, $G$ is warm at every query and
\State \hspace{1.5em}\phantom{\textbf{invariant}} accurate at the queries listed in (V3)
\State \hspace{1.5em}\textbf{procedures} $\textsc{Init}(y,\widehat J,\rho)$ at the first query \Comment{Algorithm~\ref{alg:random_weight_computation}}
\State \hspace{1.5em}\phantom{\textbf{procedures}} $\textsc{Update}(\widehat R)$ after each $G$-query \Comment{Algorithm~\ref{alg:weight_inverse_step}}
\Statex
\State \textbf{data structure} $\mathcal D^{\rm f}_{\rm est}$ \Comment{Used if $m^2\leq n$, Theorem~\ref{thm:estimator_interface}, Section~\ref{sec:randomized_oracles}}
\State \hspace{1.5em}\textbf{parameters} $N_{\max}$ and $\rho$
\State \hspace{1.5em}\textbf{member}
\State \hspace{3em}$\widehat Q\in\mathbb S^m_{++}$ with its Cholesky factor \Comment{proxy Hessian estimate of the phase}
\State \hspace{1.5em}\textbf{invariant} $\widehat Q$ is the estimate computed at the latest query with $a=1$
\State \hspace{1.5em}\textbf{procedure} $\textsc{Query}(y,G,\widehat J,a)$ \Comment{Algorithm~\ref{alg:fresh_query}}
\Statex
\State \textbf{data structure} $\mathcal D^{\rm m}_{\rm est}$\Comment{Used if $m^2>n$, Theorem~\ref{thm:estimator_interface}, Section~\ref{sec:oracle_maintenance}}
\State \hspace{1.5em}\textbf{parameters} $H$, $N_{\max}$ and $\rho$
\State \hspace{1.5em}\textbf{members}
\State \hspace{3em}$\widehat g\in\R^m$ and $\widehat R\in\mathbb S^m$ \Comment{maintained estimates}
\State \hspace{3em}$Z\in\mathbb S^m_{++}$ \Comment{reference metric of the phase}
\State \hspace{3em}$\widehat Q\in\mathbb S^m_{++}$ with its Cholesky factor \Comment{proxy Hessian estimate of the phase}
\State \hspace{3em}$y$, $G$ and $J$ \Comment{point, weight and whitener of the latest query}
\State \hspace{3em}$Y$ \Comment{sketched factors of the latest query}
\State \hspace{3em}$N_{\rm ep}\in\{1,\ldots,H\}$ \Comment{number of calls in the current epoch}
\State \hspace{1.5em}\textbf{invariant} $(\widehat g,\widehat R)$ is the pair returned at the latest query
\State \hspace{1.5em}\textbf{procedure} $\textsc{Query}(y,G,\widehat J,a)$ \Comment{Algorithm~\ref{alg:more_gradient}}
\end{algorithmic}
\end{algorithm}

\begin{definition}[Runs and valid queries]
\label{def:valid_queries}
A \emph{run} is a sequence of calls of \textnormal{\textsc{Center}} that share one oracle state $\mathsf S$. Each call starts at the point returned by the preceding call, and its objective vector $c$ is fixed when the call starts and may depend on the preceding history. A run may begin with one standalone query, which has the zero objective and is made at the starting point of its first call. Number the queries of a run by $j=1,2,\ldots$, and let $y_j$ and $G_j$ be the point and the weight $\mathcal D_{\rm weight}.G$ at query $j$. For a call that starts at $y^\circ$, its \emph{phase ball} is
\begin{equation}\label{eq:phase_ball}
 \mathfrak B(y^\circ):=\bigl\{y:\|y-y^\circ\|_{\nabla^2F(y^\circ)}\leq2/5\bigr\}.
\end{equation}
Query $j$ is \emph{valid} if the following three conditions hold.
\begin{itemize}
\item[(V1)] The call containing query $j$ starts at a point $y^\circ$ with $\lambda_c(y^\circ)\leq1/10$, and $y_j\in\mathfrak B(y^\circ)$. For the standalone query, (V1) requires only that $y_j$ is the starting point of the first call.
\item[(V2)] The weight $G_j$ is warm at $y_j$.
\item[(V3)] The weight $G_j$ is accurate at $y_j$ if query $j$ is a $y$-query, the standalone query, or the last $G$-query after a $y$-step.
\end{itemize}
\end{definition}

By Lemma~\ref{lem:sc_local_estimates}, every phase ball is feasible. In a run, the first query of a call has the same point as the last query of the preceding call, or as the standalone query. The objective of a query matters only for the anchors of $\mathcal D_{\rm slack}$. Because the standalone query has the zero objective, an anchor that it sets is continued only by calls with the zero objective. The three theorems below are proved in the sections listed in Algorithm~\ref{alg:data_structures}.

\begin{theorem}[Slack maintenance data structure]
\label{thm:slack_interface}
The deterministic data structure $\mathcal D_{\rm slack}$ supports the single procedure $\textnormal{\textsc{Whiten}}(y,c)$. Given a feasible point $y$ and the objective $c$ of the current call, this procedure forms $S(y)$ and provides the actions $X\mapsto\widehat JX$ and $X\mapsto\widehat J^\top X$ of a matrix $\widehat J$. It reports failure only if a Cholesky factorization of $S(y)$ fails. The following guarantees hold.
\begin{itemize}
\item[(a)] If $m^4>n$, the procedure factors $S(y)=LL^\top$ at every query and returns $\widehat J=L^{-1}$.
\item[(b)] Let $m^4\leq n$, and suppose that $F+\theta c^\top y$ has a minimizer for every objective $c$ of a run and every $\theta\in[1/2,2]$. If queries $1,\ldots,j$ of the run are valid, then query $j$ returns $\widehat J=\mathsf KJ$, where $J$ is an exact whitener of $S(y_j)$, the matrix $\mathsf K$ is symmetric, and $\|\mathsf K-I\|_{\rm op}\leq\varepsilon_{\rm w}$.
\item[(c)] Let $m^4\leq n$. Consider a run whose objectives are $tc^\circ$ for a fixed $c^\circ\neq0$, with multipliers $t$ that are monotone in $[t_{\min},t_{\max}]$, $t_{\min}>0$. Then $\mathcal D_{\rm slack}$ factors the slack at most $1+\lceil\log_2(t_{\max}/t_{\min})\rceil$ times. A run with zero objective factors it once.
\item[(d)] Forming $S(y)$ costs $O(mn^2)$, and each factorization costs $O(n^\omega)$. Applying either action to an $n\times q$ block costs $O(\Tmat(n,n,q))$ in Part~(a) and $O(d_{\rm w}\Tmat(n,n,q))$ in Part~(b).
\end{itemize}
\end{theorem}

The weight maintenance data structure must control both the bias caused by an approximate weight and the change in its accuracy after a point update. The next theorem collects these bounds together with the guarantees for initializing and updating the weight.

\begin{theorem}[Weight properties and maintenance]
\label{thm:weight_interface}
Let $C_{\rm pert}:=3$ and $c_{\rm pert}:=1/100$. The following three properties hold at every feasible $y$ and every weight $G\succ0$.
\begin{itemize}
\item[(a)] \emph{Warm weights.} If $\Delta_y(G)\leq1/(100s)$, then $\frac34G_*(y)\preceq G\preceq\frac54G_*(y)$, $\frac12I\preceq R_G(y)\preceq2I$ and $\tr[\mathcal L_y(G)^s]\leq2m$.
\item[(b)] \emph{Weight bias.} If $\Delta_y(G)\leq c_{\rm pert}/s$, then
\begin{gather}
 \|g_G(y)-\nabla F(y)\|_{\nabla^2F(y)^{-1}}
 \leq C_{\rm pert}\sqrt{\alpha s\Delta_y(G)},
 \label{eq:single_loop_weight_bias}\\
 \bigl(1-C_{\rm pert}\sqrt{s\Delta_y(G)}\bigr)^2Q_*(y)\preceq Q_G(y)
 \preceq\bigl(1+C_{\rm pert}\sqrt{s\Delta_y(G)}\bigr)^2Q_*(y).
 \label{eq:weight_metric_bias}
\end{gather}
\item[(c)] \emph{Motion of the gap.} Let $y(u)$, $u\in[0,1]$, be a piecewise $C^1$ feasible path with $F$-local length $d_F:=\int_0^1\|y'(u)\|_{\nabla^2F(y(u))}\,\d u$. If $\sqrt{\Delta_{y(0)}(G)}+C_{\rm pert}\sqrt{s/\alpha}\,d_F<\sqrt{c_{\rm pert}/s}$, then
\begin{equation}\label{eq:weight_gap_motion}
 \Bigl|\sqrt{\Delta_{y(1)}(G)}-\sqrt{\Delta_{y(0)}(G)}\Bigr|
 \leq C_{\rm pert}\sqrt{s/\alpha}\,d_F.
\end{equation}
\end{itemize}
The data structure $\mathcal D_{\rm weight}$ stores one weight and supports the procedures \textnormal{\textsc{Update}} and \textnormal{\textsc{Init}}, with the following guarantees.
\begin{itemize}
\item[(d)] \emph{Weight update.} Suppose that $\mathcal D_{\rm weight}.G=G$ with $\Delta_y(G)\leq1/(100s)$. Given a symmetric matrix $\widehat R$ satisfying $\|\widehat R-R_G(y)\|_F\leq\tau$, the procedure $\textnormal{\textsc{Update}}(\widehat R)$ succeeds, and the new weight $G^+$ satisfies $(1-\gamma_G)G\preceq G^+\preceq(1+\gamma_G)G$ and
\begin{equation}\label{eq:weight_correction_recurrence}
 \Delta_y(G^+)\leq(1-\gamma_G/16)\Delta_y(G)+6\gamma_G\tau^2.
\end{equation}
The procedure uses $O(m^\omega)$ arithmetic operations.
\item[(e)] \emph{Initialization.} Given a feasible point $y$, an exact whitener $\widehat J$ of $S(y)$, and $0<\rho<1/2$, the procedure $\textnormal{\textsc{Init}}(y,\widehat J,\rho)$ fixes the coefficient basis and computes a weight that is accurate at $y$ with probability at least $1-\rho$. It uses $\widetilde O(n^\omega+\mathcal C(m))$ arithmetic operations on every run.
\end{itemize}
\end{theorem}

The remaining data structure estimates the fixed-weight gradient, stationarity matrix and proxy Hessian from the whitening actions and the current weight. We state the accuracy and cost guarantees used by the path-following analysis.

\begin{theorem}[Gradient/Hessian estimation data structure]
\label{thm:estimator_interface}
Fix an integer $1\leq H\leq\sqrt m$, a query cap $N_{\max}\geq1$ and a failure budget $0<\rho<1/10$. The data structure $\mathcal D_{\rm est}$ is $\mathcal D^{\rm f}_{\rm est}$ if $m^2\leq n$ and $\mathcal D^{\rm m}_{\rm est}$ if $m^2>n$. It is created with these parameters, where only $\mathcal D^{\rm m}_{\rm est}$ uses $H$. At query $j$ of a run, its single procedure $\textnormal{\textsc{Query}}(y_j,G_j,\widehat J_j,a_j)$ takes the query state, the whitening actions returned by $\mathcal D_{\rm slack}$, and the flag $a_j$, which equals one at the first query of a call. It returns a vector $\widehat g_j$, a symmetric matrix $\widehat R_j$, and a positive definite matrix $\widehat Q_j$ with its Cholesky factor. The following guarantees hold.
\begin{itemize}
\item[(a)] \emph{Accuracy.} There is an event of probability at least $1-\rho$ on which the following holds for every $j$. Suppose that queries $1,\ldots,j$ are valid and that each of them received $\widehat J=\mathsf KJ$, where $J$ is an exact whitener, $\mathsf K$ is symmetric, $\|\mathsf K-I\|_{\rm op}\leq\varepsilon_{\rm w}$, and $\mathsf K=I$ if $m^2>n$. Then
\begin{equation}\label{eq:estimator_accuracy}
 \|\widehat g_j-g_{G_j}(y_j)\|_{\nabla^2F(y_j)^{-1}}\leq\varepsilon_{\rm g}/4,
 \qquad
 \|\widehat R_j-R_{G_j}(y_j)\|_F\leq\tau.
\end{equation}
If, in addition, $a_j=1$, then
\begin{equation}\label{eq:estimator_metric}
 0.9\,Q_{G'}(y_j)\preceq\widehat Q_j\preceq1.1\,Q_{G'}(y_j),
\end{equation}
where $G'=G_j$, or $G'=G_{j-1}$ and $y_{j-1}=y_j$.
\item[(b)] \emph{Cost.} On every run, every prefix of $N\leq N_{\max}$ queries costs $\widetilde O(\mathcal C(m)+N\mathcal B_*(H))$ operations, including the applications of the whitening actions.
\end{itemize}
\end{theorem}

The guarantees of $\mathcal D_{\rm slack}$ and $\mathcal D_{\rm est}$ apply at valid queries, while the validity of a query depends on the accuracy of the preceding answers. Section~\ref{sec:path_following} resolves this dependence by an induction over the queries of a run (Theorem~\ref{thm:centering_correctness}). The data structure $\mathcal D_{\rm est}$ returns $\widehat Q_j$ at every query, and \textsc{Center} uses it only when $a_j=1$.

\subsection{Parameters, cost functions and failure budgets}
\label{sec:parameters}
The path-following parameters are
\begin{equation}\label{eq:more_centering_parameters}
 \begin{gathered}
 \delta_{\rm c}:=\frac1{100},\qquad \gamma_t:=\frac{\delta_{\rm c}}{4\sqrt\nu},
 \qquad \mu:=\frac18,\qquad e_0:=\frac1{1600},
 \qquad \varepsilon_{\rm g}:=\frac{e_0}{\sqrt{32s}},\\
 E:=\frac{\delta_{\rm c}^2}{16},\qquad d_0:=0.0054,\qquad
 N_y:=\Bigl\lceil\frac1{\mu\gamma_y}\log\frac{4d_0}E\Bigr\rceil,
 \qquad N_{\rm c}:=(N_G+1)N_y.
 \end{gathered}
\end{equation}
The weight-correction parameters are
\begin{equation}\label{eq:single_loop_parameters}
 \begin{gathered}
 \gamma_G:=\frac1{1000s},\qquad \gamma_y:=\frac{c_y}s,\qquad
 \eta_{\rm w}:=\frac{c_{\rm w}}s,\qquad
 \eta:=\frac{\eta_0}{\alpha s^2},\qquad
 \tau:=\frac{\sqrt\eta}{32},\\
 N_G:=\left\lceil\frac{16}{\gamma_G}
                 \log\frac{4\eta_{\rm w}}\eta\right\rceil,
 \end{gathered}
\end{equation}
where $c_{\rm w},\eta_0,c_y$ are sufficiently small absolute constants, chosen in this order. First choose $c_{\rm w}\leq c_{\rm pert}/4$, with $c_{\rm pert}$ from Theorem~\ref{thm:weight_interface}. Then choose $\eta_0$ and $c_y$ such that $\eta\leq\eta_{\rm w}/16$, $c_y\leq1/200$ and
\begin{equation}\label{eq:alternating_accuracy_parameters}
 C_{\rm pert}\sqrt{\alpha s\eta}\leq\frac{\varepsilon_{\rm g}}4,
 \qquad
 \sqrt\eta+20C_{\rm pert}c_y\sqrt{s/\alpha}
 \leq\tfrac12\sqrt{\eta_{\rm w}}.
\end{equation}
These choices are possible because $\alpha s\eta=\eta_0/s$ and $\sqrt{s/\alpha}=10^{-3}s^{-1/2}$. Consequently $\eta=\Theta(s^{-4})$, $\gamma_y=\Theta(1/s)$, $N_G=O(s\log s)$, $N_y=O(s)$ and $N_{\rm c}=O(s^2\log s)$. All tolerances and counts are polylogarithmic in the dimensions. The weight gaps enter only the analysis, and the algorithm uses only the step sizes and the loop lengths.

The tolerances form a chain, which explains their sizes,
\begin{align*}
 \delta_{\rm c}=\Theta(1)\ \Longrightarrow\
 \varepsilon_{\rm g}=\Theta(s^{-1/2})\ \Longrightarrow\
 \eta=\Theta(s^{-4})\ \Longrightarrow\
 \tau=\Theta(s^{-2}).
\end{align*}
A $y$-step reaches the centering tolerance $\delta_{\rm c}$ when its gradient error in the norm of $P^{-1}$ is at most the constant $e_0$. The phase preconditioner is within a factor $O(s)$ of the Hessian, which turns this requirement into the tolerance $\varepsilon_{\rm g}=e_0/\sqrt{32s}$ in the dual norm of the Hessian (Lemma~\ref{lem:random_direct_centering}). The weight bias~\eqref{eq:single_loop_weight_bias} stays below $\varepsilon_{\rm g}/4$ when the gap is at most $\eta$, because $\alpha=\Theta(s^2)$. The $G$-steps keep the gap below $\eta$ when the stationarity estimates have error $\tau=\Theta(\sqrt\eta)$, because the fixed point of~\eqref{eq:weight_correction_recurrence} is $96\tau^2=3\eta/32$. Finally, the warm gap $\eta_{\rm w}=\Theta(1/s)$ keeps the spectral bounds of Theorem~\ref{thm:weight_interface}(a), and the step size $\gamma_y=\Theta(1/s)$ keeps the gap below $\eta_{\rm w}$ along a $y$-step.

The slack maintenance data structure uses the window and the polynomial degree
\begin{equation}\label{eq:whitening_degree}
 \Upsilon:=4m\sqrt\alpha,\qquad
 d_{\rm w}:=\Bigl\lceil\frac{\Upsilon+1}2
            \log\frac\Upsilon{\varepsilon_{\rm w}}\Bigr\rceil
 =\widetilde O(m),
\end{equation}
where $\varepsilon_{\rm w}$ is the whitening tolerance fixed in~\eqref{eq:bias_budget}. This tolerance is inverse-polynomial in $2mn$, which gives the last bound. Table~\ref{tab:parameters} summarizes the roles of the parameters.

\begin{table}[!ht]
\centering\small
\begin{tabular}{@{}lll@{}}
\toprule
Parameters & Role & Defined in\\
\midrule
$s,\alpha,\nu$ & trace power, scaling and parameter of the barrier & Definition~\ref{def:ml_s_choice}\\
$\delta_{\rm c},\gamma_t$ & centering tolerance and relative multiplier step & \eqref{eq:more_centering_parameters}\\
$\varepsilon_{\rm g}$ & gradient accuracy needed by a $y$-step & \eqref{eq:more_centering_parameters}\\
$\eta,\eta_{\rm w}$ & weight gaps of accurate and warm weights & \eqref{eq:single_loop_parameters}\\
$\tau$ & accuracy of the stationarity estimates & \eqref{eq:single_loop_parameters}\\
$\gamma_y,\gamma_G$ & sizes of $y$-steps and $G$-steps & \eqref{eq:single_loop_parameters}\\
$N_y,N_G,N_{\rm c}$ & $y$-steps per call, $G$-steps per $y$-step, queries per call & \eqref{eq:more_centering_parameters}, \eqref{eq:single_loop_parameters}\\
$\Upsilon,d_{\rm w},\varepsilon_{\rm w}$ & slack window, polynomial degree, whitening tolerance & \eqref{eq:whitening_degree}, \eqref{eq:bias_budget}\\
$H,N_{\max},\gamma$ & epoch length, query cap, failure budget & inputs\\
\bottomrule
\end{tabular}
\caption{Parameters shared by the algorithm and its data structures.}
\label{tab:parameters}
\end{table}

The running times use the following cost functions. Besides $\mathcal C(q)$ of Definition~\ref{def_more_block_cost}, let
\begin{align}
 \mathcal C_{\rm sep}(q)&:=\Tmat(n,n,mq)+\Tmat(m,m,nq)+mnq,
 \label{eq_more_separate_cost}\\
 \mathcal B(H)&:=n^\omega+\mathcal C(H\sqrt m)+\frac{\mathcal C(m)}H,
 \label{eq_more_batched_cost}\\
 \mathcal B_{\rm imp}&:=mn^2+m^\omega+\Tmat(mn,n,m)
                         +d_{\rm w}\mathcal C_{\rm sep}(m),
 \label{eq:implicit_batched_cost}
\end{align}
and define the per-query charge
\begin{equation}\label{eq:effective_query_cost}
 \mathcal B_*(H):=
 \begin{cases}
  \mathcal B_{\rm imp},&m^4\leq n,\\
  n^\omega+\mathcal C(m),&m^2\leq n<m^4,\\
  \mathcal B(H),&m^2>n.
 \end{cases}
\end{equation}
Thus $\mathcal C(q)=\Tmat(mn,n,q)+\mathcal C_{\rm sep}(q)$, where the first term applies the stacked coefficients. In $\mathcal D^{\rm m}_{\rm est}$, an epoch of $H$ phases computes one full estimate at width $\widetilde O(m)$ and corrections at width $\widetilde O(H\sqrt m)$, and $\mathcal B(H)$ is the resulting amortized cost of a query. The charge $\mathcal B_{\rm imp}$ bounds a query of $\mathcal D^{\rm f}_{\rm est}$ with polynomial whitening, in which the degree $d_{\rm w}$ multiplies only the products with the whitening actions. The first two cases of $\mathcal B_*(H)$ use fresh estimates, and only the third depends on $H$. The last two cases include a slack factorization at every query, while the first case leaves the anchor factorizations to be counted separately.

Table~\ref{tab:procedures} lists every procedure in the algorithm and its data structures, along with their callers, per-call cost, and the result that justifies the cost. Preprocessing, recovery and the complete solver of Sections~\ref{sec:preprocessing}--\ref{sec:solver_wrapup} appear only as callers.

\begin{table}[!ht]
\centering\small
\begin{tabularx}{\textwidth}{@{}>{\raggedright\arraybackslash}p{4.1cm}>{\raggedright\arraybackslash}p{3.6cm}>{\raggedright\arraybackslash}X>{\raggedright\arraybackslash}p{2.4cm}@{}}
\toprule
Procedure & Called by & Cost per call & Proved in\\
\midrule
\textsc{AmortizedSDP} & the complete solver & $N_{\rm path}+1$ calls of \textsc{Center} & Theorem~\ref{thm_more_amortized_sdp}\\
\textsc{InitializeOracleState} & \textsc{AmortizedSDP}, preprocessing & $O(1)$ & Algorithm~\ref{alg:maintained_center}\\
\textsc{Center} & \textsc{AmortizedSDP}, preprocessing & $N_{\rm c}$ queries and $O(N_{\rm c}m^\omega)$ & Lemma~\ref{lem:run_cost}\\
\textsc{BarrierOracle} & \textsc{Center}, preprocessing & one \textsc{Whiten}, one \textsc{Query}, and \textsc{Init} at the first query & Algorithm~\ref{alg:maintained_center}\\
\midrule
$\mathcal D_{\rm slack}.\textsc{Whiten}$ & \textsc{BarrierOracle} & $O(mn^2)$, and $O(n^\omega)$ if it factors & Theorem~\ref{thm:slack_interface}(d)\\
\midrule
$\mathcal D_{\rm weight}.\textsc{Init}$ & \textsc{BarrierOracle} & $\widetilde O(n^\omega+\mathcal C(m))$ & Theorem~\ref{thm:weight_interface}(e)\\
\mbox{\textsc{RandomizedWeights}} & $\mathcal D_{\rm weight}.\textsc{Init}$ & $\widetilde O(\eta^{-1}(n^\omega+\mathcal C(m)))$ & Lemma~\ref{lem:weight_initialization}\\
$\mathcal D_{\rm weight}.\textsc{Update}$ & \textsc{Center} & $O(m^\omega)$ & Theorem~\ref{thm:weight_interface}(d)\\
\textsc{UpdateWeight} & \mbox{\textsc{RandomizedWeights}} and $\mathcal D_{\rm weight}.\textsc{Update}$ & $O(m^\omega)$ & Algorithm~\ref{alg:weight_inverse_step}\\
\midrule
$\mathcal D^{\rm f}_{\rm est}.\textsc{Query}$ & \textsc{BarrierOracle} & $\widetilde O(n^\omega+\mathcal C(m))$, or $\widetilde O(\mathcal B_{\rm imp})$ with polynomial whitening & Section~\ref{sec:fresh_mode}, Lemma~\ref{lem:implicit_query_cost}\\
$\mathcal D^{\rm m}_{\rm est}.\textsc{Query}$ & \textsc{BarrierOracle} & $\widetilde O(\mathcal B(H))$ amortized, and $\widetilde O(\mathcal C(m))$ once & Lemma~\ref{lem_more_gradient_tracking}(b)\\
\textsc{ProxyHessian} & both \textsc{Query} procedures, recovery & $\widetilde O(n^\omega+\mathcal C(q))$ for sketched factors of width $q$ & Lemma~\ref{lem_more_metric}(b)\\
\bottomrule
\end{tabularx}
\caption{Procedures of the algorithm and of its data structures, grouped by owner. The actions returned by \textsc{Whiten} are paid when $\mathcal D_{\rm est}$ applies them, at the rates of Theorem~\ref{thm:slack_interface}(d).}
\label{tab:procedures}
\end{table}

\paragraph{Failure budgets.} \textsc{InitializeOracleState} receives a failure budget $\gamma$. The initialization $\mathcal D_{\rm weight}.\textsc{Init}$ receives $\gamma/3$, and the data structure $\mathcal D_{\rm est}$ receives $\rho=2\gamma/3$, which Sections~\ref{sec:randomized_oracles} and~\ref{sec:oracle_maintenance} divide among its random events. The data structure $\mathcal D_{\rm slack}$ is deterministic. Every routine reports failure if a required positive definite factorization, solve or stated check fails. Query counts, widths, trial counts and polynomial degrees are fixed in advance, so every cost bound also covers failed executions.
\section{Correctness of path following}
\label{sec:path_following}
This section proves Part~(a) of Theorem~\ref{thm_more_amortized_sdp} from the guarantees of Theorems~\ref{thm:slack_interface}, \ref{thm:weight_interface} and~\ref{thm:estimator_interface}. We first recall the central path and the certificate for the objective error. We then analyze one centering phase with accurate oracle answers, and finally show by induction over the queries of a run that the answers are accurate. Throughout, $\phi_c=F+c^\top y$ and $\lambda_c$ are the penalized objective and its decrement from Section~\ref{sec:query_quantities}. For the multiplier $t$ of the algorithm, the objective is $c=tb$, and $y$ is centered to tolerance $\delta$ if $\lambda_{tb}(y)\leq\delta$. Lemma~\ref{lem:sc_local_estimates} collects the local estimates for self-concordant functions that we use.

\subsection{The central path and the output certificate}
We first verify that the penalized objective has a unique minimizer at every positive multiplier. These minimizers define the central path followed by the algorithm.

\begin{lemma}[Existence of the central path]
\label{lem:ml_central_minimizer}
Assume primal and dual strict feasibility. For every $t>0$, the function $\phi_{tb}$ has a unique minimizer $y_t$ in $\Omega$.
\end{lemma}
\begin{proof}
We show that $\phi_{tb}$ tends to infinity both at the finite boundary of $\Omega$ and as $\|y\|_2\to\infty$, so that its sublevel sets are compact. Fix a primal feasible $X_0\succ0$, and let $A(h):=\sum_ih_iA_i$. Since the $A_i$ are linearly independent, there is a constant $c>0$ with $\|A(h)\|_F^2\geq c\|h\|_2^2$. Moreover, $\|S^{-1/2}XS^{-1/2}\|_F\geq\|S\|_{\rm op}^{-1}\|X\|_F$ for every symmetric $X$. Together with $h^\top\mathcal L_y^*(I)h=\|S(y)^{-1/2}A(h)S(y)^{-1/2}\|_F^2$, this gives
\begin{align*}
 \mathcal L_y^*(I)\succeq c\|S(y)\|_{\rm op}^{-2}I_m.
\end{align*}
Since $X_0$ is primal feasible, $b_i=\langle A_i,X_0\rangle$, and therefore
\begin{align*}
 b^\top y=\tr[X_0S(y)]+\langle C,X_0\rangle
 \geq\lambda_{\min}(X_0)\|S(y)\|_{\rm op}+\langle C,X_0\rangle,
\end{align*}
where the inequality uses $S(y)\succ0$. We evaluate~\eqref{eq:ml_weight_maximum} at $W=I$ and combine these two bounds. This gives
\begin{align*}
 \phi_{tb}(y)\geq
 t\lambda_{\min}(X_0)\|S(y)\|_{\rm op}
 -\alpha m\log\|S(y)\|_{\rm op}+c_t
\end{align*}
for a constant $c_t$ that depends only on the data and on $t$. The right side tends to infinity with $\|S(y)\|_{\rm op}$, which tends to infinity with $\|y\|_2$ because $A$ is injective. At the finite boundary points of $\Omega$, the barrier $F$ diverges by Theorem~\ref{thm:ml_matrix_sc}. Hence every sublevel set of $\phi_{tb}$ is compact and contained in $\Omega$, and $\phi_{tb}$ has a minimizer. The minimizer is unique because the Hessian of $\phi_{tb}$ is positive definite.
\end{proof}

The gradient matrix gives a primal certificate at each central point. The following identities also express its constraint residual at a general feasible point in terms of the gradient of the penalized objective.

\begin{lemma}[Central primal certificate]
\label{lem:ml_central_certificate}
At a feasible $y$, let $M_y$ be the matrix of Definition~\ref{def:ml_certificate_matrices}, and let
\begin{align*}
 X_t(y):=\frac\alpha t S(y)^{-1/2}M_yS(y)^{-1/2}.
\end{align*}
Then $X_t(y)\succ0$ and
\begin{equation}\label{eq:ml_certificate_residual}
 \bigl(\langle A_i,X_t(y)\rangle-b_i\bigr)_{i=1}^m
       =-\frac1t\nabla\phi_{tb}(y),\qquad
 \langle S(y),X_t(y)\rangle=\frac\nu t.
\end{equation}
In particular, at the minimizer $y_t$ of Lemma~\ref{lem:ml_central_minimizer}, the matrix $X_t(y_t)$ is primal feasible and has duality gap $\nu/t$.
\end{lemma}
\begin{proof}
Positive definiteness follows from $M_y\succ0$. The gradient formula $\nabla F(y)_i=-\alpha\tr[M_yB_i]$ of~\eqref{eq:ml_gradient} and $\tr[M_y]=m$ give
\begin{align*}
 \langle A_i,X_t(y)\rangle=-\nabla F(y)_i/t,\qquad
 \langle S(y),X_t(y)\rangle=\alpha m/t=\nu/t.
\end{align*}
Since $\nabla\phi_{tb}(y)=\nabla F(y)+tb$, the first identity is~\eqref{eq:ml_certificate_residual}. At $y_t$, the gradient $\nabla\phi_{tb}(y_t)$ vanishes, so the residual vanishes and the complementarity $\langle S(y_t),X_t(y_t)\rangle$ equals the duality gap.
\end{proof}

\subsection{Centrality and objective accuracy}
Combining Lemma~\ref{lem:ml_central_certificate} with the local self-concordance estimates converts a small decrement into an objective-error bound. The same estimates control the increase in decrement when the multiplier is increased and hence determine the required number of centering phases.

\begin{corollary}[Centrality, objective accuracy, and phase count]
\label{cor:ml_iterations}
Assume primal and dual strict feasibility, and let $t_+:=(1+\gamma_t)t$ with $\gamma_t=\delta_{\rm c}/(4\sqrt\nu)$. If $\lambda_{tb}(y)\leq\delta_{\rm c}$, then
\begin{equation}\label{eq:central_accuracy}
 b^\top y-\mathrm{OPT}\leq\frac{2\nu}t,\qquad
 \lambda_{t_+b}(y)<\frac1{10}.
\end{equation}
Consequently, consider a centering procedure that restores decrement at most $\delta_{\rm c}$ after each increase of the multiplier, and start from a point with decrement at most $\delta_{\rm c}$ at $t_0>0$. Then the dual objective error is at most $\epsilon$ after
\begin{align*}
 O\Bigl(\sqrt\nu\log\Bigl(2+\frac\nu{t_0\epsilon}\Bigr)\Bigr)
 =O\Bigl(\sqrt m\log(2mn)
       \log\Bigl(2+\frac\nu{t_0\epsilon}\Bigr)\Bigr)
\end{align*}
centering phases.
\end{corollary}
\begin{proof}
We bound the size of the objective in the local norm and then compare $y$ with the central point $y_t$. Since $F$ is a $\nu$-self-concordant barrier, $\|\nabla F(y)\|_{\nabla^2F(y)^{-1}}\leq\sqrt\nu$, and the triangle inequality gives
\begin{align*}
 \|tb\|_{\nabla^2F(y)^{-1}}
 \leq\|\nabla F(y)\|_{\nabla^2F(y)^{-1}}+\lambda_{tb}(y)
 \leq\sqrt\nu+\delta_{\rm c}.
\end{align*}
The increase of the multiplier adds $\gamma_ttb$ to the gradient of the penalized objective. Therefore
\begin{align*}
 \lambda_{t_+b}(y)
 \leq\lambda_{tb}(y)+\gamma_t\|tb\|_{\nabla^2F(y)^{-1}}
 \leq\delta_{\rm c}+\frac{\delta_{\rm c}}{4\sqrt\nu}(\sqrt\nu+\delta_{\rm c})
 <\frac1{10}.
\end{align*}
For the objective error, the primal feasible matrix $X_t(y_t)$ of Lemma~\ref{lem:ml_central_certificate} shows that $b^\top y_t-\mathrm{OPT}\leq\nu/t$. Moreover, \eqref{eq:sc_center_estimates} gives $\|y-y_t\|_{\nabla^2F(y)}\leq\delta_{\rm c}/(1-\delta_{\rm c})$. Hence the Cauchy--Schwarz inequality in the local norm gives
\begin{align*}
 b^\top y-\mathrm{OPT}
 &\leq b^\top y_t-\mathrm{OPT}
 +\frac1t\|tb\|_{\nabla^2F(y)^{-1}}\|y-y_t\|_{\nabla^2F(y)}\\
 &\leq\frac1t\Bigl(\nu+(\sqrt\nu+\delta_{\rm c})
       \frac{\delta_{\rm c}}{1-\delta_{\rm c}}\Bigr)
 \leq\frac{2\nu}t.
\end{align*}
It therefore suffices to reach $t\geq2\nu/\epsilon$. Each phase multiplies $t$ by $1+\gamma_t$, and $\log(1+\gamma_t)\geq\gamma_t/2$, so $O(\sqrt\nu\log(2+\nu/(t_0\epsilon)))$ phases suffice. Finally, $\sqrt\nu=1000s\sqrt m=O(\sqrt m\log(2mn))$ by Definition~\ref{def:ml_s_choice}.
\end{proof}

\subsection{Damped \texorpdfstring{$y$}{y}-steps with a fixed preconditioner}
\label{sec:random_direct_centering}
A constant-factor approximation to $\alpha Q_y$ has condition number $O(s)$ relative to the Hessian, and self-concordance preserves this comparison in a phase ball. The next lemma analyzes one $y$-step with such a preconditioner. The $G$-steps take place between $y$-steps and are analyzed in Section~\ref{sec:weight_correction_analysis}.

\begin{lemma}[Damped $y$-steps with a fixed preconditioner]
\label{lem:random_direct_centering}
Let $\lambda_c(y^\circ)\leq1/10$, let $y_c$ be the minimizer of $\phi_c$, and let $d(y):=\phi_c(y)-\phi_c(y_c)$ be the objective gap. Let $P$ satisfy
\begin{equation}\label{eq:phase_metric_sandwich}
 \tfrac12P\preceq\nabla^2F(y^\circ)\preceq(4s+2)P,
\end{equation}
and let $0<\gamma_y\leq1/(32s)$.
\begin{itemize}
\item[(a)] The minimizer $y_c$ exists, $d(y^\circ)<d_0$, and every $y\in\mathfrak B(y^\circ)$ with $d(y)\leq E/2$ satisfies $\lambda_c(y)\leq\delta_{\rm c}$.
\item[(b)] Let $y\in\mathfrak B(y^\circ)$ satisfy $d(y)\leq d_0$, and let $\|\widehat g-\nabla F(y)\|_{\nabla^2F(y)^{-1}}\leq\varepsilon_{\rm g}$. Then the segment from $y$ to $y^+:=y-\gamma_yP^{-1}(\widehat g+c)$ lies in $\mathfrak B(y^\circ)$ and has $F$-local length less than $20s\gamma_y$. Moreover,
\begin{equation}\label{eq:random_centering_contraction}
 d(y^+)\leq(1-\mu\gamma_y)d(y)+\gamma_ye_0^2/2\leq d_0.
\end{equation}
\item[(c)] Consider $N_y$ steps that start at $y^\circ$ and use estimates as in Part~(b). Then the final point $y$ satisfies $\lambda_c(y)\leq\delta_{\rm c}$.
\end{itemize}
\end{lemma}
\begin{proof}
We compare $\phi_c$ with quadratic models on the ball $\mathfrak B:=\mathfrak B(y^\circ)$, where the Hessian is comparable to $P$, and use a lower bound on the gap at the boundary of $\mathfrak B$ to keep every step inside it.

\emph{Geometry of the ball.} Since $\lambda_c(y^\circ)\leq1/10<1$, Lemma~\ref{lem:sc_local_estimates} shows that $y_c$ exists, and~\eqref{eq:sc_center_estimates} gives
\begin{equation}\label{eq:phase_center_distance}
 \|y_c-y^\circ\|_{\nabla^2F(y^\circ)}\leq\frac19,\qquad
 d(y^\circ)\leq-\frac1{10}-\log\frac9{10}<d_0.
\end{equation}
Every $y\in\mathfrak B$ has local distance at most $2/5$ from $y^\circ$, so~\eqref{eq:sc_hessian_comparison} gives $(3/5)^2\nabla^2F(y^\circ)\preceq\nabla^2F(y)\preceq(5/3)^2\nabla^2F(y^\circ)$. Combining this with~\eqref{eq:phase_metric_sandwich} and $4s+2\leq5s$ yields
\begin{equation}\label{eq:ball_metric_comparison}
 \mu P\preceq\nabla^2F(y)\preceq32sP\qquad(y\in\mathfrak B).
\end{equation}
Let $y\in\partial\mathfrak B$. The segment from $y_c$ to $y$ lies in $\mathfrak B$, and its length in the norm of $\nabla^2F(y^\circ)$ is at least $2/5-1/9$ by~\eqref{eq:phase_center_distance}. Since $\nabla\phi_c(y_c)=0$, Taylor's formula and the lower Hessian bound on $\mathfrak B$ give
\begin{equation}\label{eq:ball_boundary_gap}
 d(y)\geq\tfrac12(3/5)^2(2/5-1/9)^2>0.015\qquad(y\in\partial\mathfrak B).
\end{equation}

\emph{Part~(a).} The first two claims are~\eqref{eq:phase_center_distance}. Let $y\in\mathfrak B$ satisfy $d(y)\leq E/2$, and let $\lambda:=\lambda_c(y)$. The lower bound $\lambda-\log(1+\lambda)\leq d(y)$ of~\eqref{eq:sc_center_estimates} first implies $\lambda<1$, since $1-\log2>E/2$. For $0\leq\lambda<1$, we have $\lambda-\log(1+\lambda)\geq\lambda^2/(2(1+\lambda))\geq\lambda^2/4$, and so $\lambda^2\leq2E=\delta_{\rm c}^2/8$. Thus $\lambda\leq\delta_{\rm c}$, which proves Part~(a).

\emph{Part~(b).} Write $g:=\nabla\phi_c(y)$ and $\xi:=\widehat g-\nabla F(y)$, so that $\widehat g+c=g+\xi$. Since $P^{-1}\preceq32s\nabla^2F(y)^{-1}$ by~\eqref{eq:ball_metric_comparison},
\begin{equation}\label{eq:step_noise_metric}
 \|\xi\|_{P^{-1}}\leq\sqrt{32s}\,\varepsilon_{\rm g}=e_0.
\end{equation}
Let $0\leq u\leq1$ be such that the segment from $y$ to $y+u(y^+-y)$ lies in $\mathfrak B$. On this segment, $\nabla^2F\preceq32sP\preceq\gamma_y^{-1}P$, so the change of $\phi_c$ is at most the quadratic model $ug^\top(y^+-y)+u^2\|y^+-y\|_P^2/(2\gamma_y)$. Substituting $y^+-y=-\gamma_yP^{-1}(g+\xi)$ and completing the square in this model gives
\begin{equation}\label{eq:random_fractional_descent}
 \phi_c(y+u(y^+-y))-\phi_c(y)
 \leq-\frac{u(2-u)\gamma_y}{2}
   \Bigl\|g+\frac{1-u}{2-u}\xi\Bigr\|_{P^{-1}}^2
   +\frac{u\gamma_y}{2(2-u)}\|\xi\|_{P^{-1}}^2.
\end{equation}

We first show that the segment stays in $\mathfrak B$. Suppose that it leaves $\mathfrak B$, and let $u\leq1$ be the fraction at which it first reaches $\partial\mathfrak B$. Then~\eqref{eq:random_fractional_descent} and~\eqref{eq:step_noise_metric} bound the gap at the boundary point by $d_0+\gamma_ye_0^2/2<0.006$, which contradicts~\eqref{eq:ball_boundary_gap}. Hence the whole segment lies in $\mathfrak B$.

Next we prove the contraction. Since $y,y_c\in\mathfrak B$ and $\nabla^2F\succeq\mu P$ on $\mathfrak B$, strong convexity gives $\|g\|_{P^{-1}}^2\geq2\mu d(y)$. \eqref{eq:random_fractional_descent} with $u=1$ and~\eqref{eq:step_noise_metric} then give the first inequality in~\eqref{eq:random_centering_contraction}. The fixed point of this recurrence is $e_0^2/(2\mu)=E/4<d_0$, which gives the second inequality.

Finally we bound the length of the step. Integrating the Hessian from $y_c$ to $y$ inside $\mathfrak B$, where $\nabla^2F\preceq(25/9)\nabla^2F(y^\circ)$, gives
\begin{align*}
 \|g\|_{\nabla^2F(y^\circ)^{-1}}
 \leq\tfrac{25}9\|y-y_c\|_{\nabla^2F(y^\circ)}
 \leq\tfrac{25}9\bigl(\tfrac25+\tfrac19\bigr)<1.42.
\end{align*}
Since $\nabla^2F(y^\circ)\preceq(4s+2)P\preceq5sP$, both $\|v\|_{P^{-1}}\leq\sqrt{5s}\|v\|_{\nabla^2F(y^\circ)^{-1}}$ and $\|P^{-1}v\|_{\nabla^2F(y^\circ)}\leq\sqrt{5s}\|v\|_{P^{-1}}$ hold for every $v$. Together with~\eqref{eq:step_noise_metric}, these bounds give
\begin{align*}
 \|\widehat g+c\|_{P^{-1}}
 &\leq1.42\sqrt{5s}+e_0\leq4\sqrt s,\\
 \|y^+-y\|_{\nabla^2F(y^\circ)}
 &\leq\gamma_y\sqrt{5s}\,\|\widehat g+c\|_{P^{-1}}<9s\gamma_y.
\end{align*}
The segment lies in $\mathfrak B$, where the local norm is at most $5/3$ times the norm of $\nabla^2F(y^\circ)$. Its $F$-local length is therefore less than $15s\gamma_y<20s\gamma_y$. This proves Part~(b).

\emph{Part~(c).} By Part~(b) and induction, every iterate $y_k$ lies in $\mathfrak B$ and satisfies $d(y_k)\leq d_0$. Iterating~\eqref{eq:random_centering_contraction} gives
\begin{align*}
 d(y_k)\leq(1-\mu\gamma_y)^kd(y^\circ)+\frac E4
 \leq e^{-\mu\gamma_yk}d_0+\frac E4.
\end{align*}
The choice of $N_y$ in~\eqref{eq:more_centering_parameters} makes the first term at most $E/4$ for $k=N_y$, so $d(y_{N_y})\leq E/2$. Part~(a) now gives $\lambda_c(y_{N_y})\leq\delta_{\rm c}$, which proves Part~(c).
\end{proof}

\subsection{Restoring the weight by \texorpdfstring{$G$}{G}-steps}
\label{sec:weight_correction_analysis}
\label{sec:single_loop_tracking}
A call of \textsc{Center} takes each $y$-step with an accurate weight. The $y$-step leaves a warm weight at the new point, and the following $N_G$ $G$-steps at the fixed point restore accuracy. We first show that the proxy Hessian estimate of $\mathcal D_{\rm est}$ gives an admissible phase preconditioner.

\begin{lemma}[The phase preconditioner]
\label{lem:phase_metric}
Let $y$ be feasible, let $G'$ be accurate at $y$, and let $0.9\,Q_{G'}(y)\preceq\widehat Q\preceq1.1\,Q_{G'}(y)$. Then $P:=\alpha\widehat Q$ satisfies~\eqref{eq:phase_metric_sandwich} with $y^\circ=y$.
\end{lemma}
\begin{proof}
We compare $\widehat Q$ with $Q_*(y)$ through $Q_{G'}(y)$ and then apply the Hessian comparison of Lemma~\ref{lem:ml_gradient_certificates}. Since $\Delta_y(G')\leq\eta\leq c_{\rm pert}/s$, \eqref{eq:weight_metric_bias} applies. Moreover, \eqref{eq:alternating_accuracy_parameters} gives $C_{\rm pert}\sqrt{s\eta}\leq C_{\rm pert}\sqrt{\alpha s\eta} \leq\varepsilon_{\rm g}/4<10^{-3}$, so that $0.99\,Q_*(y)\preceq Q_{G'}(y)\preceq1.01\,Q_*(y)$. Combined with the assumption on $\widehat Q$, this gives
\begin{align*}
 0.89\,Q_*(y)\preceq\widehat Q\preceq1.12\,Q_*(y).
\end{align*}
Since $F=\alpha\Phi_p$, \eqref{eq:ml_metric_comparison} gives $\alpha Q_*(y)\preceq\nabla^2F(y)\preceq(2s+1)\alpha Q_*(y)$. Therefore $P\preceq1.12\nabla^2F(y)\preceq2\nabla^2F(y)$ and $\nabla^2F(y)\preceq(2s+1)P/0.89\preceq(4s+2)P$.
\end{proof}

We now combine the damped $y$-step analysis with the weight-update guarantee of Theorem~\ref{thm:weight_interface}. The next theorem shows that the intervening $G$-steps restore weight accuracy and keep every query valid throughout one centering call.

\begin{theorem}[Centering with $G$-steps]
\label{thm:single_loop_centering}
Consider one call of \textnormal{\textsc{Center}} with objective $c$ in a run, and assume that all queries of the run before this call are valid. Suppose that the call starts at $y^\circ$ with $\lambda_c(y^\circ)\leq1/10$ and with a weight that is accurate at $y^\circ$, and that its phase preconditioner $P$ satisfies~\eqref{eq:phase_metric_sandwich}. Suppose also that every query of the call returns estimates with
\begin{equation}\label{eq:single_loop_oracle_accuracy}
 \|\widehat g-g_G(y)\|_{\nabla^2F(y)^{-1}}\leq\varepsilon_{\rm g}/4,
 \qquad
 \|\widehat R-R_G(y)\|_F\leq\tau
\end{equation}
whenever this query and all preceding queries of the call are valid, where $(y,G)$ is its query state. Then every query of the call is valid, and every point-update segment lies in $\mathfrak B(y^\circ)$ and has $F$-local length less than $20s\gamma_y$. The call returns a feasible point $y$ with $\lambda_c(y)\leq\delta_{\rm c}$. At this point, both the final weight and the weight of the last query have gap less than $\eta/2$.
\end{theorem}
\begin{proof}
We follow the loops of the call and show that a $y$-step from an accurate weight leaves a warm weight, and that the following $G$-steps restore accuracy. Let $y_k$ be the point of the $y$-query of loop $k$, for $0\leq k<N_y$, and let $d$ be the objective gap of Lemma~\ref{lem:random_direct_centering}. We prove by induction over $k$ that the $y$-query of loop $k$ is valid, that the queries before it in the call are valid, and that $d(y_k)\leq d_0$. For $k=0$, the $y$-query is valid by the assumptions on $y^\circ$ and its weight, and $d(y^\circ)<d_0$ by Lemma~\ref{lem:random_direct_centering}(a).

\emph{The $y$-step.} Suppose that the claim holds for loop $k$, and let $G$ be the weight of its $y$-query. Since $G$ is accurate, \eqref{eq:single_loop_weight_bias}, \eqref{eq:single_loop_oracle_accuracy} and~\eqref{eq:alternating_accuracy_parameters} give
\begin{align*}
 \|\widehat g-\nabla F(y_k)\|_{\nabla^2F(y_k)^{-1}}
 \leq\frac{\varepsilon_{\rm g}}4+C_{\rm pert}\sqrt{\alpha s\eta}
 \leq\frac{\varepsilon_{\rm g}}2.
\end{align*}
Moreover, $\gamma_y=c_y/s\leq1/(32s)$ because $c_y\leq1/200$. Lemma~\ref{lem:random_direct_centering}(b) therefore keeps the segment from $y_k$ to $y_{k+1}$ in $\mathfrak B(y^\circ)$, bounds its $F$-local length by $20s\gamma_y=20c_y$, and gives $d(y_{k+1})\leq d_0$. The weight $G$ does not change along this segment, and~\eqref{eq:alternating_accuracy_parameters} gives
\begin{align*}
 \sqrt{\Delta_{y_k}(G)}+20C_{\rm pert}c_y\sqrt{s/\alpha}
 \leq\tfrac12\sqrt{\eta_{\rm w}}<\sqrt{c_{\rm pert}/s},
\end{align*}
where the last inequality uses $c_{\rm w}\leq c_{\rm pert}/4$. Thus Theorem~\ref{thm:weight_interface}(c) applies to the segment and gives
\begin{equation}\label{eq:warm_after_step}
 \Delta_{y_{k+1}}(G)\leq\tfrac14\eta_{\rm w}.
\end{equation}

\emph{The $G$-steps.} Let $\Delta_\ell$ be the gap at $y_{k+1}$ after $\ell$ $G$-steps, so that $\Delta_0\leq\eta_{\rm w}/4$ by~\eqref{eq:warm_after_step}. Suppose that the $G$-queries with gaps $\Delta_0,\ldots,\Delta_{\ell-1}$ are valid. By~\eqref{eq:single_loop_oracle_accuracy}, their stationarity estimates have error at most $\tau$. Each of these gaps is at most $\eta_{\rm w}\leq1/(100s)$, so Theorem~\ref{thm:weight_interface}(d) applies to each update. Iterating~\eqref{eq:weight_correction_recurrence}, whose fixed point is $96\tau^2=3\eta/32$, gives with $a:=\gamma_G/16$
\begin{equation}\label{eq:weight_correction_iterate}
 \Delta_\ell\leq e^{-a\ell}\Delta_0+\frac{3\eta}{32}
 \leq e^{-a\ell}\frac{\eta_{\rm w}}4+\frac{3\eta}{32}\leq\eta_{\rm w}.
\end{equation}
Hence the next $G$-query satisfies (V1) and (V2). The choice of $N_G$ in~\eqref{eq:single_loop_parameters} gives $e^{-aN_G}\eta_{\rm w}\leq\eta/4$, and so~\eqref{eq:weight_correction_iterate} yields
\begin{equation}\label{eq:weight_correction_accuracy}
 \Delta_{N_G}\leq\frac\eta{16}+\frac{3\eta}{32}<\frac\eta2,\qquad
 \Delta_{N_G-1}\leq e^a\frac\eta{16}+\frac{3\eta}{32}<\frac\eta2,
\end{equation}
where the second bound uses $a\leq1/16000$. The last $G$-query has gap $\Delta_{N_G-1}$, so it satisfies (V3), and all $G$-queries of loop $k$ are valid. The $y$-query of loop $k+1$ has the point $y_{k+1}$ and the weight with gap $\Delta_{N_G}<\eta/2$, so it is valid as well. This completes the induction.

\emph{The returned point.} The returned point is $y_{N_y}$, and it lies in $\mathfrak B(y^\circ)$, which is feasible. Lemma~\ref{lem:random_direct_centering}(c) gives $\lambda_c(y_{N_y})\leq\delta_{\rm c}$. The $G$-steps of the last loop take place at $y_{N_y}$, and~\eqref{eq:weight_correction_accuracy} for this loop gives the two gap bounds.
\end{proof}

\subsection{Correctness of a run}
\label{sec:run_correctness}
We now combine the guarantees of the three data structures with Theorem~\ref{thm:single_loop_centering}. The induction runs over the calls of a run, and each step uses accurate answers at valid queries to show that the next queries are valid.

\begin{theorem}[Correctness of a run]
\label{thm:centering_correctness}
Let $0<\gamma<1/10$, let $1\leq H\leq\sqrt m$ and $N_{\max}\geq1$ be integers, and let $\mathsf S:=\textnormal{\textsc{InitializeOracleState}}(\mathcal P,H,N_{\max},\gamma)$. Consider a run with at most $N_{\max}$ queries that uses $\mathsf S$, and suppose the following.
\begin{itemize}
\item[(i)] The starting point is feasible. Each call starts with $\lambda_c(y)\leq1/10$ whenever all preceding calls returned feasible points with decrement at most $\delta_{\rm c}$ and the standalone query, if any, returned $\widehat g$ with $\|\widehat g-\nabla F(y)\|_{\nabla^2F(y)^{-1}}\leq\varepsilon_{\rm g}/2$.
\item[(ii)] For every objective $c$ of the run and every $\theta\in[1/2,2]$, the function $\phi_{\theta c}$ has a minimizer.
\end{itemize}
Then, with probability at least $1-\gamma$, every query of the run is valid, the standalone query, if any, returns an estimate as in~(i), and every call returns a feasible point with decrement at most $\delta_{\rm c}$, at which the final weight has gap less than $\eta/2$.
\end{theorem}
\begin{proof}
We work on the intersection $\mathcal E$ of the success events of the weight initialization and of $\mathcal D_{\rm est}$, and prove by induction over the calls that all queries are valid. Theorem~\ref{thm:weight_interface}(e) with budget $\gamma/3$ and Theorem~\ref{thm:estimator_interface}(a) with budget $2\gamma/3$ give $\Pr[\mathcal E]\geq1-\gamma$.

\emph{Whiteners.} On every run, the first query factors the slack, because no anchor is stored. If $m^4>n$, Theorem~\ref{thm:slack_interface}(a) gives exact whiteners at all queries, and this includes every case with $m^2>n$. Otherwise $m^4\leq n$, and Theorem~\ref{thm:slack_interface}(b), with assumption~(ii), shows that each query $j$ receives a whitener with $\|\mathsf K-I\|_{\rm op}\leq\varepsilon_{\rm w}$ whenever queries $1,\ldots,j$ are valid. Hence, on $\mathcal E$, \eqref{eq:estimator_accuracy} and~\eqref{eq:estimator_metric} hold at every query $j$ for which queries $1,\ldots,j$ are valid.

\emph{The first query.} The first query computes the weight from an exact whitener, and on $\mathcal E$ this weight is accurate. If the run starts with a standalone query, this query is therefore valid. Its estimate satisfies~\eqref{eq:estimator_accuracy}, and~\eqref{eq:single_loop_weight_bias} and~\eqref{eq:alternating_accuracy_parameters} give $\|\widehat g-\nabla F(y)\|_{\nabla^2F(y)^{-1}}\leq\varepsilon_{\rm g}/4 +C_{\rm pert}\sqrt{\alpha s\eta}\leq\varepsilon_{\rm g}/2$.

\emph{One call.} Consider a call, and assume that all earlier queries are valid, that all earlier calls returned feasible points with decrement at most $\delta_{\rm c}$ and final weight gap less than $\eta/2$, and that the standalone query, if any, satisfied~(i). Then assumption~(i) shows that the call starts at a point $y^\circ$ with $\lambda_c(y^\circ)\leq1/10$. Its starting weight is the initial weight or the final weight of the preceding call, and in both cases it is accurate at $y^\circ$. Hence the first query $j$ of the call is valid. On $\mathcal E$, its proxy Hessian estimate satisfies~\eqref{eq:estimator_metric} with $G'=G_j$, or with $G'=G_{j-1}$ and $y_{j-1}=y^\circ$. In the second case, query $j-1$ is the standalone query or the last query of the preceding call, and its weight is accurate at $y^\circ$ by the initialization or by Theorem~\ref{thm:single_loop_centering}. Thus Lemma~\ref{lem:phase_metric} shows that the phase preconditioner satisfies~\eqref{eq:phase_metric_sandwich}. On $\mathcal E$, every query of the call satisfies~\eqref{eq:single_loop_oracle_accuracy} whenever it and all preceding queries are valid. Theorem~\ref{thm:single_loop_centering} therefore shows that all queries of the call are valid and that the call returns a feasible point with decrement at most $\delta_{\rm c}$, at which the final weight and the weight of the last query have gap less than $\eta/2$. This completes the induction.
\end{proof}

\subsection{Proof of Part (a) of the main algorithmic result}
\label{sec:central_start_correctness}
We apply Theorem~\ref{thm:centering_correctness} to the multiplier schedule of Algorithm~\ref{alg:more_sdp}. Together with Corollary~\ref{cor:ml_iterations}, this gives the objective guarantee and the bounds on the numbers of steps and queries from a central start.

\begin{lemma}[Correctness from a central start]
\label{lem:central_start_correctness}
Under the assumptions of Theorem~\ref{thm_more_amortized_sdp}, Algorithm~\ref{alg:more_sdp} satisfies Part~(a) and the bounds on the numbers of $y$-steps and oracle queries in Part~(b) of that theorem.
\end{lemma}
\begin{proof}
We apply Theorem~\ref{thm:centering_correctness} to the run of Algorithm~\ref{alg:more_sdp}, which has no standalone query and calls \textsc{Center} at the multipliers $t_j=t_0(1+\gamma_t)^j$ for $0\leq j\leq N_{\rm path}$. The run makes $(N_{\rm path}+1)N_{\rm c}=N_{\rm tot}$ queries, which equals the cap. The supplied point is feasible with $\lambda_{t_0b}(y_0)\leq1/10$. If the call at $t_j$ returns a point with decrement at most $\delta_{\rm c}$, then Corollary~\ref{cor:ml_iterations} shows that this point has decrement less than $1/10$ at $t_{j+1}$. This proves assumption~(i). Each objective $\theta t_jb$ in assumption~(ii) has $\theta t_j>0$, so Lemma~\ref{lem:ml_central_minimizer} supplies its minimizer. Theorem~\ref{thm:centering_correctness} with $\gamma=\beta$ thus shows that, with probability at least $1-\beta$, the returned point $y$ is feasible and satisfies $\lambda_{t_{N_{\rm path}}b}(y)\leq\delta_{\rm c}$.

By its definition, $N_{\rm path}$ is the least nonnegative integer with $t_0(1+\gamma_t)^{N_{\rm path}}\geq2\nu/\epsilon$, and so
\begin{align*}
 \max\{t_0,2\nu/\epsilon\}\leq t_{N_{\rm path}}
 \leq(1+\gamma_t)\max\{t_0,2\nu/\epsilon\}.
\end{align*}
Since $y$ is feasible, weak duality gives $b^\top y\geq\operatorname{OPT}$, and~\eqref{eq:central_accuracy} gives $b^\top y-\operatorname{OPT}\leq2\nu/t_{N_{\rm path}}\leq\epsilon$. This proves Part~(a).

By Corollary~\ref{cor:ml_iterations}, there are $N_{\rm path}+1=O(\sqrt\nu\log(2+\nu/(t_0\epsilon)))$ phases. Each phase makes $N_y=O(s)$ $y$-steps and $N_{\rm c}=O(s^2\log s)$ queries by~\eqref{eq:more_centering_parameters}. Since $\sqrt\nu=1000s\sqrt m$, this proves both counts.
\end{proof}
\section{Gradient/Hessian estimation}
\label{sec:randomized_oracles}
This section constructs the Gram estimator used by both implementations $\mathcal D^{\rm f}_{\rm est}$ and $\mathcal D^{\rm m}_{\rm est}$ of the gradient/Hessian estimation data structure, and it proves Theorem~\ref{thm:estimator_interface} for $\mathcal D^{\rm f}_{\rm est}$, which the algorithm uses when $m^2\leq n$. At a query state, one sketch of a common family of matrix factors gives both the gradient estimate $\widehat g$ and the stationarity estimate $\widehat R$. A compression of the same sketched factors gives the proxy Hessian estimate $\widehat Q$. Section~\ref{sec:random_implicit_operations} implements the block operations. Sections~\ref{sec:common_gram} and~\ref{sec_more_metric} analyze the Gram estimator and the proxy Hessian estimate at an exactly whitened state, and Section~\ref{sec:implicit_accuracy} bounds the effect of an inexact whitener. Section~\ref{sec:fresh_mode} assembles $\mathcal D^{\rm f}_{\rm est}$.

Throughout this section, $1\leq m\leq n$, the query state $(y,G)$ is fixed, the quantities of Definition~\ref{def:random_implicit_notation} are evaluated at this state, and we drop their arguments. Costs use the arithmetic model of Section~\ref{sec:arithmetic_model}. The Gram estimator requires only the mass bound $\tr[\mathcal L_y(G)^s]\leq3m$, and the proxy Hessian estimate also requires $I/2\preceq R\preceq2I$. By Theorem~\ref{thm:weight_interface}(a), both hold at warm weights.

\subsection{Block operations}
\label{sec:random_implicit_operations}
The implementation never forms the normalized constraint matrices. It applies them, through a whitener and a factor of the weight, to thin blocks of vectors at the cost $\mathcal C(q)$ of Definition~\ref{def_more_block_cost}.

\begin{lemma}[Implicit multiplication and coefficient normalization]
\label{lem:random_implicit_blocks}
Let $0<\rho<1/2$.
\begin{itemize}
\item[(a)] Given the $A_i$ and a feasible slack $S$, one can construct an implicit coefficient basis satisfying
\begin{equation}\label{eq:random_approx_coefficient_gram}
 \tfrac12I\preceq(\tr[B_iB_j])_{ij}\preceq2I
\end{equation}
with probability at least $1-\rho$, in $\widetilde O(n^\omega+\mathcal C(m))$ operations.
\item[(b)] Given an exact whitener, applying the coefficients, their adjoints, $\mathcal V$, $\mathcal V^\top$ or $\mathcal L_y(G)$ to $q$ columns, including the resulting Gram products and coefficient pairings, costs $O(n^\omega+\mathcal C(q))$. A degree-$d$ polynomial in $\mathcal L_y(G)$ uses $d$ such applications.
\end{itemize}
\end{lemma}
\begin{proof}
We first build a normalizing change of basis from a Gaussian sketch, and then implement each operation as a short sequence of rectangular products.

\emph{Part~(a).} Form $S$ and $J$ with $JSJ^\top=I$. Draw an $n\times k$ matrix $\Omega$ with independent $\mathcal N(0,1/k)$ entries, where $k=\lceil C(m+\log(2/\rho))\rceil$, and form
\begin{align*}
 Y_i:=JA_iJ^\top\Omega,\qquad
 \widehat H_0:=(\langle Y_i,Y_j\rangle)_{ij},\qquad
 H_0:=(\tr[(JA_iJ^\top)(JA_jJ^\top)])_{ij}.
\end{align*}
The coefficients are linearly independent, so $H_0\succ0$. Lemma~\ref{lem:gaussian_right_embedding} gives $\frac34H_0\preceq\widehat H_0\preceq\frac54H_0$ with probability at least $1-\rho$. Choose $T$ with $T^\top\widehat H_0T=I$. Then $\frac45I\preceq T^\top H_0T\preceq\frac43I$, which proves~\eqref{eq:random_approx_coefficient_gram} for $B_i=\sum_jT_{ji}JA_jJ^\top$. The orthogonal matrix $JS^{1/2}$ identifies these coefficients with the normalized constraint matrices of Definition~\ref{def:ml_potential} in the basis given by $T$. The procedure reports failure if $\widehat H_0$ is singular, and otherwise it stores $J$, $T$ and the input.

\emph{Part~(b).} Computing all $V_aX$ from an $n\times q$ block $X$ consists of three products. It applies the stacked $A_i$ to $J^\top X$, applies $J$ to the concatenated outputs, and mixes the coefficient index by the factor of the weight in the stored basis. These products cost
\begin{align*}
 \Tmat(mn,n,q),\qquad \Tmat(n,n,mq),\qquad \Tmat(m,m,nq),
\end{align*}
respectively. By Lemma~\ref{lem:matrix_transposition}, the transposed circuit has the same cost, and composing the two circuits applies $\mathcal L_y(G)=\mathcal V\mathcal V^\top$. Gram products are covered by $\Tmat(n,n,mq)$ using permutation symmetry. For gradient pairings, apply $J^\top$ to the output-factor blocks before forming their Gram matrix, and then pair it with the original coefficients in $O(mn^2)$ operations. Factorizations cost $O(n^\omega+m^\omega)$.

The same circuits evaluate the sketch in Part~(a). Splitting its $k$ columns into $O(k/m)$ blocks gives total cost $\widetilde O(n^\omega+\mathcal C(m))$. In particular, by permutation symmetry and $m\leq n$, its feature Gram product satisfies $\Tmat(m,nm,m)=\Tmat(nm,m,m)\leq\Tmat(mn,n,m)$.
\end{proof}

\subsection{One sketch for both estimates}
\label{sec:common_gram}
The gradient matrix $M$ and the stationarity matrix $R$ are Gram matrices of explicit matrix factors. We sketch these factors once and read both estimates off the sketched factors. Since $s$ is even, the factor
\begin{align*}
 \mathcal K:=\mathcal L_y(G)^{s/2-1}\mathcal V\in\R^{n\times mn}
 \qquad\text{satisfies}\qquad \mathcal K\mathcal K^\top=\mathcal L_y(G)^{s-1}
\end{align*}
and uses only integer powers of $\mathcal L_y(G)$. Define
\begin{equation}\label{eq:common_gram_factors}
 \begin{gathered}
 F_0:=2^{-1/2}\mathcal L_y(G)^{s/2}\in\R^{n\times n},\qquad
 F_a:=2^{-1/2}V_a\mathcal K\in\R^{n\times mn}\quad(1\leq a\leq m),\\
 \mathsf F:=[F_1^\top\ \cdots\ F_m^\top]^\top\in\R^{mn\times mn}.
 \end{gathered}
\end{equation}
These factors satisfy
\begin{equation}\label{eq:common_gram_identities}
 \sum_{a=0}^mF_aF_a^\top=M,\qquad
 R_{ab}=2\langle F_a,F_b\rangle\quad(a,b\geq1),\qquad
 \mathsf F^\top\mathsf F=\tfrac12\mathcal K^\top \mathcal L_y(G)\mathcal K.
\end{equation}
The first identity follows from $\mathcal K\mathcal K^\top=\mathcal L_y(G)^{s-1}$, the second from $R_{ab}=\tr[\mathcal L_y(G)^{s-1}V_aV_b]$, and the third from $\sum_aV_a^2=\mathcal L_y(G)$. The factors are the block rows of an even power of the symmetric dilation of $\mathcal V$. Indeed, after padding $F_0$ on the right and each $F_a$, $a\geq1$, on the left by zero columns,
\begin{equation}\label{eq:common_gram_dilation}
 \mathcal F_G:=\begin{pmatrix}F_0&0\\0&\mathsf F\end{pmatrix}
 =\frac1{\sqrt2}
   \begin{pmatrix}0&\mathcal V\\\mathcal V^\top&0\end{pmatrix}^{s}.
\end{equation}
This identity relates the movement of the factors to the trace-power energy of Section~\ref{sec:ml_exponential_curvature}. We apply the explicit blocks above and sketch $F_0$ separately from $F_1,\ldots,F_m$, because their input spaces differ. If $\tr[\mathcal L_y(G)^s]\leq3m$, the factors satisfy the mass bounds
\begin{equation}\label{eq:common_gram_trace}
 \sum_{a=0}^m\|F_a\|_F^2=\tr[\mathcal L_y(G)^s],\qquad
 \|\mathcal K\|_F^2=\tr[\mathcal L_y(G)^{s-1}]
 \leq n^{1/s}(3m)^{(s-1)/s}=O(m).
\end{equation}
The inequality is H\"older's inequality, and the last bound uses $n^{1/s}\leq e^{1/2}$, which holds because $s\geq2\log(2mn)$.

The two families need low-rank projections that capture their large singular directions. One range finder supplies both.

\begin{lemma}[Projections from one range finder]
\label{lem:common_gram_deflation}
Suppose $\tr[\mathcal L_y(G)^s]\leq3m$, $r\geq1$ and $0<\rho<1/2$. In $\widetilde O(n^\omega+\mathcal C(r))$ operations one can compute orthogonal projectors $\Pi_0=U_0U_0^\top$ on $\R^n$ and $\Pi_1=U_1U_1^\top$ on $\R^{mn}$, each of rank $\widetilde O(r)$, such that, with probability at least $1-\rho$,
\begin{equation}\label{eq:common_gram_tail}
 \|F_0N_0\|_{\rm op}^2\leq Cm/r,\qquad
 \|\mathsf FN_1\|_{\rm op}^2\leq Cm/r,\qquad
 \|\mathcal KN_1\|_{\rm op}^2\leq Cm/r,\qquad N_b:=I-\Pi_b.
\end{equation}
For two states, the same two projectors give these bounds at both states with the same probability and cost up to constants. If $2r+2\geq n$, use the exact factors
\begin{equation}\label{eq:common_gram_exact_factors}
 F_0^{\rm ex}:=2^{-1/2}\mathcal L_y(G)^{s/2},\qquad
 F_a^{\rm ex}:=2^{-1/2}V_a\mathcal L_y(G)^{s/2-1}\operatorname{chol}(\mathcal L_y(G)).
\end{equation}
Their Gram matrices give $M$ and $R$ exactly, at the same cost bound.
\end{lemma}
\begin{proof}
We find a right projection for $\mathcal K$ by a randomized power method and transfer it to the other two factors through the range of $\mathcal K$. Let $2r+2<n$. Apply Lemma~\ref{lem:gaussian_power_method} to $\mathcal K^\top$ with target rank $r$, oversampling $r+2$ and $\lceil\log(2mn)\rceil$ power steps. By~\eqref{eq:common_gram_trace}, the squared singular values of $\mathcal K$ sum to $O(m)$, so $\sigma_{r+1}^2\leq Cm/(r+1)$. The power-method bound gives $\E[\|\mathcal K(I-U_1U_1^\top)\|_{\rm op}]\leq C\sqrt{m/r}$, and Markov's inequality gives this bound up to a factor $2$ with probability at least $1/2$. Taking the span of $O(\log(1/\rho))$ independent ranges boosts the success probability to $1-\rho$, because enlarging the range decreases the residual.

Let $U_0$ be an orthonormal basis for the range of $\mathcal KU_1$. Then $(I-\Pi_0)\mathcal K\Pi_1=0$, and hence
\begin{align*}
 \|(I-\Pi_0)\mathcal K\|_{\rm op}\leq\|\mathcal K(I-\Pi_1)\|_{\rm op}.
\end{align*}
The mass bound gives $\|\mathcal L_y(G)\|_{\rm op}^s\leq\tr[\mathcal L_y(G)^s]\leq3m$ and therefore $\mathcal L_y(G)\preceq3I$, since $s\geq\log(3m)$. Consequently $\mathcal L_y(G)^s\preceq3\mathcal K\mathcal K^\top$ and $\mathsf F^\top\mathsf F\preceq\frac32\mathcal K^\top \mathcal K$, so that
\begin{align*}
 \|F_0N_0\|_{\rm op}^2
 =\tfrac12\|N_0\mathcal L_y(G)^sN_0\|_{\rm op}
 \leq\tfrac32\|N_0\mathcal K\|_{\rm op}^2,\qquad
 \|\mathsf FN_1\|_{\rm op}^2\leq\tfrac32\|\mathcal KN_1\|_{\rm op}^2.
\end{align*}
These inequalities prove~\eqref{eq:common_gram_tail}. For two states, take $U_1$ to span the right ranges of both states, with budget $\rho/2$ for each, and take $U_0$ to span $\mathcal K^{(0)}U_1$ and $\mathcal K^{(1)}U_1$.

Each step of the range finder applies $\mathcal K$ and $\mathcal K^\top$ to thin blocks, using integer powers of $\mathcal L_y(G)$ and the operations of Lemma~\ref{lem:random_implicit_blocks}. Forming $\mathcal KU_1$ and the thin orthogonalizations have the same cost. In particular, $\Tmat(mn,r,r)\leq\Tmat(mn,n,r)$ for $r\leq n$, and logarithmic factors in the widths are handled by blocking. In the exact case, $\mathcal L_y(G)\succ0$ by Lemma~\ref{lem:ml_sos}, and the identities follow from $\operatorname{chol}(\mathcal L_y(G))\operatorname{chol}(\mathcal L_y(G))^\top=\mathcal L_y(G)$. Forming $\mathcal L_y(G)$ and the exact factors costs $\widetilde O(n^\omega+\mathcal C(n))=\widetilde O(n^\omega+\mathcal C(r))$.
\end{proof}

Following Hutch++~\cite{mmmw21}, each sketch computes the projected part of a factor exactly and applies Gaussian vectors only to its residual. For an orthogonal projector $\Pi=UU^\top$, let $N=I-\Pi$ and draw independent standard Gaussian vectors $z_1,\ldots,z_k$ in its input space. The Hutch++ sketching matrix
\begin{equation}\label{eq:random_head_tail}
 \Xi:=[U,Nz_1/\sqrt k,\ldots,Nz_k/\sqrt k],\qquad
 \mathcal E:=\Xi\Xi^\top=\Pi+\frac1k\sum_jNz_jz_j^\top N
\end{equation}
satisfies $\E[\mathcal E\mid\Pi]=I$. For symmetric kernels $H_i$, Lemma~\ref{lem:gaussian_quadratic_forms} gives
\begin{equation}\label{eq:random_gaussian_variance}
 \E[\sum_i\bigl(\tr[(\mathcal E-I)H_i]\bigr)^2]
 =\frac2k\sum_i\|NH_iN\|_F^2.
\end{equation}
A nonsymmetric kernel may be replaced by its symmetric part, whose Frobenius norm is no larger. \eqref{eq:random_gaussian_variance} also holds after any fixed linear map of the output vector. The errors of our estimators are quadratic polynomials in Gaussian vectors, and their higher moments follow from their variances.

\begin{lemma}[Moments of quadratic Gaussian errors]
\label{lem:gaussian_moments}
Let $\xi$ be a centered random vector in a finite-dimensional inner product space whose coordinates are quadratic polynomials in a standard Gaussian vector, and let $\|\cdot\|$ be a Euclidean norm on this space. Then, for every $p\geq2$ and an absolute constant $C$,
\begin{equation}\label{eq:conditional_gaussian_moment}
 \bigl(\E[\|\xi\|^p]\bigr)^{1/p}
 \leq Cp\bigl(\E[\|\xi\|^2]\bigr)^{1/2}.
\end{equation}
The same bound holds conditionally on a sigma-field $\mathfrak F^-$ that determines the norm and the coefficients of the polynomials and is independent of the Gaussian vector.
\end{lemma}
\begin{proof}
We prove the bound for each coordinate in an orthonormal basis of the norm and then combine the coordinates by Minkowski's inequality. For a scalar coordinate, write $\xi^{(i)}=Q+L$, where $Q=z^\top Az-\tr[A]$ with $A=A^\top$ and $L=b^\top z$ for the standard Gaussian vector $z$. Integrating the tail bound of Lemma~\ref{lem:gaussian_quadratic_forms} gives $\|Q\|_{L^p}\leq Cp\|Q\|_{L^2}$, while the Gaussian moment bound gives $\|L\|_{L^p}\leq C\sqrt p\|L\|_{L^2}$. Gaussian symmetry gives $\E[QL]=0$, so $\|\xi^{(i)}\|_{L^2}^2=\|Q\|_{L^2}^2+\|L\|_{L^2}^2$. The triangle and Cauchy--Schwarz inequalities therefore give, after increasing $C$, $\|\xi^{(i)}\|_{L^p}\leq Cp\|\xi^{(i)}\|_{L^2}$. Minkowski's inequality in $L^{p/2}$ then gives
\begin{align*}
 \bigl(\E[\|\xi\|^p]\bigr)^{2/p}
 =\Bigl\|\sum_i(\xi^{(i)})^2\Bigr\|_{L^{p/2}}
 \leq\sum_i\|\xi^{(i)}\|_{L^p}^2
 \leq C^2p^2\sum_i\E[(\xi^{(i)})^2]
 =C^2p^2\,\E[\|\xi\|^2].
\end{align*}
The same calculation holds conditionally on $\mathfrak F^-$.
\end{proof}

Apply~\eqref{eq:random_head_tail} to $\Pi_0$ and $\Pi_1$ with independent residual probes, and call the results $\Xi_0$ and $\Xi_1$. Our estimator is
\begin{equation}\label{eq:common_gram_estimator}
 \begin{gathered}
 Y_0:=F_0\Xi_0,\qquad Y_a:=F_a\Xi_1\quad(1\leq a\leq m),
 \qquad \widehat M:=\sum_{a=0}^mY_aY_a^\top,\\
 \widehat g_i:=-\alpha\tr[\widehat MB_i],\qquad
 \widehat R_{ab}:=2\langle Y_a,Y_b\rangle\quad(1\leq a,b\leq m).
 \end{gathered}
\end{equation}
The outputs $\widehat M$ and $\widehat R$ are positive semidefinite. By~\eqref{eq:common_gram_identities} and $\E[\mathcal E_b\mid\Pi_b]=I$, their conditional means are $M$ and $R$, and $\E[\widehat g]=g$. To evaluate the factors, form $\mathcal L_y(G)^{s/2}\Xi_0$ and $\mathcal K\Xi_1=\mathcal L_y(G)^{s/2-1}(\mathcal V\Xi_1)$, and then apply all $V_a$ to the latter block. After $\mathcal V\Xi_1$ is formed, every matrix operation acts on blocks with $n$ rows. We write $C_a:=V_a\mathcal K=\sqrt2F_a$ for the stationarity factors. All expectations below condition on the state and on both projections.

\begin{lemma}[Covariance of the Gram estimator]
\label{lem:common_gram_covariance}
Let $F_0$ be a matrix with $n$ rows, and let $F_1,\ldots,F_m$ have $n$ rows and a common number of columns. Let
\begin{gather*}
 M=\sum_{a=0}^mF_aF_a^\top,\qquad Q[h,h]=\sum_{a=0}^m\|B_hF_a\|_F^2,\quad
 \theta=\max\{\|F_0N_0\|_{\rm op}^2,\|\mathsf FN_1\|_{\rm op}^2\}.
\end{gather*}
For the independent sketching matrices $\Xi_0,\Xi_1$ in~\eqref{eq:common_gram_estimator},
\begin{equation}\label{eq:common_gram_covariance}
 \operatorname{Cov}(\widehat g\mid\Pi_0,\Pi_1)
 \preceq(2\alpha^2\theta/k)Q.
\end{equation}
For the factors of~\eqref{eq:common_gram_factors}, these are the matrices $M$ and $Q$ of Definition~\ref{def:random_implicit_notation}, and Lemma~\ref{lem:common_gram_deflation} gives $\theta\leq Cm/r$.
\end{lemma}
\begin{proof}
We apply the variance identity~\eqref{eq:random_gaussian_variance} to each family in a fixed direction $h$. For either family, let $\mathsf H$ be its stacked factor and $\Theta=I\otimes B_h$ on its output space. Its contribution to $\widehat g^\top h$ is $-\alpha\tr[\mathcal E\,\mathsf H^\top\Theta\mathsf H]$, so its kernel is $\mathsf H^\top\Theta\mathsf H$, and
\begin{align*}
 \|N\mathsf H^\top\Theta\mathsf HN\|_F^2
 \leq\|\mathsf HN\|_{\rm op}^2\|\Theta\mathsf HN\|_F^2
 \leq\theta\|\Theta\mathsf H\|_F^2.
\end{align*}
\eqref{eq:random_gaussian_variance} bounds the variance of each contribution by $2\alpha^2\theta/k$ times $\|\Theta\mathsf H\|_F^2$. The errors of the two families are independent, so their variances add, and the two squared norms $\|\Theta\mathsf H\|_F^2$ sum to $Q[h,h]$.
\end{proof}

To bound the stationarity error, we also need operator-norm and trace bounds on the residual stationarity factors. The next lemma supplies these bounds.

\begin{lemma}[Residual stationarity mass]
\label{lem:power_compression}
For an orthogonal projector $N$ on $\R^{mn}$, let $A_N:=\sum_aV_a\mathcal KN\mathcal K^\top V_a$. Then
\begin{equation}\label{eq:stationarity_compression}
 A_N\preceq\|\mathcal KN\|_{\rm op}^2\mathcal L_y(G),\qquad
 \tr[A_N]\leq\tr[\mathcal L_y(G)^s].
\end{equation}
Under $\tr[\mathcal L_y(G)^s]\leq3m$ and~\eqref{eq:common_gram_tail}, $\|A_{N_1}\|_{\rm op}\leq Cm/r$ and $\tr[A_{N_1}]=O(m)$.
\end{lemma}
\begin{proof}
The first bound follows from $\mathcal KN\mathcal K^\top\preceq\|\mathcal KN\|_{\rm op}^2I$ and $\sum_aV_a^2=\mathcal L_y(G)$. The second follows from $\mathcal KN\mathcal K^\top\preceq \mathcal K\mathcal K^\top=\mathcal L_y(G)^{s-1}$ and $\tr[\sum_aV_a\mathcal L_y(G)^{s-1}V_a]=\tr[\mathcal L_y(G)^{s-1}\mathcal L_y(G)]=\tr[\mathcal L_y(G)^s]$. The last claim also uses $\|\mathcal L_y(G)\|_{\rm op}\leq3$.
\end{proof}

We now combine the covariance bound of Lemma~\ref{lem:common_gram_covariance} with the residual mass bound of Lemma~\ref{lem:power_compression}. This gives the accuracy and cost of computing the gradient and stationarity estimates from the same sketch.

\begin{lemma}[Accuracy of a full estimate]
\label{sec_more_projected_gradient}
\label{lem_more_projected_gradient}
Suppose $\tr[\mathcal L_y(G)^s]\leq3m$, and let the requested projection rank $r$ and the residual count $k$ both be $\widetilde\Theta(q)$. Condition on the state and on projections that satisfy~\eqref{eq:common_gram_tail}.
\begin{itemize}
\item[(a)] The error $e=\widehat g-g$ is a centered quadratic Gaussian polynomial, with
\begin{equation}\label{eq_more_absolute_covariance}
 \operatorname{Cov}(e)\preceq\widetilde O(\alpha^2m/q^2)Q,\qquad
 \E[\|e\|_{G/\alpha}^2]\leq\widetilde O(m^2/q^2).
\end{equation}
Every matrix $X$ fixed before sketching also satisfies $\E[\|Xe\|^2]\leq\widetilde O(\alpha^2m/q^2)\tr[XQX^\top]$.
\item[(b)] The stationarity estimate has mean $R$ and satisfies $\E[\|\widehat R-R\|_F^2]\leq Cm^2/(rk)$ for all $r,k\geq1$. Forming it from the sketched factors costs $\widetilde O(\Tmat(m,m,nq))$.
\end{itemize}
The total cost is $\widetilde O(n^\omega+\mathcal C(q))$.
\end{lemma}
\begin{proof}
Both parts follow from the variance identity~\eqref{eq:random_gaussian_variance} with the residual bounds~\eqref{eq:common_gram_tail}. Part~(a) follows from Lemma~\ref{lem:common_gram_covariance} by taking traces after congruences. In particular,
\begin{align*}
 \tr[GQ]=\sum_a\tr[MV_a^2]=\tr[M\mathcal L_y(G)]\leq\|\mathcal L_y(G)\|_{\rm op}\tr[M]\leq9m,
\end{align*}
because $\tr[M]=\tr[\mathcal L_y(G)^s]\leq3m$ by~\eqref{eq:common_gram_trace} and $\|\mathcal L_y(G)\|_{\rm op}\leq3$. Taking $X=(G/\alpha)^{1/2}$ gives the norm bound, where the factor $\alpha$ is polylogarithmic and suppressed.

For Part~(b), the mean follows from $\E[\widehat R_{ab}]=\langle C_a,C_b\rangle=R_{ab}$. The variance identity \eqref{eq:random_gaussian_variance} and Lemma~\ref{lem:power_compression} give
\begin{align*}
 \E[\|\widehat R-R\|_F^2]
 \leq\frac2k\sum_{a,b}\|N_1C_a^\top C_bN_1\|_F^2
 =\frac2k\tr[A_{N_1}^2]
 \leq\frac2k\|A_{N_1}\|_{\rm op}\tr[A_{N_1}]
 \leq\frac{Cm^2}{rk}.
\end{align*}
Vectorizing the sketched stationarity factors forms $\widehat R$ by a product of shape $(m,nq,m)$, which has the same cost as $(m,m,nq)$ by permutation symmetry. The projections and the factor evaluations have the claimed cost by Lemmas~\ref{lem:common_gram_deflation} and~\ref{lem:random_implicit_blocks}. The full-factor case is exact.
\end{proof}

\subsection{Sketching the proxy Hessian}
\label{sec_more_metric}
The factors sketched by the Gram estimator also determine an estimate of the proxy Hessian. We keep these factors and compress their output columns by an independent Gaussian sketching matrix when an estimate is needed. Both the sketched factors and their compression are analyzed by the Gram concentration of Lemma~\ref{lem:gaussian_gram_concentration}. Its variance bounds are natural in the coordinates of the weight, while centering needs accuracy relative to the proxy Hessian $Q$. The next lemma compares these two norms.

\begin{lemma}[Spectral bounds on the proxy Hessian]
\label{lem:metric_comparison}
Suppose that $I/2\preceq R\preceq2I$ and $\tr[\mathcal L_y(G)^s]\leq3m$. There are absolute constants $c_0,C_0>0$ such that
\begin{equation}\label{eq_more_geometry}
 c_0G^{-1}\preceq Q\preceq C_0mG^{-1},\qquad
 \tr[GQ]\leq C_0m,\qquad \|\mathcal L_y(G)\|_{\mathrm{op}}\leq C_0.
\end{equation}
In particular, if $h_1,\ldots,h_m$ are $Q$-orthonormal, then the directions $T_i:=B_{h_i}$ satisfy $\sum_iT_i^2\preceq C_0I$.
\end{lemma}
\begin{proof}
We prove the lower bound on $Q$ by H\"older's inequality and the upper bounds by the mass bound. After the coefficient change by $L$, we may assume that $G=I$, so that $V_a=B_a$, $\mathcal L_y(G)=\sum_aB_a^2$ and $R=\mathcal L_y^*(\mathcal L_y(G)^{s-1})$. All assertions of the lemma are invariant under this change. Let $h$ be a unit vector and $T:=B_h$. The Cauchy--Schwarz inequality gives $T^2\preceq \mathcal L_y(G)$, and so $\tr[T^2]\leq\tr[\mathcal L_y(G)]$. Since $h^\top Rh=\tr[\mathcal L_y(G)^{s-1}T^2]$ and $R\succeq I/2$, H\"older's inequality for the diagonal entries of $T^2$ in an eigenbasis of $\mathcal L_y(G)$ gives
\begin{align*}
 \tfrac12\leq\tr[\mathcal L_y(G)^{s-1}T^2]
 \leq\tr[\mathcal L_y(G)^sT^2]^{(s-1)/s}\tr[T^2]^{1/s}.
\end{align*}
We solve for $\tr[\mathcal L_y(G)^sT^2]$ and use $M\succeq\frac12\mathcal L_y(G)^s$. This gives
\begin{align*}
 Q[h,h]=\tr[MT^2]\geq\tfrac12\tr[\mathcal L_y(G)^sT^2]
 \geq2^{-1-s/(s-1)}\tr[\mathcal L_y(G)]^{-1/(s-1)}\geq c_0,
\end{align*}
where the last step uses H\"older's inequality $\tr[\mathcal L_y(G)]\leq n^{1-1/s}\tr[\mathcal L_y(G)^s]^{1/s}\leq n^{1-1/s}(3m)^{1/s}$ and the choice of $s$ in Definition~\ref{def:ml_s_choice}. This proves $Q\succeq c_0I=c_0G^{-1}$, and we may take $c_0\leq1$.

For the upper bounds, $\|\mathcal L_y(G)\|_{\rm op}^s\leq\tr[\mathcal L_y(G)^s]\leq3m$ and $s\geq\log(3m)$ give $\|\mathcal L_y(G)\|_{\rm op}\leq3$. Also, $\tr[\sum_aV_a\mathcal L_y(G)^{s-1}V_a]=\tr[\mathcal L_y(G)^{s-1}\mathcal L_y(G)]=\tr[\mathcal L_y(G)^s]$, so $\tr[M]=\tr[\mathcal L_y(G)^s]\leq3m$. Since $M\succeq0$,
\begin{align*}
 \tr[GQ]=\sum_a\tr[MB_a^2]=\tr[M\mathcal L_y(G)]\leq\|\mathcal L_y(G)\|_{\rm op}\tr[M]\leq9m.
\end{align*}
As $Q\succeq0$, this also gives $Q\preceq\tr[Q]\,I\preceq9m\,G^{-1}$. Finally, $Q$-orthonormal vectors satisfy $\sum_ih_ih_i^\top=Q^{-1}\preceq c_0^{-1}G$, and the positive map $\mathcal L_y$ gives
\begin{align*}
 \sum_iT_i^2=\mathcal L_y(Q^{-1})\preceq c_0^{-1}\mathcal L_y(G)
 \preceq3c_0^{-1}I.
\end{align*}
Thus the lemma holds for every $C_0\geq9/c_0$.
\end{proof}

The compression uses the following bound on a row covariance, which holds even when $M$ is singular.

\begin{lemma}[Row covariance in a weighted matrix subspace]
\label{lem:weighted_row_covariance}
Let $M\succeq0$, and let symmetric matrices $T_1,\ldots,T_m$ satisfy $\tr[MT_iT_j]=\delta_{ij}$ and $\sum_iT_i^2\preceq\beta I$. Then
\begin{equation}\label{eq_more_second_marginal}
 \sum_iT_iMT_i\preceq\sqrt{\beta\tr[M]}\,I.
\end{equation}
\end{lemma}
\begin{proof}
We bound the quadratic form of the left side at a unit vector $x$ by the Cauchy--Schwarz inequality, after showing that a certain compression of $M$ is at most $I$. Write $M=YY^\top$, and let $\Lambda_x:=\sum_iT_ixx^\top T_i\succeq0$. The matrices $T_iY$ are orthonormal in the Frobenius inner product, so Bessel's inequality, applied to $xv^\top$, gives
\begin{align*}
 v^\top Y^\top\Lambda_xYv
 =\sum_i\langle T_iY,xv^\top\rangle^2\leq\|v\|^2.
\end{align*}
It follows that $Y^\top\Lambda_xY\preceq I$. The matrix $\Lambda_x^{1/2}M\Lambda_x^{1/2}$ has the same nonzero eigenvalues, so it is also at most $I$. Consequently,
\begin{align*}
 \tr[M\Lambda_x^2]
 &=\tr[(\Lambda_x^{1/2}M\Lambda_x^{1/2})\Lambda_x]\leq\tr[\Lambda_x]\leq\beta,\\
 \Bigl(x^\top\sum_iT_iMT_ix\Bigr)^2
 &=(\tr[M\Lambda_x])^2\leq\tr[M]\tr[M\Lambda_x^2]\leq\beta\tr[M],
\end{align*}
where $\tr[\Lambda_x]=x^\top\sum_iT_i^2x\leq\beta$, and the second line applies the Cauchy--Schwarz inequality to $M^{1/2}$ and $\Lambda_xM^{1/2}$. Taking the supremum over unit vectors $x$ proves the claim.
\end{proof}

By Lemma~\ref{lem:weighted_row_covariance}, the compression needs only width $\widetilde O(\sqrt m)$.

\begin{lemma}[A proxy Hessian estimate from the sketched factors]
\label{lem_more_metric}
Suppose $I/2\preceq R\preceq2I$ and $\tr[\mathcal L_y(G)^s]\leq3m$, and let $0<\rho<1/2$ and $\epsilon:=1/40$.
\begin{itemize}
\item[(a)] Consider the Gram estimator with rank and probe widths $\widetilde\Theta(q)$, where $q\geq\sqrt m$, with sufficiently large logarithmic factors in $mn/\rho$. Keep its sketched factors
\begin{equation}\label{eq:metric_two_stage_factor}
 \begin{gathered}
 Y:=[Y_0\ Y_1\ \cdots\ Y_m]
   =[F_0\Xi_0\ F_1\Xi_1\ \cdots\ F_m\Xi_1],\\
 Q_Y[h,h]:=\tr[YY^\top B_h^2].
 \end{gathered}
\end{equation}
Conditional on the state and on projections that satisfy~\eqref{eq:common_gram_tail}, with probability at least $1-\rho$,
\begin{equation}\label{eq:raw_metric_accuracy}
 (1-\epsilon)Q\preceq Q_Y\preceq(1+\epsilon)Q,
 \qquad \|Y\|_F^2\leq\tfrac32\tr[\mathcal L_y(G)^s].
\end{equation}
For the exact factors, these conclusions hold deterministically with $Q_Y=Q$ and $\|Y\|_F^2\leq3m$.
\item[(b)] Let $Y$ have $O(mq)$ columns, $q\geq\sqrt m$, satisfy~\eqref{eq:raw_metric_accuracy}, and satisfy $\|Y\|_F^2\leq8m$. Algorithm~\ref{alg:comparison_metric} uses a new Gaussian sketching matrix of width $\widetilde O(\sqrt m)$ and returns $\widehat Q$ and its Cholesky factor. Conditional on $Y$, with probability at least $1-\rho$,
\begin{equation}\label{eq:compressed_metric_accuracy}
 (1-\epsilon)Q_Y\preceq\widehat Q\preceq(1+\epsilon)Q_Y,
 \qquad \tfrac9{10}Q\preceq\widehat Q\preceq\tfrac{11}{10}Q.
\end{equation}
The cost is $\widetilde O(n^\omega+\mathcal C(q))$.
\end{itemize}
The two parts have separate failure budgets, and the budget of the range finder is additional. Part~(b) draws only the Gaussian sketching matrix.
\end{lemma}
\begin{proof}
Both parts apply the Gram concentration of Lemma~\ref{lem:gaussian_gram_concentration} in coordinates where the target matrix is the identity. Its two hypotheses are a bound on the mean Gram matrix of one copy and a bound on its row covariance.

\emph{Part~(a).} Fix $Q$-orthonormal directions $T_1,\ldots,T_m$. Lemma~\ref{lem:metric_comparison} gives $\sum_iT_i^2\preceq C_0I$, and $\tr[M]=\tr[\mathcal L_y(G)^s]\leq3m$. Stack the two independent residuals by writing
\begin{gather*}
 \mathsf H:=\operatorname{diag}(F_0N_0,\mathsf FN_1),\qquad
 \Theta_i:=I_{m+1}\otimes T_i,\\
 X(z):=[\Theta_1\mathsf Hz\ \cdots\ \Theta_m\mathsf Hz],
\end{gather*}
where $z=(z_0,z_1)$ is standard Gaussian in $\R^n\oplus\R^{mn}$. In these coordinates, $Q=I$ is the sum of the positive semidefinite contributions of the two projected subspaces and of $\E[X(z)^\top X(z)]$. The matrix $Q_Y$ has the same projected contributions and the empirical average of the Gram matrices $X(z_j)^\top X(z_j)$. Thus
\begin{align*}
 \E[X(z)^\top X(z)]\preceq I,\qquad
 \E[X(z)X(z)^\top]
 =\sum_i\Theta_i\mathsf H\mathsf H^\top\Theta_i
 \preceq\|\mathsf H\|_{\rm op}^2\sum_i\Theta_i^2
 \preceq C(m/r)I,
\end{align*}
where the last step uses~\eqref{eq:common_gram_tail}. Lemma~\ref{lem:gaussian_gram_concentration} gives the comparison with $Q$ with failure probability $\rho/2$ as soon as $k\geq C(1+m/r)\log^2(Cmn/\rho)$. The stated widths suffice because $q^2\geq m$.

For the mass bound, write $T_b$ for the residual Gram matrix of family $b\in\{0,1\}$. It is positive semidefinite, so $\|T_b\|_F,\|T_b\|_{\rm op}\leq\tr[T_b]$. Lemma~\ref{lem:gaussian_quadratic_forms} gives
\begin{align*}
 \Pr\left[\frac1k\sum_jz_{b,j}^\top T_bz_{b,j}
                  >\tfrac32\tr[T_b]\right]\leq2e^{-k/32}.
\end{align*}
When $T_b=0$, the residual is zero. For $k\geq32\log(8/\rho)$, both residual masses are at most $3/2$ times their means, with failure probability at most $\rho/2$. The projected contributions are deterministic and nonnegative, so $\|Y\|_F^2\leq\frac32\tr[M] =\frac32\tr[\mathcal L_y(G)^s]$ on this event. The claims for the exact factors follow from~\eqref{eq:common_gram_exact_factors}. This proves Part~(a).

\emph{Part~(b).} Hold $Y$ fixed. If it has $d_Y$ columns, draw an independent $d_Y\times k_2$ Gaussian sketching matrix $\Omega$ with entries $\mathcal N(0,1/k_2)$, and form
\begin{equation}\label{eq:metric_two_stage_output}
 \widetilde Y:=Y\Omega,\qquad
 \widehat Q_{ij}:=\langle B_i\widetilde Y,B_j\widetilde Y\rangle.
\end{equation}
For $Q_Y$-orthonormal directions $\widetilde T_1,\ldots,\widetilde T_m$, the first bound of~\eqref{eq:raw_metric_accuracy} and Lemma~\ref{lem:metric_comparison} give
\begin{align*}
 \sum_i\widetilde T_i^2=\mathcal L_y(Q_Y^{-1})
 \preceq(1-\epsilon)^{-1}\mathcal L_y(Q^{-1})\preceq2C_0I.
\end{align*}
Let $\widetilde X(z):=[\widetilde T_1Yz\ \cdots\ \widetilde T_mYz]$ for a standard Gaussian $z\in\R^{d_Y}$. Then $\E[\widetilde X(z)^\top\widetilde X(z)]=I$, and Lemma~\ref{lem:weighted_row_covariance}, applied to $YY^\top$, gives
\begin{align*}
 \E[\widetilde X(z)\widetilde X(z)^\top]
 =\sum_i\widetilde T_iYY^\top\widetilde T_i
 \preceq\sqrt{2C_0\tr[YY^\top]}\,I\preceq C\sqrt m\,I.
\end{align*}
Lemma~\ref{lem:gaussian_gram_concentration} therefore gives the comparison with $Q_Y$ with failure probability $\rho$ and $k_2=\widetilde O(\sqrt m)$. The comparison with $Q$ follows from $(1-1/40)^2>0.9$ and $(1+1/40)^2<1.1$.

The compression and the final feature Gram matrix cost
\begin{align*}
 \Tmat(n,mq,k_2)\leq\widetilde O(\Tmat(n,n,mq)),\qquad
 \Tmat(m,nk_2,m)=O(\Tmat(m,m,nk_2)),
\end{align*}
respectively, by padding, blocking and permutation symmetry. Lemma~\ref{lem:random_implicit_blocks} applies all $B_i$ to $\widetilde Y$ at width $\widetilde O(\sqrt m)$. The mass check reads $O(nmq)$ entries, and the final factorization costs $O(m^\omega)$. All costs are $\widetilde O(n^\omega+\mathcal C(q))$, and the algorithm forms neither $YY^\top$ nor $Q_Y$. This proves Part~(b).
\end{proof}

\begin{algorithm}[!ht]
\caption{Proxy Hessian estimate from stored sketched factors}
\label{alg:comparison_metric}
\begin{algorithmic}[1]
\Procedure{ProxyHessian}{$B_1,\ldots,B_m,Y,\rho$} \Comment{Lemma~\ref{lem_more_metric}(b)}
 \State Report failure if $\|Y\|_F^2>8m$. \Comment{$O(nmq)$}
 \State $k_2\gets \lceil C\sqrt m\log^2(Cmn/\rho)\rceil$.
 \State Draw an independent Gaussian sketching matrix $\Omega$ with $k_2$ columns and entries $\mathcal N(0,1/k_2)$.
 \State Compute $\widetilde Y:=Y\Omega$ and all $B_i\widetilde Y$. \Comment{$\widetilde O(n^\omega+\mathcal C(q))$}
 \State Form $\widehat Q$ by~\eqref{eq:metric_two_stage_output}, and factor it or report failure. \Comment{$O(\Tmat(m,m,nk_2)+m^\omega)$}
 \State \Return $\widehat Q$ and its Cholesky factor.
\EndProcedure
\end{algorithmic}
\end{algorithm}

\subsection{Estimates with an inexact whitener}
\label{sec:implicit_accuracy}
When $m^4\leq n$, $\mathcal D_{\rm slack}$ returns $\widehat J=\mathsf KJ$ for an exact whitener $J$ and a symmetric $\mathsf K$ close to the identity (Theorem~\ref{thm:slack_interface}(b)). The data structure $\mathcal D^{\rm f}_{\rm est}$ then sketches the factors of slightly perturbed coefficients. We bound the resulting bias of one estimate and verify that the variance bounds persist. Throughout this subsection, the weight $G$ is warm, so Theorem~\ref{thm:weight_interface}(a) and Lemma~\ref{lem:metric_comparison} give $\|\mathcal L_y(G)\|_{\rm op}\leq C_0$ and $\tr[\mathcal L_y(G)^s]\leq2m$ for an absolute constant $C_0\geq1$. All factors use the Cholesky factor of $G$ in the fixed coefficient basis, and all matrix powers have integer exponents and are evaluated exactly.

\begin{lemma}[Whitening defects perturb the common factors]
\label{lem:whitening_perturbation}
Suppose that $\widehat J=\mathsf KJ$, where $J$ is an exact whitener, $\mathsf K$ is symmetric, and $\|\mathsf K-I\|_{\rm op}\leq\varepsilon_{\rm w}\leq1/(10s)$. Define the implemented coefficients and factors
\begin{align*}
 \widehat B_i&:=\mathsf KB_i\mathsf K,\qquad
 \widehat V_a:=\mathsf KV_a\mathsf K,\qquad
 \widehat{\mathcal L}_y(G):=\sum_a\widehat V_a^2,\\
 \widehat{\mathcal V}&:=[\widehat V_1\ \cdots\ \widehat V_m],\qquad
 \widehat{\mathcal K}:=\widehat{\mathcal L}_y(G)^{(s-2)/2}\widehat{\mathcal V},\\
 \widehat F_0&:=2^{-1/2}\widehat{\mathcal L}_y(G)^{s/2},\qquad
 \widehat F_a:=2^{-1/2}\widehat V_a\widehat{\mathcal K},\qquad
 \widehat C_a:=\widehat V_a\widehat{\mathcal K}\quad(a\geq1),
\end{align*}
and let $\widetilde F_a:=\mathsf K\widehat F_a$ for $0\leq a\leq m$. Then the following statements hold.
\begin{itemize}
\item[(a)] For every ordered eigenvalue,
\begin{align*}
 (1-\varepsilon_{\rm w})^4\lambda_i(\mathcal L_y(G))
 \leq\lambda_i(\widehat{\mathcal L}_y(G))
 \leq(1+\varepsilon_{\rm w})^4\lambda_i(\mathcal L_y(G)).
\end{align*}
Consequently $\|\widehat{\mathcal L}_y(G)\|_{\rm op}\leq2C_0$, $\tr[\widehat{\mathcal L}_y(G)^s]\leq3m$, and $\|\widehat{\mathcal L}_y(G)-\mathcal L_y(G)\|_{\rm op}\leq6C_0\varepsilon_{\rm w}$.
\item[(b)] Let $C_{\rm wh}:=40s(2C_0)^s$, which is polynomial in $2mn$. The factors $F_a$ and $C_a=V_a\mathcal K$ of~\eqref{eq:common_gram_factors} at the exact whitener satisfy
\begin{equation}\label{eq:whitening_factor_perturbation}
 \Bigl(\sum_{a=0}^m\|\widetilde F_a-F_a\|_F^2\Bigr)^{1/2},
 \quad
 \Bigl(\sum_a\|\widehat C_a-C_a\|_F^2\Bigr)^{1/2}
 \leq C_{\rm wh}n\varepsilon_{\rm w}.
\end{equation}
\item[(c)] One sketch of the implemented factors $\widehat F_a$ gives both outputs of~\eqref{eq:common_gram_estimator}. In the exact frame $J$, its gradient is the Gram estimator with the factors $\widetilde F_a$, and its stationarity matrix is the Gram estimator with the factors $\widehat C_a$.
\end{itemize}
\end{lemma}
\begin{proof}
All three parts compare the implemented quantities with the exact ones through the bound on $\mathsf K-I$.

\emph{Part~(a).} Since $(1-\varepsilon_{\rm w})^2I\preceq\mathsf K^2\preceq(1+\varepsilon_{\rm w})^2I$,
\begin{align*}
 (1-\varepsilon_{\rm w})^2\mathsf K\mathcal L_y(G)\mathsf K
 \preceq\widehat{\mathcal L}_y(G)
 =\mathsf K\Bigl(\sum_aV_a\mathsf K^2V_a\Bigr)\mathsf K
 \preceq(1+\varepsilon_{\rm w})^2\mathsf K\mathcal L_y(G)\mathsf K.
\end{align*}
The eigenvalues of $\mathsf K\mathcal L_y(G)\mathsf K$ equal those of $\mathcal L_y(G)^{1/2}\mathsf K^2\mathcal L_y(G)^{1/2}$, which lies between $(1\pm\varepsilon_{\rm w})^2\mathcal L_y(G)$. Eigenvalue monotonicity proves the claimed bounds. They include the trace bound, because $(1+\varepsilon_{\rm w})^{4s}\leq e^{2/5}<3/2$ and $\tr[\mathcal L_y(G)^s]\leq2m$. For the difference, we expand
\begin{align*}
 \widehat{\mathcal L}_y(G)-\mathcal L_y(G)=\mathsf K\Bigl(\sum_aV_a(\mathsf K^2-I)V_a\Bigr)\mathsf K
 +(\mathsf K-I)\mathcal L_y(G)\mathsf K+\mathcal L_y(G)(\mathsf K-I).
\end{align*}
The first term has norm at most $(1+\varepsilon_{\rm w})^2(2\varepsilon_{\rm w}+\varepsilon_{\rm w}^2)C_0$, and the triangle inequality gives $\|\widehat{\mathcal L}_y(G)-\mathcal L_y(G)\|_{\rm op}\leq6C_0\varepsilon_{\rm w}$.

\emph{Part~(b).} The identity $\widehat V_a-V_a=(\mathsf K-I)V_a\mathsf K+V_a(\mathsf K-I)$ gives
\begin{align*}
 \|\widehat{\mathcal V}-\mathcal V\|_F^2
 \leq9C_0n\varepsilon_{\rm w}^2,\qquad
 \|\mathcal V\|_{\rm op}^2\leq C_0,\qquad
 \|\widehat{\mathcal V}\|_{\rm op}^2\leq2C_0.
\end{align*}
For every integer $1\leq j\leq s/2$, the telescoping identity $\widehat{\mathcal L}_y(G)^j-\mathcal L_y(G)^j=\sum_{i=0}^{j-1}\widehat{\mathcal L}_y(G)^i(\widehat{\mathcal L}_y(G)-\mathcal L_y(G))\mathcal L_y(G)^{j-1-i}$ and Part~(a) give
\begin{align*}
 \|\widehat{\mathcal L}_y(G)^j-\mathcal L_y(G)^j\|_F
 \leq6C_0s(2C_0)^{s/2}\sqrt n\,\varepsilon_{\rm w}.
\end{align*}
The difference is zero for $j=0$. For $a\geq1$, expand
\begin{align*}
 \widehat C_a-C_a
 &=(\widehat V_a-V_a)\widehat{\mathcal K}
   +V_a(\widehat{\mathcal L}_y(G)^{(s-2)/2}-\mathcal L_y(G)^{(s-2)/2})\widehat{\mathcal V}\\
 &\qquad+V_a\mathcal L_y(G)^{(s-2)/2}(\widehat{\mathcal V}-\mathcal V).
\end{align*}
Apply the preceding bounds to the three terms, using $\|\widehat{\mathcal K}\|_{\rm op}\leq(2C_0)^{(s-1)/2}$ and $\sum_a\|V_aX\|_F^2\leq C_0\|X\|_F^2$. Their combined Frobenius norm is at most $20s(2C_0)^s\sqrt n\,\varepsilon_{\rm w}$. The extra left factor $\mathsf K$ in $\widetilde F_a$ contributes at most $\varepsilon_{\rm w}(\tr\widehat{\mathcal L}_y(G)^s)^{1/2}\leq\sqrt{3m}\,\varepsilon_{\rm w}$. For $a=0$, expand $\mathsf K\widehat{\mathcal L}_y(G)^{s/2}-\mathcal L_y(G)^{s/2} =(\mathsf K-I)\widehat{\mathcal L}_y(G)^{s/2}+\widehat{\mathcal L}_y(G)^{s/2}-\mathcal L_y(G)^{s/2}$. Combining these estimates with $m\leq n$ proves~\eqref{eq:whitening_factor_perturbation}.

\emph{Part~(c).} Let $\mathcal E_b=\Xi_b\Xi_b^\top$, $b\in\{0,1\}$, where $\Xi_0$ and $\Xi_1$ are the two independent sketching matrices of~\eqref{eq:random_head_tail}, and let $b(0)=0$ and $b(a)=1$ for $a\geq1$. For every factor, cyclicity of the trace gives
\begin{align*}
 \tr[\widehat B_i\widehat F_a\mathcal E_{b(a)}\widehat F_a^\top]
 =\tr[B_i\widetilde F_a\mathcal E_{b(a)}\widetilde F_a^\top].
\end{align*}
The stationarity output pairs the implemented factors directly, since $2\langle\widehat F_a\Xi_1,\widehat F_b\Xi_1\rangle =\langle\widehat C_a\Xi_1,\widehat C_b\Xi_1\rangle$ for $a,b\geq1$. Thus the left congruence by $\mathsf K$ enters the gradient only. The full-factor branch has the same target matrices and no sketching error.
\end{proof}

For the following estimates, define
\begin{equation}\label{eq:whitening_tolerance}
 \varepsilon_{\rm wh}:=C_{\rm wh}n\varepsilon_{\rm w}.
\end{equation}

The factor perturbation bound of Lemma~\ref{lem:whitening_perturbation} controls the bias introduced by the inexact whitener. We next show that the variance and proxy Hessian guarantees also persist at the stated tolerance.

\begin{lemma}[Bias and variance under an inexact whitener]
\label{lem:whitening_consequences}
Under the hypotheses of Lemma~\ref{lem:whitening_perturbation}, assume $\varepsilon_{\rm wh}\leq c$ for a sufficiently small absolute constant $c\leq1$. Let $\tilde g$ and $\widetilde R$ be the conditional means of an implemented estimate at the state, and let $(g,R)$ be its exact pair.
\begin{itemize}
\item[(a)] For an absolute constant $C$,
\begin{equation}\label{eq:whitening_bias}
 \|\tilde g-g\|_{G/\alpha}
 \leq C\sqrt\alpha(\sqrt m\,\varepsilon_{\rm wh}+\varepsilon_{\rm wh}^2),
 \qquad
 \|\widetilde R-R\|_F\leq C\sqrt m\,\varepsilon_{\rm wh}.
\end{equation}
\item[(b)] Let $\widehat Y$ be the sketched factors of an implemented estimate. The proxy Hessian and mass bounds for sketched factors of Lemma~\ref{lem_more_metric}(a), followed by the independent compression of Lemma~\ref{lem_more_metric}(b), give a proxy Hessian estimate with $0.9\,Q\preceq\widehat Q\preceq1.1\,Q$, at the same asymptotic cost.
\item[(c)] Conditional on the two projections, an implemented estimate satisfies the variance bounds of Lemma~\ref{lem_more_projected_gradient} up to absolute constants.
\end{itemize}
\end{lemma}
\begin{proof}
We express the means and the sketched factors through the factors of Lemma~\ref{lem:whitening_perturbation}(c) and compare them with the exact factors.

\emph{Part~(a).} By Lemma~\ref{lem:whitening_perturbation}(c), the conditional means are $\tilde g_i=-\alpha\tr[B_i\sum_a\widetilde F_a\widetilde F_a^\top]$ and $\widetilde R_{ab}=\langle\widehat C_a,\widehat C_b\rangle$. Fix $h$ with $h^\top G^{-1}h=1$. Lemma~\ref{lem:ml_sos} gives $B_h^2\preceq \mathcal L_y(G)\preceq C_0I$, and $\sum_a\|B_hF_a\|_F^2=Q[h,h]\leq\tr[M\mathcal L_y(G)]\leq C_0\tr[\mathcal L_y(G)^s]\leq2C_0m$. Expanding $\widetilde F_a\widetilde F_a^\top-F_aF_a^\top$ and applying the Cauchy--Schwarz inequality over $a$, together with~\eqref{eq:whitening_factor_perturbation}, gives
\begin{align*}
 \Bigl|\tr\Bigl[B_h\sum_a(\widetilde F_a\widetilde F_a^\top
       -F_aF_a^\top)\Bigr]\Bigr|
 \leq2\sqrt{2C_0m}\,\varepsilon_{\rm wh}+\sqrt{C_0}\,\varepsilon_{\rm wh}^2.
\end{align*}
Multiplying by $\alpha$ and taking the supremum over $h$ proves the gradient bound, since the dual of $\|\cdot\|_{G/\alpha}$ is $(\alpha h^\top G^{-1}h)^{1/2}$. In the same way, $\sum_a\|C_a\|_F^2=\tr[\mathcal L_y(G)^s]\leq2m$ and $\varepsilon_{\rm wh}\leq1$ give $\|\widetilde R-R\|_F\leq2\sqrt{2m}\,\varepsilon_{\rm wh} +\varepsilon_{\rm wh}^2\leq C\sqrt m\,\varepsilon_{\rm wh}$.

\emph{Part~(b).} The compression forms the vectors $\widehat B_h\widehat Y\Omega=\mathsf KB_h(\mathsf K\widehat Y)\Omega$, where $\widehat Y$ collects the sketched implemented factors. Its inner factor $\mathsf K\widehat Y$ is the sketch of the factors $\widetilde F_a$, and the outer $\mathsf K$ changes each norm by a factor in $[1-\varepsilon_{\rm w},1+\varepsilon_{\rm w}]$. For $h^\top G^{-1}h=1$, the factor displacement changes $(\sum_a\|B_hF_a\|_F^2)^{1/2}=Q[h,h]^{1/2}$ by at most $\sqrt{C_0}\,\varepsilon_{\rm wh}$, while $Q[h,h]\geq c_0$ by Lemma~\ref{lem:metric_comparison}. Hence the target of the sketch is within the factors $(1\pm C\varepsilon_{\rm wh})^2(1\pm\varepsilon_{\rm w})^2$ of $Q$. For sufficiently small $c$, this target is at least $(c_0/2)G^{-1}$. Together with $\widehat{\mathcal L}_y(G)\preceq3I$ and $\tr[\widehat{\mathcal L}_y(G)^s]\leq3m$, this supplies the geometric bounds used in the proof of Lemma~\ref{lem_more_metric}. Its proxy Hessian and mass bounds for sketched factors and its conditional compression guarantee therefore apply after the absolute constants of the sketch are increased, and they imply the stated comparison. The matrix $Q_Y$ uses the same sketched factors as the gradient and stationarity estimates, and only the Gaussian sketching matrix $\Omega$ is drawn afterwards.

\emph{Part~(c).} Apply Lemma~\ref{lem:common_gram_deflation} to the implemented coefficients, whose mass bound is Lemma~\ref{lem:whitening_perturbation}(a). On its success event, the two residual projectors satisfy
\begin{align*}
 \|\widehat F_0N_0\|_{\rm op}^2\leq Cm/r,\qquad
 \|\widehat{\mathsf F}N_1\|_{\rm op}^2\leq Cm/r,\qquad
 \|\widehat{\mathcal K}N_1\|_{\rm op}^2\leq Cm/r.
\end{align*}
The same factor bounds hold after left multiplication by $\mathsf K$, since $\|\mathsf K\|_{\rm op}\leq1.1$. For the native stationarity factors,
\begin{align*}
 \Bigl\|\sum_a\widehat C_aN_1\widehat C_a^\top\Bigr\|_{\rm op}
 &\leq\|\widehat{\mathcal K}N_1\|_{\rm op}^2\|\widehat{\mathcal L}_y(G)\|_{\rm op}\leq Cm/r,\\
 \tr\Bigl[\sum_a\widehat C_aN_1\widehat C_a^\top\Bigr]
 &\leq\tr[\widehat{\mathcal L}_y(G)\widehat{\mathcal K}\widehat{\mathcal K}^\top]
 =\tr[\widehat{\mathcal L}_y(G)^s]\leq3m.
\end{align*}
The two gradient errors are independent, so their covariances add. The variance proof of Lemma~\ref{lem_more_projected_gradient} therefore applies to both outputs with the same bounds up to constants.
\end{proof}

\subsection{The data structure \texorpdfstring{$\mathcal D^{\rm f}_{\rm est}$}{with fresh estimates}}
\label{sec:fresh_mode}
When $m^2\leq n$, the algorithm uses $\mathcal D^{\rm f}_{\rm est}$, which draws a new full estimate at every query, independently of all earlier sketches. It stores a single member, the proxy Hessian estimate $\mathcal D^{\rm f}_{\rm est}.\widehat Q$ of the phase, which it recomputes at every query with $a=1$. We first fix the tolerances, the width and the whitening tolerance, and then prove both parts of Theorem~\ref{thm:estimator_interface} for $\mathcal D^{\rm f}_{\rm est}$.

The estimates are accurate in a norm defined by the weight. The next lemma converts this norm to the dual norm of the Hessian.

\begin{lemma}[The Hessian norm at warm states]
\label{lem:warm_norm_comparison}
Every weight $G$ that is warm at a feasible point $y$ satisfies $\nabla^2F(y)^{-1}\preceq(2/c_0)G/\alpha$, where $c_0$ is the constant of Lemma~\ref{lem:metric_comparison}.
\end{lemma}
\begin{proof}
We pass from the Hessian to the proxy Hessian at the optimal weight, then to the proxy Hessian at $G$, and finally to $G^{-1}$. Lemma~\ref{lem:ml_gradient_certificates} gives $\nabla^2F(y)\succeq\alpha Q_*(y)$. Since $\eta_{\rm w}\leq c_{\rm pert}/s$ and $C_{\rm pert}\sqrt{s\eta_{\rm w}}=3\sqrt{c_{\rm w}}\leq3/20$, \eqref{eq:weight_metric_bias} gives $Q_*(y)\succeq(1+3/20)^{-2}Q_G(y)\succeq Q_G(y)/2$. Theorem~\ref{thm:weight_interface}(a) supplies the hypotheses of Lemma~\ref{lem:metric_comparison}, which gives $Q_G(y)\succeq c_0G^{-1}$. Hence $\nabla^2F(y)\succeq(\alpha c_0/2)G^{-1}$.
\end{proof}

Fix an absolute constant $C_Z\geq2/c_0$, large enough for Corollary~\ref{cor:phase_reference_metric} in Section~\ref{sec:oracle_maintenance}. For a vector $v$, a symmetric matrix $X$ and $Z\succ0$, define the gradient tolerance and the joint norm
\begin{equation}\label{eq:joint_norm}
 \delta_g:=\frac{\varepsilon_{\rm g}}{4\sqrt{C_Z}},\qquad
 \|(v,X)\|_Z^2:=\delta_g^{-2}\|v\|_{Z/\alpha}^2+\tau^{-2}\|X\|_F^2,
\end{equation}
where $\|v\|_{Z/\alpha}:=(v^\top Zv/\alpha)^{1/2}$. If $G$ is warm and $\|(\widehat g-g_G,\widehat R-R_G)\|_G\leq1$, then Lemma~\ref{lem:warm_norm_comparison} gives $\|\widehat g-g_G\|_{\nabla^2F^{-1}}\leq\sqrt{C_Z}\,\delta_g =\varepsilon_{\rm g}/4$ and $\|\widehat R-R_G\|_F\leq\tau$, which is~\eqref{eq:estimator_accuracy}.

By Lemmas~\ref{lem_more_projected_gradient} and~\ref{lem:whitening_consequences}(c), there is a polylogarithmic $V_0\geq1$ with the following property. At a warm state and on the success event of the projections, the centered error $\xi$ of a full estimate of width $q$, exact or implemented, satisfies
\begin{equation}\label{eq:full_estimate_variance}
 \E[\|\xi\|_G^2]\leq V_0m^2/q^2.
\end{equation}
The number $V_0$ absorbs the tolerances, the factor $s$ and all absolute constants. Let $\rho$ be the failure budget of $\mathcal D^{\rm f}_{\rm est}$, and let
\begin{equation}\label{eq:fresh_width}
 p:=2+\Bigl\lceil\log\frac{100(N_{\max}+1)}\rho\Bigr\rceil,\qquad
 q_{\rm fresh}:=\bigl\lceil C_{\rm f}\,p\sqrt{V_0}\,m\bigr\rceil,\qquad
 \rho_{\rm met}:=\frac\rho{4N_{\max}},
\end{equation}
where $C_{\rm f}$ is a sufficiently large absolute constant. The requested projection rank and the probe count are $q_{\rm fresh}$, multiplied by logarithmic factors in $mn/\rho$ that are large enough for Lemma~\ref{lem_more_metric}(a) with budget $\rho_{\rm met}$. Thus $q_{\rm fresh}=\widetilde O(m)$.

Finally, we fix the whitening tolerance. Let $C$ be the constant of Lemma~\ref{lem:whitening_consequences}(a), and choose the largest $\varepsilon_{\rm w}\leq\min\{1/100,1/(10s)\}$ such that $\varepsilon_{\rm wh}\leq c$ and
\begin{equation}\label{eq:bias_budget}
 C\sqrt\alpha(\sqrt m\,\varepsilon_{\rm wh}
                 +\varepsilon_{\rm wh}^2)\leq\delta_g/4,
 \qquad C\sqrt m\,\varepsilon_{\rm wh}\leq\tau/4.
\end{equation}
This tolerance is inverse-polynomial in $2mn$, because $C_{\rm wh}$ is polynomial and $\delta_g^{-1}$ and $\tau^{-1}$ are polylogarithmic.

\begin{algorithm}[!ht]
\caption{The procedure $\mathcal D^{\rm f}_{\rm est}.\textsc{Query}$}
\label{alg:fresh_query}
\begin{algorithmic}[1]
\Procedure{FreshQuery}{$y,G,\widehat J,a$} \Comment{$\mathcal D^{\rm f}_{\rm est}.\textsc{Query}$}
 \State Apply the coefficients through $\widehat J$ and the Cholesky factor of $G$ (Lemma~\ref{lem:random_implicit_blocks}).
 \State Compute the projections of Lemma~\ref{lem:common_gram_deflation} with failure budget $\rho/(4N_{\max})$, or use the exact factors~\eqref{eq:common_gram_exact_factors} if $2q_{\rm fresh}+2\geq n$. \Comment{$\widetilde O(n^\omega+\mathcal C(q_{\rm fresh}))$}
 \State Draw new residual probes, and compute $Y$, $\widehat g$ and $\widehat R$ by~\eqref{eq:metric_two_stage_factor} and~\eqref{eq:common_gram_estimator}. \Comment{$\widetilde O(n^\omega+\mathcal C(q_{\rm fresh}))$}
 \If{$a=1$}
  \State $\mathcal D^{\rm f}_{\rm est}.\widehat Q\gets\Call{ProxyHessian}{\widehat B_1,\ldots,\widehat B_m,Y,\rho_{\rm met}}$. \Comment{$\widetilde O(n^\omega+\mathcal C(q_{\rm fresh}))$}
 \EndIf
 \State \Return $(\widehat g,\widehat R,\mathcal D^{\rm f}_{\rm est}.\widehat Q)$.
\EndProcedure
\end{algorithmic}
\end{algorithm}

In Algorithm~\ref{alg:fresh_query}, $\widehat B_i=\widehat JA_i\widehat J^\top$ are the implemented coefficients, which equal the exact coefficients when $\widehat J$ is an exact whitener.

\begin{proof}[Proof of Theorem~\ref{thm:estimator_interface} for $m^2\leq n$]
We bound the failure probability of each query conditionally on the preceding history, and then take a union bound over the at most $N_{\max}$ queries.

\emph{The estimates at one query.} Fix $j\leq N_{\max}$, and condition on the history before the Gaussian sketching matrices of query $j$ are drawn. Suppose that queries $1,\ldots,j$ are valid and that query $j$ received $\widehat J=\mathsf KJ$ with $\|\mathsf K-I\|_{\rm op}\leq\varepsilon_{\rm w}$. Otherwise Part~(a) makes no claim at query $j$. The weight $G_j$ is warm, so the mass bound of Lemma~\ref{lem:whitening_perturbation}(a) holds. With failure probability at most $\rho/(4N_{\max})$, the projections satisfy~\eqref{eq:common_gram_tail} for the implemented factors. On this event, let $\xi$ be the difference between the estimate and its conditional mean $(\tilde g,\widetilde R)$. The vector $\xi$ is a centered quadratic polynomial in the probes, so~\eqref{eq:full_estimate_variance} and~\eqref{eq:conditional_gaussian_moment} and Markov's inequality give
\begin{align*}
 \Pr\bigl[\|\xi\|_G>\tfrac12\bigr]
 \leq\Bigl(\frac{2Cp\sqrt{V_0}\,m}{q_{\rm fresh}}\Bigr)^p
 \leq e^{-p}\leq\frac\rho{4N_{\max}},
\end{align*}
where $C$ is the constant of~\eqref{eq:conditional_gaussian_moment} and $C_{\rm f}\geq2eC$. Lemma~\ref{lem:whitening_consequences}(a) and~\eqref{eq:bias_budget} bound the bias by $\|(\tilde g-g,\widetilde R-R)\|_G\leq(1/16+1/16)^{1/2}=\sqrt2/4$. Hence $\|(\widehat g-g,\widehat R-R)\|_G\leq1/2+\sqrt2/4<1$, and the remark after~\eqref{eq:joint_norm} gives~\eqref{eq:estimator_accuracy}.

\emph{The proxy Hessian.} If $a_j=1$, then query $j$ is a $y$-query or the standalone query, so $G_j$ is accurate by (V3). Lemma~\ref{lem_more_metric} for an exact whitener, or Lemma~\ref{lem:whitening_consequences}(b) for an inexact one, gives~\eqref{eq:estimator_metric} with $G'=G_j$. The matrix $Q_Y$ and its compression each fail with probability at most $\rho_{\rm met}$.

\emph{Union bound.} Summing over the at most $N_{\max}$ queries, the failure probabilities are at most $N_{\max}(2\rho/(4N_{\max})+2\rho_{\rm met})=\rho$. On the complementary event, every query $j$ for which queries $1,\ldots,j$ are valid and received admissible whiteners satisfies~\eqref{eq:estimator_accuracy} and~\eqref{eq:estimator_metric}. This proves Part~(a) for $\mathcal D^{\rm f}_{\rm est}$.

\emph{Cost.} With exact whitening actions, Lemmas~\ref{lem:common_gram_deflation} and~\ref{lem_more_projected_gradient} bound a query by $\widetilde O(n^\omega+\mathcal C(q_{\rm fresh})) =\widetilde O(n^\omega+\mathcal C(m))$. Here we split the columns into blocks of width $m$, which costs a logarithmic factor. A proxy Hessian compression has the same cost by Lemma~\ref{lem_more_metric}(b). These bounds hold on every run, because the widths and trial counts are fixed. Thus every prefix of $N$ queries costs $\widetilde O(N\mathcal B_*(H))$ when $m^4>n$. With polynomial whitening actions, Theorem~\ref{thm:implicit_amortized_cost} gives the same bound. This proves Part~(b) for $\mathcal D^{\rm f}_{\rm est}$.
\end{proof}

The following bound specializes the fresh charge to current multiplication exponents.

\begin{lemma}[Fresh full estimates in the thin regime]
\label{lem:fresh_query_cost}
If $m^2\leq n$, then under Fact~\ref{fact:rectangular_thresholds}, for every fixed $\zeta>0$,
\begin{equation}\label{eq:fresh_query_cost}
 \mathcal C(m)\leq O_\zeta(mn^2m^\zeta).
\end{equation}
Consequently, a fresh query with an exact slack factorization costs $\widetilde O_\zeta(n^{\omega+\zeta}+n^2m^{1+\zeta})$.
\end{lemma}
\begin{proof}
We bound the four terms of $\mathcal C(m)$ by blocking. Since $m\leq\sqrt n$, the threshold $\omega(1,1,1/2)=2$ gives $\Tmat(n,n,m)\leq O_\zeta(n^2m^\zeta)$ by blocking both large dimensions into blocks of side $m^2$. Split $\Tmat(mn,n,m)$ into $m$ products of this shape, and $\Tmat(n,n,m^2)$ into $m$ blocks of $m$ columns. The mixing term satisfies $\Tmat(m,m,nm)\leq O_\zeta(nm^{\omega+\zeta})\leq O_\zeta(n^2m^{1+\zeta})$, because $m^{\omega-1}\leq n$. The linear term $nm^2$ has the same bound. Logarithmic width factors are absorbed by column blocking.
\end{proof}
\section{The weight maintenance data structure}
\label{sec:cost_per_iter}
This section implements $\mathcal D_{\rm weight}$ and proves Theorem~\ref{thm:weight_interface}. The data structure stores a single member, the weight $\mathcal D_{\rm weight}.G$, which we write $G$ in this section. One damped inverse update serves both procedures. $\textsc{Init}$ applies it repeatedly from a fixed starting weight, and $\textsc{Update}$ applies it once per $G$-query. We analyze the update in Section~\ref{sec:cost_weight_iteration} and the initialization in Section~\ref{sec:weight_initialization}. Sections~\ref{sec:weight_sublevel} and~\ref{sec:random_sharp_weight_accuracy} relate the weight gap to the matrix factors, the gradient bias and the proxy Hessian, and Section~\ref{sec:weight_movement} bounds the movement of the weight.

Until Section~\ref{sec:random_sharp_weight_accuracy}, the point $y$ and a coefficient basis are fixed, and we write $f=f_y$, $\mathcal L=\mathcal L_y$, $\mathcal L^*=\mathcal L_y^*$ and $\Delta=\Delta_y$. We use Definition~\ref{def:random_implicit_notation} with any $L$ satisfying $G=LL^\top$, so that $R=L^\top\mathcal L^*(\mathcal L(G)^{s-1})L$. Stationarity of the optimal weight gives $\mathcal L^*(\mathcal L(G_*)^{s-1})=G_*^{-1}$. The proofs are invariant under fixed changes of the coefficient basis, and the implementation uses the Cholesky factor.

\subsection{One inverse update for all weight computations}
\label{sec:cost_weight_iteration}
For an estimate $\widehat R$ of $R$, define
\begin{equation}\label{eq:weight_inverse_update}
 \gamma_G:=\frac1{1000s},\qquad
 G^+:=L\bigl((1-\gamma_G)I+2\gamma_G(I+\widehat R)^{-1}\bigr)L^\top.
\end{equation}
The update damps the inverse of the average of $G^{-1}$ and $L^{-\top}\widehat RL^{-1}$, an estimate of $\mathcal L^*(\mathcal L(G)^{s-1})$. It fixes $G$ when $\widehat R=I$, and it bounds the relative movement of $G$ even for large stationarity residuals. The next lemma measures the progress of the update by the residual energy
\begin{equation}\label{eq:weight_residual_energy}
 \mathcal R(\widehat R):=
 \tr[(\widehat R-I)^2(\widehat R+I)^{-1}].
\end{equation}

\begin{lemma}[Descent and mass control for the inverse update]
\label{lem:weight_inverse_step}
Suppose $\widehat R\succeq0$ and $\|\widehat R-R\|_F\leq\vartheta\leq1/8$. Then the update~\eqref{eq:weight_inverse_update} satisfies
\begin{align}
 (1-\gamma_G)G&\preceq G^+\preceq(1+\gamma_G)G,
 \label{eq:weight_relative_step}\\
 f(G^+)-f(G)&\leq-\frac{\gamma_G}4\mathcal R(\widehat R)
                         +5\gamma_G\vartheta^2.
 \label{eq:weight_inverse_descent}
\end{align}
It also preserves the mass bound $\tr[\mathcal L(G)^s]\leq3m$.
\end{lemma}
\begin{proof}
We first derive an upper model of $f$ along multiplicative updates, and then evaluate it at the update. For symmetric $X$ with $\|X\|_{\rm op}\leq1/(100s)$, the model is
\begin{equation}\label{eq:weight_upper_model}
 f(Le^{X}L^\top)-f(G)\leq\tr[(R-I)X]+2s\tr[RX^2].
\end{equation}
To prove it, let $\mathcal V=[V_1\ \cdots\ V_m]$ and $\mathsf B(t)=\mathcal V\exp(t(X\otimes I_n)/2)$. Then $\mathcal L(Le^{tX}L^\top)=\mathsf B(t)\mathsf B(t)^\top$. For $u(t)=s^{-1}\tr[\mathcal L(Le^{tX}L^\top)^s]$, cyclicity and $(\mathcal V^\top\mathcal V)^s=\mathcal V^\top \mathcal L(G)^{s-1}\mathcal V$ give $u'(0)=\tr[RX]$ and $\tr[(\mathcal V^\top\mathcal V)^s(X^2\otimes I_n)]=\tr[RX^2]$. Lemma~\ref{lem:ml_exponential_curvature}, with left generator $0$ and right generator $(X\otimes I_n)/2$, therefore gives
\begin{align*}
 0\leq u''(0)\leq2s\tr[RX^2],\qquad
 |u'''(t)|\leq20s\|X\|_{\rm op}u''(t).
\end{align*}
Thus $u''(t)\leq e^{1/5}u''(0)$ on $[0,1]$, and the Taylor remainder of $u$ at $t=1$ is at most $2s\tr[RX^2]$. Subtracting the change $\tr[X]$ of the log-determinant proves~\eqref{eq:weight_upper_model}.

For the update, put
\begin{align*}
 A=(\widehat R-I)(\widehat R+I)^{-1},\qquad
 X=\log(I-\gamma_GA),\qquad d^2=\mathcal R(\widehat R).
\end{align*}
Then $G^+=Le^{X}L^\top$ and $\|A\|_{\rm op}\leq1$, which proves~\eqref{eq:weight_relative_step}. Comparing eigenvalues gives $\|A\|_F^2\leq d^2$ and $\tr[\widehat RA^2]\leq d^2$. The scalar bound $|\log(1-x)|\leq2|x|$ for $|x|\leq1/2$ consequently gives
\begin{align*}
 \|X\|_{\rm op}\leq2\gamma_G,\qquad
 \|X\|_F\leq2\gamma_Gd,\qquad
 \|X^2\|_F\leq4\gamma_G^2d,\\
 \tr[RX^2]\leq\tr[\widehat RX^2]+\vartheta\|X^2\|_F
             \leq4\gamma_G^2(d^2+\vartheta d).
\end{align*}
For $r\geq0$ and $a=(r-1)/(r+1)$,
\begin{align*}
 (r-1)\log(1-\gamma_Ga)
 \leq-\frac{\gamma_G}{1+\gamma_G}\frac{(r-1)^2}{r+1}.
\end{align*}
This follows from $\log(1-x)\leq-x$ for $r\geq1$ and from $\log(1+x)\geq x/(1+x)$ for $r\leq1$. Summing this bound over the eigenvalues of $\widehat R$ bounds $\tr[(\widehat R-I)X]$, and $|\tr[(R-\widehat R)X]|\leq2\gamma_G\vartheta d$. The upper model~\eqref{eq:weight_upper_model} therefore yields
\begin{align*}
 f(G^+)-f(G)
 &\leq-\frac{\gamma_G}{1+\gamma_G}d^2+2\gamma_G\vartheta d
          +8s\gamma_G^2(d^2+\vartheta d)\\
 &\leq-\tfrac23\gamma_Gd^2+\tfrac52\gamma_G\vartheta d
 \leq-\tfrac14\gamma_Gd^2+5\gamma_G\vartheta^2,
\end{align*}
where the last two steps use $s\gamma_G=1/1000$ and Young's inequality. This proves~\eqref{eq:weight_inverse_descent}.

For the mass, the same Taylor bound before subtracting $\tr[X]$ gives
\begin{align*}
 \tr[\mathcal L(G^+)^s]\leq\tr[R]+s\tr[RX]+2s^2\tr[RX^2].
\end{align*}
We use $X\preceq-\gamma_GA$, $X^2\preceq4\gamma_G^2I$ and
\begin{align*}
 \tr[R(I+\widehat R)^{-1}]
 \leq\tr[\widehat R(I+\widehat R)^{-1}]
       +\vartheta\|(I+\widehat R)^{-1}\|_F
 \leq m+\vartheta\sqrt m.
\end{align*}
Since $-A=2(I+\widehat R)^{-1}-I$ and $R\succeq0$, these bounds give
\begin{align*}
 \tr[RX]\leq-\gamma_G\tr[RA]
 \leq2\gamma_G(m+\vartheta\sqrt m)-\gamma_G\tr[R],\qquad
 \tr[RX^2]\leq4\gamma_G^2\tr[R].
\end{align*}
The resulting bound is
\begin{align*}
 (1-s\gamma_G+8s^2\gamma_G^2)\tr[R]+2s\gamma_G(m+\vartheta\sqrt m).
\end{align*}
It increases with $\tr[R]=\tr[\mathcal L(G)^s]$, and at $\tr[R]=3m$ it equals $3m+s\gamma_Gm(-1+24s\gamma_G+2\vartheta/\sqrt m)<3m$.
\end{proof}

\begin{algorithm}[!ht]
\caption{The damped inverse update and the procedure $\mathcal D_{\rm weight}.\textsc{Update}$}
\label{alg:weight_inverse_step}
\begin{algorithmic}[1]
\Procedure{UpdateWeight}{$G,\widehat R$}
 \State Compute $L=\operatorname{chol}(G)$ and $\gamma_G=1/(1000s)$. \Comment{$O(m^\omega)$}
 \State Check $I+\widehat R\succ0$ and $(1-\gamma_G)I+2\gamma_G(I+\widehat R)^{-1}\succ0$. \Comment{$O(m^\omega)$}
 \State Report failure if a check or a factorization fails.
 \State \Return $L\bigl((1-\gamma_G)I+2\gamma_G(I+\widehat R)^{-1}\bigr)L^\top$. \Comment{$O(m^\omega)$}
\EndProcedure
\Statex
\Procedure{Update}{$\widehat R$} \Comment{$\mathcal D_{\rm weight}.\textsc{Update}$}
 \State $\mathcal D_{\rm weight}.G\gets\Call{UpdateWeight}{\mathcal D_{\rm weight}.G,\widehat R}$. \Comment{$O(m^\omega)$}
\EndProcedure
\end{algorithmic}
\end{algorithm}

The descent in~\eqref{eq:weight_inverse_descent} is measured by the residual energy, and the next lemma converts it into a decrease of the weight gap.

\begin{lemma}[An objective-gap certificate]
\label{lem:weight_gap_certificate}
The weight gap satisfies
\begin{equation}\label{eq:weight_gap_certificate}
 0\leq\Delta(G)\leq\tr[R]-m-\log\det R.
\end{equation}
If $R\succeq aI$ with $0<a\leq1/2$, put $C_{\rm gap}:=2+6\log(1/a)$. Every $\widehat R\succeq0$ with $\|\widehat R-R\|_F\leq\vartheta$ satisfies
\begin{equation}\label{eq:weight_gap_residual}
 \Delta(G)\leq2C_{\rm gap}\mathcal R(\widehat R)+8C_{\rm gap}\vartheta^2.
\end{equation}
For $R\succeq I/2$, one may take $C_{\rm gap}=2$.
\end{lemma}
\begin{proof}
We bound $f(G_*)$ from below by convexity, compare the resulting scalar functions with the residual energy, and finally pass from $R$ to $\widehat R$. Convexity of the trace power gives, for every $G'\succ0$,
\begin{align*}
 f(G')\geq\tr[L^{-\top}RL^{-1}G']-\log\det G'-m-\frac{s-1}{s}\tr[\mathcal L(G)^s],
\end{align*}
because $L^{-\top}RL^{-1}=\mathcal L^*(\mathcal L(G)^{s-1})$ is the gradient of the trace-power term at $G$. The first three terms have minimum $\log\det R-\log\det G$, attained at $G'=LR^{-1}L^\top$. Subtracting this lower bound on $f(G_*)$ from $f(G)$ proves~\eqref{eq:weight_gap_certificate}, since $\tr[R]=\tr[\mathcal L(G)^s]$.

For $r\geq1/2$, the bound $r-1-\log r\leq2(r-1)^2/(r+1)$ follows by integration over $[r,1]$ when $r\leq1$, and from $\log r\geq1-1/r$ when $r\geq1$. For $a\leq r<1/2$, the left side is at most $\log(1/a)$ and $(r-1)^2/(r+1)$ is at least $1/6$. Summing over the eigenvalues of $R$ gives $\Delta(G)\leq C_{\rm gap}\mathcal R(R)$.

Finally, the function $\psi(r)=(r-1)/\sqrt{r+1}$ satisfies $|\psi'|\leq2$ on $[0,\infty)$. For eigenpairs $(r_i,u_i)$ of $R$ and $(\hat r_j,v_j)$ of $\widehat R$,
\begin{align*}
 \|\psi(R)-\psi(\widehat R)\|_F^2
 =\sum_{i,j}|u_i^\top v_j|^2|\psi(r_i)-\psi(\hat r_j)|^2
 \leq4\|R-\widehat R\|_F^2.
\end{align*}
Since $\mathcal R(R)=\|\psi(R)\|_F^2$, this implies $\mathcal R(R)\leq2\mathcal R(\widehat R)+8\vartheta^2$, which proves~\eqref{eq:weight_gap_residual}.
\end{proof}

\subsection{Computing a weight from scratch}
\label{sec:weight_initialization}
\label{sec:random_stationarity}
Initialization uses the same Gram estimator as $\mathcal D_{\rm est}$. Its mass hypothesis holds before the weight is accurate, and its stationarity output is positive semidefinite on every run.

\begin{lemma}[Stationarity accuracy from the Gram estimator]
\label{lem:random_fast_stationarity}
If $\tr[\mathcal L(G)^s]\leq3m$ and $0<\vartheta,\rho<1/2$, the Gram estimator of Section~\ref{sec:common_gram} returns $\widehat R\succeq0$ with $\Pr[\|\widehat R-R\|_F\leq\vartheta]\geq1-\rho$ in $\widetilde O(\vartheta^{-2}(n^\omega+\mathcal C(m)))$ operations. The suppressed factors are polynomial in $s$ and logarithmic in $1/\rho$ and $1/\vartheta$.
\end{lemma}
\begin{proof}
We combine the variance bound of Lemma~\ref{lem_more_projected_gradient}(b) with the moment bound~\eqref{eq:conditional_gaussian_moment}. Take $r=m$, $p=2+\lceil\log(2/\rho)\rceil$ and $k=\lceil C'p^2m/\vartheta^2\rceil$ for a sufficiently large absolute $C'$. Compute $\Pi_1$ with failure probability $\rho/2$, and then draw $k$ new residual probes. The Gram output is positive semidefinite and has conditional mean $R$. Conditional on a successful projection, Lemma~\ref{lem_more_projected_gradient}(b) and~\eqref{eq:conditional_gaussian_moment} give
\begin{align*}
 \E[\|\widehat R-R\|_F^2\mid\Pi_1]\leq C_0m/k,\qquad
 \bigl(\E[\|\widehat R-R\|_F^p\mid\Pi_1]\bigr)^{1/p}
 \leq C_1p\sqrt{m/k}\leq\vartheta/e.
\end{align*}
Markov's inequality bounds the probability of a sketching error greater than $\vartheta$ by $e^{-p}\leq\rho/2$. If $2m+2\geq n$, the exact factors of~\eqref{eq:common_gram_exact_factors} give no sketching error. Splitting the projection columns and probes into $\widetilde O(\vartheta^{-2})$ blocks of at most $m$ columns gives the cost by Lemma~\ref{lem:random_implicit_blocks}. The stationarity Gram product of each block has shape $(m,nm,m)$.
\end{proof}

We now combine the inverse-update descent and gap certificate with the stationarity estimator of Lemma~\ref{lem:random_fast_stationarity}. This gives an initialization procedure that computes an accurate weight while preserving the mass bound.

\begin{lemma}[Weight computation from scratch]
\label{lem:weight_initialization}
Assume that the coefficient basis satisfies~\eqref{eq:random_approx_coefficient_gram}, and let $0<\eta<1/100$ and $0<\rho<1/2$.
\begin{itemize}
\item[(a)] Algorithm~\ref{alg:random_weight_computation} returns a weight with $\Delta(G)\leq\eta$ with probability at least $1-\rho$, in $\widetilde O(\eta^{-1}(n^\omega+\mathcal C(m)))$ operations. For inverse-polylogarithmic $\eta$, this is $\widetilde O(n^\omega+\mathcal C(m))$.
\item[(b)] On the same event, every intermediate weight satisfies $\tr[\mathcal L(G)^s]\leq3m$.
\end{itemize}
\end{lemma}
\begin{proof}
Start at $G_0=I/(4m)$ and apply the inverse update with new positive semidefinite estimates. We first assume that every estimate has error at most $\vartheta\leq1/8$, derive spectral bounds along the iteration, and then obtain a geometric contraction of the gap.

\emph{Spectral bounds.} The normalization \eqref{eq:random_approx_coefficient_gram} gives $\tr[\mathcal L(G_0)]\leq1/2$, and hence initial mass at most $2^{-s}$. Lemma~\ref{lem:weight_inverse_step} preserves the mass bound. Since $\tr[\mathcal L(G)]\geq\tr[G]/2$ by~\eqref{eq:random_approx_coefficient_gram} and $\tr[\mathcal L(G)]\leq n^{1-1/s}(3m)^{1/s}$ by H\"older's inequality, we have $G\preceq b_0I$ for $b_0=2n^{1-1/s}(3m)^{1/s}$. Also $\mathcal L^*(\mathcal L(G)^{s-1})\preceq6mI$, because for unit $h$,
\begin{align*}
 h^\top \mathcal L^*(\mathcal L(G)^{s-1})h=\tr[\mathcal L(G)^{s-1}B_h^2]
 \leq2(3m)^{(s-1)/s}\leq6m.
\end{align*}
Thus $I+\widehat R\preceq6mL^\top L+(1+\vartheta)I$, since $R=L^\top\mathcal L^*(\mathcal L(G)^{s-1})L$. Inverting this bound and using $L(6mL^\top L+cI)^{-1}L^\top=G(6mG+cI)^{-1}$ yields
\begin{align*}
 G^+\succeq(1-\gamma_G)G+
       2\gamma_GG((1+\vartheta)I+6mG)^{-1}.
\end{align*}
The scalar function on the right is increasing, and at $a_0=1/(32m)$ it is at least $a_0(1+11\gamma_G/21)$, because $1+\vartheta+6ma_0\leq21/16$. Induction gives
\begin{equation}\label{eq:random_weight_spectral_region}
 a_0I\preceq G\preceq b_0I,\qquad \tr[\mathcal L(G)^s]\leq3m.
\end{equation}
For unit $h$, put $T=B_h$. Then $T^2\preceq \mathcal L(G)/a_0$ by Lemma~\ref{lem:ml_sos}, and $\tr[T^2]\geq1/2$ by~\eqref{eq:random_approx_coefficient_gram}. In an eigenbasis $(\lambda_j,v_j)$ of $\mathcal L(G)$, put $t_j=v_j^\top T^2v_j\leq\lambda_j/a_0$. H\"older's inequality gives
\begin{align*}
 h^\top \mathcal L^*(\mathcal L(G)^{s-1})h=\sum_j\lambda_j^{s-1}t_j
 \geq a_0^{s-1}\sum_jt_j^s
 \geq a_0^{s-1}n^{1-s}2^{-s}.
\end{align*}
Consequently $R\succeq a_RI$ for $a_R=a_0^sn^{1-s}2^{-s}$, with $\log(1/a_R)=O(s\log(2mn))$.

\emph{Contraction.} Write $\Delta_j=\Delta(G_j)$ and $C_{\rm gap}=2+6\log(1/a_R)$. Lemma~\ref{lem:weight_gap_certificate} bounds $\Delta_j$ by $2C_{\rm gap}\mathcal R(\widehat R)+8C_{\rm gap}\vartheta^2$, and Lemma~\ref{lem:weight_inverse_step} decreases $f$ by $\gamma_G\mathcal R(\widehat R)/4$ up to $5\gamma_G\vartheta^2$. Combining the two bounds gives
\begin{equation}\label{eq:random_weight_contraction}
 \Delta_{j+1}\leq(1-\gamma_G/(8C_{\rm gap}))\Delta_j+6\gamma_G\vartheta^2.
\end{equation}
The optimal weight also satisfies $G_*\preceq b_0I$, since $\tr[\mathcal L(G_*)^s]=m$. Thus $f(G_0)\leq s^{-1}2^{-s}+m\log(4m)-m$ and $f(G_*)\geq s^{-1}m-m\log b_0-m$, which give $\Delta_0\leq\Delta_{\max}:=m\log(4mb_0)$. Choose
\begin{equation}\label{eq:random_weight_initial_parameters}
 \vartheta^2:=\frac{\eta}{384C_{\rm gap}},\qquad
 N_0:=\left\lceil\frac{8C_{\rm gap}}{\gamma_G}
                  \log\frac{4\Delta_{\max}}{\eta}\right\rceil.
\end{equation}
The fixed point of~\eqref{eq:random_weight_contraction} is $48C_{\rm gap}\vartheta^2=\eta/8$, and the initial contribution after $N_0$ steps is at most $\eta/4$. Hence $\Delta_{N_0}<\eta$.

\emph{Probability and cost.} Each call of the Gram estimator receives failure budget $\rho/N_0$. Conditional on the success of the preceding calls, \eqref{eq:random_weight_spectral_region} supplies the mass hypothesis of Lemma~\ref{lem:random_fast_stationarity}, and a union bound proves both claims. Finally, $C_{\rm gap}$ and $N_0$ are polylogarithmic and $\vartheta^{-2}=\widetilde O(\eta^{-1})$. The sketching cost follows from Lemma~\ref{lem:random_fast_stationarity}, and the $O(m^\omega)$ cost of each inverse update is absorbed.
\end{proof}

The procedure $\mathcal D_{\rm weight}.\textsc{Init}$ first normalizes the coefficient basis by Lemma~\ref{lem:random_implicit_blocks}(a), in $\widetilde O(n^\omega+\mathcal C(m))$ operations, and keeps this basis fixed during the computation. If the normalized basis is $B'_i=\sum_jT_{ji}B_j$ and the computed weight is $G'$, the data structure stores $G:=TG'T^\top$. Indeed, $\mathcal L_y'(G')=\mathcal L_y(TG'T^\top)$ for the map $\mathcal L_y'$ of the normalized basis, and therefore
\begin{equation}\label{eq:weight_basis_change}
 f_y(G)=f'_y(G')-2\log|\det T|,
\end{equation}
where $f'_y$ is the weight objective in the normalized basis. Thus the conversion preserves the gap and the mass. All later queries use the original coefficient basis and the Cholesky factor of $G$.

\begin{algorithm}[!ht]
\caption{Weight initialization with the damped inverse update}
\label{alg:random_weight_computation}
\begin{algorithmic}[1]
\Procedure{RandomizedWeights}{$B_1,\ldots,B_m,\eta,\rho$}
 \State $G\gets I/(4m)$, and take $\vartheta,N_0$ from~\eqref{eq:random_weight_initial_parameters}.
 \For{$j=1,\ldots,N_0$}
  \State Estimate $R$ by Lemma~\ref{lem:random_fast_stationarity} with tolerance $\vartheta$ and failure budget $\rho/N_0$. \Comment{$\widetilde O(\vartheta^{-2}(n^\omega+\mathcal C(m)))$}
  \State $G\gets\Call{UpdateWeight}{G,\widehat R}$. \Comment{$O(m^\omega)$}
 \EndFor
 \State \Return $G$.
\EndProcedure
\Statex
\Procedure{Init}{$y,\widehat J,\rho$} \Comment{$\mathcal D_{\rm weight}.\textsc{Init}$}
 \State Normalize the coefficient basis at the exact whitener $\widehat J$ by Lemma~\ref{lem:random_implicit_blocks}(a) with failure budget $\rho/2$, and let $T$ be the change of basis. \Comment{$\widetilde O(n^\omega+\mathcal C(m))$}
 \State $G'\gets\Call{RandomizedWeights}{B'_1,\ldots,B'_m,\eta,\rho/2}$. \Comment{$\widetilde O(\eta^{-1}(n^\omega+\mathcal C(m)))$}
 \State $\mathcal D_{\rm weight}.G\gets TG'T^\top$. \Comment{$O(m^\omega)$}
\EndProcedure
\end{algorithmic}
\end{algorithm}

\subsection{Weight accuracy and the common factors}
\label{sec:weight_sublevel}
The weight gap controls the distance between the factors of Section~\ref{sec:common_gram} at the current weight and at the optimal weight. Let
\begin{equation}\label{eq:weight_dilation_factors}
 \mathcal Z_G:=\begin{pmatrix}0&\mathcal V\\\mathcal V^\top&0\end{pmatrix},
 \qquad \mathcal F_G=2^{-1/2}\mathcal Z_G^s
\end{equation}
be the dilation and the stacked factor of~\eqref{eq:common_gram_dilation}, formed with the Cholesky factor of $G$. Since $s$ is even, the consecutive $n$-row blocks of $\mathcal F_G$ are $[F_0,0]$ and $[0,F_a]$, $1\leq a\leq m$. Thus the sum of the squared distances between corresponding factors equals the squared distance between the stacked factors. The following bound holds without any assumption on the size of the gap.

\begin{lemma}[Weight gap and factor distance]
\label{lem:weight_factor_distance}
For every $G\succ0$, the Cholesky-based factors satisfy
\begin{equation}\label{eq:weight_sublevel_factors}
 \|\mathcal F_G-\mathcal F_{G_*}\|_F^2\leq s\Delta(G).
\end{equation}
\end{lemma}
\begin{proof}
We prove a scalar inequality for even powers, apply it to the two dilations, and identify the linear term through stationarity. For real symmetric matrices $A,B$ and even $s$,
\begin{equation}\label{eq:even_power_distance}
 \|A^s-B^s\|_F^2
 \leq\tr[A^{2s}-B^{2s}-2sB^{2s-1}(A-B)].
\end{equation}
Indeed, for scalars $a,b$ with $b\neq0$, the right side minus the left side is $2b^{2s}[(a/b)^s-1-s(a/b-1)]\geq0$ by convexity of the even power, and for $b=0$ it is zero. Apply this inequality to the eigenvalues of $A$ and $B$, multiply by the squared overlaps of their eigenvectors, and sum to obtain~\eqref{eq:even_power_distance}.

Normalize the coefficient basis by $L_*:=\operatorname{chol}(G_*)$. The normalized weight has Cholesky factor $L_*^{-1}L$, so this change leaves the matrices $V_a$ and $\mathcal F_G$ unchanged. We may therefore assume $G_*=I$. Stationarity gives $\tr[\mathcal L(G_*)^s]=m$ and $\mathcal L^*(\mathcal L(G_*)^{s-1})=I$, and hence
\begin{align*}
 \tr[\mathcal Z_*^{2s-1}(\mathcal Z_G-\mathcal Z_*)]
 =2\tr[L-I].
\end{align*}
Since $\tr[\mathcal Z_G^{2s}]=2\tr[\mathcal L(G)^s]$, applying~\eqref{eq:even_power_distance} with $A=\mathcal Z_G$ and $B=\mathcal Z_*$ gives
\begin{align*}
 \|\mathcal F_G-\mathcal F_*\|_F^2
 &\leq\tr[\mathcal L(G)^s]-m-2s\tr[L-I]\\
 &=s\bigl(\Delta(G)+\log\det G-2\tr[L-I]\bigr)
 \leq s\Delta(G).
\end{align*}
The equality uses $\Delta(G)=s^{-1}(\tr[\mathcal L(G)^s]-m)-\log\det G$, which holds because $f(I)=s^{-1}m-m$, and the last inequality follows from $\log\det G=2\sum_i\log L_{ii}\leq2\sum_i(L_{ii}-1)$.
\end{proof}

A sufficiently small weight gap also controls the spectra of the weight and stationarity matrices. The following bounds establish the warm-weight guarantees used by the data structures.

\begin{lemma}[Geometry of a weight sublevel set]
\label{lem:weight_sublevel}
Suppose that $\Delta(G)\leq1/(100s)$. In the coefficient basis normalized by $\operatorname{chol}(G_*)$, so that $G_*=I$, we have
\begin{equation}\label{eq:weight_sublevel_region}
 \tfrac34I\preceq G\preceq\tfrac54I,\qquad
 \tfrac12I\preceq R\preceq2I,\qquad \tr[\mathcal L(G)^s]\leq2m.
\end{equation}
\end{lemma}
\begin{proof}
We split the gap into two nonnegative terms. Stationarity at $I$ gives the decomposition
\begin{equation}\label{eq:weight_bregman_gap}
 \Delta(G)=s^{-1}\tr[\mathcal L(G)^s-\mathcal L(G_*)^s-s\mathcal L(G_*)^{s-1}(\mathcal L(G)-\mathcal L(G_*))]
              +\tr[G-I]-\log\det G.
\end{equation}
Both terms are nonnegative by convexity. Since $x-1-\log x>1/100$ outside $[3/4,5/4]$, the second term gives the bounds on $G$.

The matrices $\sqrt2F_a(G_*)$, $a\geq1$, have Gram matrix $R_*=I$. Lemma~\ref{lem:weight_factor_distance} bounds the sum of their squared perturbations by $2s\Delta(G)$. For every unit vector $h\in\R^m$, the triangle inequality therefore gives
\begin{align*}
 \bigl|\sqrt{h^\top Rh}-1\bigr|
 \leq\sqrt{2s\Delta(G)}\leq\sqrt2/10.
\end{align*}
This proves the bounds on $R$. The mass bound follows from $\tr[\mathcal L(G)^s]=\tr[R]$.
\end{proof}

\subsection{Weight accuracy, gradient bias and the proxy Hessian}
\label{sec:random_sharp_weight_accuracy}
We now let the point vary. Write $M_G,g_G,Q_G$ for the matrices in Definition~\ref{def:random_implicit_notation}, and use the subscript $*$ for their values at $G_*$. In particular, $g_*=\nabla F(y)$ and $Q_*=Q_y$. The factor-distance bound of Lemma~\ref{lem:weight_factor_distance} controls both oracle biases.

\begin{lemma}[A dimension-free weight perturbation bound]
\label{lem:random_sharp_weight_accuracy}
Suppose that $\Delta_y(G)\leq1/(100s)$, and let $e:=\sqrt{s\Delta_y(G)}$.
\begin{itemize}
\item[(a)] The proxy Hessians satisfy
\begin{equation}\label{eq:random_metric_weight_error}
 (1-3e)^2Q_*\preceq Q_G\preceq(1+3e)^2Q_*.
\end{equation}
\item[(b)] The fixed-weight gradient satisfies
\begin{equation}\label{eq:random_gradient_weight_error}
 \|g_G-\nabla F(y)\|_{(\alpha Q_*)^{-1}}
 \leq3\sqrt{\alpha s\Delta_y(G)}.
\end{equation}
The same bound holds in the norm of $\nabla^2F(y)^{-1}$.
\end{itemize}
\end{lemma}
\begin{proof}
Both parts express a quantity as a quadratic form in the stacked factor $\mathcal F_G$ and compare it with its value at $\mathcal F_*$. By Lemma~\ref{lem:weight_factor_distance}, $\|\mathcal F_G-\mathcal F_*\|_F\leq e$. For $h\in\R^m$, let $\Theta_h=I_{m+1}\otimes B_h$. Lemma~\ref{lem:ml_rounding} gives
\begin{equation}\label{eq:random_operator_metric}
 \|\Theta_h\|_{\rm op}=\|B_h\|_{\rm op}\leq3\sqrt{Q_*[h,h]}.
\end{equation}
Since $Q_G[h,h]=\|\Theta_h\mathcal F_G\|_F^2$, the reverse triangle inequality gives
\begin{align*}
 \bigl|\sqrt{Q_G[h,h]}-\sqrt{Q_*[h,h]}\bigr|
 \leq\|\Theta_h\|_{\rm op}\|\mathcal F_G-\mathcal F_*\|_F
 \leq3e\sqrt{Q_*[h,h]}.
\end{align*}
This proves Part~(a).

For Part~(b), expand the two quadratic forms $\tr[B_hM_G]=\langle\Theta_h\mathcal F_G,\mathcal F_G\rangle$ and $\tr[B_hM_*]=\langle\Theta_h\mathcal F_*,\mathcal F_*\rangle$ around $\mathcal F_*$. The Cauchy--Schwarz inequality, $\|\Theta_h\mathcal F_*\|_F^2=Q_*[h,h]$ and~\eqref{eq:random_operator_metric} give
\begin{align*}
 |\tr[B_h(M_G-M_*)]|\leq(2e+3e^2)\sqrt{Q_*[h,h]}.
\end{align*}
Multiplying by $\alpha$ and taking the supremum over $\alpha Q_*[h,h]\leq1$ yields
\begin{align*}
 \|g_G-g_*\|_{(\alpha Q_*)^{-1}}
 \leq\sqrt\alpha\,(2e+3e^2)\leq3\sqrt\alpha\,e,
\end{align*}
because $e\leq1/10$. Finally, $\nabla^2F(y)\succeq\alpha Q_*$ by Lemma~\ref{lem:ml_gradient_certificates}, which gives the bound in the norm of $\nabla^2F(y)^{-1}$.
\end{proof}

The bias bound also controls how fast the weight gap changes when the point moves and the weight stays fixed.

\begin{lemma}[Variation of the square root of the weight gap]
\label{lem:weight_gap_motion}
Fix $G$ in the coefficient coordinates, let $y(u)$, $u\in[0,1]$, be a piecewise $C^1$ feasible path, and let $d_F$ be its $F$-local length. If $\Delta_{y(u)}(G)\leq c_{\rm pert}/s$ for all $u$, then~\eqref{eq:weight_gap_motion} holds. The assumption on the path holds automatically if $\sqrt{\Delta_{y(0)}(G)}+C_{\rm pert}\sqrt{s/\alpha}\,d_F<\sqrt{c_{\rm pert}/s}$.
\end{lemma}
\begin{proof}
We differentiate the weight gap along the path and bound the derivative by the bias bound of Lemma~\ref{lem:random_sharp_weight_accuracy}(b). The coefficient basis stays fixed. At a point where the slack is normalized to $I$, let $h\in\R^m$ and $T:=B_h$. The normalized constraint matrices have derivative $-\frac12(TB_i+B_iT)$ in the direction $h$, and so
\begin{align*}
 \nabla_y\mathcal L_y(G)[h]&=-\tfrac12(T\mathcal L_y(G)+\mathcal L_y(G)T)-\sum_aV_aTV_a,\\
 \nabla_yf_y(G)[h]&=\tr[\mathcal L_y(G)^{s-1}\nabla_y\mathcal L_y(G)[h]]
                  =-2\tr[M_GT]=\frac2\alpha g_G^\top h,
\end{align*}
where $V_a:=\sum_iL_{ia}B_i$ for $G=LL^\top$. In particular, $g_G$ is the gradient of $\alpha f_y(G)/2$ with respect to $y$. Since $\min_Gf_y(G)=2\Phi_p(y)=2F(y)/\alpha$, it follows that
\begin{align*}
 \frac{\d}{\d u}\Delta_{y(u)}(G)
 =\frac2\alpha\bigl(g_G-\nabla F(y(u))\bigr)^\top y'(u).
\end{align*}
By duality and Lemma~\ref{lem:random_sharp_weight_accuracy}(b), with $C_{\rm pert}=3$,
\begin{align*}
 \Bigl|\frac{\d}{\d u}\Delta_{y(u)}(G)\Bigr|
 \leq2C_{\rm pert}\sqrt{s/\alpha}\,\sqrt{\Delta_{y(u)}(G)}\,
 \|y'(u)\|_{\nabla^2F(y(u))}.
\end{align*}
Dividing by $2\sqrt{\Delta_{y(u)}(G)}$ and integrating proves~\eqref{eq:weight_gap_motion}. Where the gap vanishes, we apply the same argument to $\sqrt{\Delta_{y(u)}(G)+\varepsilon}$ and let $\varepsilon\downarrow0$.

For the last assertion, suppose that $\Delta_{y(u)}(G)=c_{\rm pert}/s$ for some $u$, and let $u_0$ be the first such $u$. Applying~\eqref{eq:weight_gap_motion} to the path on $[0,u_0]$ gives $\sqrt{\Delta_{y(u_0)}(G)}<\sqrt{c_{\rm pert}/s}$, a contradiction. Therefore $\Delta_{y(u)}(G)<c_{\rm pert}/s$ for all $u$.
\end{proof}

\subsection{Movement of the weight}
\label{sec:weight_movement}
The data structure $\mathcal D^{\rm m}_{\rm est}$ (Section~\ref{sec:oracle_maintenance}) compares sketched factors at consecutive query states. The next lemma shows that the weights of two such states are close. Part~(a) bounds the velocity of the optimal weight along a path, and Part~(b) integrates it between two points of a phase ball.

\begin{lemma}[Movement of the weight]
\label{lem:weight_movement}
\begin{itemize}
\item[(a)] Along a smooth feasible path $y(u)$, with $G_*$ and $\nabla^2\Phi_p$ evaluated at $y(u)$,
\begin{equation}\label{eq_more_weight_movement}
 \|G_*^{-1/2}(\partial_uG_*)G_*^{-1/2}\|_F
 \leq5\sqrt{(\partial_uy)^\top\nabla^2\Phi_p\,\partial_uy}.
\end{equation}
\item[(b)] Let $y_0$ and $y_1$ be two points of a common phase ball $\mathfrak B(y^\circ)$, with weights $G_0$ and $G_1$ that are warm at $y_0$ and $y_1$, respectively. Write $\Delta_0:=\Delta_{y_0}(G_0)$ and $\Delta_1:=\Delta_{y_1}(G_1)$, and let $d_F$ be the $F$-local length of the segment from $y_0$ to $y_1$. Then
\begin{equation}\label{eq:weight_movement_bound}
 \|\log(G_0^{-1/2}G_1G_0^{-1/2})\|_F
 \leq\frac{5d_F}{\sqrt\alpha}+C(\sqrt{\Delta_0}+\sqrt{\Delta_1}),
\end{equation}
and the right side is at most a small absolute constant.
\end{itemize}
\end{lemma}
\begin{proof}
Part~(a) restates Lemma~\ref{lem:ml_directional_response}. Part~(b) integrates Part~(a) and then passes from the optimal weights to the actual weights at the two endpoints.

\emph{Part~(a).} By the chain rule, $\partial_uG_*=\nabla G_*(y(u))[\partial_uy]$. Hence Part~(a) is the bound of Lemma~\ref{lem:ml_directional_response} in arbitrary coordinates, with $h:=\partial_uy$ and $\sqrt6<5$.

\emph{Part~(b).} The distance $\|\log(G^{-1/2}G'G^{-1/2})\|_F$ between positive definite matrices satisfies the triangle inequality, and it is at most the length $\int\|G^{-1/2}(\partial_uG)G^{-1/2}\|_F\,\d u$ of every path from $G$ to $G'$. Let $y(u):=y_0+u(y_1-y_0)$ for $u\in[0,1]$. Since $F=\alpha\Phi_p$, Part~(a) bounds the length of the path $G_*(y(u))$ by $5d_F/\sqrt\alpha$. At each endpoint, the log-determinant term in~\eqref{eq:weight_bregman_gap} and the spectral bounds of Lemma~\ref{lem:weight_sublevel} give $\|\log(G_*(y_0)^{-1/2}G_0G_*(y_0)^{-1/2})\|_F\leq C\sqrt{\Delta_0}$, and similarly at $y_1$. The triangle inequality gives~\eqref{eq:weight_movement_bound}. The segment lies in $\mathfrak B(y^\circ)$, whose diameter is $4/5$ in the norm of $\nabla^2F(y^\circ)$, and on which the local norm is at most $5/3$ times this norm by~\eqref{eq:sc_hessian_comparison}. Thus $d_F\leq4/3$. The choice of $\alpha$ and the bound $\Delta_i\leq c_{\rm w}/s$ then make the right side small.
\end{proof}

\subsection{Proof of the weight maintenance guarantee}
\begin{proof}[Proof of Theorem~\ref{thm:weight_interface}]
We collect the lemmas of this section. Part~(a) is Lemma~\ref{lem:weight_sublevel}, whose bounds on $G$ in the basis with $G_*=I$ read $\frac34G_*\preceq G\preceq\frac54G_*$ in the original basis. Part~(b) is Lemma~\ref{lem:random_sharp_weight_accuracy} with $C_{\rm pert}=3$, since $\Delta_y(G)\leq c_{\rm pert}/s=1/(100s)$. Part~(c) is Lemma~\ref{lem:weight_gap_motion}.

For Part~(d), Part~(a) gives $R\succeq I/2$, and $\tau\leq1/8$, so $\widehat R\succeq(1/2-\tau)I\succ0$ and the checks of Algorithm~\ref{alg:weight_inverse_step} pass. Lemma~\ref{lem:weight_inverse_step} with $\vartheta=\tau$ gives~\eqref{eq:weight_relative_step} and the descent $f(G^+)-f(G)\leq-\gamma_G\mathcal R(\widehat R)/4+5\gamma_G\tau^2$. Lemma~\ref{lem:weight_gap_certificate} with $C_{\rm gap}=2$ gives $\mathcal R(\widehat R)\geq\Delta(G)/4-4\tau^2$. Substituting this bound proves~\eqref{eq:weight_correction_recurrence}. The update uses one Cholesky factorization, one inversion of an $m\times m$ matrix and two products, which cost $O(m^\omega)$.

For Part~(e), the normalization fails with probability at most $\rho/2$, and Lemma~\ref{lem:weight_initialization} with budget $\rho/2$ computes $G'$ with gap at most $\eta$ in the normalized basis. By~\eqref{eq:weight_basis_change}, $G=TG'T^\top$ has the same gap in the original basis. The cost follows from Lemmas~\ref{lem:random_implicit_blocks}(a) and~\ref{lem:weight_initialization}(a), since $\eta^{-1}=\Theta(s^4)$ is polylogarithmic. The widths and the iteration count are fixed, so the cost bound holds on every run.
\end{proof}
\section{Gradient/Hessian estimate maintenance}
\label{sec_more_amortized}
\label{sec:oracle_maintenance}
\label{sec:local_estimates}
This section implements $\mathcal D^{\rm m}_{\rm est}$ and proves Theorem~\ref{thm:estimator_interface} when $m^2>n$. In this regime a full estimate at every query is too expensive. The data structure therefore maintains one pair of estimates across queries. An \emph{epoch} consists of at most $H$ calls of \textsc{Center} and starts with a full estimate. Every other query adds a correction of fixed width $\widetilde O(H\sqrt m)$, which estimates the change of the pair between two consecutive query states with coupled sketches. The proxy Hessian estimate of a call is compressed from the sketched factors of the current query at an epoch start, and from those of the preceding query otherwise.

Section~\ref{sec_more_local_geometry} aligns the factors of consecutive states, and Section~\ref{sec_more_gradient_increments} bounds the variance of a coupled correction. The gradient must be accurate in a norm that follows the weight. Section~\ref{sec_more_reference_metric} constructs reference metrics for this norm, and Section~\ref{sec:predictable_contractions} proves a maximal inequality for errors transported between such norms. Sections~\ref{sec_more_gradient_tracking}--\ref{sec:maintained_proof} combine these tools, and Section~\ref{sec:epoch_choice} chooses the epoch length.

The data structure $\mathcal D^{\rm m}_{\rm est}$ stores the maintained estimates $\mathcal D^{\rm m}_{\rm est}.\widehat g$ and $\mathcal D^{\rm m}_{\rm est}.\widehat R$, the reference metric $\mathcal D^{\rm m}_{\rm est}.Z$, the proxy Hessian estimate $\mathcal D^{\rm m}_{\rm est}.\widehat Q$ of the phase, the point $\mathcal D^{\rm m}_{\rm est}.y$, the weight $\mathcal D^{\rm m}_{\rm est}.G$, the whitener $\mathcal D^{\rm m}_{\rm est}.J$ and the sketched factors $\mathcal D^{\rm m}_{\rm est}.Y$ of the preceding query, and the number $\mathcal D^{\rm m}_{\rm est}.N_{\rm ep}$ of calls in the current epoch. Inside this section we write $\widehat g,\widehat R,Z,\widehat Q,Y,N_{\rm ep}$ for these members, and we index them by the query number $j$ when they change between queries. We keep the prefix for the preceding point, weight and whitener, which differ from the arguments of the current query. We write $w_j:=(g_{G_j}(y_j),R_{G_j}(y_j))$ for the target pair at query $j$ and $\widehat w_j:=(\widehat g_j,\widehat R_j)$ for the maintained pair.

Since $m^2>n$ implies $m^4>n$, every slack in this regime is factored exactly by Theorem~\ref{thm:slack_interface}(a). All probability statements condition on the query states and on common projections that satisfy Lemma~\ref{lem:common_gram_deflation}, before the new residual probes are drawn. They therefore apply conditionally on an adaptive history. The factors are those of~\eqref{eq:common_gram_factors}. The factor $F_0$ has $n$ columns, and $F_1,\ldots,F_m$ have $mn$ columns. The two families use independent residual probes, and within each family the two endpoints of a correction share their probes.

\subsection{Aligned factors at consecutive queries}
\label{sec_more_local_geometry}
The factors depend on the whitening coordinates. Before estimating a difference, we align the coordinates at the two states.

\begin{lemma}[Aligned factors and their implementation]
\label{lem_more_local_geometry}
\label{lem_more_aligned_whitener}
Let $y_0,y_1,G_0,G_1,\Delta_0,\Delta_1$ and $d_F$ be as in Lemma~\ref{lem:weight_movement}(b). Let $J_0S(y_0)J_0^\top=I$ and $\Sigma_{01}:=J_0S(y_1)J_0^\top$. Then $I/2\preceq\Sigma_{01}\preceq2I$, and $J_1:=\Sigma_{01}^{-1/2}J_0$ is an exact whitener of $S(y_1)$. Form the factors of~\eqref{eq:common_gram_factors} with these whiteners and the Cholesky factors of the weights. They satisfy
\begin{equation}\label{eq_more_factor_movement}
 \sum_{a=0}^m\|F_a^{(1)}-F_a^{(0)}\|_F^2
 \leq Cs(d_F^2/\alpha+\Delta_0+\Delta_1).
\end{equation}
The same bound holds with $\Sigma_{01}^{-1/2}F_a^{(1)}$ in place of $F_a^{(1)}$. In $\widetilde O(n^\omega)$ operations one can compute an exact whitener $\widehat J_1=OJ_1$, where $O$ is orthogonal and $\|O-I\|_{\rm op}\leq\epsilon_{\rm al}:=c_1/\sqrt m$. Using $\widehat J_1$ adds at most $Cm\epsilon_{\rm al}^2$ to either movement bound. The same conclusions hold for the stationarity factors $C_a=\sqrt2F_a$, $a\geq1$, up to an absolute factor.
\end{lemma}
\begin{proof}
We first bound the velocity of the stacked factor along the segment from $y_0$ to $y_1$ by the curvature bounds of Part~\ref{part:barrier}, and then integrate it. Normalize the slack at $y_0$, write $\widetilde S(u)=I+uT_0$ and $J(u)=(I+uT_0)^{-1/2}$, and use the optimal weight and its Cholesky factor $L$ along this segment. Put $v:=(\partial_uy)^\top\nabla^2\Phi_p\,\partial_uy$. Use the dilation $\mathcal Z$ and the stacked factor $\mathcal F=2^{-1/2}\mathcal Z^s$ of~\eqref{eq:weight_dilation_factors}.

\emph{The symmetric part of the velocity.} At each $u$, normalize the weight by its Cholesky factor, and put $T:=J(u)T_0J(u)$ and $V:=L^{-1}(\partial_uG_*)L^{-\top}$. The slack and weight derivatives agree with those of the exponential curves in Section~\ref{sec:ml_exponential_curvature}. For the slack this uses the commutation of $J(u)$ and $T_0$, and for the weight we first use the symmetric tangent $V/2$. If $\mathcal T$ is the resulting tangent of the dilation, Lemma~\ref{lem:ml_symmetric_curvature} gives
\begin{align*}
 \|\mathrm D(\mathcal Z^s)[\mathcal T]\|_F^2
 &\leq s\sum_{a=0}^{s-1}\|\mathcal Z^a\mathcal T\mathcal Z^{s-1-a}\|_F^2
 =\tfrac s2u_{T,V}''(0)=s(v-Q(T))\leq sv.
\end{align*}

\emph{The rotation of the Cholesky factor.} The derivative of the Cholesky factor satisfies $\partial_uL=L\mathsf X$ with $\mathsf X+\mathsf X^\top=V$. Thus $\mathsf A:=\mathsf X-V/2$ is skew-symmetric, and $\|\mathsf A\|_F^2\leq\|V\|_F^2/4\leq3v/2$ by Lemma~\ref{lem:ml_directional_response}. Set $\mathcal A:=\operatorname{diag}(0,\mathsf A\otimes I_n)$. Its contribution to $\partial_u(\mathcal Z^s)$ is the commutator $[\mathcal Z^s,\mathcal A]$, and stationarity gives
\begin{align*}
 \|[\mathcal Z^s,\mathcal A]\|_F^2
 \leq4\tr[(\mathcal V^\top\mathcal V)^s(\mathsf A\mathsf A^\top\otimes I_n)]
 =4\|\mathsf A\|_F^2\leq6v.
\end{align*}
Since $\mathcal F=2^{-1/2}\mathcal Z^s$, the squared triangle inequality yields
\begin{equation}\label{eq:factor_velocity}
 \sum_{a=0}^m\|\partial_uF_a\|_F^2
 =\|\partial_u\mathcal F\|_F^2\leq(s+6)v.
\end{equation}

\emph{Integration.} Lemma~\ref{lem:ml_rounding} gives $\|T\|_{\rm op}\leq3\sqrt v$. Integrating $|\partial_u\log(1+u\lambda)|\leq3\sqrt v$ for every eigenvalue $\lambda$ of $T_0$ gives $|\log(1+u\lambda)|\leq3d_F/\sqrt\alpha\leq\log2$, and so $I/2\preceq I+uT_0\preceq2I$. In the initial slack frame the derivative of the factor is $J(\partial_uF_a-TF_a/2)$. Its squared norm, summed over $a$, is also $O(sv)$, because $\|J\|_{\rm op}\leq\sqrt2$ and $\sum_a\|TF_a\|_F^2=\tr[MT^2]\leq v$. Integrating~\eqref{eq:factor_velocity} with $\int\sqrt v\,\d u=d_F/\sqrt\alpha$ gives squared movement at most $Csd_F^2/\alpha$ for either version of the factors. \eqref{eq:weight_sublevel_factors} contributes $Cs\Delta_i$ at endpoint $i$, which proves~\eqref{eq_more_factor_movement}.

\emph{The aligned whitener.} Use Lemma~\ref{lem:inverse_sqrt_polynomial} to choose a polynomial $\chi_2$ such that $E_{01}:=\chi_2(\Sigma_{01})\Sigma_{01}^{1/2}$ satisfies $\|E_{01}-I\|_{\rm op}\leq c_2\epsilon_{\rm al}/\sqrt n$. Form $L_{\rm al}:=\operatorname{chol}(E_{01}^2)$ and $\widehat J_1:=L_{\rm al}^{-1}\chi_2(\Sigma_{01})J_0$. Only the polynomial matrix $E_{01}^2=\chi_2(\Sigma_{01})\Sigma_{01}\chi_2(\Sigma_{01})$ is formed. Then $\widehat J_1S(y_1)\widehat J_1^\top=I$, and $O:=L_{\rm al}^{-1}E_{01}$ is orthogonal. Integrating the derivative of the Cholesky factor along $I+t(E_{01}^2-I)$ gives $\|L_{\rm al}-I\|_F\leq C\|E_{01}^2-I\|_F\leq Cc_2\epsilon_{\rm al}$, and so
\begin{align*}
 \|O-I\|_{\rm op}\leq
 \|L_{\rm al}^{-1}\|_{\rm op}(\|E_{01}-I\|_{\rm op}+\|L_{\rm al}-I\|_F)
 \leq\epsilon_{\rm al}
\end{align*}
for small $c_2$. The polynomial degree is $O(\log(2mn))$, and the cost is $\widetilde O(n^\omega)$.

The new factors are $OF_0^{(1)}O^\top$ and $OF_a^{(1)}(I_m\otimes O^\top)$ for $a\geq1$. In the initial slack frame, the left map $O$ is replaced by $\Sigma_{01}^{-1/2}$. Their total squared norm is $O(m)$, so the rotation adds $Cm\epsilon_{\rm al}^2$ to the movement bounds. Exact whitening preserves the targets $(g,R)$.
\end{proof}

\subsection{Coupled increments}
\label{sec_more_gradient_increments}
\label{sec:shared_pair_sampler}
By Lemma~\ref{lem_more_local_geometry}, the factors of the Gram estimator move only slightly between consecutive queries. If both endpoints use the same projections and the same Gaussian vectors, the sketching errors of the two estimates nearly cancel in their difference. The same pair of sketches also gives an estimate at each endpoint, so the transport correction below needs no further sketch. Each stationarity estimate uses the factors in the exact slack whitening of its own state. When the gradients of two states are compared in the slack frame of the first state, as in the proof below, the additional left whitening map applies to the gradient estimates only.

\begin{lemma}[Coupled increments]
\label{lem_more_gradient_increments}
\label{lem:paired_estimates}
\label{lem:shared_gradient_stationarity}
Let $y_0$ and $y_1$ be consecutive query points with weights $G_0$ and $G_1$ as in Lemma~\ref{lem:weight_movement}(b), and let $Z\succ0$ satisfy $Z\preceq CG_0$ and $Z\preceq CG_1$. For $q\geq1$, evaluate the estimates of~\eqref{eq:common_gram_estimator} at both states with common projections of requested rank $r$ and the same $k$ residual probes in each of the two input spaces, where $r,k=\widetilde\Theta(q)$ and the probes are drawn after $Z$ is fixed. Write $(g_u,R_u)$ for the targets at endpoint $u$. Condition on the two states and on projections that satisfy~\eqref{eq:common_gram_tail}.
\begin{itemize}
\item[(a)] The difference $\widehat d:=\widehat g_1-\widehat g_0$ has conditional mean $g_1-g_0$, and
\begin{equation}\label{eq_more_increment_variance}
 \E[\|\widehat d-(g_1-g_0)\|_{Z/\alpha}^2]
 \leq\widetilde O(m/q^2).
\end{equation}
\item[(b)] The centered error $\widehat d-(g_1-g_0)$ is a quadratic polynomial in the Gaussian probes with values in $\R^m$.
\item[(c)] Let $Q_0$ and $Q_1$ be the proxy Hessians at the two states. Each endpoint estimate satisfies
\begin{align*}
 \operatorname{Cov}(\widehat g_u-g_u)
 \preceq\widetilde O(\alpha^2m/q^2)Q_u,\qquad u\in\{0,1\}.
\end{align*}
\item[(d)] Let $C_{u,a}$, $u\in\{0,1\}$, be the stationarity factors of the two states, formed with the whiteners $J_0$ and $\widehat J_1$ of Lemma~\ref{lem_more_aligned_whitener}. Then
\begin{equation}\label{eq_more_stationarity_variance}
 \E[\|\widehat R_1-\widehat R_0-(R_1-R_0)\|_F^2]
 \leq\frac{Cm}{rk}\sum_a\|C_{1,a}-C_{0,a}\|_F^2
 =\widetilde O(m/(rk)).
\end{equation}
\end{itemize}
Both pairs of estimates together cost $\widetilde O(n^\omega+\mathcal C(q))$ operations.
\end{lemma}
\begin{proof}
We write the two gradient errors in one slack frame, so that the variance identity~\eqref{eq:random_gaussian_variance} applies to their difference with a kernel that is small by Lemma~\ref{lem_more_local_geometry}. With the exact factors, both endpoint estimates are exact and all centered errors vanish. Otherwise, use the factors and projections of~\eqref{eq:common_gram_estimator}.

\emph{One frame.} Let $J_{01}:=\widehat J_1J_0^{-1}$, so that $B_i^{(1)}=J_{01}B_i^{(0)}J_{01}^\top$ and $\|J_{01}\|_{\rm op},\|J_{01}^{-1}\|_{\rm op}\leq\sqrt2$. Represent the second gradient in the first slack frame by replacing $F_a^{(1)}$ with $J_{01}^\top F_a^{(1)}$. This leaves its sketched value unchanged, because $\tr[J_{01}^\top\widehat M_1J_{01}B_i^{(0)}]=\tr[\widehat M_1B_i^{(1)}]$. In this frame let $\mathsf F_0^{(u)}=F_0^{(u)}$, and let $\mathsf F_1^{(u)}$ stack $F_1^{(u)},\ldots,F_m^{(u)}$ vertically. Lemma~\ref{lem_more_aligned_whitener} and Lemma~\ref{lem:weight_movement}(b) give
\begin{align*}
 \sum_{b=0}^1\|\mathsf F_b^{(1)}-\mathsf F_b^{(0)}\|_F^2
 \leq Cs(d_F^2/\alpha+\Delta_0+\Delta_1)+Cm\epsilon_{\rm al}^2
 =\widetilde O(1).
\end{align*}
For the common residual projectors $N_0,N_1$, write $\mathsf H_b^{(u)}=\mathsf F_b^{(u)}N_b$. \eqref{eq:common_gram_tail} and the bound on $J_{01}$ imply $\|\mathsf H_b^{(u)}\|_{\rm op}^2\leq Cm/r$.

\emph{Gradient.} Let $T_i:=B_{Z^{1/2}e_i}^{(0)}$, and take $\Theta_{0,i}:=T_i$ and $\Theta_{1,i}:=I_m\otimes T_i$. Since $\sum_iT_i^2=\mathcal L_{y_0}(Z)\preceq C\mathcal L_{y_0}(G_0)\preceq CI$, we have $\sum_i\Theta_{b,i}^2\preceq CI$ for either family. The residual quadratic kernel of family $b$ for the $i$th coordinate of the difference is
\begin{align*}
 &\mathsf H_b^{(1)\top}\Theta_{b,i}\mathsf H_b^{(1)}
       -\mathsf H_b^{(0)\top}\Theta_{b,i}\mathsf H_b^{(0)}\\
 &\quad=(\mathsf H_b^{(1)}-\mathsf H_b^{(0)})^\top
                    \Theta_{b,i}\mathsf H_b^{(1)}
   +\mathsf H_b^{(0)\top}\Theta_{b,i}
                    (\mathsf H_b^{(1)}-\mathsf H_b^{(0)}).
\end{align*}
Its squared Frobenius norm, summed over $i$, is at most $C\max_u\|\mathsf H_b^{(u)}\|_{\rm op}^2 \|\mathsf F_b^{(1)}-\mathsf F_b^{(0)}\|_F^2$. The two probe families are independent, so their variances add. The variance identity~\eqref{eq:random_gaussian_variance} therefore gives
\begin{align*}
 \E\|\widehat d-(g_1-g_0)\|_{Z/\alpha}^2
 \leq\frac{C\alpha m}{rk}\sum_{b=0}^1
       \|\mathsf F_b^{(1)}-\mathsf F_b^{(0)}\|_F^2
 =\widetilde O(m/q^2).
\end{align*}
Each endpoint estimate is conditionally unbiased, and its centered error is quadratic in the two independent Gaussian blocks. This proves Parts~(a) and~(b). Lemma~\ref{lem:common_gram_covariance} proves Part~(c).

\emph{Stationarity.} Put $A_{N_1,u}:=\sum_aC_{u,a}N_1C_{u,a}^\top$. Lemma~\ref{lem:power_compression} gives $\|A_{N_1,u}\|_{\rm op}\leq Cm/r$. We expand
\begin{align*}
 C_{1,a}^\top C_{1,b}-C_{0,a}^\top C_{0,b}
 =(C_{1,a}-C_{0,a})^\top C_{1,b}
                  +C_{0,a}^\top(C_{1,b}-C_{0,b})
\end{align*}
and sum the squared Frobenius norms of the compressed kernels. The sum is at most $2(\|A_{N_1,0}\|_{\rm op}+\|A_{N_1,1}\|_{\rm op}) \sum_a\|C_{1,a}-C_{0,a}\|_F^2$. The variance identity~\eqref{eq:random_gaussian_variance} proves Part~(d), and the movement bound for the stationarity factors follows from $C_a=\sqrt2F_a$. The Cholesky factor fixes the coefficient index, and an orthogonal change of exact slack whitening leaves the stationarity target unchanged.

\emph{Cost.} Both endpoints share $\Xi_0$ and $\Xi_1$. The applications of integer powers and the $m$ coefficient operations cost $\widetilde O(n^\omega+\mathcal C(q))$ by Lemma~\ref{lem:random_implicit_blocks}(b). The Gram products, including the stationarity product of cost $\widetilde O(\Tmat(m,m,nq))$, and the $O(mn^2)$ coefficient pairings fit this bound.
\end{proof}

Every estimate in~\eqref{eq:common_gram_estimator} has the conditional mean $(g,R)$ of its state, for every width and every pair of projections. Hence a state has the same target whenever it is sketched again. With the exact factors of~\eqref{eq:common_gram_exact_factors}, both sketching errors vanish.

\subsection{A reference metric with a bounded refresh charge}
\label{sec_more_reference_metric}
The norm in which the gradient must be accurate changes from query to query. To keep this change from amplifying old errors, we construct reference metrics $Z_t$ and transport maps $\Gamma_t$ with $\|\Gamma_tv\|_{Z_t}\leq\|v\|_{Z_{t-1}}$ for every vector $v$. We clip the preceding reference from above relative to the current weight and from below relative to the proxy Hessian estimate. The transported estimate misses the component $(I-\Gamma_t)g_t$ of the new gradient. A trace $r_t$, the refresh charge, controls the variance of its new estimate, and the log-determinant divergence bounds the total refresh charge of an epoch. We first record the bounds supplied by the phase preconditioner.

\begin{lemma}[Bounds at the start of a phase]
\label{lem:lagged_metric}
Let query $j$ start a call at the point $y:=y_j$. Let $G_-:=G_j$ at an epoch start and $G_-:=G_{j-1}$ otherwise, and suppose in the latter case that $y_{j-1}=y$. Suppose that $G_-$ and $G_j$ are accurate at $y$ and that $0.9\,Q_{G_-}(y)\preceq\widehat Q\preceq1.1\,Q_{G_-}(y)$. Let $P:=\alpha\widehat Q$ and $\underline Z:=c_Z\widehat Q^{-1}$ for a sufficiently small absolute constant $c_Z$. Then
\begin{equation}\label{eq:lagged_metric_bounds}
 \tfrac12P\preceq\nabla^2F(y)\preceq(4s+2)P,\qquad
 cG_j/m\preceq \underline Z\preceq G_j/100,\qquad Q_{G_j}(y)\preceq Cs\,\widehat Q.
\end{equation}
\end{lemma}
\begin{proof}
The first pair of bounds is Lemma~\ref{lem:phase_metric}. Both weights are accurate at the same point, so Theorem~\ref{thm:weight_interface}(a) gives $\frac34G_*(y)\preceq G_-,G_j\preceq\frac54G_*(y)$, and hence $\frac35G_-\preceq G_j\preceq\frac53G_-$. The bounds $c_0G_-^{-1}\preceq Q_{G_-}(y)\preceq C_0mG_-^{-1}$ of Lemma~\ref{lem:metric_comparison} and the assumption on $\widehat Q$ therefore prove the bounds on $\underline Z$ for small $c_Z$. Finally, \eqref{eq:weight_metric_bias} and Lemma~\ref{lem:ml_gradient_certificates} give $Q_{G_j}(y)\preceq1.01\,Q_*(y)\preceq1.01\nabla^2F(y)/\alpha$, and the upper Hessian bound gives $Q_{G_j}(y)\preceq1.01(4s+2)\widehat Q$.
\end{proof}

Using the bounds of Lemma~\ref{lem:lagged_metric}, we now construct the reference metrics and transport maps by spectral clipping. The construction also controls the accumulated refresh charge.

\begin{lemma}[Reference metrics by spectral clipping]
\label{lem_more_reference_metric}
Let $G_t$ be a sequence of accurate weights, and let $\widehat Q_t$ be proxy Hessian estimates that satisfy the assumptions of Lemma~\ref{lem:lagged_metric}. Suppose that $\ell_t:=\|\log(G_{t-1}^{-1/2}G_tG_{t-1}^{-1/2})\|_F$ is bounded by an absolute constant. The index $t$ counts reference updates.
\begin{itemize}
\item[(a)] Let $c_Z>0$ be the constant of Lemma~\ref{lem:lagged_metric}, and let $\underline Z_t:=c_Z\widehat Q_t^{-1}$, so that
\begin{equation}\label{eq_more_lower_metric}
 cG_t/m\preceq \underline Z_t\preceq G_t/100.
\end{equation}
Starting an epoch with $Z_0:=G_0$, one can compute matrices $Z_t^-$ and $Z_t$ and maps $\Gamma_t$ such that
\begin{equation}\label{eq_more_reference_sandwich}
 \underline Z_t\preceq Z_t\preceq3G_t,\qquad Z_t^-\preceq2G_t,\qquad
 \Gamma_t^\top Z_t\Gamma_t\preceq Z_t^-\preceq Z_{t-1}.
\end{equation}
\item[(b)] The refresh charge $r_t:=\tr[Z_t(I-\Gamma_t)\underline Z_t^{-1}(I-\Gamma_t)^\top]$ satisfies $0\leq r_t<2m$. For every prefix of $N\leq H$ reference updates in an epoch,
\begin{equation}\label{eq_more_rank_budget}
 \sum_{t=1}^{N}r_t\leq C\Bigl[1+\Bigl(\sum_{t=1}^{N}\ell_t\Bigr)^2\Bigr]
 =\widetilde O(1+N^2).
\end{equation}
\item[(c)] Let $\zeta$ be a centered random vector with conditional covariance at most $\sigma^2\underline Z_t^{-1}$, where the metric and the maps are fixed before $\zeta$ is drawn. Then $\E[\|(I-\Gamma_t)\zeta\|_{Z_t}^2]\leq\sigma^2r_t$.
\item[(d)] The corrections and maps cost $\widetilde O(m^\omega)$ arithmetic operations per update, using scalar square roots and comparisons. Their spectral computations require only constant absolute accuracy, and the logarithm appears only in the analysis.
\end{itemize}
\end{lemma}
\begin{proof}
We construct the two clippings, show that the first one does not increase the log-determinant divergence from the current weight and the second one decreases it by half the refresh charge, and bound the increase of this divergence between updates. Lemma~\ref{lem:lagged_metric} gives~\eqref{eq_more_lower_metric}. Suppress the index $t$ during one update, and fix Cholesky factors $G=X_GX_G^\top$, $Z_{t-1}=X_ZX_Z^\top$ and $\underline Z=\underline X\,\underline X^\top$.

\emph{Clipping from above.} Apply Lemma~\ref{lem:approximate_diagonalization} to $X_Z^\top G^{-1}X_Z$, and write $U^\top X_Z^\top G^{-1}X_ZU=\operatorname{diag}(d)+E_1$, where $\|E_1\|_{\rm op}\leq1/4$ and $E_1$ has zero diagonal. Define
\begin{align*}
 t_i:=\min\{1,1/d_i\},\qquad
 T_-:=U\operatorname{diag}(t_i)U^\top,\qquad Z^-:=X_ZT_-X_Z^\top.
\end{align*}
Since $0\prec T_-\preceq I$, we have $Z^-\preceq Z_{t-1}$. Moreover, $T_-^{1/2}X_Z^\top G^{-1}X_ZT_-^{1/2}\preceq(5/4)I$, because the clipped diagonal is at most one and $\|E_1\|_{\rm op}\leq1/4$. This matrix has the same eigenvalues as $X_G^{-1}Z^-X_G^{-\top}$, which gives $Z^-\preceq(5/4)G$.

\emph{Clipping from below.} Write $\widetilde Z^-:=\underline X^{-1}Z^-\underline X^{-\top}$, and apply the same approximate diagonalization to it. With its new $U,d,E_1$, define
\begin{align*}
 \mathsf D:=U\operatorname{diag}((5/4-d_i)_+)U^\top,\qquad
 \widetilde Z:=\widetilde Z^-+\mathsf D,\qquad
 Z:=\underline X\widetilde Z\underline X^\top,\qquad \Gamma:=Z^{-1}Z^-.
\end{align*}
Then $\widetilde Z\succeq I$ and $0\preceq\mathsf D\preceq(5/4)I$. Hence
\begin{align*}
 \underline Z\preceq Z\preceq\tfrac54G+\tfrac54\underline Z\preceq2G,\qquad
 \Gamma^\top Z\Gamma=Z^-Z^{-1}Z^-\preceq Z^-,
\end{align*}
where the last inequality uses $Z^-\preceq Z$. This proves Part~(a).

\emph{The divergence.} For Part~(b), use the log-determinant divergence
\begin{equation}\label{eq_more_rank_potential}
 \Psi(Z,G):=\tr[G^{-1}Z]-\log\det(G^{-1}Z)-m,
\end{equation}
which is the Bregman divergence of $-\log\det$ and is therefore nonnegative. The upper clipping does not increase it. Indeed, writing $\phi(x)=x-\log x-1$ and using the zero diagonal of $E_1$,
\begin{align*}
 \Psi(Z^-,G)-\Psi(Z_{t-1},G)
 =\sum_i[d_i(t_i-1)-\log t_i]
 =-\sum_{d_i>1}\phi(d_i)\leq0.
\end{align*}
For the lower clipping, $I-\Gamma=\underline X^{-\top}\widetilde Z^{-1}\mathsf D\underline X^\top$, so cyclicity gives $r=\tr[\mathsf D\widetilde Z^{-1}\mathsf D]\leq25m/16<2m$. Let $\Pi_{\mathsf D}$ project onto the support of $\mathsf D$. In the diagonalizing coordinates, $\Pi_{\mathsf D}\widetilde Z\Pi_{\mathsf D} =(5/4)\Pi_{\mathsf D}+\Pi_{\mathsf D}E_1\Pi_{\mathsf D}\preceq(3/2)\Pi_{\mathsf D}$. The Schur-complement inequality gives $\Pi_{\mathsf D}\widetilde Z^{-1}\Pi_{\mathsf D}\succeq(2/3)\Pi_{\mathsf D}$. Consequently,
\begin{align*}
 \tr[\mathsf D]\leq\tfrac32\tr[\widetilde Z^{-1}\mathsf D],\qquad
 r\leq\tfrac54\tr[\widetilde Z^{-1}\mathsf D].
\end{align*}
Using $\log(1+x)\geq x/(1+x)$ and $\underline X^\top G^{-1}\underline X\preceq I/100$, we obtain
\begin{align*}
 \Psi(Z^-,G)-\Psi(Z,G)
 \geq\tr[\widetilde Z^{-1}\mathsf D]-\tfrac1{100}\tr[\mathsf D]
 \geq\tfrac{197}{200}\tr[\widetilde Z^{-1}\mathsf D]\geq\tfrac12r.
\end{align*}

\emph{Between updates.} Hold $Z$ fixed and move $G$ along its geodesic. The bound on consecutive weight distances and $Z\preceq2G_{t-1}$ give $G^{-1/2}ZG^{-1/2}\preceq CI$ along it. Differentiating $\Psi$ and using $(x-1)^2\leq C\phi(x)$ on this interval gives
\begin{equation}\label{eq_more_potential_growth}
 \mathrm d\Psi
 =\tr[(I-G^{-1/2}ZG^{-1/2})G^{-1/2}(\mathrm dG)G^{-1/2}]
 \leq C\sqrt\Psi\,\mathrm d\ell_G,\qquad
 \mathrm d\ell_G:=\|G^{-1/2}(\mathrm dG)G^{-1/2}\|_F.
\end{equation}
The divergence starts at zero, because $Z_0=G_0$. Its square root increases by at most $C\,\mathrm d\ell_G$ between updates and does not increase at updates. Integrating and accounting for the decreases of at least $r_t/2$ gives $\sum_{t=1}^{N}r_t\leq C(\sum_{t=1}^{N}\ell_t)^2$. Regularizing $\Psi$ by a positive constant and then taking the limit justifies the argument where $\Psi$ vanishes. This proves~\eqref{eq_more_rank_budget}. Part~(c) follows from $\E\|(I-\Gamma_t)\zeta\|_{Z_t}^2 =\tr[Z_t(I-\Gamma_t)\operatorname{Cov}(\zeta)(I-\Gamma_t)^\top] \leq\sigma^2r_t$.

\emph{Cost.} The matrix $X_Z^\top G^{-1}X_Z$ has norm $O(1)$ at every valid update, while $\|\widetilde Z^-\|_{\rm op}=O(m)$ follows from $Z^-\preceq(5/4)G$ and $\underline Z\succeq cG/m$. Thus $B=Cm^2$ bounds both Frobenius norms, which suffices for both calls to Lemma~\ref{lem:approximate_diagonalization}. Together with the remaining factorizations and matrix products, this gives Part~(d).
\end{proof}

We apply Lemma~\ref{lem_more_reference_metric} at the start of centering calls and keep the reference metric fixed within each call. The next corollary gives the resulting norm comparisons and refresh budget over an epoch.

\begin{corollary}[Reference metrics within a centering phase]
\label{cor:phase_reference_metric}
Suppose that queries $1,\ldots,j$ of a run are valid and that the proxy Hessian estimates of the calls that start at these queries satisfy the assumptions of Lemma~\ref{lem:lagged_metric}. At the start of an epoch, set $Z=G$. At a later call start in the same epoch, apply Lemma~\ref{lem_more_reference_metric} with the current weight and the new proxy Hessian estimate. At all other queries, keep $Z$ fixed and set $\Gamma=I$. Then, for the constant $C_Z$ of Section~\ref{sec:fresh_mode},
\begin{equation}\label{eq:phase_reference_sandwich}
 Z_i\preceq C_ZG_i,\qquad
 \nabla^2F(y_i)^{-1}\preceq C_ZZ_i/\alpha,\qquad
 \Gamma_i^\top Z_i\Gamma_i\preceq Z_{i-1}
\end{equation}
for $i\leq j$, where the third inequality holds at queries that do not start an epoch. Set $r_i=0$ at queries with identity transport, and use the charge of Lemma~\ref{lem_more_reference_metric} at call starts. Over every prefix of an epoch,
\begin{equation}\label{eq:phase_refresh_budget}
 \sum_ir_i\leq CH^2.
\end{equation}
The charges are used only in the analysis.
\end{corollary}
\begin{proof}
We compare every query with the start of its call, where the reference is constructed. Let $\mathcal H_\circ:=\nabla^2F(y^\circ)$ at the start $y^\circ$ of the call that contains query $i$. By~\eqref{eq:sc_hessian_comparison}, $(9/25)\mathcal H_\circ\preceq\nabla^2F(y)\preceq(25/9)\mathcal H_\circ$ on the phase ball, and any segment in this ball has $F$-local length at most $(5/3)(4/5)=4/3$. By Lemma~\ref{lem:weight_movement}(b), the weights of any two queries of the call have geodesic distance $O(\alpha^{-1/2}+\sqrt{\eta_{\rm w}})$ and are comparable by an absolute factor. Thus a reference constructed at the start of the call satisfies $Z\preceq3G_\circ\preceq C_ZG_i$ throughout the call. Clipping gives $Z\succeq c_Z\widehat Q^{-1}$, and at an epoch start $Z=G_\circ$ satisfies the same lower bound by~\eqref{eq:lagged_metric_bounds}. The lower Hessian bound in~\eqref{eq:lagged_metric_bounds} gives $c_Z\widehat Q^{-1}\succeq(c_Z\alpha/2)\mathcal H_\circ^{-1}$, and so $\nabla^2F(y_i)^{-1}\preceq(25/9)\mathcal H_\circ^{-1} \preceq(50/(9c_Z))Z_i/\alpha$. This gives the second bound of~\eqref{eq:phase_reference_sandwich} for large $C_Z$.

Apply Lemma~\ref{lem_more_reference_metric} to the first weight of the epoch and to the weights at its later call starts. Consecutive states of this subsequence lie in one phase ball, so their geodesic distance is bounded. The lemma gives the contraction and~\eqref{eq:phase_refresh_budget}, since at most $H$ reference updates occur during an epoch. Identity transport preserves the norm at all other queries.
\end{proof}

\subsection{Concentration under predictable contractions}
\label{sec:predictable_contractions}
The reference norm changes with the history of the algorithm. The following maximal inequality controls the accumulated error through its total conditional variance, and it allows the variance proxy to vary from query to query.

\begin{lemma}[Predictable contractions]
\label{lem:predictable_contractions}
Let $Z_j\succ0$ and, for $j\geq1$, $\Gamma_j$ be measurable with respect to a sigma-field $\mathfrak F_j^-$ that contains the preceding history and all auxiliary randomness revealed before $\xi_j$. Suppose that
\begin{align*}
 x_0=\xi_0,\qquad x_j=\Gamma_jx_{j-1}+\xi_j,\qquad
 \Gamma_j^\top Z_j\Gamma_j\preceq Z_{j-1},\qquad
 \E[\xi_j\mid\mathfrak F_j^-]=0.
\end{align*}
Let $v_j\geq0$ be $\mathfrak F_j^-$-measurable. Assume that, for every integer $k\geq2$,
\begin{align*}
 \E[\|\xi_j\|_{Z_j}^k\mid\mathfrak F_j^-]
 \leq k!C_0^k v_j^{k/2},\qquad
 \sum_{j=0}^{N} v_j\leq V
\end{align*}
almost surely, where $C_0$ is an absolute constant and $V>0$ is deterministic. Then, for every $a>0$,
\begin{equation}\label{eq:predictable_contraction_moment}
 \Pr\!\left[\max_{0\leq j\leq N}\|x_j\|_{Z_j}\geq a\right]
 \leq2\exp\!\left(\frac12-\frac{a}{4C_0\sqrt V}\right).
\end{equation}
\end{lemma}
\begin{proof}
We construct an exponential supermartingale from the function $\cosh(\lambda\|x\|_2)$ and stop it at the first exceedance. Let $\widetilde x_j:=Z_j^{1/2}x_j$, $\widetilde\xi_j:=Z_j^{1/2}\xi_j$, and $A_j:=Z_j^{1/2}\Gamma_jZ_{j-1}^{-1/2}$ for $j\geq1$. The maps $A_j$ are predictable contractions in the Euclidean norm, and the transformed noises remain conditionally centered. It suffices to prove the claim in these coordinates.

For $f(x)=\cosh(\lambda\|x\|_2)$, the radial and tangential eigenvalues of the Hessian at $x\neq0$ are $\lambda^2\cosh(\lambda\|x\|_2)$ and $\lambda\sinh(\lambda\|x\|_2)/\|x\|_2$, respectively. The inequality $\tanh u\leq u$ gives $\nabla^2f(x)\preceq\lambda^2f(x)I$, and at zero the same bound follows from $\nabla^2f(0)=\lambda^2I$. Since $f(x+tz)\leq f(x)e^{\lambda t\|z\|_2}$ for $0\leq t\leq1$, Taylor's integral formula gives
\begin{align*}
 f(x+z)\leq f(x)+\langle\nabla f(x),z\rangle
       +\frac{\lambda^2}{2}f(x)\|z\|_2^2e^{\lambda\|z\|_2}.
\end{align*}
For $\lambda C_0\sqrt V\leq1/2$, the conditional moment assumption and the exponential series imply
\begin{align*}
 \E[\|\widetilde\xi_j\|_2^2e^{\lambda\|\widetilde\xi_j\|_2}
       \mid\mathfrak F_j^-]
 &\leq C_0^2v_j
       \sum_{k\geq0}(k+1)(k+2)(\lambda C_0\sqrt{v_j})^k\\
 &=\frac{2C_0^2v_j}{(1-\lambda C_0\sqrt{v_j})^3}
 \leq16C_0^2v_j.
\end{align*}
Conditional centering removes the linear term in Taylor's formula, and $f(A_jx)\leq f(x)$. Therefore
\begin{align*}
 \mathcal M_j:=\exp\!\left(-8C_0^2\lambda^2\sum_{i=0}^jv_i\right)
                 \cosh(\lambda\|\widetilde x_j\|_2)
\end{align*}
is a nonnegative supermartingale, starting with expectation at most one. This includes the initial innovation if the preceding vector is taken to be zero. Stop at the first $j$ with $\|\widetilde x_j\|_2\geq a$, or at $N$ if there is none. On the former event, the stopped value is at least $e^{-8C_0^2\lambda^2V}\cosh(\lambda a)$. Consequently,
\begin{align*}
 \Pr\!\left[\max_j\|\widetilde x_j\|_2\geq a\right]
 \leq2e^{8C_0^2\lambda^2V-\lambda a}.
\end{align*}
The choice $\lambda=(4C_0\sqrt V)^{-1}$ proves the claim. The conditional moment bound justifies the exponential series and the integrability used above. The matrix square roots are used only in this analysis.
\end{proof}

The innovations of our estimator are quadratic polynomials in Gaussian probes. For an integer $k\geq2$, the bound $k^k\leq e^kk!$ converts the moment bound~\eqref{eq:conditional_gaussian_moment} with $p=k$ into the factorial moment hypothesis of Lemma~\ref{lem:predictable_contractions}.

\subsection{Joint maintenance within an epoch}
\label{sec_more_gradient_tracking}
We now describe the procedure $\mathcal D^{\rm m}_{\rm est}.\textsc{Query}$ and state its accuracy. An epoch begins with a full estimate at the first query of its first call. Every subsequent query uses a correction of the same fixed width. The reference norm is constant within each call and changes only at call starts. A standalone query counts as a call with one query, and all other calls have $N_{\rm c}=(N_G+1)N_y=O(s^2\log s)$ queries. Thus an epoch has at most $HN_{\rm c}$ queries.

We measure the pair in the joint norm of~\eqref{eq:joint_norm} with $Z=Z_j$, and we transport it by
\begin{equation}\label{eq:paired_norm_transport}
 \mathfrak P_j(v,X):=(\Gamma_jv,X).
\end{equation}
By Corollary~\ref{cor:phase_reference_metric}, $\mathfrak P_j$ is a contraction from the joint norm of $Z_{j-1}$ to that of $Z_j$, even when its construction depends on preceding estimation errors.

At query $j$, the flag $a_j$ records whether the query starts a call. If it does and the epoch continues, the reference is updated by Lemma~\ref{lem_more_reference_metric}. Otherwise the reference stays fixed and the transport is the identity. An epoch reset takes precedence and sets $Z_j=G_j$. After the endpoints and the transport are fixed, new common projections and shared residual probes give the correction
\begin{equation}\label{eq_more_projected_refresh}
 \widehat v_j:=(\widehat g_j^{\rm pair}-\Gamma_j\widehat g_{j-1}^{\rm pair},
                 \widehat R_j^{\rm pair}-\widehat R_{j-1}^{\rm pair}),
 \qquad v_j:=w_j-\mathfrak P_jw_{j-1}.
\end{equation}
The estimate at the preceding endpoint in this expression comes from a new sketch, distinct from the stored estimate $\widehat w_{j-1}$. Both endpoint estimates are unbiased for their targets, and therefore
\begin{equation}\label{eq:one_level_pair_recursion}
 \widehat w_j=\mathfrak P_j\widehat w_{j-1}+\widehat v_j,\qquad
 w_j=\mathfrak P_jw_{j-1}+v_j.
\end{equation}

\emph{Variance of a correction.} Let $\mathfrak F_j^-$ be the sigma-field generated by the history before the residual probes of query $j$. Write $e_j^{\rm pair}:=\widehat g_j^{\rm pair}-g_{G_j}(y_j)$, and similarly at the preceding state. The gradient error of the correction is
\begin{align*}
 \Gamma_j(e_j^{\rm pair}-e_{j-1}^{\rm pair})
                    +(I-\Gamma_j)e_j^{\rm pair}.
\end{align*}
At a call start, the first term is measured in $Z_j^-$, which is at most both $2G_j$ and $Z_{j-1}\preceq C_ZG_{j-1}$, so Lemma~\ref{lem_more_gradient_increments}(a) applies to it. The second term uses the current endpoint. By Lemma~\ref{lem:lagged_metric}, the proxy Hessian at this endpoint is at most $Cs$ times its estimate $\widehat Q$ of the phase, so Lemma~\ref{lem_more_gradient_increments}(c) and Lemma~\ref{lem_more_reference_metric}(c) apply to it. Together they give
\begin{equation}\label{eq:correction_variance}
 \E\bigl[\|\widehat v_j-v_j\|_{Z_j}^2\bigm|\mathfrak F_j^-\bigr]
 \leq V_0m(1+r_j)/q^2
\end{equation}
at width $q$, after enlarging the polylogarithmic constant $V_0$ of~\eqref{eq:full_estimate_variance}. At an interior query, the same bound follows with $\Gamma_j=I$ and $r_j=0$ from $Z_j=Z_{j-1}\preceq C_ZG_j,C_ZG_{j-1}$. The stationarity variance is bounded by Lemma~\ref{lem_more_gradient_increments}(d). A full estimate at width $q$ has joint variance at most $V_0m^2/q^2$ by~\eqref{eq:full_estimate_variance}. At a query that adds a correction, the centered error is a quadratic polynomial in the probes, conditional on the states, the stored proxy Hessian estimate, the transport and the projections. At a full estimate, we condition only on the state and the projections. The norm $Z=G$ is already fixed, and the proxy Hessian estimate is computed from the sketched factors afterwards.

\emph{Widths.} With $p$ from~\eqref{eq:fresh_width}, we choose the fixed parameters
\begin{equation}\label{eq_more_refresh_width}
 \Lambda:=C_*p^2\sqrt{N_{\rm c}V_0},\qquad
 q_{\rm ref}:=\lceil\Lambda m\rceil,\qquad
 q_{\rm corr}:=\lceil\Lambda H\sqrt m\rceil,
\end{equation}
where $C_*$ is a sufficiently large absolute constant. The requested projection rank and the residual probe count are of the indicated width, up to logarithmic factors. We also choose $C_*$ large enough that the sketched factors at every current endpoint satisfy Lemma~\ref{lem_more_metric}(a) with budget $\rho_{\rm met}$. Each range finder at a state receives failure budget $\rho/[80(N_{\max}+1)]$, including its internal repetitions. If the requested rank $r$ satisfies $2r+2\geq n$, the exact factors of~\eqref{eq:common_gram_exact_factors} replace the sketched factors, and the residual variance is zero.

\begin{algorithm}[!ht]
\caption{The procedure $\mathcal D^{\rm m}_{\rm est}.\textsc{Query}$}
\label{alg:more_gradient}
\begin{algorithmic}[1]
\Procedure{MaintainedQuery}{$y,G,\widehat J,a$} \Comment{$\mathcal D^{\rm m}_{\rm est}.\textsc{Query}$}
\If{this is the first query, or $a=1$ and $N_{\rm ep}=H$}
  \State $N_{\rm ep}\gets1$, $(\mathcal D^{\rm m}_{\rm est}.y,\mathcal D^{\rm m}_{\rm est}.G,\mathcal D^{\rm m}_{\rm est}.J)\gets(y,G,\widehat J)$, and $Z\gets G$.
  \State Compute a full estimate $(\widehat g,\widehat R)$ at width $q_{\rm ref}$ with the whitener $\widehat J$, and store its sketched factors in $Y$. \Comment{$\widetilde O(n^\omega+\mathcal C(m))$}
  \State $\widehat Q\gets\Call{ProxyHessian}{B_1,\ldots,B_m,Y,\rho_{\rm met}}$. \Comment{$\widetilde O(n^\omega+\mathcal C(m))$}
  \State \Return $(\widehat g,\widehat R,\widehat Q)$.
\EndIf
\State Compute the aligned whitener $\widehat J_1$ of $S(y)$ from $\mathcal D^{\rm m}_{\rm est}.J$ by Lemma~\ref{lem_more_aligned_whitener}. \Comment{$\widetilde O(n^\omega)$}
\If{$a=1$}
  \State $N_{\rm ep}\gets N_{\rm ep}+1$.
  \State $\widehat Q\gets\Call{ProxyHessian}{B_1^-,\ldots,B_m^-,Y,\rho_{\rm met}}$, where $B_i^-:=(\mathcal D^{\rm m}_{\rm est}.J)A_i(\mathcal D^{\rm m}_{\rm est}.J)^\top$. \Comment{$\widetilde O(n^\omega+\mathcal C(q))$ at the width $q$ of $Y$}
  \State Update $Z$ and compute $\Gamma$ by Lemma~\ref{lem_more_reference_metric} with $G$ and $\widehat Q$. \Comment{$\widetilde O(m^\omega)$}
\Else
  \State Keep $Z$ and $\widehat Q$, and $\Gamma\gets I$.
\EndIf
\State Fix both endpoint states and the transport before sketching.
  \label{line:more_gradient_freeze}
\State Draw new common projections for both states and shared probes, at width $q_{\rm corr}$. \Comment{$\widetilde O(n^\omega+\mathcal C(H\sqrt m))$}
\State Evaluate both endpoint pairs with the whiteners $\mathcal D^{\rm m}_{\rm est}.J$ and $\widehat J_1$, and form $\widehat v$ by~\eqref{eq_more_projected_refresh}. \Comment{$\widetilde O(n^\omega+\mathcal C(H\sqrt m))$}
\State $(\widehat g,\widehat R)\gets\mathfrak P(\widehat g,\widehat R)+\widehat v$, $(\mathcal D^{\rm m}_{\rm est}.y,\mathcal D^{\rm m}_{\rm est}.G,\mathcal D^{\rm m}_{\rm est}.J)\gets(y,G,\widehat J_1)$, and store the sketched factors of the current endpoint in $Y$.
\State \Return $(\widehat g,\widehat R,\widehat Q)$.\label{line:more_gradient_return}
\EndProcedure
\end{algorithmic}
\end{algorithm}

An epoch lasts $H$ calls, with no reset inside a call, and the final epoch may end earlier. At an epoch start, the full estimate is sketched before its proxy Hessian estimate is computed. At other call starts, the sketched factors of the preceding query supply the proxy Hessian estimate before the new correction is sketched. Because the call starts at the point of the preceding query, the coefficients $B_i^-$ are those of an exact whitener of $S(y)$. These orders ensure that each transport is fixed before its correction probes are drawn.

\begin{lemma}[Joint oracle maintenance]
\label{lem_more_gradient_tracking}
\label{lem:maintained_pair}
Let $m^2>n$, $1\leq H\leq\sqrt m$, $N_{\max}\geq1$, and let $\rho<1/10$ be the failure budget of $\mathcal D^{\rm m}_{\rm est}$.
\begin{itemize}
\item[(a)] There is an event of probability at least $1-\rho/2$ on which the following holds for every $j$. Suppose that queries $1,\ldots,j$ are valid, and that the proxy Hessian estimates of the calls that start at these queries satisfy the assumptions of Lemma~\ref{lem:lagged_metric}. Then
\begin{equation}\label{eq_more_tracked_error}
 \|\widehat w_j-w_j\|_{Z_j}\leq1.
\end{equation}
In particular, $\|\widehat R_j-R_j\|_F\leq\tau$ and $\|\widehat g_j-g_{G_j}(y_j)\|_{Z_j/\alpha}\leq\delta_g$.
\item[(b)] On every run, including a failed run, a prefix of $N\leq N_{\max}$ queries costs
\begin{equation}\label{eq_more_gradient_total}
 \widetilde O\bigl(\mathcal C(m)+N\mathcal B(H)\bigr).
\end{equation}
This includes the proxy Hessian estimates and all full estimates, and it also holds for a prefix that ends inside an epoch.
\end{itemize}
\end{lemma}
\begin{proof}[Proof of Part~(a)]
We apply Lemma~\ref{lem:predictable_contractions} to the error of the maintained pair in each epoch and stop the process before the first query at which its local hypotheses fail.

\emph{The recursion.} Consider an epoch, conditioned on the history before its first query. Let $\mathfrak F_j^-$ contain both endpoint states and all auxiliary randomness revealed before the residual probes of query $j$. At an epoch start this excludes the proxy Hessian estimate, which is computed from the full estimate afterwards. \eqref{eq:one_level_pair_recursion} gives
\begin{align*}
 \widehat w_j-w_j=\mathfrak P_j(\widehat w_{j-1}-w_{j-1})+\xi_j,\qquad
 \xi_j:=\widehat v_j-v_j,\qquad
 \E[\xi_j\mid\mathfrak F_j^-]=0.
\end{align*}
For the initial full estimate, let $\xi_0$ be its error. By~\eqref{eq:full_estimate_variance}, \eqref{eq:correction_variance} and~\eqref{eq_more_refresh_width}, conditional variance proxies for these innovations are
\begin{align*}
 \sigma_0^2:=V_0/\Lambda^2,\qquad
 \sigma_j^2:=\frac{V_0(1+r_j)}{\Lambda^2H^2}\quad(j\geq1).
\end{align*}
Corollary~\ref{cor:phase_reference_metric} implies, for every prefix with $N\leq HN_{\rm c}$ corrections,
\begin{align*}
 \sum_{j=0}^{N} \sigma_j^2
 \leq\frac{V_0}{\Lambda^2}
       \left(1+\frac{HN_{\rm c}+CH^2}{H^2}\right)
 \leq \frac{C_1N_{\rm c}V_0}{\Lambda^2}
 =\frac{C_1}{C_*^2p^4}.
\end{align*}
\eqref{eq:conditional_gaussian_moment} and the remark after Lemma~\ref{lem:predictable_contractions} supply the factorial conditional moments. Applying that lemma in the joint norms with $a=1/2$ bounds the probability that any error in this epoch exceeds $1/2$ by
\begin{align*}
 2\exp(1/2-cC_*p^2)\leq\rho/[8(N_{\max}+1)]
\end{align*}
for a sufficiently large absolute $C_*$. The bound covers every partial epoch and holds conditionally on its starting history.

\emph{Stopping.} To justify the local hypotheses without conditioning on future success, stop an auxiliary process before the first residual draw at which validity, the assumptions of Lemma~\ref{lem:lagged_metric}, or the guarantee of the projections fails. Thereafter, freeze the error and the norm, using identity transport, zero innovation and zero variance proxy. An invalid epoch start has zero error in a fixed positive definite norm. At a valid epoch start, the full estimate is sketched and its error is bounded in $Z=G$ before the proxy Hessian estimate is computed. A failure of this later estimate is charged to the events of Lemma~\ref{lem_more_metric} and stops the auxiliary process before its next residual draw. Thus every stopping decision for an innovation is made before the probes of that innovation. Every valid prefix has the variance bound above, so the stopped process satisfies Lemma~\ref{lem:predictable_contractions} on every history, and it agrees with the actual error until it stops.

\emph{Union bound.} At most $N_{\max}$ epochs start, so a union bound gives total failure probability less than $\rho/8$ for the errors. The budgets of the range finders, for at most $N_{\max}$ full-estimate states and $2N_{\max}$ coupled endpoints, sum to less than $\rho/8$. On the complement of these events, the actual error is at most $1/2$ at every query whose predecessors satisfy the local hypotheses, which proves Part~(a). Every call has the same target at a given state, so~\eqref{eq:one_level_pair_recursion} telescopes without a deterministic sketching bias.
\end{proof}

The maintained estimate itself need not be unbiased, because later transports may depend on earlier errors. The proof requires conditional centering only for each new correction.

\subsection{Cost of \texorpdfstring{$\mathcal D^{\rm m}_{\rm est}$}{the maintained data structure}}
\label{sec:oracle_cost_analysis}
\label{sec:oracle_single_point_cost}
\label{sec:oracle_gradient_cost}
\label{sec:operation_caps}
The sketch widths in~\eqref{eq_more_refresh_width} are fixed before the execution. This gives the same cost bound on successful and failed runs.

\begin{proof}[Proof of Lemma~\ref{lem_more_gradient_tracking}(b)]
We charge the corrections, the full estimates and the compressions of the proxy Hessian estimates separately. Every correction costs $\widetilde O(n^\omega+\mathcal C(H\sqrt m))$ by Lemma~\ref{lem_more_gradient_increments}. Its shared coefficient product is charged by $\Tmat(mn,n,H\sqrt m)$, and its separate products by $\mathcal C_{\rm sep}(H\sqrt m)$. Every completed epoch contains $H$ nonempty calls and therefore at least $H$ queries. A prefix of $N$ queries consequently has at most $1+N/H$ full estimates, whose total cost is
\begin{align*}
 \widetilde O\bigl(\mathcal C(m)+N\mathcal C(m)/H+Nn^\omega\bigr).
\end{align*}
The sketched factors of width $q$ have $\widetilde O(mq)$ columns, and Lemma~\ref{lem_more_metric}(b) compresses them in $\widetilde O(n^\omega+\mathcal C(q))$ operations. The stored sketched factors of a query are compressed at most twice, once at its own epoch start, if applicable, and once at the next call start. Charging these costs to the original sketches preserves these bounds, including for prefixes. Reference updates cost $\widetilde O(m^\omega)$, aligned whiteners cost $\widetilde O(n^\omega)$, weight factorizations and solves cost $O(m^\omega)$, and slack formation and coefficient pairings cost $O(mn^2)$. Since $m^\omega\leq mn^2\leq\mathcal C(q)$, all these costs fit the correction charge.

The widths, the repetitions of the range finder and the power counts are fixed. Thus this accounting holds even if a random estimate is inaccurate, and a failed factorization only terminates the query early. Charging a partially executed query in full also proves the prefix bound~\eqref{eq_more_gradient_total}.
\end{proof}

\subsection{Proof of the gradient/Hessian estimation guarantee for \texorpdfstring{$\mathcal D^{\rm m}_{\rm est}$}{the maintained data structure}}
\label{sec:maintained_proof}
\label{sec:maintained_center_interface}
\label{sec:maintained_query_cost}
\begin{proof}[Proof of Theorem~\ref{thm:estimator_interface} for $m^2>n$]
We combine the event of Lemma~\ref{lem_more_gradient_tracking}(a) with the events of the proxy Hessian estimates. The pair event has failure probability at most $\rho/2$. The event for $Q_Y$ and the mass of Lemma~\ref{lem_more_metric}(a) for the sketched factors at each current endpoint, and the compression event of Lemma~\ref{lem_more_metric}(b) for each call start, receive budget $\rho_{\rm met}=\rho/(4N_{\max})$. There are at most $N_{\max}$ events of each kind, so their total failure probability is at most $\rho/2$. Let $\mathcal E$ be the intersection of all these success events, which has probability at least $1-\rho$.

\emph{Accuracy.} Work on $\mathcal E$, and suppose that queries $1,\ldots,j$ are valid. We first show by induction that the proxy Hessian estimates of the call starts up to $j$ satisfy the assumptions of Lemma~\ref{lem:lagged_metric}. At a call start $i\leq j$, the estimate is compressed from the sketched factors of query $i$ at an epoch start, and from those of query $i-1$ otherwise. In the second case, $y_{i-1}=y_i$ and query $i-1$ is the last $G$-query of the preceding call or the standalone query. In both cases, the state of these sketched factors is valid and its weight is accurate by (V3), so on $\mathcal E$ Lemma~\ref{lem_more_metric} gives~\eqref{eq:estimator_metric} with $G'=G_i$ or $G'=G_{i-1}$. This is the assumption of Lemma~\ref{lem:lagged_metric}. Lemma~\ref{lem_more_gradient_tracking}(a) then gives $\|\widehat w_j-w_j\|_{Z_j}\leq1$, and~\eqref{eq:phase_reference_sandwich} gives
\begin{align*}
 \|\widehat g_j-g_{G_j}(y_j)\|_{\nabla^2F(y_j)^{-1}}
 \leq\sqrt{C_Z}\,\|\widehat g_j-g_{G_j}(y_j)\|_{Z_j/\alpha}
 \leq\sqrt{C_Z}\,\delta_g=\varepsilon_{\rm g}/4.
\end{align*}
Together with $\|\widehat R_j-R_j\|_F\leq\tau$, this proves~\eqref{eq:estimator_accuracy}, and Part~(a) follows.

\emph{Cost.} Lemma~\ref{lem_more_gradient_tracking}(b) is Part~(b), since $\mathcal B_*(H)=\mathcal B(H)$ when $m^2>n$.
\end{proof}

\subsection{Choosing the epoch length}
\label{sec:epoch_choice}
Longer epochs reduce the cost of full estimates but increase the width of corrections. The following choice balances the two costs.

\begin{corollary}[Choosing the epoch length]
\label{cor:more_epoch_choice}
Under the multiplication bounds in Fact~\ref{fact:rectangular_thresholds}, let $\omega$ be a square exponent bound with $2<\omega\leq9/4$, and choose $H:=\lfloor m^{1/6}\rfloor$. For every fixed $\zeta>0$,
\begin{equation}\label{eq_more_explicit_cost}
 \mathcal B(H)=\widetilde O_\zeta
       (n^{\omega+\zeta}+n^2m^{1+\zeta}).
\end{equation}
\end{corollary}
\begin{proof}
We bound the three terms of $\mathcal B(H)$ by splitting each product into blocks whose shapes appear in Fact~\ref{fact:rectangular_thresholds}. The chosen length satisfies $\tfrac12m^{1/6}\leq H\leq m^{1/6}$ and $H\sqrt m\leq m^{2/3}$. We use $\omega(2,1,2/3)=3$, $\omega(2,1,1)\leq25/8$ and $\omega\leq9/4$ from Fact~\ref{fact:rectangular_thresholds}, allowing fixed positive exponent losses whose sum is at most $\zeta$.

\emph{The shared product of a correction.} Split the $mn$ rows and the $n$ inner coordinates into blocks of sizes $m^2$ and $m$, respectively. Then
\begin{align*}
 \Tmat(mn,n,H\sqrt m)
 \leq O((n/m)^2)\Tmat(m^2,m,\lceil m^{2/3}\rceil)
 \leq O_\zeta(n^2m^{1+\zeta}).
\end{align*}

\emph{Full estimates.} The same blocking at width $m$, and permutation symmetry for the other two products, give
\begin{align*}
 \Tmat(mn,n,m),\ \Tmat(n,n,m^2)
 &\leq O_\zeta(n^2m^{9/8+\zeta}),\\
 \Tmat(m,m,nm)&\leq O_\zeta(nm^{17/8+\zeta})
              \leq O_\zeta(n^2m^{9/8+\zeta}).
\end{align*}
Thus $\mathcal C(m)/H=O_\zeta(n^2m^{23/24+\zeta})$.

\emph{The separate products of a correction.} Take $q=m^{2/3}$ and an admissible square exponent $\tau_\omega:=\omega+\zeta'$ with $0<\zeta'\leq\min\{\zeta,1/100\}$. Splitting the wide dimensions into square blocks gives
\begin{align*}
 \mathcal C_{\rm sep}(q)
 &\leq O_{\zeta'}(n^{\tau_\omega}+mqn^{\tau_\omega-1}+nqm^{\tau_\omega-1}+mnq)\\
 &\leq O_{\zeta'}(n^{\tau_\omega}+n^2m^{\tau_\omega-4/3})
 \leq O_{\zeta'}(n^{\omega+\zeta'}+n^2m^{11/12+\zeta'}).
\end{align*}
The middle bound uses $m\leq n$, $2<\tau_\omega<3$ and $n^{\tau_\omega-3}\leq m^{\tau_\omega-3}$. Adding these bounds and $n^\omega$ proves~\eqref{eq_more_explicit_cost}, and rounding affects only constants.
\end{proof}

The same bound holds for every epoch length $H$ with $m^{1/8}\leq H\leq m^{1/6}$. Then $H\sqrt m\leq m^{2/3}$, the full estimates cost $O_\zeta(n^2m^{9/8+\zeta}/H)\leq O_\zeta(n^2m^{1+\zeta})$, and the separate products cost $O_\zeta(n^{\omega+\zeta}+n^2m^{3/4+\zeta}H) \leq O_\zeta(n^{\omega+\zeta}+n^2m^{11/12+\zeta})$.
\section{The slack maintenance data structure}
\label{sec:reused_factorization}
This section implements $\mathcal D_{\rm slack}$ and proves Theorem~\ref{thm:slack_interface}. If $m^4>n$, the data structure factors the slack at every query, and this factorization is part of the per-query charge $\mathcal B_*(H)$. If $m^4\leq n$, the data structure reuses an \emph{anchor} factorization while the objective multiplier stays within a factor of two of the multiplier of the anchor. A polynomial of degree $d_{\rm w}=\widetilde O(m)$ supplies the whitening actions between anchors. This regime uses $\mathcal D^{\rm f}_{\rm est}$, so a whitening error affects only one estimate, and Section~\ref{sec:implicit_accuracy} bounds its effect.

The data structure stores the anchor objective $\mathcal D_{\rm slack}.c$ and the Cholesky factor $\mathcal D_{\rm slack}.L$ of the slack at the query that set the anchor. Section~\ref{sec:slack_window} shows that the slack stays spectrally close to the anchor slack, and Section~\ref{sec:polynomial_whitening} turns this window into whitening actions. Section~\ref{sec:slack_interface_proof} proves Theorem~\ref{thm:slack_interface}, and Section~\ref{sec:reused_running_time} bounds the cost of a fresh query with these actions.

\subsection{One factorization for a range of path parameters}
\label{sec:slack_window}
Along a run of calls, write the objective as $tc^\circ$ with a monotone multiplier $t$. The next comparison uses the Bregman divergence of $F$ in both directions.

\begin{lemma}[Slack comparison within a factor-two range]
\label{thm:slack_window}
\label{lem:central_bregman}
\label{lem:query_near_center}
Suppose that $\phi_{tc^\circ}$ has a minimizer $x(t)$ for every $t$ in an interval $\mathcal I\subseteq(0,\infty)$. Consider two calls of \textnormal{\textsc{Center}} with objectives $t_{\rm anc}c^\circ$ and $tc^\circ$, where $t_{\rm anc},t\in\mathcal I$ and $t_{\rm anc}/2\leq t\leq2t_{\rm anc}$, and suppose that each call starts with decrement at most $1/10$. If $y_{\rm anc}$ and $y$ lie in the phase balls of these calls, then
\begin{equation}\label{eq:slack_window}
 \Upsilon^{-1}S(y_{\rm anc})\preceq S(y)\preceq\Upsilon S(y_{\rm anc}),
\end{equation}
with $\Upsilon=4m\sqrt\alpha=4000sm$ from~\eqref{eq:whitening_degree}. For the zero objective, the same conclusion holds with the factor $1.01$ in place of $\Upsilon$ whenever the minimizer exists.
\end{lemma}
\begin{proof}
We first compare the two central points through the Bregman divergence $\Delta_F(z,y):=F(z)-F(y)-\nabla F(y)^\top(z-y)$, then convert divergence into local distance, and finally compare each query with the central point of its call. Exchanging the roles of the two calls if necessary, we may assume $t\geq t_{\rm anc}$.

\emph{Central points.} For $0\leq u\leq\log(t/t_{\rm anc})\leq\log2$, let $x_u:=x(t_{\rm anc}e^u)$. Differentiating the optimality condition $\nabla F(x_u)=-t_{\rm anc}e^uc^\circ$ gives $\partial_ux_u=\nabla^2F(x_u)^{-1}\nabla F(x_u)$, where the implicit function theorem applies because $\nabla^2F(x_u)\succ0$. Hence
\begin{align*}
 \partial_uF(x_u)=\|\nabla F(x_u)\|_{\nabla^2F(x_u)^{-1}}^2\leq\nu,
 \qquad
 \partial_u(c^{\circ\top}x_u)=-\frac{\partial_uF(x_u)}{t_{\rm anc}e^u},
\end{align*}
where the bound is the gradient bound of a $\nu$-self-concordant barrier. Substituting these identities into the two divergences gives
\begin{align}
 \Delta_F(x_u,x_0)
 &=\int_0^u(1-e^{-v})\,\partial_vF(x_v)\,\d v
 \leq\nu(u+e^{-u}-1)<\nu/4,
 \label{eq:forward_bregman}\\
 \Delta_F(x_0,x_u)
 &=\int_0^u(e^{u-v}-1)\,\partial_vF(x_v)\,\d v
 \leq\nu(e^u-1-u)<\nu/3.
 \label{eq:reverse_bregman}
\end{align}

\emph{From divergence to distance.} Let $d:=\|z-y\|_{\nabla^2F(y)}$. The function $\psi(v):=(z-y)^\top\nabla^2F(y+v(z-y))(z-y)$ satisfies $|(\psi^{-1/2})'|\leq1$ by self-concordance, so $\psi(v)\geq d^2/(1+vd)^2$ and
\begin{align*}
 \Delta_F(z,y)=\int_0^1(1-v)\psi(v)\,\d v
 \geq\int_0^1\frac{(1-v)d^2}{(1+vd)^2}\,\d v
 =d-\log(1+d)\geq d/2-1.
\end{align*}
By~\eqref{eq:forward_bregman} and~\eqref{eq:reverse_bregman} and $\nu=\alpha m\geq6$, the distance between the two central points is at most $2\nu/3+2\leq\nu$ in the local norm at either of them. By Lemma~\ref{lem:ml_rounding} with $h=z-y$, a local distance $d$ at $y$ gives $\|S(y)^{-1/2}S(z)S(y)^{-1/2}-I\|_{\rm op}\leq3d/\sqrt\alpha$. Applied at the two central points, with $3\nu/\sqrt\alpha=3m\sqrt\alpha$, this gives
\begin{align*}
 S(x(t))\preceq(1+3m\sqrt\alpha)S(x(t_{\rm anc})),\qquad
 S(x(t_{\rm anc}))\preceq(1+3m\sqrt\alpha)S(x(t)).
\end{align*}
Both directions are needed because $3m\sqrt\alpha>1$.

\emph{Queries.} For a call that starts at $y^\circ$ with central point $x_c$ and decrement at most $1/10$, \eqref{eq:sc_center_estimates} gives $\|y^\circ-x_c\|_{\nabla^2F(y^\circ)}\leq1/9$, and~\eqref{eq:sc_hessian_comparison} gives $\nabla^2F(x_c)\preceq(9/8)^2\nabla^2F(y^\circ)$. Hence every query $y$ in the phase ball satisfies
\begin{equation}\label{eq:query_near_center}
 \|y-x_c\|_{\nabla^2F(x_c)}
 \leq\tfrac98\bigl(\tfrac25+\tfrac19\bigr)=\tfrac{23}{40}<\tfrac35,\qquad
 (1-10^{-3})S(x_c)\preceq S(y)\preceq(1+10^{-3})S(x_c),
\end{equation}
where the second claim follows from Lemma~\ref{lem:ml_rounding} and $9/(5\sqrt\alpha)\leq10^{-3}$. Combining the comparisons of the central points and of the queries in either direction gives the factor
\begin{align*}
 \frac{1+10^{-3}}{1-10^{-3}}(1+3m\sqrt\alpha)
 \leq4m\sqrt\alpha=\Upsilon.
\end{align*}
For the zero objective, the central points coincide, so only the two query comparisons are needed, and their ratio is less than $1.01$.
\end{proof}

A query with objective $c$ \emph{continues} the anchor if $c=\theta\,\mathcal D_{\rm slack}.c$ for some $\theta\in[1/2,2]$, or if $c=\mathcal D_{\rm slack}.c=0$. Otherwise it sets a new anchor. The objective is fixed during a call of \textsc{Center}, so this $O(m)$ test can fail only at the first query of a call, and anchors change only between calls. Lemma~\ref{thm:slack_window} therefore supplies the spectral window of every anchor. Anchor changes leave the weight unchanged.

\begin{lemma}[Number of anchors]
\label{lem:anchor_count}
Consider a run of calls with objectives $tc^\circ$, where $c^\circ\neq0$ and the multipliers $t$ are monotone in $[t_{\min},t_{\max}]$ with $t_{\min}>0$. The run sets at most $1+\lceil\log_2(t_{\max}/t_{\min})\rceil$ anchors, and a run with zero objective sets one anchor. In particular, Algorithm~\ref{alg:more_sdp} sets at most
\begin{equation}\label{eq:anchors_core}
 A_{\rm core}:=3+\log_2\left(1+\frac{2\nu}{t_0\epsilon}\right)
 =\widetilde O(L_\epsilon)
\end{equation}
anchors, and this bound is fixed before any Gaussian sketching matrix is drawn.
\end{lemma}
\begin{proof}
The multipliers of successive anchors in a monotone run differ by more than a factor of two, so the range of the multipliers bounds the number of anchors. The schedule of Algorithm~\ref{alg:more_sdp} ends at a multiplier at most $(1+\gamma_t)\max\{t_0,2\nu/\epsilon\}$, and $1+\gamma_t\leq2$. Substitution gives~\eqref{eq:anchors_core}. A run with zero objective continues its first anchor throughout.
\end{proof}

\subsection{Polynomial whitening and its implementation}
\label{sec:polynomial_whitening}
Let $J_{\rm anc}:=(\mathcal D_{\rm slack}.L)^{-1}$ be the whitener of the anchor slack. At a query with slack $S$, the slack relative to the anchor and the exact whitener aligned with the anchor are
\begin{equation}\label{eq:relative_slack}
 \Sigma:=J_{\rm anc}SJ_{\rm anc}^\top,\qquad
 J:=\Sigma^{-1/2}J_{\rm anc}.
\end{equation}
Inside the window of Lemma~\ref{thm:slack_window}, the spectrum of $\Sigma$ lies in $[\Upsilon^{-1},\Upsilon]$. Using the relative approximation of Lemma~\ref{lem:inverse_sqrt_polynomial}, define
\begin{equation}\label{eq:implicit_whitener}
 \widehat J:=\chi_\Upsilon(\Sigma)J_{\rm anc},
 \qquad \widehat J^\top=J_{\rm anc}^\top\chi_\Upsilon(\Sigma).
\end{equation}
The evaluation of $\chi_\Upsilon(\Sigma)X$ by Horner's rule uses only the action
\begin{equation}\label{eq:relative_slack_action}
 \Sigma X=J_{\rm anc}\bigl(S(J_{\rm anc}^\top X)\bigr),
\end{equation}
which consists of two triangular solves with $\mathcal D_{\rm slack}.L$ and one product with the stored slack. At the query that sets the anchor, $\Sigma=I$, and the data structure uses $J_{\rm anc}$ itself.

\begin{lemma}[Cost and quality of polynomial whitening]
\label{lem:implicit_whitening_cost}
Applying either action in~\eqref{eq:implicit_whitener} to a block with $n$ rows and $q$ columns costs $O(d_{\rm w}\Tmat(n,n,q))$, where $d_{\rm w}=\widetilde O(m)$ is the degree in~\eqref{eq:whitening_degree}. If the spectrum of $\Sigma$ lies in $[\Upsilon^{-1},\Upsilon]$, then $\widehat J=\mathsf KJ$ for the exact whitener $J$ of~\eqref{eq:relative_slack}, where $\mathsf K$ is symmetric and $\|\mathsf K-I\|_{\rm op}\leq\varepsilon_{\rm w}$.
\end{lemma}
\begin{proof}
Each application of $\Sigma$ in~\eqref{eq:relative_slack_action} costs $O(\Tmat(n,n,q))$ by Lemma~\ref{lem:fast_factorizations}, and the evaluation of Lemma~\ref{lem:inverse_sqrt_polynomial} with $\Upsilon=4m\sqrt\alpha$ uses $d_{\rm w}$ of them. The same lemma gives $\widehat J=\chi_\Upsilon(\Sigma)\Sigma^{1/2}J=\mathsf KJ$ with the stated properties of $\mathsf K$.
\end{proof}

\begin{algorithm}[!ht]
\caption{The procedure $\mathcal D_{\rm slack}.\textsc{Whiten}$}
\label{alg:implicit_slack}
\label{alg:apply_whitening}
\begin{algorithmic}[1]
\Procedure{Whiten}{$y,c$} \Comment{$\mathcal D_{\rm slack}.\textsc{Whiten}$}
\State Form and store $S\gets\sum_i y_iA_i-C$. \Comment{$O(mn^2)$}
\If{$m^4>n$, or no anchor is stored, or $c$ does not continue the anchor} \Comment{$O(m)$}
  \State Factor $S=LL^\top$, reporting failure if this factorization fails. \Comment{$O(n^\omega)$}
  \If{$m^4\leq n$}
    \State $(\mathcal D_{\rm slack}.c,\mathcal D_{\rm slack}.L)\gets(c,L)$.
  \EndIf
  \State \Return the exact actions $X\mapsto L^{-1}X$ and $X\mapsto L^{-\top}X$. \Comment{$O(\Tmat(n,n,q))$ per $n\times q$ block}
\EndIf
\State Let $\Sigma$ act by~\eqref{eq:relative_slack_action}. \Comment{$O(\Tmat(n,n,q))$ per $n\times q$ block}
\State \Return $X\mapsto\chi_\Upsilon(\Sigma)(J_{\rm anc}X)$ and $X\mapsto J_{\rm anc}^\top(\chi_\Upsilon(\Sigma)X)$. \Comment{$O(d_{\rm w}\Tmat(n,n,q))$ per block}
\EndProcedure
\end{algorithmic}
\end{algorithm}

The polynomial has a fixed degree and is defined for every symmetric slack, so the actions are also defined on histories outside the window. A factorization or solve that is undefined reports failure.

\subsection{Proof of the slack maintenance guarantee}
\label{sec:slack_interface_proof}
\begin{proof}[Proof of Theorem~\ref{thm:slack_interface}]
Parts~(a), (c) and~(d) follow from the algorithm and the lemmas above, and Part~(b) follows from the window of Lemma~\ref{thm:slack_window}.

\emph{Part~(a).} If $m^4>n$, the first branch of Algorithm~\ref{alg:implicit_slack} runs at every query.

\emph{Part~(b).} Let $m^4\leq n$, and let queries $1,\ldots,j$ be valid. If query $j$ sets an anchor, it returns the exact whitener $L^{-1}$, and $\mathsf K=I$. Otherwise, the anchor was set at an earlier query $j_{\rm a}$, with objective $c_{\rm a}:=\mathcal D_{\rm slack}.c$, and query $j$ belongs to a call whose objective continues this anchor. This objective is $\theta c_{\rm a}$ with $\theta\in[1/2,2]$, or $c_{\rm a}=0$ and the objective is zero. By (V1), this call starts with decrement at most $1/10$, and its phase ball contains $y_j$.

We next find a call with objective $c_{\rm a}$ that starts with decrement at most $1/10$ and whose phase ball contains $y_{j_{\rm a}}$. If query $j_{\rm a}$ is the first query of a call, this call has these properties by (V1). Otherwise, query $j_{\rm a}$ is the standalone query, and $c_{\rm a}=0$ by Definition~\ref{def:valid_queries}. After the standalone query, anchors change only at first queries of calls, and the anchor set by query $j_{\rm a}$ is still stored at query $j$. Hence the first call continued this anchor, and so its objective is zero. This call starts at $y_{j_{\rm a}}$, with decrement at most $1/10$ by (V1), so it has the required properties.

Apply Lemma~\ref{thm:slack_window} to these two calls, with $c^\circ=c_{\rm a}$, $t_{\rm anc}=1$ and $t=\theta$ if $c_{\rm a}\neq0$, and in its zero-objective case if $c_{\rm a}=0$. Its hypothesis on the minimizers is the assumption of Part~(b). It shows that the spectrum of $\Sigma$ lies in $[\Upsilon^{-1},\Upsilon]$, and Lemma~\ref{lem:implicit_whitening_cost} gives $\widehat J=\mathsf KJ$ with $\|\mathsf K-I\|_{\rm op}\leq\varepsilon_{\rm w}$.

\emph{Part~(c).} This is Lemma~\ref{lem:anchor_count}.

\emph{Part~(d).} Forming $S$ uses $m$ scaled additions of $n\times n$ matrices, which cost $O(mn^2)$. A Cholesky factorization costs $O(n^\omega)$, and an exact action costs $O(\Tmat(n,n,q))$, by Lemma~\ref{lem:fast_factorizations}. The cost of a polynomial action is Lemma~\ref{lem:implicit_whitening_cost}.
\end{proof}

\subsection{Cost of fresh queries with polynomial whitening}
\label{sec:reused_running_time}
We now bound the cost of $\mathcal D^{\rm f}_{\rm est}.\textsc{Query}$ when its whitening actions are polynomial. The degree $d_{\rm w}$ multiplies only the products with the whitening actions, while the product with the stacked coefficients is computed once.

\begin{lemma}[Cost of one query with polynomial whitening]
\label{lem:implicit_query_cost}
Excluding the anchor factorization, a fresh query of width $q$ costs
\begin{equation}\label{eq:implicit_query_cost}
 \widetilde O\bigl(mn^2+m^\omega+\Tmat(mn,n,q)
                     +d_{\rm w}\mathcal C_{\rm sep}(q)\bigr).
\end{equation}
The suppressed factors include power steps and trials of the range finder.
\end{lemma}
\begin{proof}
We follow the circuit of Lemma~\ref{lem:random_implicit_blocks}(b) with $\widehat J$ in place of the exact whitener. Applying $\widehat J^\top$, the stacked coefficients and $\widehat J$ costs $d_{\rm w}\Tmat(n,n,q)$, $\Tmat(mn,n,q)$ and $d_{\rm w}\Tmat(n,n,mq)$, respectively, by Lemma~\ref{lem:implicit_whitening_cost}. Coefficient mixing costs $\Tmat(m,m,nq)$, and transposition gives the cost of the adjoint. Gram products lie in $\mathcal C_{\rm sep}(q)$. For the gradient, concatenate the sketched factors in a matrix $Y$ with $\widetilde O(mq)$ columns. Form $\widehat J^\top Y$ and its Gram matrix, and then pair it with every $A_i$. These steps cost $\widetilde O(d_{\rm w}\Tmat(n,n,mq))$, $\widetilde O(\Tmat(n,n,mq))$ and $O(mn^2)$. Slack formation costs $O(mn^2)$, and the weight and preconditioner factorizations cost $O(m^\omega)$.

The requested rank $r$ of a Hutch++ projection may exceed $m$. If $2r+2\geq n$, the query uses the exact factors of~\eqref{eq:common_gram_exact_factors}. Otherwise the two ranges have $n$ and $mn$ rows and $\widetilde O(r)$ columns. Orthogonalizing these blocks and factoring their Gram matrices cost $\widetilde O(\Tmat(mn,r,r))$, where $r\leq n$ and $r=\widetilde O(q)$. This is at most $\widetilde O(\Tmat(mn,n,q))$, which is already charged.
\end{proof}

We now sum the query costs of Lemma~\ref{lem:implicit_query_cost} and include initialization and anchor factorizations. This gives a bound for every prefix of queries, including prefixes of failed runs.

\begin{theorem}[Cost with polynomial whitening]
\label{thm:implicit_amortized_cost}
Suppose $m^4\leq n$. Every prefix of $N$ queries that sets $A$ anchors costs
\begin{equation}\label{eq:implicit_prefix_cost}
 \widetilde O\bigl(An^\omega+N\mathcal B_{\rm imp}\bigr),
\end{equation}
with $\mathcal B_{\rm imp}$ from~\eqref{eq:implicit_batched_cost}. On every run, the same bound holds with a deterministic upper bound on the number of anchors in place of $A$. In particular, Theorem~\ref{thm:estimator_interface}(b) holds when $m^4\leq n$.
\end{theorem}
\begin{proof}
Every query has width $q_{\rm fresh}=\widetilde O(m)$, so Lemma~\ref{lem:implicit_query_cost} gives the charge $\widetilde O(\mathcal B_{\rm imp})$, and the proxy Hessian compression at a call start has no larger width. Initializing the weight uses the exact factor of the first anchor and costs $\widetilde O(n^\omega+\mathcal C(m))$, where $\mathcal C(m)\leq\mathcal B_{\rm imp}$ because $d_{\rm w}\geq1$. The $A$ exact factorizations cost $O(An^\omega)$. Fixed widths, polynomial degrees and trial counts give the same accounting for failed runs and partially executed queries.
\end{proof}

To simplify this cost when $m^4\leq n$, we first bound the separated block products using the rectangular multiplication estimates of Fact~\ref{fact:rectangular_thresholds}.

\begin{lemma}[Thin-product bound]
\label{lem:thin_products}
Suppose $m^4\leq n$. Under Fact~\ref{fact:rectangular_thresholds}, for every fixed $\zeta>0$ and $1\leq q\leq m$,
\begin{equation}\label{eq:thin_products}
 \mathcal C_{\rm sep}(q)\leq O_\zeta(n^2m^\zeta).
\end{equation}
\end{lemma}
\begin{proof}
We bound the three terms of $\mathcal C_{\rm sep}(q)$ by blocking. By monotonicity it suffices to take $q=m$, and for $m=1$ every term is $O(n^2)$. Otherwise, split both dimensions $n$ of $\Tmat(n,n,m^2)$ into blocks of size $m^4$. Each of the $O(n^2/m^8)$ products has shape $(m^4,m^4,m^2)$, so the bound $\omega(1,1,1/2)=2$ gives $\Tmat(n,n,m^2)=O_\zeta(n^2m^\zeta)$. For the coefficient mixing, splitting the last dimension into blocks of size $m$ gives $\Tmat(m,m,nm)=O_\zeta(nm^{\omega+\zeta})\leq O_\zeta(n^2m^\zeta)$, because $\omega\leq9/4<4$ and $m^4\leq n$. The remaining linear term is $nm^2\leq n^2$. Rounding and the change of the base dimension are absorbed by the exponent slack.
\end{proof}

Combining Lemma~\ref{lem:thin_products} with the coefficient-product bound of Lemma~\ref{lem:fresh_query_cost} and the polynomial degree gives the following explicit per-query cost.

\begin{corollary}[Cost of a fresh query with polynomial whitening]
\label{cor:reused_periods}
Suppose $m^4\leq n$. With the whitening tolerance of~\eqref{eq:bias_budget} and under Fact~\ref{fact:rectangular_thresholds}, for every fixed $\zeta>0$,
\begin{equation}\label{eq:implicit_explicit_cost}
 \mathcal B_{\rm imp}=\widetilde O_\zeta(n^2m^{1+\zeta}).
\end{equation}
\end{corollary}
\begin{proof}
Lemma~\ref{lem:fresh_query_cost} bounds the coefficient product $\Tmat(mn,n,m)\leq\mathcal C(m)$ by $O_\zeta(n^2m^{1+\zeta})$. Lemma~\ref{lem:thin_products} and $d_{\rm w}=\widetilde O(m)$ give $d_{\rm w}\mathcal C_{\rm sep}(m)=\widetilde O_\zeta(n^2m^{1+\zeta})$. The remaining term satisfies $m^\omega\leq mn^2$.
\end{proof}
\section{Proof of the main algorithmic result}
\label{sec:randomized_solver}
Lemma~\ref{lem:central_start_correctness} proves Part~(a) of Theorem~\ref{thm_more_amortized_sdp} and the bounds on the numbers of $y$-steps and queries. This section bounds the arithmetic cost. We first add the costs of the three data structures over a run, and then specialize the bound to the current multiplication exponents.

\subsection{The cost of a run}
\label{sec:running_time}
\label{sec_more_total}
\begin{lemma}[Cost of a run]
\label{lem:run_cost}
Let $\mathsf S:=\textnormal{\textsc{InitializeOracleState}}(\mathcal P,H,N_{\max},\gamma)$ with $1\leq H\leq\sqrt m$, and consider a run that uses $\mathsf S$. On every run, including a failed one, every prefix of $N\leq N_{\max}$ queries costs
\begin{equation}\label{eq:run_cost}
 \widetilde O\bigl(An^\omega+\mathcal C(m)+N\mathcal B_*(H)\bigr).
\end{equation}
Here $A$ is a deterministic upper bound on the number of anchors set during the prefix if $m^4\leq n$, and $A=0$ otherwise.
\end{lemma}
\begin{proof}
We add the costs of the three data structures and of the remaining work of \textsc{Center}, using that every term of the form $mn^2$ or $m^\omega$ is at most $\mathcal B_*(H)$. Indeed, $\mathcal B_{\rm imp}\geq mn^2$ by definition, and $\mathcal C(q)\geq\Tmat(mn,n,q)\geq mn^2$ because a product reads its input. Moreover, $m^\omega\leq mn^2$.

\emph{Slack.} By Theorem~\ref{thm:slack_interface}(d), forming the slack costs $O(mn^2)$ per query. If $m^4>n$, every query makes one factorization of cost $O(n^\omega)$, which is part of $\mathcal B_*(H)$. Otherwise the prefix makes at most $A\geq1$ factorizations, of total cost $O(An^\omega)$.

\emph{Weight.} By Theorem~\ref{thm:weight_interface}(e), the initialization costs $\widetilde O(n^\omega+\mathcal C(m))$. Its term $n^\omega$ is covered by $An^\omega$ if $m^4\leq n$, and by $\mathcal B_*(H)$ otherwise. By Theorem~\ref{thm:weight_interface}(d), each update costs $O(m^\omega)$.

\emph{Gradient and Hessian estimation.} By Theorem~\ref{thm:estimator_interface}(b), $\mathcal D_{\rm est}$ costs $\widetilde O(\mathcal C(m)+N\mathcal B_*(H))$, including the applications of the whitening actions.

\emph{Center.} Each call scales the Cholesky factor of $\widehat Q$ to obtain that of its phase preconditioner, and each step query solves with this factor and updates the point. These steps cost $O(m^2)$ per query.

The widths, trial counts, polynomial degrees and query cap are fixed in advance, and a failed check only terminates the run. Charging a partially executed query in full therefore proves~\eqref{eq:run_cost} on every run.
\end{proof}

Table~\ref{tab:query_cost_terms} lists the work of one query term by term in the three cases of~\eqref{eq:effective_query_cost}, together with the technique that controls each term and the bounds that the proof of Part~(c) uses.

\begin{table}[!ht]
\centering\small
\begin{tabularx}{\textwidth}{@{}>{\raggedright\arraybackslash}X>{\raggedright\arraybackslash}p{3.5cm}lccc@{}}
\toprule
 & & & \multicolumn{3}{c}{Cost per query}\\
\cmidrule(l){4-6}
Term & Technique & Section & $m^4\leq n$ & $m^2\leq n<m^4$ & $m^2>n$\\
\midrule
slack formation and pairings with the $A_i$ & & & \multicolumn{3}{c}{$mn^2$}\\
\addlinespace
slack factorization & reuse between anchors & \ref{sec:reused_factorization} & at anchors & $n^\omega$ & $n^\omega$\\
\addlinespace
sketch width $q$ & Hutch++ sketches, coupled increments & \ref{sec:common_gram}, \ref{sec:oracle_maintenance} & $\widetilde O(m)$ & $\widetilde O(m)$ & $\widetilde O(H\sqrt m)$\\
\addlinespace
product with the stacked $A_i$ & block operations & \ref{sec:random_implicit_operations} & \multicolumn{3}{c}{$\Tmat(mn,n,q)$}\\
\addlinespace
whitening and mixing of $m$ blocks & exact or polynomial whitener & \ref{sec:polynomial_whitening} & $d_{\rm w}\mathcal C_{\rm sep}(q)$ & $\mathcal C_{\rm sep}(q)$ & $\mathcal C_{\rm sep}(q)$\\
\addlinespace
$G$-step and $y$-step & damped inverse update & \ref{sec:cost_weight_iteration} & \multicolumn{3}{c}{$m^\omega$}\\
\addlinespace
alignment and reference update & spectral clipping & \ref{sec_more_reference_metric} & -- & -- & $n^\omega$\\
\midrule
\multicolumn{3}{@{}l}{charge $\mathcal B_*(H)$ in~\eqref{eq:effective_query_cost}} & $\mathcal B_{\rm imp}$ & $n^\omega+\mathcal C(m)$ & $\mathcal B(H)$\\
\addlinespace
\multicolumn{3}{@{}l}{bound on the charge, in $\widetilde O_\zeta(\cdot)$} & $n^2m^{1+\zeta}$ & \multicolumn{2}{c}{$n^{\omega+\zeta}+n^2m^{1+\zeta}$}\\
\bottomrule
\end{tabularx}
\caption{The work of one query, term by term, in the three cases of~\eqref{eq:effective_query_cost}. In the third case, the full estimate of width $\widetilde O(m)$ that starts an epoch is amortized over the at least $H$ queries of the epoch. The last row uses Fact~\ref{fact:rectangular_thresholds}, with $H=\lfloor m^{1/6}\rfloor$ in the third case, and its bounds are proved in Corollary~\ref{cor:reused_periods}, Lemma~\ref{lem:fresh_query_cost}, and Corollary~\ref{cor:more_epoch_choice}, respectively.}
\label{tab:query_cost_terms}
\end{table}

\subsection{Proof of the main theorem}
\begin{proof}[Proof of Theorem~\ref{thm_more_amortized_sdp}]
Part~(a) and the counts in Part~(b) are Lemma~\ref{lem:central_start_correctness}. In particular, the run makes
\begin{align*}
 N_{\rm tot}=O\left(s^3\log s\sqrt m\,
                  \log\left(2+\frac\nu{t_0\epsilon}\right)\right)
 =\widetilde O(\sqrt mL_\epsilon)
\end{align*}
queries, because $\log(2+\nu/(t_0\epsilon))\leq\log\alpha+L_\epsilon$.

\emph{Part~(b).} If $m^4\leq n$, Lemma~\ref{lem:anchor_count} gives the deterministic bound $A_{\rm core}=\widetilde O(L_\epsilon)$ on the number of anchors. Lemma~\ref{lem:run_cost} with $N_{\max}=N_{\rm tot}$ therefore gives the prefix bound \eqref{eq_more_prefix_cost}. In every regime, $\mathcal C(m)\leq\sqrt m\,\mathcal B_*(H)$. Indeed, $\mathcal C(m)\leq\mathcal B_*(H)$ in the first two cases of~\eqref{eq:effective_query_cost}, and $\mathcal C(m)/\sqrt m\leq\mathcal C(m)/H\leq\mathcal B(H)$ in the third. Substituting $N=N_{\rm tot}$ proves~\eqref{eq_more_sdp_total}.

\emph{Part~(c).} We bound $\mathcal B_*(H)$ in each case of~\eqref{eq:effective_query_cost}. If $m^4\leq n$, Corollary~\ref{cor:reused_periods} gives $\mathcal B_{\rm imp}=\widetilde O_\zeta(n^2m^{1+\zeta})$. If $m^2\leq n<m^4$, Lemma~\ref{lem:fresh_query_cost} gives $\mathcal C(m)=O_\zeta(n^2m^{1+\zeta})$. If $m^2>n$, Corollary~\ref{cor:more_epoch_choice} gives $\mathcal B(H)=\widetilde O_\zeta(n^{\omega+\zeta}+n^2m^{1+\zeta})$. In the last two cases, the slack factorizations add $n^\omega$ per query, and $n<m^4$. Thus, with $0<\zeta'\leq\zeta/4$ for their factorization costs and $\omega\leq9/4$,
\begin{align*}
 \sqrt m\,n^{\omega+\zeta'}
 \leq n^2\sqrt m\,n^{1/4+\zeta'}
 \leq n^2m^{3/2+4\zeta'}
 \leq n^2m^{3/2+\zeta}.
\end{align*}
The coefficient costs also fit $\widetilde O_\zeta(n^2m^{3/2+\zeta})$ after multiplication by $\sqrt m$. Substituting these bounds into~\eqref{eq_more_sdp_total} proves~\eqref{eq:reused_core_explicit}.
\end{proof}
\section{Preprocessing the initial point}
\label{sec:preprocessing}
We construct a near-central starting point from the SDP data and the feasibility bounds $B,\rho$ by a two-phase method. An \emph{auxiliary initialization SDP}, which penalizes a scalar shift of the slack, has a bounded domain with an explicit interior point. \emph{Auxiliary path following} toward the barrier center of this SDP (reverse path following) then leads to the desired central path. Preprocessing returns only the point and a positive multiplier, and the core initializes its own oracle state.

\subsection{Constructing a bounded pencil}
\label{sec:bounded_pencils}
\label{sec:complete_solver_description}
Assume that a dual optimum $y_*$ and a strictly feasible point $y_{\rm f}$ satisfy
\begin{align*}
 \|y_*\|_\infty,\|y_{\rm f}\|_\infty\leq B,
 \qquad S(y_{\rm f})\succeq\rho I,
 \qquad B\geq1,\quad\rho>0.
\end{align*}
The algorithm receives only the bounds $B$ and $\rho$. For the prescribed accuracy $0<\epsilon\leq1$, define
\begin{equation}\label{eq:penalty_parameters}
 a:=\Bigl(\sum_i\|A_i\|_F^2\Bigr)^{1/2},\quad z_{\max}:=2(\|C\|_F+\rho+1),\quad K:=1+\frac{4B\|b\|_1+1}\rho,\quad \delta:=\frac\epsilon{16K}.
\end{equation}
We measure input scale by
\begin{equation}\label{eq:complete_kappa}
 \kappa:=2+m+n+B+a+\|C\|_F+\|b\|_1+\rho+\rho^{-1}+\epsilon^{-1}.
\end{equation}
We introduce a scalar shift $z$ and the augmented slack
\begin{equation}\label{eq:penalty_pencil}
 \overline S(y,z)
 \mathrel{:=}\operatorname{diag}(S(y)+zI,2B\mathbf1+y,2B\mathbf1-y,
                         z_{\max}-z,\delta+z),
\end{equation}
denoted by $\overline{\mathcal P}$, with objective $c:=(b,K)$ and explicit interior point $u_0:=(0,z_{\max}/2)$. We call the SDP with this pencil and objective the auxiliary initialization SDP. The penalty forces a near-optimal shift to be negative, and its magnitude is controlled by $K\delta=\epsilon/16$ and $\delta\leq\min\{\rho/16,1/16\}$.

\begin{lemma}[The auxiliary initialization SDP]
\label{lem:penalty_pencil}
The pencil $\overline{\mathcal P}$ has $\bar m:=m+1$ variables and order $\bar n:=n+2m+2$. Its coefficient matrices are linearly independent and have trivial common kernel, and its primal and dual problems are strictly feasible. The Euclidean ball of radius
\begin{align*}
 r:=\min\{B/2,z_{\max}/4,1/[2(a+1)]\}
\end{align*}
about $u_0$ is feasible. The diameter of the domain and the width of the objective on the domain are at most
\begin{align*}
 R:=4B\sqrt m+2z_{\max},\qquad
 W_c:=4B\|b\|_1+K(z_{\max}+\delta).
\end{align*}
Moreover, every feasible $u=(y,z)$ with
\begin{align*}
 c^\top u-\inf_{\overline S\succ0}c^\top u\leq\epsilon/64
\end{align*}
satisfies
\begin{align*}
 -\delta<z\leq-\delta/4,\qquad
 S(y)\succ0,\qquad
 0\leq b^\top y-\mathrm{OPT}\leq3\epsilon/64.
\end{align*}
\end{lemma}
\begin{proof}
\emph{Interior point and bounds.} At $u_0$, the leading block of $\overline S$ satisfies
\begin{align*}
 -C+(z_{\max}/2)I\succeq(\rho+1)I.
\end{align*}
For $u=(y,z)$ with $\|u-u_0\|_2\leq r$, the leading block changes by at most
\begin{align*}
 \Bigl\|\sum_iy_iA_i+(z-z_{\max}/2)I\Bigr\|_{\rm op}
 \leq a\|y\|_2+|z-z_{\max}/2|
 \leq(a+1)\|u-u_0\|_2\leq1/2,
\end{align*}
The choices $r\leq B/2$ and $r\leq z_{\max}/4$ also keep the scalar blocks positive. The bounds $|y_i|<2B$ and $-\delta<z<z_{\max}$, with $\delta<z_{\max}$, give $R$ and $W_c$.

\emph{Structure.} The scalar blocks distinguish every variable, so the augmented coefficient matrices are linearly independent. In a vector of their common kernel, the scalar blocks force the scalar coordinates to vanish, and the identity block in the coefficient of $z$ forces the remaining coordinates to vanish. For primal strict feasibility, choose any positive definite leading primal block. Each variable has a pair of scalar coefficients with opposite signs, and for any required difference $d$, the positive numbers $1+\max\{d,0\}$ and $1+\max\{-d,0\}$ have difference $d$. Choosing these pairs to satisfy the coefficient equations gives a positive definite augmented primal matrix.

\emph{Nonnegative shifts.} Suppose that a feasible point has $z\geq0$. The point $y':=(\rho y+zy_{\rm f})/(\rho+z)$ is feasible for the original problem, because
\begin{align*}
 S(y')=\frac{\rho S(y)+zS(y_{\rm f})}{\rho+z}
 \succeq\frac{\rho(S(y)+zI)}{\rho+z}\succ0.
\end{align*}
The optimality of $y_*$ and $b^\top y_{\rm f}-\mathrm{OPT}\leq2B\|b\|_1$ therefore give
\begin{align*}
 \mathrm{OPT}\leq\frac{\rho b^\top y+zb^\top y_{\rm f}}{\rho+z},\qquad
 b^\top y\geq\mathrm{OPT}-\frac{2B\|b\|_1}\rho z,
 \qquad
 c^\top u\geq\mathrm{OPT}+\Bigl(K-\frac{2B\|b\|_1}\rho\Bigr)z\geq\mathrm{OPT}.
\end{align*}

\emph{The auxiliary SDP has a value below $\mathrm{OPT}$.} The comparison point
\begin{align*}
 y_\delta:=(1-\delta/\rho)y_*+(\delta/\rho)y_{\rm f}
\end{align*}
satisfies $\|y_\delta\|_\infty\leq B$ and $S(y_\delta)\succeq\delta I$. Hence $(y_\delta,-\delta)$ lies in the closure of the augmented domain, and since $K\geq4B\|b\|_1/\rho$,
\begin{equation*}
 c^\top(y_\delta,-\delta)
 \leq\mathrm{OPT}+\frac{2B\|b\|_1}\rho\delta-K\delta
 \leq\mathrm{OPT}-K\delta/2
 =\mathrm{OPT}-\epsilon/32.
\end{equation*}
Convex combinations with $u_0$ approach this point from the interior. A feasible point with the stated objective error therefore satisfies
\begin{align*}
 c^\top u\leq\mathrm{OPT}-\epsilon/64<\mathrm{OPT},
\end{align*}
and the preceding paragraph excludes $z\geq0$.

\emph{Accuracy.} Now $S(y)\succ-zI\succ0$, so $y$ is feasible for the original problem and $b^\top y\geq\mathrm{OPT}$. Using this bound once for $z$, and then $z>-\delta$ for the original objective, gives
\begin{align*}
 Kz&=c^\top u-b^\top y\leq-\epsilon/64=-K\delta/4,\\
 b^\top y&=c^\top u-Kz
 \leq\mathrm{OPT}-\epsilon/64+K\delta
 =\mathrm{OPT}+3\epsilon/64.
\end{align*}
\end{proof}

\subsection{The preprocessing algorithm}
\label{sec:bounded_stage_description}
For this subsection and the next two, $(m,n)$ are the dimensions of a bounded pencil $\mathcal P$. The routine receives a target objective $c$, an interior ball $(u_0,r)$ and a diameter bound $R$, and uses the centering parameters of~\eqref{eq:more_centering_parameters} and the slack access of Algorithm~\ref{alg:implicit_slack}. Its first gradient estimate $g_0$ makes $u_0$ nearly central for the auxiliary objective $c_0=-g_0$. The routine then follows the central path of $F+\tau c_0^\top u$ toward the barrier center by decreasing $\tau$ from one, and it introduces $t_0c$ at the end. It uses the parameters
\begin{equation}\label{eq:bounded_stage_schedule}
\begin{gathered}
 H:=\lfloor m^{1/6}\rfloor,\qquad
 G_{\max}:=2(\nu+1)/r,\qquad
 t_0:=\begin{cases}
 \min\{1,\delta_{\rm c}/(R\|c\|_2)\},&c\ne0,\\
 1,&c=0,
 \end{cases}\\
 N_{\rm aux}:=\left\lceil\frac{\max\{0,\log(RG_{\max}/\delta_{\rm c})\}}
                            {-\log(1-\gamma_t)}\right\rceil,
 \qquad N_{\rm pre}:=(N_{\rm aux}+3)N_{\rm c}+1,
\end{gathered}
\end{equation}
which are fixed before any Gaussian sketching matrix is drawn.
\begin{algorithm}[!ht]
\caption{Preprocessing a bounded pencil to obtain a central start}
\label{alg:bounded_stage}
\begin{algorithmic}[1]
\Procedure{PreprocessBoundedStart}{$\mathcal P,c,u_0,r,R,\gamma$} \Comment{Lemma~\ref{lem:bounded_stage}}
\State Compute $H,G_{\max},t_0,N_{\rm aux},N_{\rm pre}$ from~\eqref{eq:bounded_stage_schedule}.
\State $\mathsf S\gets\Call{InitializeOracleState}{\mathcal P,H,N_{\rm pre},\gamma}$. \Comment{Section~\ref{sec:centering_description}}
\State $(g_0,\cdot,\cdot)\gets\Call{BarrierOracle}{\mathcal P,0,u_0,1;\mathsf S}$.\Comment{Initial phase query}
\State Report failure if $\|g_0\|_2>G_{\max}$, and $c_0\gets-g_0$ and $\tau\gets1$.
\State $u\gets\Call{Center}{\mathcal P,c_0,u_0;\mathsf S}$. \Comment{Algorithm~\ref{alg:maintained_center}}
\For{$j=1,\ldots,N_{\rm aux}$} \Comment{auxiliary path following}
  \State $\tau\gets(1-\gamma_t)\tau$.
  \State $u\gets\Call{Center}{\mathcal P,\tau c_0,u;\mathsf S}$.
\EndFor
\State $u\gets\Call{Center}{\mathcal P,0,u;\mathsf S}$.
\If{$c\ne0$}
  \State $u\gets\Call{Center}{\mathcal P,t_0c,u;\mathsf S}$.
\EndIf
\State Discard $\mathsf S$ and \Return $(u,t_0)$.
\EndProcedure
\end{algorithmic}
\end{algorithm}

For the original input, apply this routine to $\overline{\mathcal P}$ and $c=(b,K)$, with $u_0,r,R$ from Lemma~\ref{lem:penalty_pencil}. Return these pencil data and bounds with the computed point for use in recovery.

\subsection{Correctness of preprocessing}
\label{sec:bounded_stage}
\label{sec:cost_euclidean_control}
To analyze preprocessing, we need to compare local norms with the Euclidean bounds supplied for the feasible region. The next lemma bounds the barrier Hessian near the central path and the gradient at the center of the interior ball.

\begin{lemma}[A Euclidean bound near the central path]
\label{lem:cost_bounded_metric}
Let $F$ be a $\nu$-self-concordant barrier on a bounded convex domain $\Omega\subseteq\R^m$ that contains a Euclidean ball of radius $r>0$, and let $W_c:=\sup_{u,v\in\Omega}|c^\top(u-v)|$. If $u\in\Omega$ and $t\geq0$ satisfy $\lambda_{tc}(u)<1$, then
\begin{equation}\label{eq:cost_bounded_metric}
 \sqrt{\lambda_{\max}(\nabla^2F(u))}\leq
 \frac{2\nu+2+tW_c}{(1-\lambda_{tc}(u))r}.
\end{equation}
At the center $u_0$ of the interior ball, one also has
\begin{align*}
 \|\nabla F(u_0)\|_2\leq\nu/r,\qquad
 \sqrt{\lambda_{\max}(\nabla^2F(u_0))}\leq(3\nu+2)/r.
\end{align*}
\end{lemma}
\begin{proof}
For $x,z\in\Omega$, the barrier gradient inequality states
\begin{align*}
 \nabla F(x)^\top(z-x)\leq\nu.
\end{align*}
Indeed, the derivative $\phi(v)$ of $F(x+v(z-x))$ satisfies $\phi'\geq\phi^2/\nu$. If $\phi(0)>\nu$, integration would force a blow-up of $\phi$ before $v=1$, which contradicts the smoothness of $F$ on the feasible segment.

Let $\ell:=\|z-u\|_{\nabla^2F(u)}$. The self-concordant lower Hessian bound along the segment gives
\begin{align*}
 \bigl(\nabla F(u+(z-u)/2)-\nabla F(u)\bigr)^\top(z-u)
 \geq\int_0^{1/2}\frac{\ell^2}{(1+v\ell)^2}\,\d v
 =\frac{\ell^2}{\ell+2}.
\end{align*}
The gradient inequality at the midpoint toward $z$ bounds $\nabla F(u+(z-u)/2)^\top(z-u)$ by $2\nu$. Since $\ell^2/(\ell+2)\geq\ell-2$, we obtain
\begin{equation}\label{eq:cost_barrier_radius}
 \|z-u\|_{\nabla^2F(u)}\leq2\nu+2-\nabla F(u)^\top(z-u).
\end{equation}
Writing $\nabla F(u)=(\nabla F(u)+tc)-tc$ gives
\begin{align*}
 (1-\lambda_{tc}(u))\|z-u\|_{\nabla^2F(u)}\leq2\nu+2+tW_c.
\end{align*}
For the endpoints $z_\pm$ of any diameter of the interior ball, $\|z_+-z_-\|_{\nabla^2F(u)}\leq\|z_+-u\|_{\nabla^2F(u)}+ \|z_--u\|_{\nabla^2F(u)}$; dividing by $2r$ proves~\eqref{eq:cost_bounded_metric}.

At the center $u_0$, the gradient inequality in two opposite radial directions gives $\|\nabla F(u_0)\|_2\leq\nu/r$. For each radius $d$ of length $r$, this bounds $-\nabla F(u_0)^\top d$ by $\nu$, and~\eqref{eq:cost_barrier_radius} gives $\|d\|_{\nabla^2F(u_0)}\leq3\nu+2$, which is the remaining Hessian bound. If the ball is open, we take limits from smaller radii.
\end{proof}

We now apply Theorem~\ref{thm:centering_correctness} to the sequence of objectives in Algorithm~\ref{alg:bounded_stage}. Lemma~\ref{lem:cost_bounded_metric} controls the initial gradient estimate, while the diameter bound permits the auxiliary objective to be removed and the target objective to be introduced.

\begin{lemma}[A central starting point for a bounded pencil]
\label{lem:bounded_stage}
Let $1\leq m\leq n$, and let $\mathcal P$ satisfy the structural assumptions of the barrier. Suppose that its domain $\Omega$ has diameter at most $R$ and contains the Euclidean ball of radius $r$ about $u_0$. Let $c$ be a target objective for which the primal problem is strictly feasible, and let $0<\gamma<1/10$. With probability at least $1-\gamma$, Algorithm~\ref{alg:bounded_stage} returns $u_{\rm c}\in\Omega$ and the multiplier $t_0>0$ of~\eqref{eq:bounded_stage_schedule} with
\begin{align*}
 \lambda_{t_0c}(u_{\rm c})
 =\|\nabla F(u_{\rm c})+t_0c\|_{\nabla^2F(u_{\rm c})^{-1}}
 \leq\delta_{\rm c}.
\end{align*}
It uses at most $N_{\rm pre}$ oracle queries.
\end{lemma}
\begin{proof}
The initial query and all calls share one state with failure budget $\gamma$. Their objectives depend only on preceding outputs, and the bounded domain ensures that every $F+c'^\top u$ has a minimizer. Thus Theorem~\ref{thm:centering_correctness} applies once we verify its assumption~(i), and we work on its success event.

\emph{The auxiliary objective.} By Theorem~\ref{thm:centering_correctness}, the initial query gives
\begin{equation*}
 \|g_0-\nabla F(u_0)\|_{\nabla^2F(u_0)^{-1}}
 \leq\varepsilon_{\rm g}<\delta_{\rm c},\qquad
 \lambda_{c_0}(u_0)\leq\delta_{\rm c}.
\end{equation*}
Thus the first centering call is admissible. Moreover, Lemma~\ref{lem:cost_bounded_metric} bounds the Euclidean norm of the same estimate,
\begin{align*}
 \|g_0\|_2\leq\|\nabla F(u_0)\|_2
       +\delta_{\rm c}\sqrt{\lambda_{\max}(\nabla^2F(u_0))}
 \leq\frac\nu r+\frac{3\nu+2}{100r}
 \leq G_{\max},
\end{align*}
so the explicit check of $\|g_0\|_2$ passes.

\emph{Decreasing the objective.} At a point with $\lambda_{\tau c_0}(u)\leq\delta_{\rm c}$, the barrier gradient bound gives $\|\tau c_0\|_{\nabla^2F(u)^{-1}}\leq\sqrt\nu+\delta_{\rm c}$. Consequently
\begin{align*}
 \lambda_{(1-\gamma_t)\tau c_0}(u)
 \leq\delta_{\rm c}+\gamma_t(\sqrt\nu+\delta_{\rm c})<1/10.
\end{align*}

\emph{Removing the auxiliary objective.} The unit Dikin ellipsoid at $u$ lies in the domain, whose diameter is at most $R$. Hence
\begin{align*}
 \nabla^2F(u)^{-1}\preceq R^2I,
 \qquad
 \|v\|_{\nabla^2F(u)^{-1}}\leq R\|v\|_2.
\end{align*}
At the end of the loop, the choice of $N_{\rm aux}$ gives
\begin{align*}
 \tau R\|c_0\|_2
 \leq(1-\gamma_t)^{N_{\rm aux}}RG_{\max}\leq\delta_{\rm c}.
\end{align*}
Removing $\tau c_0$ therefore leaves decrement at most $2\delta_{\rm c}<1/10$, and the next call returns a point with decrement at most $\delta_{\rm c}$ for $F$.

\emph{Introducing the target.} Introducing $t_0c$ changes the gradient by a vector of local dual norm at most $t_0R\|c\|_2\leq\delta_{\rm c}$. The last call therefore also starts with decrement at most $2\delta_{\rm c}$ and returns the claimed bound. If $c=0$, the point returned for $F$ already has this property with $t_0=1$.

There are at most $N_{\rm aux}+3$ centering calls after the initial query, each with $N_{\rm c}$ queries, so the routine uses at most $N_{\rm pre}$ queries.
\end{proof}

\subsection{Running time and the preprocessing guarantee}
\label{sec:bounded_stage_cost}
It remains to bound the work of the preprocessing schedule, including initialization and slack factorizations. We apply Lemma~\ref{lem:run_cost} and count the anchors created as the auxiliary objective decreases.

\begin{lemma}[Work of the preprocessing routine]
\label{lem:bounded_stage_work}
With the parameters of~\eqref{eq:bounded_stage_schedule}, Algorithm~\ref{alg:bounded_stage} has arithmetic cost
\begin{align*}
 \widetilde O\bigl(A_{\rm pre}n^\omega+\mathcal C(m)
                    +N_{\rm pre}\mathcal B_*(H)\bigr)
\end{align*}
on every run, where $\mathcal B_*$ is defined in~\eqref{eq:effective_query_cost} and
\begin{equation}\label{eq:anchors_preprocessing}
 A_{\rm pre}:=6+\log_2(2+RG_{\max}/\delta_{\rm c}).
\end{equation}
\end{lemma}
\begin{proof}
If $m^4>n$, Lemma~\ref{lem:run_cost} with $A=0$ gives the claim. Otherwise the decreasing schedule satisfies
\begin{align*}
 N_{\rm aux}\log_2\frac1{1-\gamma_t}
 \leq\log_2(2+RG_{\max}/\delta_{\rm c})+1.
\end{align*}
Lemma~\ref{lem:anchor_count} bounds its anchors by $3+\log_2(2+RG_{\max}/\delta_{\rm c})$ (one if $c_0=0$). The initial query, removing the auxiliary objective and introducing $t_0c$ add at most three anchors. Thus the deterministic bound $A_{\rm pre}$ and Lemma~\ref{lem:run_cost} give the claim, including initialization and failed runs.
\end{proof}

Combining the auxiliary construction of Lemma~\ref{lem:penalty_pencil} with the correctness and work bounds of Lemmas~\ref{lem:bounded_stage} and~\ref{lem:bounded_stage_work} gives the complete preprocessing guarantee.

\begin{theorem}[Preprocessing guarantee]
\label{thm:preprocessing}
Let $1\leq m\leq n$, and let $A_1,\ldots,A_m,C$ be real symmetric $n\times n$ matrices, with the $A_i$ linearly independent and $\bigcap_i\ker A_i=\{0\}$. Let $b\in\R^m$. Suppose that the supplied numbers $B\geq1$ and $\rho>0$ certify the existence of a dual optimum $y_*$ and a point $y_{\rm f}$ with
\begin{align*}
 \|y_*\|_\infty,\|y_{\rm f}\|_\infty\leq B,
 \qquad S(y_{\rm f})\succeq\rho I.
\end{align*}
Given $0<\epsilon\leq1$ and $0<\gamma<1/10$, preprocessing forms the pencil $\overline{\mathcal P}$ of the auxiliary initialization SDP and applies Algorithm~\ref{alg:bounded_stage} as described above. With probability at least $1-\gamma$, it returns
\begin{align*}
 (\overline{\mathcal P},c,u_{\rm c},t_0,r,R,W_c),
\end{align*}
where $r,R,W_c$ are the bounds of Lemma~\ref{lem:penalty_pencil}, $t_0>0$, and
\begin{align*}
 \overline S(u_{\rm c})\succ0,
 \qquad
 \|\nabla F(u_{\rm c})+t_0c\|_{\nabla^2F(u_{\rm c})^{-1}}
 \leq1/100.
\end{align*}
Here $F$ is the matrix Lee--Sidford barrier of $\overline{\mathcal P}$, which has $\bar m$ variables and order $\bar n$, satisfies the structural assumptions, and has strictly feasible primal and dual problems.

Under Fact~\ref{fact:rectangular_thresholds} and a square exponent bound $2<\omega\leq9/4$, for every fixed $\zeta>0$, the work on every run is
\begin{align*}
 \widetilde O_\zeta\left(
 (n^{\omega+\zeta}+n^2m^{3/2+\zeta})
 \log^{O(1)}\kappa\right),
\end{align*}
where $\kappa$ is the input scale in~\eqref{eq:complete_kappa}. The suppressed factors include logarithms of $1/\gamma$.
\end{theorem}
\begin{proof}
Lemmas~\ref{lem:penalty_pencil} and~\ref{lem:bounded_stage}, applied in dimensions $(\bar m,\bar n)$, give the output guarantee. These dimensions satisfy $\bar m=O(m)$, $\bar n=O(n)$ and $\bar m\leq\bar n$. The quantities $r^{-1},R,W_c,\|c\|_2$ and $t_0^{-1}=\max\{1,R\|c\|_2/\delta_{\rm c}\}$ are polynomially bounded in $\kappa$. Since $-\log(1-\gamma_t)\geq\gamma_t$,
\begin{align*}
 N_{\rm aux}&=O\bigl(1+\sqrt\nu\log(2+RG_{\max})\bigr),\\
 N_{\rm pre}&=O\bigl(s^2\log s[1+\sqrt\nu\log(2+RG_{\max})]\bigr)
 =\widetilde O(\sqrt m\log\kappa),
 \qquad A_{\rm pre}=O(\log\kappa).
\end{align*}
Substitute these bounds into Lemma~\ref{lem:bounded_stage_work} and apply the dimension-regime estimates in the proof of Theorem~\ref{thm_more_amortized_sdp}(c). They give the claimed cost for each branch; initialization is absorbed by $\mathcal C(\bar m)\leq\sqrt{\bar m}\,\mathcal B_*(H)$. Constructing the augmented coefficients and their bounds costs $O(mn^2)$. Finally, replace $(\bar m,\bar n)$ by $(O(m),O(n))$. All bounds apply on failed runs as well.
\end{proof}
\section{Postprocessing}
\label{sec:postprocessing}
From the final near-central point, recovery draws a fresh positive semidefinite gradient representative and corrects its moments by a modified Newton iteration. The iteration moves along positive congruences of the representative, which keep it positive semidefinite, and its residual contracts geometrically, so high accuracy requires only logarithmically many iterations.

\subsection{Algorithm and output guarantee}
\label{sec:recovery_description}
\begin{theorem}[Postprocessing a near-central dual point]
\label{thm:postprocessing}
Let $1\leq m\leq n$, and let $\overline{\mathcal P}$ be the pencil of the auxiliary initialization SDP in~\eqref{eq:penalty_pencil}, constructed under the assumptions of Lemma~\ref{lem:penalty_pencil} with $0<\epsilon\leq1$. Let $c:=(b,K)$, let $r$ and $W_c$ be the bounds of Lemma~\ref{lem:penalty_pencil}, and let $\nu=\alpha\bar m$ and $\gamma_t=\delta_{\rm c}/(4\sqrt\nu)$ be the parameters of Theorem~\ref{thm_more_amortized_sdp} for $\overline{\mathcal P}$. Suppose that the input point $u=(y,z)$ and the multiplier $t$ satisfy
\begin{gather*}
 \overline S(u)\succ0,\qquad z<0,\qquad
 \lambda_{tc}(u)\leq1/100,\qquad \frac{128\nu}{\epsilon}\leq t
 \leq(1+\gamma_t)\frac{128\nu}{\epsilon}.
\end{gather*}
Let $0<\gamma<1/10$. With probability at least $1-\gamma$, Algorithm~\ref{alg:more_congruence_recovery} returns an explicit matrix $X\succeq0$ with
\begin{equation*}
 \|\mathcal A(X)-b\|_1\leq7\epsilon/16,\qquad
 b^\top y-\langle C,X\rangle\leq5\epsilon/16.
\end{equation*}
This probability bound holds conditionally on any preceding history that fixes an input satisfying these assumptions, because recovery uses fresh randomness. Lemma~\ref{lem:postprocessing_work} bounds the arithmetic cost on every run.
\end{theorem}

Since $z<0$, the unchanged point $y$ is dual feasible, and the objective guarantee implies $\langle C,X\rangle\geq\operatorname{OPT}-5\epsilon/16$.

Let $J$ whiten $\overline S(u)$ exactly and $B_i=J\overline A_iJ^\top$ in the original coefficient basis. For $\overline X=(\alpha/t)J^\top MJ$, augmented primal feasibility means $\tr[B_iM]=(c_{\rm rec})_i$, where $c_{\rm rec}=(t/\alpha)c$. Lemma~\ref{lem:more_recovery_density} supplies $M_{\rm init}\succeq0$ with small moment residual and an estimate $\widehat Q$ of the proxy Hessian. Starting at $w=0$, Algorithm~\ref{alg:more_congruence_recovery} parametrizes the candidates by the congruence $M(w)=(I+B_w)M_{\rm init}(I+B_w)$, which keeps them positive semidefinite, and applies the modified Newton iteration $w\gets w-\widehat Q^{-1}f(w)/2$ to the residual $f(w)_i=\tr[B_iM(w)]-(c_{\rm rec})_i$. The Jacobian of $f$ at $w=0$ is $2Q_0$, where $Q_0[h,h]:=\tr[M_{\rm init}B_h^2]$, and the iteration uses the fixed approximation $2\widehat Q$ of this Jacobian at every iterate.

The number of iterations depends on $\widehat Q$, so recovery first checks the explicit cap
\begin{equation}\label{eq:more_recovery_metric_cap}
 \tr[\widehat Q]\leq U_Q
 :=\frac{2.2\bar m}{\alpha}
       \left(\frac{2\nu+2+tW_c}{0.99r}\right)^2
\end{equation}
and reports failure if it fails. It then chooses the accuracy and the number of iterations
\begin{equation}\label{eq:more_recovery_tolerance}
 \delta_{\rm rec}:=\min\left\{10^{-4},
 \frac{t\epsilon}{16\alpha B\sqrt{\bar m}\sqrt{\tr[\widehat Q]}}\right\},\qquad
 N_{\rm rec}:=\left\lceil
       \log_4\frac{2.24\cdot10^{-3}}{\delta_{\rm rec}}\right\rceil.
\end{equation}
\begin{algorithm}[!ht]
\caption{Primal recovery from the auxiliary initialization SDP}
\label{alg:more_congruence_recovery}
\begin{algorithmic}[1]
\Procedure{RecoverPrimal}{$\overline{\mathcal P},c,u,t,B,r,W_c,\epsilon,\gamma$} \Comment{Theorem~\ref{thm:postprocessing}}
\State Use the parameters $\bar m,\alpha,\nu$ and the coefficients $\overline A_i$ of $\overline{\mathcal P}$.
\State Factor $\overline S(u)=LL^\top$, and set $J\gets L^{-1}$.
\State Obtain $M_{\rm init}$ and $\widehat Q$ from Lemma~\ref{lem:more_recovery_density} with failure budget $\gamma$.
\State Compute $U_Q$ by~\eqref{eq:more_recovery_metric_cap}, and report failure if $\tr[\widehat Q]>U_Q$.
\State $c_{\rm rec}\gets(t/\alpha)c$ and $w\gets0$, and compute $\delta_{\rm rec},N_{\rm rec}$ by~\eqref{eq:more_recovery_tolerance}.
\For{$j=0,\ldots,N_{\rm rec}-1$} \Comment{modified Newton iteration}
  \State Form $B_w\gets J(\sum_iw_i\overline A_i)J^\top$ and $M(w)\gets(I+B_w)M_{\rm init}(I+B_w)$.
  \State Evaluate $f(w)_i\gets\tr[J\overline A_iJ^\top M(w)]-(c_{\rm rec})_i$.
  \State $w\gets w-\widehat Q^{-1}f(w)/2$.
\EndFor
\State Recompute $M(w)$ and $f(w)$.
\State Report failure if $\|f(w)\|_{\widehat Q^{-1}}>\delta_{\rm rec}$ or $\tr[M(w)]>7\bar m$.
\State Form $\overline X\gets(\alpha/t)J^\top M(w)J$, and let $X$ be its leading block.
\State \Return $X$.
\EndProcedure
\end{algorithmic}
\end{algorithm}

\subsection{A positive semidefinite gradient representative}
\label{sec:more_positive_density}
The recovery iteration needs a positive semidefinite gradient representative and a compatible approximation to its proxy Hessian. The following construction supplies both for any pencil with $m$ variables and order $n$.

\begin{lemma}[A positive semidefinite gradient representative for primal recovery]
\label{lem:more_recovery_density}
At a strictly feasible point $y$, let $J$ satisfy $JS(y)J^\top=I$, let $B_i:=JA_iJ^\top$ in the original coefficient coordinates, and write $Q_*$ for the proxy Hessian at the optimal weight. For $0<\gamma<1/10$, one can compute an explicit $M_{\rm init}\in\mathbb S_+^n$ and a factored matrix $\widehat Q\succ0$ in $\widetilde O(n^\omega+\mathcal C(m))$ operations such that, with probability at least $1-\gamma$, the following properties hold.
\begin{itemize}
\item[(a)] The vector $\widehat g_i:=-\alpha\tr[M_{\rm init}B_i]$ satisfies $\|\widehat g-\nabla F(y)\|_{\nabla^2F(y)^{-1}}\leq1/100$.
\item[(b)] The matrix $Q_0[h,h]:=\tr[M_{\rm init}B_h^2]$ satisfies $Q_*/2\preceq Q_0\preceq2Q_*$, and $0.9Q_0\preceq\widehat Q\preceq1.1Q_0$.
\item[(c)] The trace satisfies $\tr[M_{\rm init}]\leq6m$.
\end{itemize}
The suppressed factors include logarithms of $1/\gamma$.
\end{lemma}
\begin{proof}
\emph{Construction.} Using the exact whitener $J$, normalize the coefficients and compute a weight of gap at most $\eta$ by Lemmas~\ref{lem:random_implicit_blocks} and~\ref{lem:weight_initialization}. Convert the weight back to the original coefficient basis as in~\eqref{eq:weight_basis_change}. Draw a fresh full estimate of width $\widetilde O(m)$ and keep its sketched factors:
\begin{align*}
 Y_0=F_0\Xi_0,\qquad Y_a=F_a\Xi_1\quad(a\geq1),\\
 Y=[Y_0\ \cdots\ Y_m],\qquad M_{\rm init}=YY^\top.
\end{align*}
The rank and probe count are $r,k=\widetilde\Theta(m)$, with the exact-factor fallback of~\eqref{eq:common_gram_exact_factors} when $2r+2\geq n$. Compute $\widehat Q$ by $\textsc{ProxyHessian}(B_1,\ldots,B_m,Y,\gamma/4)$. The gradient represented by $M_{\rm init}$ is $\widehat g_i=-\alpha\tr[M_{\rm init}B_i]$. Allocate total failure budget $\gamma/4$ to initialization, the range finder and the full-estimate accuracy bound.

\emph{Part~(a).} The full-estimate bound of Lemma~\ref{lem_more_projected_gradient}, with the width $q_{\rm fresh}$ of~\eqref{eq:fresh_width}, gives $\|\widehat g-g_G\|_{G/\alpha}\leq\delta_g$, independently of the dimension regime. Lemma~\ref{lem:warm_norm_comparison} and~\eqref{eq:single_loop_weight_bias} therefore imply $\|\widehat g-\nabla F(y)\|_{\nabla^2F(y)^{-1}} \leq\tfrac34\varepsilon_{\rm g}<1/100$.

\emph{Part~(b).} Let $Q$ be the proxy Hessian at $G$. Condition on the weight and on projections $\Pi_0,\Pi_1$ that satisfy the residual bound of Lemma~\ref{lem:common_gram_deflation}. Lemma~\ref{lem_more_metric}(a) gives $(39/40)Q\preceq Q_0\preceq(41/40)Q$ with failure probability $\gamma/4$; the required $O(\log^2(Cmn/\gamma))$ probes fit the full-estimate width after increasing its absolute constant. Lemma~\ref{lem:random_sharp_weight_accuracy}(a) gives $0.99Q_*\preceq Q\preceq1.01Q_*$, proving $Q_*/2\preceq Q_0\preceq2Q_*$.

Condition now on $Y$ and apply the independent compression guarantee of Lemma~\ref{lem_more_metric}(b). With conditional failure probability at most $\gamma/4$, it gives $(39/40)Q_0\preceq\widehat Q\preceq(41/40)Q_0$, which implies the stated comparison. Thus recovery uses the same proxy Hessian routine as centering.

\emph{Part~(c).} On the same event used in Part~(b), Lemma~\ref{lem_more_metric}(a) gives
\begin{align*}
 \tr M_{\rm init}\leq\tfrac32\tr[\mathcal L_y(G)^s]\leq3m<6m,
\end{align*}
where the weight accuracy gives the second inequality. The full-factor branch has trace at most $2m$ deterministically. The routine checks the bound $\tr M_{\rm init}\leq6m$ explicitly. A union bound over the initialization and full-estimate event, the joint event for $Q_0$ and the mass, and the conditional compression event gives total failure probability at most $\gamma$.

\emph{Cost.} The normalization, the weight, the projections, the full estimate and the compression of $Q_0$ cost $\widetilde O(n^\omega+\mathcal C(m))$ by Lemmas~\ref{lem:random_implicit_blocks}, \ref{lem:weight_initialization} and~\ref{lem_more_metric}. The matrix $Y$ has $d_Y=\widetilde O(m^2)$ columns, so forming $M_{\rm init}=YY^\top$ costs $\Tmat(n,d_Y,n)=\widetilde O(\Tmat(n,n,m^2)) \leq\widetilde O(\mathcal C(m))$. Every width and trial count is fixed in advance, so the same cost holds on every run.
\end{proof}

\subsection{Correcting the residual by a modified Newton iteration}
\label{sec:more_primal_congruence}
We next correct the constraint residual of the positive semidefinite representative from Lemma~\ref{lem:more_recovery_density} by congruence updates. The following lemma shows that a fixed approximate proxy Hessian gives geometric convergence while preserving positive semidefiniteness and controlling the trace.

\begin{lemma}[Modified Newton iteration over positive congruences]
\label{lem:more_congruence_correction}
Let $M_{\rm init}\succeq0$, let $B_1,\ldots,B_{\bar m}$ be symmetric matrices, and write $B_w:=\sum_iw_iB_i$ and $Q_0[w,w]:=\tr[M_{\rm init}B_w^2]$. Suppose that $Q_0\succ0$ and
\begin{align*}
 0.9Q_0\preceq\widehat Q\preceq1.1Q_0,
 \qquad \|B_v\|_{\mathrm{op}}\leq6\|v\|_{Q_0}
 \quad\text{for every }v\in\R^{\bar m}.
\end{align*}
Let $c_{\rm rec}\in\R^{\bar m}$, and define
\begin{equation}\label{eq:more_recovery_iteration}
 M(w):=(I+B_w)M_{\rm init}(I+B_w),\qquad
 f(w)_i:=\tr[B_iM(w)]-(c_{\rm rec})_i,\qquad
 w_{j+1}:=w_j-\tfrac12\widehat Q^{-1}f(w_j),
\end{equation}
with $w_0:=0$. If $\|f(0)\|_{\widehat Q^{-1}}\leq10^{-4}$, then the following statements hold.
\begin{itemize}
\item[(a)] Every iterate satisfies $\|w_j\|_{\widehat Q}\leq10^{-3}$, $M(w_j)\succeq0$ and $\tr[M(w_j)]\leq1.02\tr[M_{\rm init}]$.
\item[(b)] The iterates converge to a zero of $f$, and $\|f(w_j)\|_{\widehat Q^{-1}}\leq2.24\cdot10^{-3}4^{-j}$.
\end{itemize}
\end{lemma}
\begin{proof}
The map $f$ is quadratic, and for directions $h,v$ its derivative satisfies
\begin{align*}
 h^\top\bigl(\nabla f(w)[v]-2Q_0v\bigr)
 =\tr[M_{\rm init}(B_wB_hB_v+B_vB_hB_w)].
\end{align*}
Cauchy--Schwarz and the bound on $\|B_h\|_{\rm op}$ give
\begin{align*}
 \bigl|h^\top\bigl(\nabla f(w)[v]-2Q_0v\bigr)\bigr|
 \leq2\|B_h\|_{\mathrm{op}}\|w\|_{Q_0}\|v\|_{Q_0}
 \leq12\|h\|_{Q_0}\|w\|_{Q_0}\|v\|_{Q_0}.
\end{align*}
Since $\|\cdot\|_{Q_0}\leq(0.9)^{-1/2}\|\cdot\|_{\widehat Q}$, it follows that
\begin{align*}
 \bigl\|\widehat Q^{-1/2}\bigl(\nabla f(w)-2Q_0\bigr)\widehat Q^{-1/2}\bigr\|_{\mathrm{op}}
 \leq12(0.9)^{-3/2}\|w\|_{\widehat Q}.
\end{align*}
Since $\widehat Q^{-1/2}Q_0\widehat Q^{-1/2}$ has spectrum in $[1/1.1,1/0.9]$, the derivative of $T(w)=w-\widehat Q^{-1}f(w)/2$ has $\widehat Q$-operator norm at most
\begin{align*}
 1/9+6(0.9)^{-3/2}\cdot10^{-3}<1/4
\end{align*}
on $\|w\|_{\widehat Q}\leq10^{-3}$. Also $\|T(0)\|_{\widehat Q}\leq5\cdot10^{-5}$, so $T$ maps this ball into itself and is a $1/4$-contraction there.

\emph{Part~(a).} Positive semidefiniteness follows from congruence, also when $M_{\rm init}$ is singular. On the ball, $\|B_w\|_{\rm op}\leq6\|w\|_{Q_0}\leq6\cdot10^{-3}/\sqrt{0.9}<0.0064$, and since $M_{\rm init}\succeq0$,
\begin{align*}
 \tr[M(w)]\leq\|I+B_w\|_{\mathrm{op}}^2\tr[M_{\rm init}]
 <1.02\tr[M_{\rm init}].
\end{align*}

\emph{Part~(b).} By the contraction mapping theorem, the iterates converge to a fixed point $w_*$ of the map in the ball, which is a zero of $f$, and $\|w_j-w_*\|_{\widehat Q}\leq10^{-3}4^{-j}$. On the ball,
\begin{align*}
 \bigl\|\widehat Q^{-1/2}\nabla f(w)\widehat Q^{-1/2}\bigr\|_{\rm op}
 \leq2/0.9+12(0.9)^{-3/2}\cdot10^{-3}<2.24,
\end{align*}
and integration along the segment from $w_j$ to $w_*$ proves Part~(b).
\end{proof}

\subsection{Removing the artificial scalar blocks}
\label{sec:penalty_blocks}
The correction is performed for the augmented pencil, whereas the final primal matrix is its leading block. The next lemma bounds the residual and objective error of this block in terms of the augmented residual, complementarity and penalty parameters.

\begin{lemma}[Primal recovery from the auxiliary initialization SDP]
\label{lem:penalty_recovery_identities}
Let $u=(y,z)$ satisfy $\overline S(u)\succ0$ and $z<0$, and let $\overline X\succeq0$ be an augmented primal matrix. Let $X$ be its leading block, and let $p^+,p^-\in\R^m$ and $p_{\max},p_\delta\in\R$ be its diagonal entries in the blocks of $2B\mathbf1+y$, $2B\mathbf1-y$, $z_{\max}-z$ and $\delta+z$, respectively. Write
\begin{align*}
 e:=\overline{\mathcal A}(\overline X)-c=(e_y,e_z),\qquad
 q_0:=\langle\overline S(u),\overline X\rangle.
\end{align*}
Then $X\succeq0$, the entries $p^\pm,p_{\max},p_\delta$ are nonnegative, and
\begin{align}
 B\mathbf1^\top(p^++p^-)&\leq q_0+K\delta+3B\|e\|_1,
 \label{eq:penalty_scalar_mass}\\
 \|\mathcal A(X)-b\|_1&\leq4\|e\|_1+(q_0+K\delta)/B,
 \label{eq:penalty_residual}\\
 b^\top y-\langle C,X\rangle&\leq q_0+K\delta+2B\|e\|_1.
 \label{eq:penalty_objective}
\end{align}
\end{lemma}
\begin{proof}
Only diagonal blocks enter the coefficient pairings; their positive semidefiniteness gives the sign claims. Since $z<0$, the point $y$ is feasible for the original problem, so $b^\top y_*\leq b^\top y$. At $(y_*,0)$ the leading slack is positive semidefinite and the slacks of the box constraints on $y$ are at least $B$. Since $\overline S$ is affine,
\begin{align*}
 B\mathbf1^\top(p^++p^-)
 &\leq\langle\overline S(y_*,0),\overline X\rangle
 =q_0+(y_*-y)^\top(b+e_y)-z(K+e_z)\\
 &\leq q_0+K\delta+3B\|e_y\|_1+\delta|e_z|
 \leq q_0+K\delta+3B\|e\|_1.
\end{align*}
The second line uses $b^\top y_*\leq b^\top y$, $\|y_*-y\|_\infty\leq3B$ and $-\delta<z<0$, and the last step uses $\delta\leq B$. This proves~\eqref{eq:penalty_scalar_mass}. The constraints for the coordinates of $y$ give
\begin{align*}
 \mathcal A(X)-b=e_y-p^++p^-,
 \qquad
 \|\mathcal A(X)-b\|_1\leq\|e_y\|_1+\mathbf1^\top(p^++p^-),
\end{align*}
and~\eqref{eq:penalty_scalar_mass} proves~\eqref{eq:penalty_residual}. Finally, expanding $q_0$ over the blocks gives
\begin{align*}
 b^\top y-\langle C,X\rangle
 &=q_0-Kz-y^\top e_y-ze_z
       -2B\mathbf1^\top(p^++p^-)-z_{\max}p_{\max}-\delta p_\delta\\
 &\leq q_0+K\delta+2B\|e_y\|_1+\delta|e_z|
 \leq q_0+K\delta+2B\|e\|_1,
\end{align*}
where the inequality uses the nonnegativity of the discarded entries and $\|y\|_\infty<2B$. \end{proof}

\subsection{Proof of the postprocessing guarantee}
\label{sec:recovery_proof}
\begin{proof}[Proof of Theorem~\ref{thm:postprocessing}]
\emph{The initial residual.} Condition on any history that fixes the input $(u,t)$. Lemma~\ref{lem:more_recovery_density}, applied to $\overline{\mathcal P}$ at $u$ with fresh randomness, succeeds with conditional probability at least $1-\gamma$, and we work on its success event. Since $f(0)=-(\widehat g+tc)/\alpha$, Lemma~\ref{lem:more_recovery_density}(a) and the input condition $\lambda_{tc}(u)\leq1/100$ give
\begin{align*}
 \|f(0)\|_{\widehat Q^{-1}}=\alpha^{-1}\|\widehat g+tc\|_{\widehat Q^{-1}}
 \leq\frac2{100}\sqrt{\frac{2s+1}{0.45\alpha}}
 <10^{-4}.
\end{align*}
The norm conversion uses $\widehat Q\succeq0.45Q_*$ and $\nabla^2F\preceq\alpha(2s+1)Q_*$, and the last inequality uses $\alpha=10^6s^2$ and $s\geq2$. Moreover, \eqref{eq:random_operator_metric} and $Q_0\succeq Q_*/2$ give
\begin{align*}
 \|B_v\|_{\rm op}\leq3\sqrt2\|v\|_{Q_0}<6\|v\|_{Q_0}.
\end{align*}
Together with Lemma~\ref{lem:more_recovery_density}(b), these are the hypotheses of Lemma~\ref{lem:more_congruence_correction}.

\emph{The trace check.} Since $\widehat Q\preceq1.1Q_0\preceq2.2Q_*\preceq2.2\nabla^2F(u)/\alpha$, and Lemma~\ref{lem:cost_bounded_metric} with $\lambda_{tc}(u)\leq1/100$ bounds the largest eigenvalue of $\nabla^2F(u)$, we have
\begin{align*}
 \tr[\widehat Q]\leq\frac{2.2\bar m}\alpha\lambda_{\max}(\nabla^2F(u))\leq U_Q.
\end{align*}
Hence the check of~\eqref{eq:more_recovery_metric_cap} passes.

\emph{The iteration.} After $N_{\rm rec}$ steps, Lemma~\ref{lem:more_congruence_correction} gives $\|f(w)\|_{\widehat Q^{-1}}\leq2.24\cdot10^{-3}4^{-N_{\rm rec}}\leq\delta_{\rm rec}$, and together with Lemma~\ref{lem:more_recovery_density}(c) it gives $\tr[M(w)]\leq1.02\tr[M_{\rm init}]\leq6.12\bar m<7\bar m$. Hence both final checks pass. The matrix $\overline X=(\alpha/t)J^\top M(w)J$ is positive semidefinite, and its residual $e:=\overline{\mathcal A}(\overline X)-c=(\alpha/t)f(w)$ and its complementarity satisfy
\begin{align*}
 \|e\|_1\leq\frac\alpha t\sqrt{\bar m}\sqrt{\tr[\widehat Q]}\,\delta_{\rm rec}
 \leq\frac\epsilon{16B},\qquad q_0:=\langle\overline S(u),\overline X\rangle=\frac\alpha t\tr[M(w)]
 \leq\frac{7\nu}t\leq\frac{7\epsilon}{128}<\frac\epsilon8.
\end{align*}
\emph{The scalar blocks.} The leading block $X$ of $\overline X$ is positive semidefinite, and Lemma~\ref{lem:penalty_recovery_identities} gives
\begin{align*}
 \|\mathcal A(X)-b\|_1\leq4\|e\|_1+(q_0+K\delta)/B
 \leq7\epsilon/16,\qquad b^\top y-\langle C,X\rangle\leq q_0+K\delta+2B\|e\|_1
 \leq5\epsilon/16,
\end{align*}
since $K\delta=\epsilon/16$ and $B\geq1$. These are the asserted guarantees, and they hold with conditional probability at least $1-\gamma$.
\end{proof}

\subsection{Running time}
\label{sec:recovery_cost}
We finally account for the initial estimate and the modified Newton iterations. The explicit metric cap bounds the number of iterations on every run, including runs that report failure.

\begin{lemma}[Work of postprocessing]
\label{lem:postprocessing_work}
Under the input assumptions of Theorem~\ref{thm:postprocessing}, Algorithm~\ref{alg:more_congruence_recovery} uses
\begin{align*}
 \widetilde O\bigl(\bar n^\omega+\mathcal C(\bar m)\bigr)
   +O\Bigl((\bar n^\omega+\bar m\bar n^2)\log\Bigl(2+\frac{\alpha B\sqrt{\bar mU_Q}}{t\epsilon}\Bigr)\Bigr)
\end{align*}
arithmetic operations on every run, including runs that report failure, where $\mathcal C$ is the cost function of Definition~\ref{def_more_block_cost} for $\overline{\mathcal P}$. The suppressed factors include logarithms of $1/\gamma$. For every fixed $\zeta>0$, this implies the bound
\begin{equation}\label{eq:more_recovery_cost}
 \widetilde O_\zeta\left(
 (n^{\omega+\zeta}+n^2m^{\omega-1+\zeta})
 \log^{O(1)}\kappa\right),
\end{equation}
where $\kappa$ is the input scale in~\eqref{eq:complete_kappa}. These bounds use only that $\overline S(u)\succ0$ and that $t$ lies in the stated interval, because the widths of the oracle query are fixed and the cap $U_Q$ bounds the number of iterations.
\end{lemma}
\begin{proof}
Lemma~\ref{lem:more_recovery_density} charges the oracle query by $\widetilde O(\bar n^\omega+\mathcal C(\bar m))$ on every run. A failed query or trace check terminates recovery. Otherwise the cap and~\eqref{eq:more_recovery_tolerance} give
\begin{align*}
 N_{\rm rec}
 =O\Bigl(\log\Bigl(2+\frac{\alpha B\sqrt{\bar mU_Q}}{t\epsilon}\Bigr)\Bigr).
\end{align*}
Each iteration forms and whitens one coefficient combination, applies the congruence, pairs with the coefficients and solves a factored $\widehat Q$-system, costing $O(\bar n^\omega+\bar m\bar n^2)$. The final matrix and checks have the same cost.

\emph{Original dimensions.} We have $\bar m=O(m)$ and $\bar n=O(n)$. Splitting into square blocks of side $\bar m$ gives
\begin{align*}
 \mathcal C(\bar m)\leq O_\zeta(\bar n^2\bar m^{\omega-1+\zeta}),
\end{align*}
because each term of $\mathcal C(\bar m)$ either consists of $O(\bar n^2/\bar m)$ square products of side $\bar m$ or is at most $\bar n\bar m^2$. Hence the oracle query costs $\widetilde O_\zeta(n^{\omega+\zeta}+n^2m^{\omega-1+\zeta})$. The penalty parameters make $B,r^{-1},W_c$ polynomially bounded in $\kappa$, and $128\nu/\epsilon\leq t\leq(1+\gamma_t)128\nu/\epsilon$. Thus $U_Q$ is polynomially bounded in $\kappa$ and the barrier parameters, the number of iterations is $O(\log\kappa)$, and the cost $O(n^\omega+mn^2)$ per iteration fits the same bound. This proves~\eqref{eq:more_recovery_cost}.
\end{proof}
\section{Putting the three components together}
\label{sec:solver_wrapup}
The three modules use independent oracle states and communicate only through their returned points and multipliers.

\subsection{The complete solver}
We combine the preprocessing guarantee of Theorem~\ref{thm:preprocessing}, the optimization guarantee of Theorem~\ref{thm_more_amortized_sdp} and the postprocessing guarantee of Theorem~\ref{thm:postprocessing}. The resulting solver returns both primal and dual outputs with the total cost stated below.

\begin{theorem}[Complete randomized SDP solver]
\label{thm:more_complete_sdp}
Let $1\leq m\leq n$, let $A_1,\ldots,A_m,C$ be real symmetric $n\times n$ matrices, and let $b\in\R^m$. Assume that the $A_i$ are linearly independent and $\bigcap_i\ker A_i=\{0\}$. Suppose that numbers $B\geq1$ and $\rho>0$ are supplied such that there exist a dual optimum $y_*$ and a point $y_{\rm f}$ with
\begin{align*}
 \|y_*\|_\infty\leq B,\qquad
 \|y_{\rm f}\|_\infty\leq B,\qquad S(y_{\rm f})\succeq\rho I.
\end{align*}
Let $0<\epsilon\leq1$ and $0<\beta<1/10$, and let $\kappa$ be the input scale in~\eqref{eq:complete_kappa}. Then Algorithm~\ref{alg:more_complete_sdp} has the following guarantees.
\begin{itemize}
\item[(a)] With probability at least $1-\beta$, it returns a vector $y$ and an explicit matrix $X\succeq0$ satisfying
\begin{gather*}
 S(y)\succ0,\qquad 0\leq b^\top y-\operatorname{OPT}\leq\epsilon,\\
 \sum_{i=1}^m|\langle A_i,X\rangle-b_i|\leq\epsilon,\qquad \langle C,X\rangle\geq\operatorname{OPT}-\epsilon.
\end{gather*}
\item[(b)] Under Fact~\ref{fact:rectangular_thresholds} and a square exponent bound $2<\omega\leq9/4$, for every fixed $\zeta>0$, its total arithmetic cost on every run is
\begin{align*}
 \widetilde O_\zeta(
 (n^{\omega+\zeta}+n^2m^{3/2+\zeta})
 \log^{O(1)}\kappa).
\end{align*}
The suppressed factors include logarithms of $1/\beta$. The bound includes preprocessing, optimization of one augmented pencil, and explicit primal recovery, in the arithmetic model of Section~\ref{sec:algorithm_preliminaries}.
\end{itemize}
\end{theorem}

\begin{remark}
The hypotheses of Theorem~\ref{thm:more_complete_sdp} concern only the coefficient matrices and the dual problem. The point $y_{\rm f}$ makes the dual problem strictly feasible, and $y_*$ makes its optimal value finite. Hence, the primal and dual optimal values agree~\cite[Theorem~3.1]{vb96}, and $\operatorname{OPT}=b^\top y_*$ is their common value, even when every primal-feasible matrix is singular. The algorithm applies the barrier only to the augmented pencil of Lemma~\ref{lem:penalty_pencil}, whose primal and dual problems are strictly feasible by construction, and Lemmas~\ref{lem:penalty_pencil} and~\ref{lem:penalty_recovery_identities} transfer its guarantees to the original problem through $y_*$ and $y_{\rm f}$ alone.
\end{remark}

\begin{algorithm}[!ht]
\caption{Complete randomized SDP solver}
\label{alg:more_complete_sdp}
\begin{algorithmic}[1]
\Procedure{CompleteAmortizedSDP}{$A_1,\ldots,A_m,C,b,B,\rho,\epsilon,\beta$} \Comment{Theorem~\ref{thm:more_complete_sdp}}
\State Form the pencil $\overline{\mathcal P}$, the objective $c:=(b,K)$ and the bounds $u_0,r,R,W_c$ of Lemma~\ref{lem:penalty_pencil}.
\State Use the slack access of Algorithm~\ref{alg:implicit_slack} in all oracle states, with the dimensions $(\bar m,\bar n)$.
\State $(u_{\rm c},t_0)\gets \Call{PreprocessBoundedStart}{\overline{\mathcal P},c,u_0,r,R,\beta/3}$. \Comment{Theorem~\ref{thm:preprocessing}}
\State $H\gets\lfloor\bar m^{1/6}\rfloor$.
\State $(u,t)\gets\Call{AmortizedSDP}{\overline{\mathcal P},c,u_{\rm c},t_0, \epsilon/64,\beta/3;H}$. \Comment{Theorem~\ref{thm_more_amortized_sdp}}
\State Write $u=(y,z)$, and report failure if $z>-\delta/4$. \Comment{Lemma~\ref{lem:penalty_pencil}}
\State $X\gets\Call{RecoverPrimal}{\overline{\mathcal P},c,u,t,B,r,W_c,\epsilon,\beta/3}$. \Comment{Theorem~\ref{thm:postprocessing}}
\State \Return $(y,X)$.
\EndProcedure
\end{algorithmic}
\end{algorithm}

A failure reported by any module terminates the algorithm. Within preprocessing and the core, weights persist across phases. Fresh estimates are used when $\bar m^2\leq\bar n$, and maintained estimates are used otherwise. Recovery uses fresh randomness after its input is fixed.

\subsection{Correctness of the composition}
\label{sec:more_complete_solver}
\begin{proof}[Proof of Theorem~\ref{thm:more_complete_sdp}\textup{(a)}]
Theorem~\ref{thm:preprocessing} supplies an admissible augmented pencil and $u_{\rm c},t_0$ with $\overline S(u_{\rm c})\succ0$, $\lambda_{t_0c}(u_{\rm c})\leq1/100$ and $0<t_0\leq1$. Theorem~\ref{thm_more_amortized_sdp}, at accuracy $\epsilon/64$, returns
\begin{align*}
 c^\top u-\inf_{\overline S\succ0}c^\top u\leq\epsilon/64,\qquad
 \lambda_{tc}(u)\leq1/100,\qquad
 \frac{128\nu}\epsilon\leq t\leq(1+\gamma_t)\frac{128\nu}\epsilon.
\end{align*}
Lemma~\ref{lem:penalty_pencil} then gives $-\delta<z\leq-\delta/4$, $S(y)\succ0$ and $0\leq b^\top y-\operatorname{OPT}\leq3\epsilon/64$. The sign check passes, and Theorem~\ref{thm:postprocessing} applies:
\begin{align*}
 X\succeq0,\qquad
 \|\mathcal A(X)-b\|_1\leq7\epsilon/16,\qquad
 \langle C,X\rangle\geq b^\top y-5\epsilon/16
 \geq\operatorname{OPT}-\epsilon.
\end{align*}
Each module, conditioned on any successful preceding history, has failure probability at most $\beta/3$ under its fresh randomness. Summing the probabilities that each module is the first to fail gives total failure probability at most $\beta$.
\end{proof}

\subsection{Total running time}
\label{sec:complete_running_time}
\begin{proof}[Proof of Theorem~\ref{thm:more_complete_sdp}\textup{(b)}]
The augmented dimensions satisfy $\bar m=m+1=O(m)$, $\bar n=n+2m+2=O(n)$ and $\bar m\leq\bar n$. Theorem~\ref{thm:preprocessing} gives the claimed preprocessing cost. Its deterministic multiplier satisfies $t_0^{-1}=\max\{1,R\|c\|_2/\delta_{\rm c}\}$, where $c=(b,K)\ne0$. Since $R$ and $\|c\|_2$ are polynomially bounded in $\kappa$,
\begin{align*}
 \log\left(2+\frac{64\nu}{t_0\epsilon}\right)=O(\log\kappa),
 \qquad L_\epsilon=O(\log\kappa).
\end{align*}
Theorem~\ref{thm_more_amortized_sdp}(c) therefore gives the claimed core cost, including its independent initialization and slack anchors. Lemma~\ref{lem:postprocessing_work} gives recovery cost \eqref{eq:more_recovery_cost}, which fits the same bound because $\omega-1\leq5/4<3/2$.

These costs also cover an undetected failure in an earlier module. The core schedule and sketch widths remain deterministic, so its prefix cost does not require a central input. Recovery either reports a failed slack factorization or proceeds with fixed sketch widths and the explicit cap $U_Q$ of Lemma~\ref{lem:postprocessing_work}, whose cost likewise does not require centrality. Final checks and writing $X$ are included in these bounds.
\end{proof}
\section*{AI Methodology}
\addcontentsline{toc}{section}{AI Methodology}
\begin{flushleft}
\renewcommand{\arraystretch}{1.3}
\begin{tabular}{@{}l@{\qquad}c@{}}
\textbf{Activity} & \textbf{AI used} \\
Asking the research question & \texttt{No} \\
Coming up with the approach & \texttt{Yes} \\
Proof development & \texttt{Yes} \\
Writing and exposition & \texttt{Yes} \\
Checking for bugs & \texttt{Yes} \\
Other supporting tasks & \texttt{Yes} \\
\end{tabular}
\end{flushleft}

\begin{sloppypar}
\paragraph{Details} In preparing this paper, the authors used the following AI tools: GPT-6 Astra, ChatGPT 6 Pro, Claude Fable 5.1, Claude Opus 5.5, and Gemini 4 Argon. The authors first attempted to generalize the Lee--Sidford barrier to SDPs while working on~\cite{hjstrz22}, but were unable to find a matrix analog of the potential function for LP~\eqref{eq:geo_lp_potential} with the desired properties. Around July 2026, under the authors' guided prompting, GPT-6 Astra constructed the matrix potential~\eqref{eq:geo_weight_form} for the Lee--Sidford barrier. Its initial proof was quite involved and relied on a chain of specialized matrix inequalities. While interacting with these AI tools to verify and understand the construction, the authors realized that the dual formulation of the matrix potential~\eqref{eq:geo_potential} provides a much clearer geometric picture and that the proof can be greatly simplified using matrix-exponential curves (Section~\ref{sec:geo_proof}). On the algorithmic side, the authors originated the core ideas for an efficient implementation of the SDP solver, including matrix sketching, amortized maintenance, and fast rectangular matrix multiplication (Section~\ref{sec:overview_sketch}), while using these AI tools to assist with detailed calculations. An earlier version of the algorithm had running time $\widetilde{O}(\sqrt{m}(n^\omega + mn^2))$. Through further rounds of interaction, ChatGPT 6 Pro suggested reusing the slack matrix factorization across multiple iterations to amortize the $n^\omega$ cost (Section~\ref{sec:overview_slack}). The final algorithm presented in this paper emerged from several rounds of simplification by the authors and these AI tools.

The authors also used these AI tools to assist with drafting and revising the exposition, plotting the figures, and improving the presentation of the proofs. The authors verified the correctness and originality of all content, including references, and take full responsibility for the paper.
\end{sloppypar}
\clearpage
\bibliographystyle{alpha}
\bibliography{ref}

\end{document}